\documentclass[12pt]{amsart}
\usepackage{amssymb}
\usepackage{amsbsy}
\usepackage{amscd}
\usepackage[mathscr]{eucal}
\usepackage{verbatim}
\usepackage{version}
\usepackage{diagrams}
\usepackage{pstricks}
\usepackage{pst-all}

\makeatletter
\def\cal{\mathcal}
\def\Bbb{\mathbb}

\newenvironment{pf}{\proof[\proofname]}{\endproof}
\newenvironment{pf*}[1]{\proof[#1]}{\endproof}
\makeatother

\newenvironment{aenume}{%
  \begin{enumerate}%
  }{\end{enumerate}}
\makeatletter
\@ifclasslater{amsart}{1999/11/24}{}{
\renewcommand*\subjclass[2][1991]{%
  \def\@subjclass{#2}%
  \@ifundefined{subjclassname@#1}{%
    \ClassWarning{\@classname}{Unknown edition (#1) of Mathematics
      Subject Classification; using '1991'.}%
  }{%
    \@xp\let\@xp\subjclassname\csname subjclassname@#1\endcsname
  }%
}
\renewcommand{\subjclassname}{%
  \textup{1991} Mathematics Subject Classification}
\@xp\let\csname subjclassname@1991\endcsname \subjclassname
\@namedef{subjclassname@2000}{%
  \textup{2000} Mathematics Subject Classification}
}
\makeatother
\newenvironment{NB}{
{\bf NB}. \footnotesize \color{red}
}{}
\excludeversion{NB}

\newtheorem{Theorem}[equation]{Theorem}
\newtheorem{Corollary}[equation]{Corollary}
\newtheorem{Lemma}[equation]{Lemma}
\newtheorem{Proposition}[equation]{Proposition}

\theoremstyle{definition}
\newtheorem{Definition}[equation]{Definition}

\newtheorem{Notation}[equation]{Notation}

\theoremstyle{remark}
\newtheorem{Remark}[equation]{Remark}
\newtheorem*{Claim}{Claim}

\numberwithin{equation}{section}

\newcommand{\thmref}[1]{Theorem~\ref{#1}}
\newcommand{\secref}[1]{\S\ref{#1}}
\newcommand{\lemref}[1]{Lemma~\ref{#1}}
\newcommand{\propref}[1]{Proposition~\ref{#1}}
\newcommand{\corref}[1]{Corollary~\ref{#1}}
\newcommand{\subsecref}[1]{\S\ref{#1}}

\newcommand{\defeq}{{:=}}
\newcommand{\C}{{\Bbb C}}
\newcommand{\Z}{{\Bbb Z}}
\newcommand{\Q}{{\Bbb Q}}
\newcommand{\R}{{\Bbb R}}
\newcommand{\proj}{{\Bbb P}}

\newcommand{\End}{\operatorname{End}}

\newcommand{\Ext}{\operatorname{Ext}}

\newcommand{\len}{\mathop{\text{\rm len}}\nolimits}

\newcommand{\pd}[2]{\frac{\partial#1}{\partial#2}}

\newcommand{\ve}{\varepsilon}

\newcommand{\Pic}{\operatorname{Pic}}

\newcommand{\ch}{\operatorname{ch}}
\newcommand{\Wedge}{{\textstyle \bigwedge}}
\newcommand{\Todd}{\operatorname{Todd}}

\newcommand{\Li}{\operatorname{Li}}
\newcommand{\Zin}{Z^{\text{\rm inst}}}

\newcommand{\bZin}{\widehat{Z}^{\text{\rm inst}}}
\newcommand{\bZ}{\widehat{Z}}

\newcommand{\Fper}{F^{\text{\rm pert}}}
\newcommand{\Fin}{F^{\text{\rm inst}}}

\newcommand{\Finz}{\mathcal F^{\text{\rm inst}}}
\newcommand{\q}{\mathfrak q}
\newcommand{\bbeta}{\boldsymbol\beta}

\newcommand{\rk}{\mathop{{\rm rk}}}

\newcommand{\Coeff}{\mathop{\text{\rm Coeff}}}
\def\oo{{\cal O}}
\def\I{{\cal I}}
\def\AA{{\cal A}}
\def\A{{\mathbb A}}

\def\P{{\mathbb P}}
\def\<{\langle}
\def\>{\rangle}

\def\F{{\mathcal F}}
\def\G{{\mathcal G}}
\def\E{{\mathcal E}}
\def\cc{{\mathcal C}}

\newcommand{\sn}{\operatorname{sn}}
\newcommand{\cn}{\operatorname{cn}}
\newcommand{\dn}{\operatorname{dn}}

\newcommand{\Res}{\operatornamewithlimits{Res}}

\makeatletter
\newcommand{\vechatom}{
    {\Vec{\omega}}
    \,\smash[b]{\hbox{\lower2\ex@\hbox{$\m@th\hat{\null}$}}}
}
\makeatother

\begin{document}
\title[K-theoretic  Donaldson invariants via instanton counting]
{K-theoretic  Donaldson invariants via instanton counting}
\author{Lothar G\"ottsche}
\address{International Centre for Theoretical Physics, Strada Costiera 11, 
34014 Trieste, Italy}
\email{gottsche@ictp.trieste.it}
\author{Hiraku Nakajima}
\address{Department of Mathematics, Kyoto University, Kyoto 606-8502,
Japan}
\email{nakajima@math.kyoto-u.ac.jp}
\thanks{The second author is supported by the Grant-in-aid
for Scientific Research (No.15540023, 17340005), JSPS}

\author{K\={o}ta Yoshioka}
\address{Department of Mathematics, Faculty of Science, Kobe University,
Kobe 657-8501, Japan}
\email{yoshioka@math.kobe-u.ac.jp}

\dedicatory{To Friedrich Hirzebruch on the occasion of  his eightieth birthday}



\begin{abstract} In this paper we study the holomorphic Euler characteristics
of determinant line bundles on moduli spaces of rank $2$ semistable sheaves
on an algebraic surface $X$, which can be viewed as $K$-theoretic versions of the Donaldson invariants. In particular if  $X$ is a smooth projective toric surface, we determine these invariants and their wallcrossing in terms of the $K$-theoretic version of
the Nekrasov partition function (called 5-dimensional supersymmetric Yang-Mills theory compactified on a circle in the physics literature). Using the results of \cite{NY3} 
we give an explicit generating function for the  wallcrossing of these invariants in terms of elliptic functions and 
modular forms. 
\end{abstract}

\maketitle
\tableofcontents

\begin{NB}
  I have put the TOC at the end of the introduction several times (as
  I thought it is more reasonable), but the publisher always move it
  to the top. Nov. 9, HN.  
\end{NB}

\section*{Introduction}

This paper is a sequel to \cite{GNY}.%
\begin{NB}
I have add `a'. Nov. 27, HN.
\end{NB}
In \cite{GNY} we expressed the wallcrossing terms of equivariant
Donaldson invariants for a smooth toric surface in terms of the
Nekrasov partition function, and then using the solution of the
Nekrasov conjecture \cite{NY1},\cite{NO},\cite{BE} and its refinement
\cite{NY2} we gave the wallcrossing formula for simply connected
projective surfaces with $p_g = 0$ in terms of modular forms, thus
recovering%
\begin{NB} corrected 27.8.LG\end{NB}
 the formula in \cite{G} originally proved assuming the
Kotschick-Morgan conjecture~\cite{KM}.%
\begin{NB}
I add a sentence. Nov. 27, HN.
\end{NB}
The Nekrasov partition function is defined as the generating function
of the integrals of the equivariant cohomology class $1$ on the
Uhlenbeck partial compactifications $M_0(r,n)$ of the moduli spaces
of $SU(r)$-instantons on $\mathbb R^4$ with $c_2=n$. (As $M_0(r,n)$ is
noncompact, we need a justification of the integration. See
\cite{NY1} for details.)

There is a natural $K$-theoretic counterpart of the Nekrasov partition
function, namely we replace the integration in equivariant
cohomology by the character of the coordinate ring of $M_0(r,n)$,
where we view $M_0(r,n)$ as an affine algebraic variety via the ADHM
description. The coordinate ring itself is infinite dimensional, but
the weight spaces are 
finite dimensional (see \cite{NY1}), so the 
character is well-defined. This $K$-theoretic counterpart is called
the $5$-dimensional supersymmetric Yang-Mills theory compactified on a
circle in the physics literature \cite{Nek5d},\cite{LN}.
In \cite{NY3} we proved the analogues of the results obtained
in \cite{NY1} in the $K$-theoretic version. (The approach in
\cite{NO} can be applied to the $K$-theoretic version, while it seems
difficult to generalize that of \cite{BE}.)
There is also a mathematical reason why we should study the $K$-theoretic
Nekrasov partition function. By the geometric engineering of Katz,
Klemm and Vafa \cite{KKV}, it is (after a parameter is specialized)
equal to the generating function of all genus, all degree
Gromov-Witten invariants for a certain noncompact toric Calabi-Yau
$3$-fold. (See \cite{Zhou} for a mathematically rigorous proof).%
\begin{NB}
Corrected. H.N. Nov. 27  
\end{NB}
Gromov-Witten invariants for toric Calabi-Yau $3$-folds have been
studied intensively both in mathematics and physics (see e.g.\
\cite{MNOP} and the references therein).

On the other hand, the $K$-theoretic Donaldson invariants have {\it
  not\/} been studied very much in the mathematical literature, as far as
the authors know. One of the reasons might be a lack of
motivation, as it is unlikely that there is an application to
$4$-dimensional topology. But another reason seems to lie in 
technical difficulties in defining the invariants. For example, the
dimension counting argument used in the definition of the Donaldson
invariants cannot be applied to the $K$-theoretic situation.
Instead of attacking this problem, we restrict our interest to the
case when the base $4$-manifold is a projective surface $X$. Then we
can use Gieseker-Maruyama moduli spaces of semistable sheaves and
define the $K$-theoretic Donaldson invariants as the holomorphic Euler
characteristics of the determinant line bundles. Then the
algebro-geometric techniques used in \cite{GNY} to derive the
wallcrossing formula for the ordinary Donaldson invariants can be
equally applied to the $K$-theoretic invariants. We will express the
generating function of wallcrossing terms of the $K$-theoretic
Donaldson invariants in terms of elliptic functions, which have a
power series development in terms of modular forms. Their lowest order
terms are the modular forms which occur in the wallcrossing formula in
\cite{GNY} for the usual Donaldson invariants. If the moduli spaces
are smooth of the expected dimension, it is easy to see that this
is compatible with  the Hirzebruch-Riemann-Roch formula.
Our approach is very similar to the one in \cite{GNY}, though the
final step identifying invariants with the $q$-developments of modular
forms and elliptic functions is more involved than in \cite{GNY}.
We want to remark that our final answer for  the wallcrossing formula
strongly suggests that there should exist a definition of
$K$-theoretic Donaldson invariants for any $4$-manifold with a
$Spin^c$-structure (see \subsecref{subsec:digress}).

The holomorphic Euler characteristics of determinant line bundles are
interesting algebro-geometric objects in their own right.  They are
refinements of the usual Donaldson invariants, which contain a lot of
geometrical information about the moduli spaces of stable sheaves on
$X$, their Uhlenbeck compactifications and the linear systems on them.
For instance by a result of \cite{Li} the morphism associated to
certain determinant line bundles defines a projective embedding of the
Uhlenbeck moduli spaces. The corresponding Donaldson invariants will
determine the degree of the Uhlenbeck compactification and under
suitable assumptions one would expect that the $K$-theoretic Donaldson
invariants determine its Hilbert polynomial.  

\begin{NB}
The following paragraph is editted. I emphasize the difference between
$H^0$ and the Euler characteristic. H.N. Nov.27.
\end{NB}
The $K$-theoretic Donaldson invariant is a natural $2$-dimensional
analogue of the dimension of the space of conformal blocks (nonabelian
theta functions). Another, closely related, analogue is the genuine
space of sections of a determinant line bundle, rather than the
alternating sum of cohomology groups. Its conjectural formula appeared
as four dimensional Verlinde formula in the physics literature
\cite{LMNS},\cite{LNS}.
(It is given as the space of sections, but it is not clear to the
authors whether the physical approach actually yields the space of
sections, not Euler characteristic.)
A mathematical formulation was given in \cite{Nak}, where it was
called the space of conformal blocks in $4$D WZW-Theory.
In case the base manifold is the projective plane the
strange duality conjecture of Le~Potier (see e.g. \cite{Da2}) gives a
duality between the spaces of sections of determinant line bundles for
moduli spaces of sheaves of positive rank and their analogues on
moduli spaces of pure sheaves of rank $0$.  This conjecture has been
checked in some cases in \cite{Da1},\cite{Da2}. An analogue of this
conjecture has been proved for some moduli spaces of sheaves on
K3-surfaces in \cite{OG}.
In many examples, the determinant line bundle is ample, or at least
nef and big, so we have the vanishing of higher cohomology groups. In
those cases there are no difference between the spaces of sections and
Euler characteristics. However it is not clear whether we can control
the spaces of sections in general.

\begin{NB}
I add the reference \cite{HP} because it is mentioned in \cite{LMNS}.
But I fail to understand they consider the holomorphic Euler
characteristic.  
If I understand \cite{LMNS} correctly, they say that if you replace the line bundle $L$  by $tL$ and divide by $t$ to the power the dimension of the  moduli space, and then take the limit as $t\to \infty$ one gets the  formulas of \cite{HP}. However we know that this limit is just the usual Donaldson invariants, which \cite{HP} computes the wallcrossing for. I think the reference \cite{HP}
should be removed. 8.11.LG
\end{NB}

\begin{NB}
Yes. I agree with Lothar. I deleted. Nov.9, HN
\end{NB}

\begin{NB}
There is a `kind' of formula in \cite[(6.3)]{LMNS}.  But it is written
that this is ill-defined. An ill-defined formula is not a formula (in
a mathematical sense), but it may represent something..... How should
we refer them ?  Nov. 2, H.N.
\end{NB}

\begin{NB} I do not understand their formula. I think the illdefinedness refers to the  integral and they mean that one has to regularize the integral
(but they do not explain how). On the other hand such integrals are considered
in \cite{MW}, where they explain the regularization. 8.11.LG
\end{NB}

The paper is organized as follows. In Sect.~\ref{sec:background} we
collect background material on the holomorphic Euler characteristic
of the determinant line bundle and the $K$-theoretic Nekrasov partition
function. We also explain the partition function with $5D$
Chern-Simons terms (see \cite{IMS,Ta}), which naturally appears in our
approach. We also calculate the $K$-theoretic Donaldson invariants for
$K3$ surface (see \subsecref{subsec:K3}).%
\begin{NB}
A sentence added, H.N. Nov.27  
\end{NB}
In Sect.~\ref{sec:Hilb} we express the wallcrossing terms in
terms of the holomorphic Euler characteristic of some virtual vector
bundles on the Hilbert schemes $X_2^{[n]}$ of points on two copies of
$X$. In Sect.~\ref{sec:comparison} we take $X$ a smooth projective
toric surface and express the equivariant wallcrossing terms in terms
of the $K$-theoretic Nekrasov partition function. These two sections
are parallel to \cite[Sect's.~2,3]{GNY}. Then in
Sect.~\ref{sec:modular} we take the nonequivariant limit and give the
formula of wallcrossing terms in terms of modular forms
and elliptic functions. We use the
solution of the Nekrasov conjecture and its refinement.
In particular, we determine the Hilbert series of the determinant line
bundles on $M_H^{\P^2}(0,d)$ and $M_H^{\P^2}(H,d)$ for small $d$ in
\subsecref{subsec:P2}.%
\begin{NB}
 I add the \subsecref{subsec:P2}. Nov. 6, HN
\end{NB}
In Appendix~\ref{sec:SWcurve} we explain the Seiberg-Witten curve for
the the $5$-dimensional supersymmetric Yang-Mills theory. We prove
that the Seiberg-Witten prepotential defined via the period of the
curve satisfies the contact term equation, which was also satisfied by
the nonequivariant limit of the Nekrasov partition function
\cite{NY3}. This completes our proof of Nekrasov's conjecture started
in \cite{NY3}, as the solution of the contact term equation is unique.
\begin{NB}
`up to the perturbation term.' But we do not mention it. Oct. 27, H.N.
\end{NB}

This paper is dedicated to 
Friedrich Hirzebruch, one of the founders of $K$-theory. Among the  other  subjects of this paper related his work are the Hirzebruch-Riemann-Roch theorem, modular forms and elliptic functions. The first-named author particularly wants to thank him, his teacher, for all the things he learned from him.

\subsection*{Acknowledgement}
The project started in 2004 Jan.\ when the first-named author visited Kyoto
for a workshop organized by the second and third-named authors. They are
grateful to the Kyoto University for its hospitality. 
The second-named author thanks Yuji Tachikawa and Hiroaki Kanno for
their explanations of the partition function with $5D$ Chern-Simons terms.
Part of this paper was written while the second and third-named
authors were visiting the International Centre for Theoretical
Physics, and also while the first-named author was visiting the
Institut-Mittag-Leffler. We thank both institutes for the hospitality.

\begin{NB}
Hiraku suggested  that one writes as dedication something like "To 
Prof. Hirzebruch one of the founders of K-theory". For some reason
I find it a bit impersonal, but maybe I am wrong. If I was a single author I think I should dedicate it to my teacher F. Hirzebruch, but this is obviously not possible

I am wondering if it would be possible to do the following:
The dedication is something like 
"To Friedrich Hirzebruch on his 80-th birthday", but then 
one adds a couple of sentences to the introduction, why  we find it fitting to 
contribute in particular this paper. One has to find a very good and short
formulation otherwise it sounds bad.

A first attempt (it is maybe too strong):

We dedicate this paper to 
Friedrich Hirzebruch, one of the founders of K-theory. Among the  other  subjects of this paper related his work are the Hirzebruch-Riemann-Roch formula, modular forms and elliptic functions. The first-named particularly wants to thank him as his teacher and advisor, who introduced him to these subjects.
\end{NB}
\begin{NB} slightly edited 15.11 LG\end{NB}

\section{Background Material}\label{sec:background}
We will  work over $\C$. We usually consider homology and cohomology with
rational coefficients and for a variety $Y$ we will 
 write $H_i(Y)$, and $H^i(Y)$ for $H_i(Y,\Q)$ and
 $H^i(Y,\Q)$ respectively. If $Y$ 
is projective and $\alpha \in H^*(Y)$, we denote 
$\int_Y\alpha$ its evaluation on the fundamental cycle of $Y$.
If $Y$ carries an action of a torus $T$, $\alpha$ is a $T$-equivariant class, 
and $p:X\to pt$ is the projection to a point, we denote  $\int_Y\alpha:=p_*(\alpha)\in H^*_T(pt)$.

In this whole paper $X$ will be a nonsingular projective surface over $\C$. Later 
we will specialize $X$ to a smooth projective toric surface.
For a class $\alpha\in H^*(X)$, we denote $\<\alpha\>:=\int_X\alpha$.
If $X$ is a toric surface 
 we use the same notation for the equivariant
pushforward to a point.

Let $X$ be simply connected smooth projective surface with $p_g(X)=0$. Let $H$ be an ample divisor on $X$. 
We denote by $M_H^X(r,c_1,c_2)$  the moduli space of rank $r$ torsion-free $H$-semistable 
sheaves
(in the sense of Gieseker and Maruyama) with $c_1(E)=c_1$, $c_2(E)=c_2$. 
Let $M^X_H(r, c_1,c_2)_s$ be the open 
subset of stable sheaves.

\subsection{Determinant line bundles}\label{sec:detbun}
\begin{NB}
I changed `determinant bundle' to `determinant line bundle'. Nov. 9, HN.
\end{NB}
We briefly review the determinant line bundle on the moduli space
\cite{DN},\cite{LP}, for more details we refer to \cite[Chap.~8]{HL}.

For a Noetherian scheme $Y$ we denote by $K(Y)$ and $K^0(Y)$ the Grothendieck groups of coherent sheaves and locally free sheaves on $Y$ respectively.
Then $K^0(Y)$ is a commutative ring with $1=[\oo_Y]$, with the multiplication given
by the tensor product of locally free sheaves. If $Y$ is nonsingular and 
quasiprojective, then $K(Y)=K^0(Y)$. 
In particular we have $K(X) = K^0(X)$ for the smooth projective
surface $X$. We will identify $K^0(X)$ with $K(X)$ hereafter.
If we want to distinguish a sheaf $\F$ and its class in $K(Y)$, we
denote the latter by $[\F]$. But we may also write $\F$ for the class in
$K(Y).$ 
For a proper morphism $f\colon Y_1\to Y_2$ we have the pushforward
homomorphism 
\(
   f_!\colon K(Y_1)\to K(Y_2)
\)
defined by
\(
   f_!([\F]) = \sum_i (-1)^i [R^if_*\F].
\)
When $Y_2 = \mathrm{pt}$, this is the Euler characteristic of $\F$
under the identification of $K(\mathrm{pt}) \cong\Z$:
\(
   f_!([\F]) = \chi(Y_1,\F) = \sum_i (-1)^i \dim H^i(Y_1,\F).
\)
We also have a pushforward homomorphism $K^0(Y_1)\to K^0(Y_2)$ when
$f$ is a locally complete intersection morphism.
(See \cite[\S4.4]{BFM}.)%
\begin{NB}
In \cite{Fulton} the RR was formulated in the Chow groups. There is a
cycle map $A_*(X)\to H_*(X)$, but $A^*(X)\to H^*(X)$ was not
defined (see \cite[Ex.~19.3.12]{Fulton}). In \cite{BFM} the
definition of $A^*(X)$ is given as $\varinjlim A^*(Y)$, where $f\colon
X\to Y$ with $Y$ nonsingular. We have $\varinjlim A^*(Y)\to A^*(X)$
(\cite[Ex.~17.4.9]{Fulton}). Then there is a cycle map $\varinjlim A^*(Y)\to
H^*(X)$ (\cite[Ex.~19.3.12]{Fulton}). Anyway, the RR holds for either Chow
cohomology group, as we only need axiomatic properties of the
cohomology theory. H.N. Nov.30 
\end{NB}
For any morphism $f\colon Y_1\to Y_2$ we have the pullback
homomorphism
\(
   f^*\colon K^0(Y_2)\to K^0(Y_1)
\)
defined by
\(
  f^*[\F] = [f^*\F] 
\)
for a locally free sheaf $\F$ on $Y_2$.

On $K(X)$ we have a quadratic form $(u,v)\mapsto \chi(X,u\otimes v)
\equiv \chi(u\otimes v)$. (We denote $\chi(X,u\otimes v)$ by
$\chi(u\otimes v)$ for brevity hereafter.)
We say that $u,v\in K(X)$ are numerically equivalent if $u-v$ is in the radical of 
this quadratic form, and denote $K(X)_{{\rm num}}$ the set of numerical equivalence classes. 
Let $c\in K(X)_{{\rm num}}$.
Let $\E$ be a flat family of coherent sheaves of class $c$ on $X$ parametrized by a scheme $S$,
and let $p:X\times S\to S$, $q:X\times S\to X$ be the projections.
Define $\lambda_\E:K(X)\to \Pic(S)$ as the composition
\begin{equation}\label{eq:lambdaE}
\begin{diagram}[height=.8em,width=2em,scriptlabels] 
K(X)&\rTo^{q^*}& K^0(X\times S)&
\rTo{\otimes [\E]}&K^0(X\times S)&
\rTo{p_!} &K^0(S)&\rTo^{\det}&\Pic(S),
\end{diagram}\end{equation}
(see also \cite[(2.1.10), (2.1.11)]{HL}).
The following elementary facts are important for working with these line bundles:
\begin{enumerate}
\item $\lambda_\E$ is a homomorphism, i.e. $\lambda_\E(v_1+v_2)=\lambda_\E(v_1)\otimes \lambda_{\E}(v_2)$.
\item If $\mu\in \Pic(S)$ is a line bundle, then $\lambda_{\E\otimes \mu}(v)=
\lambda_{\E}(v)\otimes \mu^{\chi(c\otimes v)}$.
\item $\lambda_\E$ is compatible with base change: if 
$\phi:S'\to S$ is a morphism, then $\lambda_{\phi^*\E}(v)=\phi^*\lambda_{\E}(v)$.
\end{enumerate}

Let $H$ be a very ample divisor on $X$. 
For a class $c\in K(X)_{{\rm num}}$ we denote by $K_c:=c^\perp=\big\{v\in K(X)\bigm|
\chi(v\otimes c)=0\big\}$.%
\begin{NB} Changed: we did not introduce the notation $(v,c)$
27.11.LG\end{NB}
We denote by $K_{c,H}:=c^\perp\cap\{1,h,h^2\}^{\perp\perp}$, where $h=[\oo_H]$.
Now let $c\in K(X)_{{\rm num}}$ be the class of an element in $M_H^X(r,c_1,c_2)$.
There are  homomorphisms
$\lambda\colon K_c\to \Pic(M_H^X(r,c_1,c_2)_s)$,  and $\lambda\colon K_{c,H}\to \Pic(M_H^X(r,c_1,c_2))$,
 such that $\lambda$ commute with the  inclusions $K_{c,H}\subset K_c$ and $\Pic(M_H^X(r,c_1,c_2))\subset \Pic(M_H^X(r,c_1,c_2)_s)$. 
Note that $\lambda_\E(v)$ is independent of the choice of the universal
family $\E$ for $v\in K_c$ by the property (2) above, and in fact, we
do not need the existence of the universal sheaf to define the map $\lambda$.
We call $H$ {\it general} with respect to 
$(r,c_1,c_2)$ if all the strictly semistable sheaves in $M_H^X(r,c_1,c_2)$
are strictly semistable with respect to all ample divisors on $X$ in a neighbourhood of $H$ (the ample cone has the topology induced from the 
Euclidean topology on $H^2(X,\R)$).%
\begin{NB} 
changed according to Kota's message 24.11.06 LG\end{NB} 
In this case any strictly semistable sheaves in $M_H^X(r,c_1,c_2)$ is
of type $1$ in the sense of \cite[0.3]{Dre}. Then the stabilizer
subgroup $\operatorname{Aut}F$ (appeared in the proof of
\cite[Theorem~8.1.5]{HL}) acts trivially on the fiber of the
determinant line bundle (on the open subscheme of the quot-scheme).
Therefore $\lambda\colon K_{c,H}\to \Pic(M_H^X(r,c_1,c_2))$ can be
extended to $K_c$.

If $\E$ is a flat family of semistable sheaves  of rank $r$ and with Chern classes 
$c_1,c_2$ on $X$ parametrized by $S$,  then  we have 
$\phi_{\E}^*(\lambda(v))=\lambda_{\E}(v)$ for all $v\in K_{c,H}$ for $\phi_\E^*:\Pic(M^X_H(r,c_1,c_2))\to \Pic(S)$ the pullback  by the classifying morphism. 
If $H$ is general with respect to $(r,c_1,c_2)$, the same statement holds
with $K_{c,H}$ replaced by $K_c$. 
If $\E$ is a flat family of stable sheaves, the same statement holds with  $K_{c,H}$, $M^X_H(r,c_1,c_2)$ replaced by  $K_{c}$, $M^X_H(r,c_1,c_2)_s$.

\subsection{$K$-theoretic Donaldson invariants}\label{backDong}
We write $M^X_H(c_1,d)$ for $M^X_H(2,c_1,c_2)$ with $d=4c_2-c_1^2-3$.
Let $v\in K_c$, where $c$ is the class of a coherent rank $2$ sheaf 
with Chern classes $c_1,c_2$. 
Assume that $H$ is general with respect to $(2,c_1,c_2)$. 
The {\em $K$-theoretic Donaldson invariant\/} of $X$ with respect to 
$v,c_1,c_2,H$ is the holomorphic Euler characteristic
$\chi(M^X_H(c_1,d),\lambda(v))$
of the line bundle $\lambda(v)$.

 
 \begin{Notation}\label{ui} We introduce the following notation that we will
often use in the paper. For $i\ge 0$, we put  $v^{(i)}:=[\ch(v)e^{c_1/2}\Todd(X)]_i$.
Thus $v^{(0)}=\rk(v)$, $v^{(1)}=c_1(v)+\frac{\rk(v)}{2}(c_1-K_X)$ and $v^{(2)}$ could be interpreted as $\chi(v\otimes \oo(c_1/2))$.
\end{Notation}

 

By the Riemann-Roch Theorem it follows that 
\begin{equation}
\label{eq:chuc} \chi(v\otimes c)=2v^{(2)}-\rk(v)(c_2-\frac{c_1^2}{4}),
\end{equation}
in particular we see that the condition $v\in K_c$ is independent of $d$ if  $\rk(v)=0$.


An important special case is the following:
Let $L$ be a line bundle on $X$. Assume that $\<c_1(L),c_1\>$ is even
(otherwise replace $L$ by $L^{\otimes 2}$).
Then for $c$ the class of a rank $2$ coherent sheaf with Chern classes 
$c_1,c_2$, we put 
\begin{equation}\label{eq:uL} v(L):=-(1-L^{-1})-\<\frac{c_1(L)}{2},c_1(L)+K_X+c_1\>[\oo_x]\in K_c.\end{equation}
Note that $v(L)$ is independent of $c_2$. 
The condition that $\<c_1(L),c_1\>$ is even implies that $v(L)\in K(X)$.
Assume that $H$ is general with respect to $(2,c_1,c_2)$. Then we denote
$\mu(L):=
\lambda(v(L))\in \Pic(M^X_H(c_1,d))$.
The {\it $K$-theoretic Donaldson invariant\/} of $X$, with respect to $L,c_1,d,H$ is 
$\chi(M^X_H(c_1,d),\mathcal O(\mu(L)))$.
The generating function is 
\begin{equation}\label{eq:Kdon}
\begin{split}
\chi_{c_1}^H(L;\Lambda)&:=\sum_{d\ge 0}\Lambda^d
\chi(M^X_H(c_1,d),\mathcal O(\mu(L))).
\end{split}
\end{equation}
If $\E$ is a flat family of coherent sheaves parametrized by $S$, we
have $c_1(\mu(L)) = (c_2(\E)-\frac14c_1(\E)^2)/PD(c_1(L))\in H^2(S)$
by the Riemann-Roch for a smooth morphism
(\cite[\S4.3]{BFM}).%
\begin{NB}
I have changed the reference. H.N. Nov. 30.   
\end{NB}
It extends to a class in $H^2(M^X_H(c_1,d))$ by the same argument for
$\mu(L)$.%
\begin{NB}
Editted by H.N. Nov.27. We do not need the smoothness of the moduli
space. The following was original:  

If we assume that $M^X_H(c_1,d)$ is nonsingular of the expected dimension,
and $\E$ is a universal sheaf on $X\times M^X_H(c_1,d)$, then  the Grothendieck-Riemann-Roch theorem implies that 
$c_1(\mu(L))=(c_2(\E)-\frac14c_1(\E)^2)/PD(c_1(L))$.
\end{NB}
This coincides with the definition of $\mu(c_1(L))$ appearing in the usual
Donaldson invariant. This is the reason why we denote the line bundle
by $\mu(L)$.
Thus it follows from the definitions and the singular Riemann-Roch
theorem \cite{BFM}%
\begin{NB}
I changed the reference. H.N. Nov.30  
\end{NB}
that $\chi(M^X_H(c_1,d),\mathcal
O(\mu(nL)))$ is a polynomial of degree $d$ in $n$, whose leading term
is the algebraic geometric version of the Donaldson invariants
$n^d\Phi^H_{c_1}( c_1(L)^d/d!)$ (in the notations of \cite{GNY}) when
$M^X_H(c_1,d)$ is of the expected dimension.
\begin{NB}
Editted by H.N. Nov.27. The following was original:  

Thus it follows from the definitions and the Hirzebruch-Riemann-Roch
theorem \cite{Hi} that $\chi(M^X_H(c_1,d),\mathcal O(\mu(nL)))$ is a
polynomial of degree $d$ in $n$, whose leading term is the algebraic
geometric version of the Donaldson invariants $n^d\Phi^H_{c_1}(
c_1(L)^d/d!)$ (in the notations of \cite{GNY}).
\end{NB}

The above argument also implies that the invariant
$\chi(M^X_H(c_1,d),\lambda(v))$ depends only on $\ch(v)\in H^*(X)$.
Therefore the invariant is well-defined on
$K(X)_{\rm hom} := K(X)/\!\!\sim$ where $v\sim v'$ if and only if $\ch(v)
= \ch(v')$.

\begin{NB}
  The above paragraph is added by H.N. Nov.27. It seems that it is not
  known whether $K(X)_{\rm hom} = K(X)_{\rm num}$ or not. See
  \cite[\S19.3]{Fulton}.
\end{NB}

\subsection{A digression on the definition of the invariants}\label{subsec:digress}

The definition of the $K$-theoretic Donaldson invariants above is only ad hoc  and will in general need to be modified, so that the invariants have
good properties and so that they might be related to 
gauge-theoretical invariants. 

We expect that for general $X$, when the moduli space $M_H^X(c_1,d)$
does not have the expected dimension, one needs to use a 
virtual structure sheaf (see \cite{Lee}) in the definition. If
$M_H^X(c_1,d)$ consists only of stable sheaves, 
the perfect obstruction theory was constructed in \cite[Th.~3.30]{Thomas}.
Then we just need to replace $\chi(M_H^X(c_1,d),\lambda(v))$ by
$\chi(M_H^X(c_1,d),\oo^{{\rm virt}}\otimes \lambda(v))$. If
$M_H^X(c_1,d)$ has the expected dimension, then by
\cite[Prop.~2]{Lee}, the virtual structure sheaf is just the usual
structure sheaf, and this definition reduces to our previous
definition.
When $M_H^X(c_1,d)$ may contain a strictly semistable sheaf, we need
to construct a perfect obstruction theory on another moduli space with
additional structures and prove that it is independent of the
additional structure as in \cite{Moc}, or use the blowup formula as in
the definition of the usual Donaldson invariants (see
\cite[\S1.1]{GNY}). See \subsecref{subsec:blowup} below for the first
step in this approach.

Let us examine the possibility to extend our definition of invariants to a
$C^\infty$ $4$-manifold $X$. To avoid a technical difficulty, we first
assume the moduli space $M_H^X(c_1,d)$ is smooth. Our definition
depends on the complex structure of $X$, and if we have a gauge
theoretic definition, it should be independent of the complex
structure, and the definition must be modified.
Our guess is to consider the index of a Dirac operator instead of the
holomorphic Euler characteristic. If $X$ is spin, then we have a
square root $K_X^{1/2}$ of $K_X$, and then $\mu(K_X)$ is a line
bundle. It is known that this is isomorphic to half of the canonical
bundle of $M^X_H(c_1,d)$ when it is smooth (see
e.g.~\cite[\S8.3]{HL}).%
\begin{NB}
Reference corrected by H.N. Nov. 27
\end{NB}
Therefore
\(
    \chi(M^X_H(c_1,d), \oo(\mu(K_X))
\) 
is equal to the index of the Dirac operator. In this special case, our
main result \corref{cor:unim} is simplified as $v^{(1)} = -K_X$. In
particular, the answer is independent of the complex structure except
the term $\sqrt{-1}^{\<\xi,K_X\>}$ which corresponds to the
orientation of the moduli space.

More generally the complex structure on $X$ and a line bundle $L$ on
$X$ induces the $Spin^c$-structure $W^+ = (\Wedge^{0,0}\oplus
\Wedge^{0,2})\otimes L$, $W^- = \Wedge^{0,1}\otimes L$ on $X$ with the
characteristic line bundle $\det W^+ = \det W^- = -K_X + 2L$.
We conjecture that it induces a $Spin^c$-structure on the moduli
space. The recipe should be somewhat similar to the definition of the
orientation of the moduli space induced from the homological
orientation on $H^0(X)\oplus H^1(X)^*\oplus H^2_+(X)$, but we do not
know how to define it in general (even on the nonsingular part of the
moduli space). However in our situation, an obvious candidate for the index is
\(
    \chi(M^X_H(c_1,d), \oo(\mu(2L)).
\)%
\begin{NB}
When $L = K_X^{1/2}$, the $Spin^c$ structure is the $Spin$ structure
considered above. Therefore $2L = K_X$ is compatible with the above.
Nov.5 HN 
\end{NB}
This means that the $Spin^c$-structure is the one given by the complex
structure twisted by the line bundle $\mu(2L)$.  The answer given in
\corref{cor:unim} is written in terms of $v^{(1)} + K_X = K_X - 2L$.
As this is the negative of the characteristic line bundle of the
$Spin^c$ structure, the candidate seems reasonable.

Now we come to discuss more technical points.
For a $C^\infty$ $4$-manifold $X$, we do not have the Gieseker-Maruyama
compactification $M^X_H(c_1,d)$ and we need to use the Uhlenbeck
compactification $N^X_H(c_1,d)$ of the moduli space of instantons
instead. We also need to use its topological $K$-homology group
$K^{\operatorname{top}}(N^X_H(c_1,d))$. When $X$ is a projective
surface, we have a homomorphism $\pi_*\colon
K^{\operatorname{top}}(M^X_H(c_1,d))\to\linebreak[3]
K^{\operatorname{top}}(N^X_H(c_1,d))$ given by $\pi\colon
M^X_H(c_1,d)\to N^X_H(c_1,d)$ and we can pushforward the virtual
structure sheaf on $M^X_H(c_1,d)$ to $N^X_H(c_1,d)$. And it can be
shown that the line bundle $\mu(L)$ is a pull-back of a line bundle
from $N^X_H(c_1,d)$ under some conditions.  (See
\subsecref{subsec:blowup}. And this assertion, at least for a
topological line bundle, is well-known in the gauge theory context.)
Therefore the invariants in \eqref{eq:Kdon} can be defined in terms of
$N^X_H(c_1,d)$ and the framework of the topological $K$-group.
However it is not clear, at least to the authors, how to define the
$K$-theoretic fundamental class $[\oo_{N^X_H(c_1,d)}]\in
K^{\operatorname{top}}(N^X_H(c_1,d))$ for an arbitrary $C^\infty$
$4$-manifold $X$ even under the assumption that $N^X_H(c_1,d)$ is of
expected dimension.

We have considered the determinant line bundle $\mu(L)$ above. This is
the case $\rk(v) = 0$. When $v\in K_c$ is suitably chosen (see
\cite[\S8.1]{HL}) with $\rk(v)=2$, the determinant line bundle
$\lambda(v)$ is ample on $M^X_H(c_1,d)$ and does {\it not\/} come from
$N^X_H(c_1,d)$. This observation seems to suggest that the invariant
can be defined {\it only\/} for a restricted class $v$ on a $C^\infty$
$4$-manifold $X$. We have discussed $\rk(v) = 0$ is sufficient for the
existence of the line bundle $\lambda(v)$ on $N^X_H(c_1,d)$ above, but
we do not know whether this is necessary.

Also we do not give the definition of the analog of $\mu(p) \in
H^4(M^X_H(c_1,d))$ where $p$ is the point class of $H_0(X)$. It may be
defined as
\begin{equation*}
   \chi(M^X_H(c_1,d),p_! (q^* v\otimes \E)),
\end{equation*}
but it is not independent of the choice of the universal bundle $\E$
in general, and may not be defined when we do not have a universal
bundle. A possible candidate, which can be defined for any $v\in
K(X)$, is given by replacing $\E$ by $\E\otimes \E^\vee$, where
${}^\vee$ is the involution on $K^0(X\times M^X_H(c_1,d))$ defined by
taking the dual of a vector bundle.
Or more generally, if we have a representation $\rho \colon
PGL(2,\C)\to GL(V)$, we may consider
\(
  \chi(M^X_H(c_1,d), p_! (q^* v\otimes\rho([\E]))),
\)
or applying $\rho$ (with an appropriate change of $PGL(2,\C)$) after
the pushforward $p_!$. 
But we do not study these `higher' invariants, and stick to
our $\chi(M^X_H(c_1,d), \lambda(v))$, which we believe most basic.

\begin{NB}
Changed according to Hiraku's suggestion 4.10.06 LG.

Actually I think that this change makes the remark in some sense less correct, as I think  this is not a correct description of the 
situation: as far as I can see Mochizuki does not directly put a virtual 
fundamental class of the moduli space of sheaves: the reason is that strictly
speaking the usual moduli space with its natural stack structure is an Artin
stack (even when stable=semistable) and one does not know how to put an obstruction theory  on that. This is reflected in the fact that the natural 
obstruction theory is given by $Hom(E,E)$ as infinitesimal automorphisms,
$Ext^1(E,E)$ as tangent and $Ext^2(E,E)$ as obstructions. 

Instead he looks at the moduli space of oriented sheaves (i.e. these are 
sheaves together with an isomorphism of the determinant to a given line bundle
if we assume $X$ is simply connected). For these if stable=semistable we get 
a Deligne Mumford stack with a perfect obstruction theory (I think).
He always studies such sheaves with nowhere vanishing $L$-sections, 
and does the obstruction theory for those (I think under his assumptions they are always Deligne Mumford stacks and (independent on whether for the moduli 
of sheaves stable=semistable). 
I must however admit that I have not been able to understand very well what 
Mochizuki is doing. 

Thus if I am not mistaken, there are two differences to what we write above.
\begin{description}
\item[(1)] He never puts a an obstruction theory directly on the moduli of sheaves, but always on a Deligne Mumford stack obtained by adding some data.
\item[(2)] It is not necessary for him in order to do this to assume that on
the  moduli of sheaves stable=semistable.
\end{description}

I think instead that if one is willing to assume that stable=semistable,
then one can define the virtual fundamental class in a simpler way,
without using the complicated constructions of Mochizuki, directly on the 
moduli space. One just has to take both as tangent and obstruction the 
$Ext^1_0(E,E)$, $Ext^2_0(E,E)$, the tracefree endomorphisms.
If I am not mistaken this is what Donaldson Thomas and MNOP do in 
order to define the Donaldson-Thomas invariants. 

The point is that if one wants to study the wallcrossing, this is not good enough 
because in the wallcrossing one comes across the case where no longer
stable=semistable.

\end{NB}
\begin{NB}
I have edited my comment according to Lothar's comment on Oct. 4. I
think that there is a hope to prove a blowup formula using J.~Li's
technique of the degeneration of the moduli space, though I need to
understand his paper with Gieseker and his papers on relative GW
invariants.

I also add a comment on invariants for smooth $4$-manifolds, and
`higher' invariants.
Oct. 18, H.N.
\end{NB}

\subsection{Blowup formula and the invariants for moduli spaces with
  strictly semistable sheaves}\label{subsec:blowup}
\begin{NB} 
I include the note of Kota, which gives the blowup formula for rational surfaces.
I for the moment only include the part which is relevant for this and I edit it
a bit to make it compatible with our notations. 
Actually one might wonder whether one should include the whole note,
because it shows also a statement about universal sheaves, that for the 
moment we just claim without proof in the paper. 6.10.LG
\end{NB}
\begin{NB}
I think that this subsection is worthwhile including (I do not consider
about the remaining part of Kota's note), but I want to make it a new
subsection so that a reader in hurry can skip it.
Oct. 17, H.N.
\end{NB}
As mentioned in the previous subsection, we give a proposal of the
definition of invariants when moduli spaces may contain strictly
semistable sheaves by using a blowup formula. We assume a kind of
smoothness of moduli spaces on the blowup. This allows us to avoid the
virtual structure sheaf. However the smoothness assumption is used much
more essentially as we use Kawamata-Viehweg vanishing theorem. If we
could prove the same vanishing theorem under the assumption that the
moduli space is of expected dimension, we could use the blowup formula
as the definition of the invariant, as is done in the context of usual
Donaldson invariants (see e.g., \cite{GNY}).
Then the invariant is integral, in contrast with the ordinary Donaldson
invariants in which we must divide by powers of $2$.%
\begin{NB}
I add a comment on the integrality. Nov.30, HN.  
\end{NB}
Moreover we also prove that the pushforward of the structure sheaf of
the Gieseker-Maruyama compactification is equal to the the structure
sheaf of the Uhlenbeck compactification. This seems an evidence of our
belief that the $K$-theoretic Donaldson invariant has a gauge
theoretic definition. The material in this subsection is technical,
so a reader in hurry can just read the statement of \corref{cor:blow}
and skip the rest.

\begin{NB}
A paragraph added. Nov. 9, HN.  
\end{NB}

Let $(X,H)$ be a polarized rational surface.  Let $\widehat X$ be the
blowup of $X$ in a point and $C$ the exceptional divisor. In the
following we always denote a class in $H^*(X,\Z)$ and its pullback by
the same letter. Write $c:=(2,c_1,c_2)$, and $M_H(c):=M_H^X(c)$.
Let $Q$ be an open subset of a suitable quot-scheme
such that $M_H(c)=Q/GL(N)$.

Let $N_H(c)$ be the Uhlenbeck compactification of
the moduli space of slope stable vector bundles on $X$.
The line bundle $\mu(2D)$ is a pull-back of a line bundle from $N_H(c)$
if $D\in \bigcap_{\xi:\langle H,\xi\rangle=0} \xi^\perp$.
In fact, the stability is the same for $H$ and 
$H+\varepsilon D$ with $D \in \cap_{\xi:\langle H,\xi\rangle=0}\xi^{\perp}$
for a sufficiently small $\varepsilon$.
Then $\mu(H+\varepsilon D)$ is nef and big and gives a map
to the Uhlenbeck compactification. In particular,
$\mu(2 D)$ is the pull-back of a line bundle on the Uhlenbeck
compactification, which we denote by the same symbol.
We further assume $H$ is general with respect to $c$, then 
we have $\{ \xi\mid \langle H,\xi\rangle=0\} = \{ 0\}$. Therefore
$\mu(2D)$ is the pull-back of a line
bundle on $N_H(c)$ for any $D$.

\begin{NB}
I change several $(\ ,\ )$ to  $\langle\ ,\ \rangle$. (Also in below.)
Nov. 2, HN
\end{NB}

\begin{NB}
The above is modified according to Kota's suggestion. But I do not
understand Kota's message well, so please check it
is OK. The followings are Kota's message on Oct.\ 19 and the
original. Oct. 26, H.N.; Kota checked it in the message on Oct. 27.
Therefore this NB will disappear eventually. Oct. 30, H.N.
\begin{verbatim}
(2)
sect. 2.3

(i)
Assume that $H$ is general with... or X is a Del Pezzo...

It may be better say that

``$\mu(2D)$ is the pull-back of a line bundle on N_H(c),
if D \in \cap_{\xi;(H,\xi)=0}\xi^{\perp}."
And
``In particular, if $H$ is general with respect to $c$, then
$\mu(2D)$ is the pull-back of a line bundle on N_H(c)."

...
Indeed ``if H is general," then stability
is the same for H and H+\varepsilon D...

will be

Indeed stability
is the same for H and H+\varepsilon D...,

since the generic means that \{\xi|(H,\xi)=0\}=\{0\}.
\end{verbatim}

Assume now that $H$ is general with respect to $c$ or
$X$ is a Del Pezzo surface and $H=-K_X$.
Then $\mu(2 K_X)$ 
\verb|\begin{NB}| For security I multiply with 2, 5.10 LG\verb|end{NB}|
defines a line bundle on $N_H(c)$.
Indeed if $H$ is general, then
stability is the same for $H$ and 
$H+\varepsilon D$ with $D \in \cap_{\xi;(H,\xi)=0}\xi^{\perp}$
for a sufficiently small $\varepsilon$.
Then $\mu(H+\varepsilon D)$ is nef and big and gives a map
to the Uhlenbeck compactification. In particular,
$\mu(2 D)$ is the pull-back of a line bundle on the Uhlenbeck
compactification, which we denote by the same symbol.
\verb+\begin{NB}+ again multiplied by $2$, 5.10. LG\verb+end{NB}+
\end{NB}

We shall study the singularities of $M_H(c)$ and $N_H(c)$.
\begin{Lemma}[\cite{B}]
Assume that $Q$ is smooth \textup(e.g.\ $\langle -K_X,H\rangle>0$\textup).
Then $M_H(c)=Q/GL(N)$ is normal and has only rational singularities.
\end{Lemma}
We next consider the singularities of
$N_H(c)$.
Replacing $N_H(c)$ by its normalization, we may assume that
$N_H(c)$ is normal.
\begin{Lemma}
Assume that $Q$ is smooth \textup(e.g.\ $\langle -K_X,H\rangle>0$\textup).

\textup{(1)} 
Then $N_H(c)$ has only rational singularities.

\textup{(2)}
${\bf R}\pi_*({\cal O}_{M_H(c)})={\cal O}_{N_H(c)}$.
\end{Lemma}

\begin{proof}
(1)
We first assume that $c$ is primitive. Then there is a resolution
of $\pi\colon M^{\alpha}_H(c) \to N_H(c)$,
where $M^{\alpha}_H(c)$ is the moduli space of  $\alpha$-twisted semi-stable sheaves for suitable $\alpha$.
Then by the Grauert-Riemenschneider vanishing theorem,
${\bf R}\pi_*(K_{M^{\alpha}_H(c)})=\pi_*(K_{M_H^{\alpha}(c)})$.
Since $K_{M_H^{\alpha}(c)} \cong \mu(2K_X)$
comes from $N_H(c)$ and $N_H(c)$ is normal,
${\bf R}\pi_*({\cal O}_{M_H^{\alpha}(c)})=
\pi_*({\cal O}_{M_H^{\alpha}(c)})={\cal O}_{N_H(c)}$.
Thus $N_H(c)$ has only rational singularities.

We next treat $M_H(c)$ with $c=(2,0,2n)$.
We set $\widehat c:=(2,C,2n)$ and $\widehat{M}_H(\widehat{c}):=
{M}^{\widehat X}_{H-\varepsilon C}(\widehat c)$.%
\begin{NB} 
I changed this, with the current notations $\widehat n=2n$ 5.10.06 LG \end{NB}
Then there is a surjective 
morphism $\widehat{\pi}\colon \widehat{M}_H(\widehat{c})\to N_H(c)$
which is
generically a ${\Bbb P}^1$-bundle.
Since $-\mu(C)$ is $\widehat{\pi}$-nef and big,
the Kawamata-Viehweg vanishing theorem implies that
$R^i \widehat{\pi}_*(K_{\widehat{M}_H(\widehat{c})}(-2\mu(C)))=0, i>0$.
By our assumption,
$\mu(K_X)$ comes from $N_H(c)$.
This implies that
$$
c_1(K_{\widehat{M}_H(\widehat{c})})=2\mu(c_1(K_{\widehat{X}}))=2\mu(c_1(K_X))+2\mu(C)
\equiv 2\mu(C) \mod \widehat{\pi}^* H^2(N_H(c),{\Bbb Q}).
$$
Hence $R^i \widehat{\pi}_*({\cal O}_{\widehat{M}_H(\widehat{c})})=0,i>0$.
Thus  ${\bf R} \widehat{\pi}_*({\cal O}_{\widehat{M}_H(\widehat{c})})
={\cal O}_{N_H(c)}$.
Then $N_H(c)$ has rational singularities by
\cite[Thm. 1]{Kov}.

(2) It is sufficient to prove the following: For a proper birational
map $f\colon Y \to Z$ of normal varieties $Y,Z$ with only rational
singularities, ${\bf R}f_*({\cal O}_Y)={\cal O}_Z$.

Proof of the claim:
Let $g\colon Y' \to Y$ be a resolution of the singularities.
Since $Y$ has only rational singularities,
$R^i g_*({\cal O}_{Y'})=0$, $i>0$.
Then $R^i f_*({\cal O}_{Y})=R^i (f \circ g)_*({\cal O}_{Y'})=0, i>0$.
Hence we get our claim.

\end{proof}

By the proof, we also get the following. 
\begin{Corollary}\label{cor:blow}
Let $\widehat{M}_H(\widehat{c})$ be the moduli space of stable
sheaves on $\widehat{X}$ such that $\widehat{c}=(2,c_1+kC,c_2)$ with
$k=0,1$.%
\begin{NB}
Corrected according to Kota's e-mail. Oct.21, H.N. 
\end{NB}
Then
$$
{\bf R} \widehat{\pi}_*({\cal O}_{\widehat{M}_H(\widehat{c})})=
{\cal O}_{N_H(c)}
={\bf R} \pi_*({\cal O}_{M_H(c)}).
$$
In particular,
$$
\chi(\widehat{M}_H(\widehat{c}),\mu(D))
=\chi(N_H(c),\mu(D))
=\chi(M_H(c),\mu(D))
$$
for any line bundle  $D$ on $X$ such that $\<D, c_1\>$ is even and 
$\<D,\xi\>=0$ for $\xi$ any class of 
type $(c_1,4c_1-c_1^2-3)$ on $\widehat X$ with $\<H,\xi\>=0$. 
\end{Corollary}
\begin{NB}
The condition $\<D, c_1\>$ is even is added.
Nov. 2. HN.
\end{NB}

\begin{Remark}\label{rem:canonical}
Since $M_H(c)$ is normal,
the dualizing sheaf $\omega_{M_H(c)}$ is reflexive.
If $H$ is a general polarization, then
${\cal O}_{M_H(c)}(2\mu(K_X))$ is a line bundle on $M$
which coincides with the dualizing sheaf
on the locus of stable sheaves $M_H(c)^s$.
If $\dim (M_H(c) \setminus M_H(c)^s) \leq \dim M_H(c)-2$,
then $\omega_{M_H(c)}={\cal O}_{M_H(c)}(2\mu(K_X))$.
\end{Remark}

\begin{NB}
I put here an older  version of the same section  where  introduce $u_d$
 Please tell
me what you think better. I think that the  version without $u_d$ is better.

\subsection{$K$-theoretic Donaldson invariants}\label{backDon}
We write $M^X_H(c_1,d)$ for $M^X_H(2,c_1,c_2)$ with $d=4c_2-c_1^2-3$.
Let $v\in K_c$, where $c$ is the class of a coherent rank $2$ sheaf 
with Chern classes $c_1,c_2$. 
Assume that $H$ is general with respect to $(2,c_1,c_2)$. 
The {\em $K$-theoretic Donaldson invariant} of $X$ with respect to 
$u,c_1,c_2,H$ is the holomorphic Euler characteristic
$\chi(M^X_H(c_1,c_2),\lambda(v))$
of the line bundle $\lambda(v)$.

Let $v\in K(X)$. 
We want to study generating functions in $d$ for the numbers 
 $\chi(M^X_H(c_1,d),\lambda(v))$.
Fix $c_2$, and put $d:=4c_2-c_1^2-3$. $c\in K(X)_{{\rm num}}$ be the  class of a 
coherent sheaf $E$ of rank $2$ with Chern classes $c_1,c_2$. 
 In general  the condition that $\chi(v\otimes c)=0$ will depend on $d$. 
 
 \begin{Notation}\label{ui} We introduce the following notation that we will
often use in the paper. For $i\ge 0$, we put  $v^{(i)}:=[\ch(v)e^{c_1/2}\Todd(X)]_i$.
Thus $v^{(0)}=\rk(v)$,  $v^{(1)}=c_1(v)+\frac{\rk(v)}{2}(c_1-K_X)$ and $v^{(2)}$ could be interpreted as $\chi(v\otimes \oo(c_1/2))$.
\end{Notation}

We assume the condition  
\begin{equation}\label{even}  \hbox{$\rk(v)$ is even and $\rk(v)\frac{c_1^2}{4}+2v^{(2)}$
 is even.} \end{equation}
 
 Then we put 
\begin{equation}\label{eq:ud}
v_d:=v-\frac{\chi(v\otimes c)}{2}[\oo_x]\in K_c.
\end{equation}

By the Riemann-Roch Theorem it follows that 
\begin{equation}
\label{eq:chuc} \chi(v\otimes c)=2v^{(2)}-\rk(v)(c_2-\frac{c_1^2}{4})
\end{equation}
So we see that the condition \eqref{even} implies that $v_d\in K_c$. 
In case $\rk(v)=0$ and $v^{(2)}\in \Z$, we see that $v'=v-v^{(2)}[\oo_x]\in K_c$ for all $d$. 

Let $v\in K(X)$ satisfying \eqref{even}. Then we define 
\begin{equation}\label{eq:Kdonu}
\begin{split}
\chi_{c_1}^H(v;\Lambda)&:=\sum_{d\ge 0}\Lambda^d  \chi(M^X_H(c_1,d),
\lambda(v_d)).
\end{split}
\end{equation}

The most important special case is the following:
Let $L$ be a line bundle on $X$. Assume that $\<c_1(L),c_1\>$ is even.
Then for $c$ the class of a rank $2$ coherent sheaf with Chern classes 
$c_1,c_2$, we put 
\begin{equation}\label{eq:uL} u(L):=-(1-L^{-1})-\<\frac{c_1(L)}{2},c_1(L)+K_X+c_1\>[\oo_x]\in K_c.\end{equation}
Note that $u(L)$ is independent of $c_2$. 
The condition that $\<c_1(L),c_1\>$ is even implies that $u(L)\in K(X)$.
Assume that $H$ is general with respect to $(2,c_1,c_2)$. Then we denote
$\mu(L):=
\lambda(v(L))\in \Pic(M^X_H(c_1,d))$, and call 
$\mu(L)$ the {\it Donaldson line bundle} associated to $L$.
The {\em $K$-theoretic Donaldson invariant} of $X$, with respect to $L,c_1,d,H$ is 
$\chi(M^X_H(c_1,d),\mu(L))$.
The generating function is 
\begin{equation}\label{eq:Kdon}
\begin{split}
\chi_{c_1}^H(L;\Lambda)&:=\sum_{d\ge 0}\Lambda^d  \chi(M^X_H(c_1,d),\mu(L)).
\end{split}
\end{equation}
Assume that $M^X_H(c_1,d)$ is nonsingular of the expected dimension,
and that $\E$ is a universal sheaf on $X\times M^X_H(c_1,d)$. 
By the Grothendieck-Riemann-Roch theorem we see that 
$c_1(\mu(L))=(c_2(\E)-c_1(\E)^2/4)/PD(c_1(L))$. 
Let $\Phi^H_{c_1}$ be the algebraic geometric version of the Donaldson invariants
(in the notations of \cite{GNY}), then it follows from the definitions and the 
Hirzebruch-Riemann-Roch theorem \cite{Hi} that 
$\chi(M^X_H(c_1,d),\mu(L))$ is a polynomial of degree $d$ in $n$, whose leading term is 
$\Phi^H_{c_1}( (nL)^d/d!)$.

\end{NB}

\subsection{$K$-theoretic invariant for $K3$ surfaces and strange
  duality}\label{subsec:K3}%
\begin{NB} This is the modified version of the section, I also change its 
introduction to fit with the change of the contents.
14.11.LG\end{NB}
Let $X$ be a projective $K3$ surface. In this subsection we calculate
the $K$-theoretic invariants for $X$ as examples.  We also give a
formula for the $K$-theoretic invariants of rank $1$ sheaves on abelian surfaces.%
\begin{NB} I remove all reference to strange duality from this paragraph. 
 I only talk about it in a remark after the result.
 24.11. LG\end{NB}
\begin{NB}
A paragraph added. Nov.9, HN
\end{NB}

For any projective algebraic surface $Y$ and $c\in K(Y)_{\rm hom}$ we denote by 
$M_H^Y(c)$ the moduli space of $H$-stable 
sheaves $E$ on $Y$ with $\ch(E)=\ch(c)$.
This is just a%
\begin{NB} changed 27.11. LG\end{NB}
change of notation, but is convenient to see the
strange duality.
We also define the discriminant by $\Delta(c)=2\rk(c)c_2(c)-(\rk(c)-1)c_1(c)^2$,
$\Delta(E)=2\rk(E)c_2(E)-(\rk(E)-1)c_1(E)^2$.%
\begin{NB}
I changed the definition slightly. 
Strictly speaking, the above may not be equivalent to the original
definition ($E$ on $Y$ with $\rk(E)=\rk(c)$, $c_1(E)=c_1(c)$,
$c_2(E)=c_2(c)$) as cohomology group may have torsions.  
H.N. Nov.27.
\end{NB}

\begin{Proposition}\label{prop:K3} Let $c\in K(X)_{\rm hom}$ with either $\rk(c)>0$ or
$\rk(c)=0$ and $c_1(c)$ nef and big. 
Assume that $M_H^X(c)$ consists only of stable sheaves. Then
for $v \in K_c$,
\begin{equation*}
\chi(M_H^X(c),\lambda(v))=\binom{\frac{\Delta(c)}{2}-\rk(c)^2+\frac{\Delta(v)}{2}-\rk(v)^2+2}{\frac{\Delta(c)}{2}-\rk(c)^2+1}.
\end{equation*}
\end{Proposition}

\begin{Corollary}\label{cor:K3dual}
Let $c,v \in K(X)_{\rm hom}$ with $\chi(v\otimes c) = 0$.%
\begin{NB}
Editted by H.N. Nov.27  
Corrected 27.11.LG
\end{NB}
Assume that both $c$ and $v$ fulfill the assumptions for $c$ in 
\propref{prop:K3}.  Then 
$\chi(M_H^X(c),\lambda(-v))=\chi(M_H^X(v),\lambda(-c))$.
\end{Corollary}

Recall that our invariant is well-defined on $K(X)_{\rm hom}$ (see
\secref{backDong}). We have $v\in K_c = \{ v\mid \chi(v\otimes c) = 0\}$%
\begin{NB} corrected 27.11 LG\end{NB}
 if and
only if $c\in K_v$, therefore the line bundles $\lambda(-v)$,
$\lambda(-c)$ exist on $M_H^X(c)$, $M_H^X(v)$ respectively.

\begin{NB}
The above paragraph is added by H.N. Nov.27.  
\end{NB}

\begin{Remark}
(1)
\corref{cor:K3dual} can be viewed as a weak version of an analogue of the strange duality conjecture, which was formulated by Le Potier for $\P^2$, 
and which is in turn an analogue of the strange duality
(level-rank duality) for moduli spaces of vector bundles on curves
(see \cite{Bea},\cite{Donagi},\cite{Nakani}).
Let $c\in K(\P^2)$  with $\rk(c)>0$ and 
$v\in K_c$ with $\rk(v)=0$ and $c_1(v)>0$ and assume  $M^{\P^2}(c)\ne\emptyset\ne
M^{\P^2}(v)$. Then  the strange duality conjecture of Le Potier (see \cite{Da1},
\cite{Da2}) predicts an explicit duality between  $H^0(M^{\P^2}(c),\lambda(-v))$ and $H^0(M^{\P^2}(v),\lambda(-c))$. It is shown in \cite{Da2} that the higher 
cohomology groups $H^i(M^{\P^2}(c),\lambda(-v))$ vanish and in the known 
cases also the higher cohomology groups $H^i(M^{\P^2}(v),\lambda(-c))$
are zero, thus one has in particular that $\chi(M^{\P^2}(c),\lambda(-v))=\chi(M^{\P^2}(v),\lambda(-c))$.%
\begin{NB} Formulation changed. 23.11.06 LG\end{NB}
Thus \corref{cor:K3dual} says that on K3 surfaces this is  true more generally  for 
$c$, $v$ of any nonnegative rank,  at least when $M^{X}_H(c)$ and $M^X_H(v)$ consist only of stable sheaves. 
It seems natural to conjecture that the condition that the moduli spaces only consist of stable sheaves can be dropped. 

In the context of Brill-Noether theory of K3 surfaces Markman proposed
to put $M^X_H(v):=M^X_H(-v^\vee)$, in case $\rk(v)$ is negative (see
\cite{Markman}). It is easy to see that if $\chi(v\otimes c)=0$, then
also $\chi(-v^{\vee}\otimes c)=0$, and $\Delta(-v^{\vee})=\Delta(v)$.
Thus with this definition \propref{prop:K3} also holds if $\rk(v)$ or
$\rk(c)$ are negative.

In \cite{Yo3} the proof of an equivalent formulation of  \propref{prop:K3} in terms of the Mukai vector is sketched. In \cite{OG} there is a short sketch%
\begin{NB} One could also write "there is a hint of the proof", what do you prefer
24.11.06 LG\end{NB}
  of the 
proof of \propref{prop:K3}. Furthermore the duality map $H^0(M^{X}_H(c),\lambda(-v))^\vee \to H^0(M^X_H(v),\lambda(-c))$ is constructed and 
it is checked in some cases that it is an isomorphism.%
\begin{NB}
Corrected ($\P^2$ to $X$) Nov.28, HN.
\end{NB}
\end{Remark}

We first recall some properties of the moduli spaces
$M_H^X(r,c_1,c_2)$.
The Mukai lattice of $X$ is 
$H^*(X,{\Bbb Z})$ with the symmetric bilinear form
\begin{equation}
  \langle {w},{w}'\rangle=\int_X (c_1 \wedge c_1'-r \wedge a' \varrho-
r'\wedge a \varrho)\,,
\end{equation}
for any ${w}=(r,c_1,a)\in H^{*}(X,{\Bbb Z})$ and ${w}'=(r',c_1',a')\in
H^{*}(X,{\Bbb Z})$. Here the notation $w=(r,c_1,a)$ means $w=r\oplus c_1
\oplus a \varrho$ with $r\in H^0(X,{\Bbb Z})$, $c_1\in H^2(X,{\Bbb Z})$, 
$a\in {\Bbb Z}$ and $\varrho\in H^4(X,{\Bbb Z})$ 
is the fundamental cohomology class of
$X$ so that $\int_X \varrho =1$.  
We define a weight 2 Hodge structure on $H^*(X,{\Bbb Z})$
by $H^{p,q}(H^*(X,{\Bbb C})):=\oplus_{i}H^{p+i,q+i}(X)$. 
We set $H^*(X,{\Bbb Z})_{\mathrm{alg}}:=H^*(X,{\Bbb Z}) \cap
H^{1,1}(H^*(X,{\Bbb C}))$.
Let $\phi:K(X) \to H^*(X,{\Bbb Z})$ be a homomorphism such that 
\begin{equation*}
\begin{split}
\phi(E):=&\left(\ch(E)\sqrt{\Todd(X)}\right)^{\vee}\\
=&
(\rk(E),-c_1(E),(c_1(E)^2)/2-c_2(E)+\rk (E)).
\end{split}
\end{equation*}
Then we see that $\phi$ is injective and the image is
$H^*(X,{\Bbb Z})_{\mathrm{alg}}$.
We set $w:=(r,c_1,(c_1^2)/2-c_2+r) \in H^*(X,{\Bbb Z})$.  
By the definition of the lattice structure,
$\phi$ induces an isomorphism
$\phi:K_c \to w^{\perp} \cap H^*(X,{\Bbb Z})_{\mathrm{alg}}$.
There is a homomorphism $\theta_w$ which makes the following diagram
commutative:
\begin{equation*}
\begin{CD}
K_c @>{\lambda}>> \Pic(M_H^X(r,c_1,c_2))\\
@V{\phi}VV @VVV \\
w^{\perp} @>{\theta_w}>> H^2(M_H^X(r,c_1,c_2),{\Bbb Z})
\end{CD}
\end{equation*}
If there is a universal family ${\cal E}$, 
then $\theta_w$ is given by 
$$
\theta_w(x)=
\left[
p_{M_H^X(r,c_1,c_2)*}\left(\ch {\cal E} \sqrt{\Todd(X)} x^{\vee} \right) 
\right]_1.
$$
For $M_H^X(c)$, the following is known (cf. \cite{K3}, \cite{Yo2}).
\begin{Theorem}\label{thm:k3def}
Let $c\in K(X)_{\rm hom}$ with $\rk(c)>0$ or $\rk(c)=0$ and $c_1(c)$ nef and big.
Assume that $M_H^X(c)$ consists only of stable sheaves.

\textup{(1)}
$M_H^X(c)$ is an irreducible symplectic manifold
which is deformation equivalent to 
$X^{[n]}$, where $n=\Delta(c)/2-(\rk(c)^2-1)$.

\textup{(2)}
If $\Delta(c)/2-(\rk(c)^2-1)>1$, then
$\theta_w$ is an isomorphism such that
$\theta_w$ preserves the Hodge structure
and the Beauville quadratic form 
$q_{M_H^X(c)}$ coincides with the quadratic form
associated to the Mukai lattice:
$\langle x^2 \rangle=q_{M_H^X(c)}(\theta_w(x))$.
If $\Delta(c)/2-(\rk(c)^2-1)=1$, then
$\theta_w$ is surjective with the kernel ${\Bbb Z}w$ 
and similar properties hold.


\end{Theorem} 
 For the Euler characteristic of an irreducible symplectic 
manifold, we can use the following result due to Fujiki 
(cf. \cite[Corollary 23.18]{GHJ}).
\begin{Theorem}\label{thm:fuj}
For an irreducible symplectic manifold $M$,
there is a polynomial $f(x) \in {\Bbb Q}[x]$
such that for all $D \in H^2(M,{\Bbb Z})$,
\begin{equation*}
\int_M e^{D}\Todd(M)=f(q_M(D)),
\end{equation*} 
where $q_M$ is the Beauville quadratic form on $H^2(M,{\Bbb Z})$.
Obviously $f(x)$ is deformation invariant.
\end{Theorem}
Thus it is sufficient to compute the Euler characteristic of $\lambda(v)$ 
for the Hilbert scheme $X^{[n]}$ of $n$ points on a $K3$ surface $X$.
In this case,
the Euler characteristic is determined by \cite{EGL} 
(cf. \cite[Example 23.19]{GHJ}).
\begin{equation}\label{eq:rank1}
\chi(X^{[n]},\lambda(v))=
\binom{\frac{q_{X^{[n]}}(\lambda(v))}{2}+2+n-1}{n}.
\end{equation}
Now let $c\in K(X)$ be general. Then
for $\phi(v)=(\rk v,-c_1(v),c_1(v)^2/2-c_2(v)+\rk(v))$ 
we have
$$
q_{M_H^X(c)}(\lambda(v))=
\langle \phi(v)^2 \rangle=-2\rk(v)(c_1(v)^2/2-c_2(v)+\rk(v))+c_1(v)^2=\Delta(v)-2\rk(v)^2.
$$
Therefore Proposition \ref{prop:K3} follows from \eqref{eq:rank1}.
\begin{NB}
If $\Delta/2-(r^2-1)=0$, then
$M_H^X(r,c_1,c_2)$ is a point and the claim holds.
\end{NB}

If $c$ is a class in $K(Y)$ for a surface $Y$ we want to momentarily 
introduce the following notation. We write 
$\overline M^Y_H(c)$ for the moduli space of $H$-semistable sheaves $E$
on $Y$ with $\rk(E)=\rk(c)$, $\det(E)=\det(c)$ and $c_2(E)=c_2(c)$,
i.e.\ the moduli space with fixed determinant.%
\begin{NB}
I have added an explanation. Nov. 27. HN.  
\end{NB}
Now let $A$ be an abelian surface, we have a formula very similar to  \propref{prop:K3}.

\begin{Remark} Let $c\in K(A)$ be a class  $\rk(c)=1$.%
\begin{NB}
Corrected by H.N. Nov.27. (Originally it was written as $\rk(c) > 0$.)    
\end{NB}
 Let $v\in K_c$. Then 
$$\chi(\overline M^A_H(c),\lambda(v))=\frac{\Delta(v)+\rk(v)^2\Delta(c)}{\Delta(v)+\Delta(c)}\binom{\frac{\Delta(v)}{2}+\frac{\Delta(c)}{2}}{\frac{\Delta(c)}{2}}.$$
\end{Remark}
\begin{pf}
Put $n:=\frac{\Delta(c)}{2}$. Then $\overline M^A_H(c)=A^{[n]}$, and if $Z\subset A\times A^{[n]}$ is the universal subscheme, then
the universal sheaf  is $\I_Z\otimes p_A^*(\det(c))$. Thus for any $v\in K_c$ we get $\lambda(v)=\lambda'(v\otimes \det(c))$, where 
$\lambda'(w)$ is the determinant bundle on $A^{[n]}=\overline M^A_H(c\otimes \det(c)^{-1})$ defined via the universal sheaf $\I_Z$.
Thus replacing $v$ by $v\otimes \det(c)$ we can assume that $\det(c)=0$. We write $r=\rk(v)$. Then we get 
\begin{equation}\label{eq:abel}\begin{split} \lambda(v)&=\det(p_{A^{[n]}!}( p_A^*(v)\otimes \I_Z))=\det(p_{A^{[n]}*}(p_A^*(v)\otimes \oo_Z))^{-1}\\&=
\det(p_{A^{[n]}*}(\oo_Z))^{\otimes (1-r)}\otimes \det(p_{A^{[n]}*}(p_A^*(\det(v))\otimes \oo_Z))^{-1}.\end{split}\end{equation}
In the last line we use that $\det(p_{A^{[n]}*} (p_A^*(v)\otimes \oo_Z))$ depends only on $\rk(v)$ and $\det(v)$, so we can replace $v$ by 
$\oo_A^{\oplus(r-1)}\oplus \det(v)$.
Thus by \cite[Theorem 5.3]{EGL} we get 
$$\chi(A^{[n]},\lambda(v))=
\frac{\frac{c_1(v)^2}{2}}{\frac{c_1(v)^2}{2}-(r^2-1)n}\binom{\frac{c_1(v)^2}{2}-(r^2-1)n}{n}.$$
Finally the condition $\chi(c\otimes v)$ gives $c_1(v)^2/2-c_2(v)=rn$, which is equivalent to 
$c_1(v)^2/2=r^2n+\frac{\Delta(v)}{2}$. The result follows.
\end{pf}

It seems natural to expect that a similar formula also holds for $\rk(c)\ge 0$ 
arbitrary. The simplest formula possible seems to be 
$$\chi(\overline M^A_H(c),\lambda(v))=\frac{\rk(c)^2\Delta(v)+\rk(v)^2\Delta(c)}{\Delta(v)+\Delta(c)}\binom{\frac{\Delta(v)}{2}+\frac{\Delta(c)}{2}}{\frac{\Delta(c)}{2}}.$$

\begin{NB} I kept in the moment my formula without calling it a conjecture. 
If you prefer you can also remove it. 
24.11.LG

I think that it is OK to keep this.
Nov. 27. HN.
\end{NB}

\begin{Remark} Let $Y$ be projective  surface and let $c,v\in K(Y)$ with $\rk(v)=\rk(c)=1$ and $\chi(c\otimes v)=0$. Then
$\chi(\overline M^Y_H(c),\lambda(-v))=\binom{\frac{\Delta(c)}{2}+\frac{\Delta(v)}{2}}{\frac{\Delta(c)}{2}},$
in particular $\chi(\overline M^Y_H(c),\lambda(-v))=\chi(\overline M^Y_H(v),\lambda(-c)).$
\end{Remark}
\begin{pf} Write $c_1(c)=L$, $c_1(v)=M$, $c_2(c)=l=\Delta(c)/2$, $c_2(v)=m=\Delta(c)/2$.
Let $Z\subset Y\times Y^{[l]}$ be the universal subscheme. 
Then the universal sheaf on $Y\times M^Y_H(c)=Y\times Y^{[l]}$ is $\I_Z\otimes p_Y^*(L)$. 
Thus 
$$\lambda(-v)=-\det(p_{Y^{[l]}!}(\I_Z\otimes p_{Y}^*(L\otimes M)))=\det(p_{Y^{[l]}*}(\oo_Z\otimes p_Y^*(L\otimes M))).$$
Thus
\cite[Lemma 5.1]{EGL}, we get $\chi(\overline M^Y_H(c),\lambda(-v))=\binom{\chi(L\otimes M)}{l}$.
By the Riemann-Roch theorem  $\chi(c\otimes v)=0$ is 
equivalent to $\frac{\Delta(c)}{2}+\frac{\Delta(v)}{2}=\chi(L\otimes M)$. The result follows.
\end{pf}

\begin{NB} I comment out the previous version of this section 14.11.LG
\subsection{$K$-theoretic invariant for $K3$ surfaces}

Let $X$ be a projective $K3$ surface. In this subsection we calculate
the $K$-theoretic invariants for $X$ as examples. The results of this
subsection are independent of the rest of the paper, and will not be
used later.
\begin{Proposition}\label{prop:K3}
Assume that $M_H^X(r,c_1,c_2)$ consists only of stable sheaves. Then
for $v \in K_c$,
\begin{equation*}
\chi(M_H^X(r,c_1,c_2),\lambda(v))=
\binom{\frac{1}{2}(\frac{\rk v}{r}c_1+c_1(v))^2+2+
\frac{r^2-(\rk v)^2}{r^2}(\frac{\Delta}{2}-r^2)}
{\frac{\Delta}{2}-r^2+1}
\end{equation*}
where $\Delta=2rc_2-(r-1)(c_1^2)$.
\end{Proposition}

\begin{Remark}
$\frac{1}{2}(\frac{\rk v}{r}c_1+c_1(v))^2+2+
\frac{r^2-(\rk v)^2}{r^2}(\frac{\Delta}{2}-r^2)$ and 
$\frac{\Delta}{2}-r^2+1$ are integers. 
\end{Remark}
We first recall some properties of the moduli spaces
$M_H^X(r,c_1,c_2)$.
The Mukai lattice of $X$ is 
$H^*(X,{\Bbb Z})$ with the symmetric bilinear form
\begin{equation}
  \langle {w},{w}'\rangle=\int_X (c_1 \wedge c_1'-r \wedge a' \varrho-
r'\wedge a \varrho)\,,
\end{equation}
for any ${w}=(r,c_1,a)\in H^{*}(X,{\Bbb Z})$ and ${w}'=(r',c_1',a')\in
H^{*}(X,{\Bbb Z})$. Here the notation $w=(r,c_1,a)$ means $w=r\oplus c_1
\oplus a \varrho$ with $r\in H^0(X,{\Bbb Z})$, $c_1\in H^2(X,{\Bbb Z})$, 
$a\in {\Bbb Z}$ and $\varrho\in H^4(X,{\Bbb Z})$ 
is the fundamental cohomology class of
$X$ so that $\int_X \varrho =1$.  
We define a weight 2 Hodge structure on $H^*(X,{\Bbb Z})$
by $H^{p,q}(H^*(X,{\Bbb C})):=\oplus_{i}H^{p+i,q+i}(X)$. 
We set $H^*(X,{\Bbb Z})_{\mathrm{alg}}:=H^*(X,{\Bbb Z}) \cap
H^{1,1}(H^*(X,{\Bbb C}))$.
Let $\phi:K(X) \to H^*(X,{\Bbb Z})$ be a homomorphism such that 
\begin{equation*}
\begin{split}
\phi(E):=&\left(\ch(E)\sqrt{\Todd(X)}\right)^{\vee}\\
=&
(\rk(E),-c_1(E),(c_1(E)^2)/2-c_2(E)+\rk (E)).
\end{split}
\end{equation*}
Then we see that $\phi$ is injective and the image is
$H^*(X,{\Bbb Z})_{\mathrm{alg}}$.
We set $w:=(r,c_1,(c_1^2)/2-c_2+r) \in H^*(X,{\Bbb Z})$.  
By the definition of the lattice structure,
$\phi$ induces an isomorphism
$\phi:K_c \to w^{\perp} \cap H^*(X,{\Bbb Z})_{\mathrm{alg}}$.
There is a homomorphism $\theta_w$ which makes the following diagram
commutative:
\begin{equation*}
\begin{CD}
K_c @>{\lambda}>> \Pic(M_H^X(r,c_1,c_2))\\
@V{\phi}VV @VVV \\
w^{\perp} @>{\theta_w}>> H^2(M_H^X(r,c_1,c_2),{\Bbb Z})
\end{CD}
\end{equation*}
If there is a universal family ${\cal E}$, 
then $\theta_w$ is given by 
$$
\theta_w(x)=
\left[
p_{M_H^X(r,c_1,c_2)*}\left(\ch {\cal E} \sqrt{\Todd(X)} x^{\vee} \right) 
\right]_1.
$$
For $M_H^X(r,c_1,c_2)$, the following is known (cf. \cite{K3}).
\begin{Theorem}
Assume that $M_H^X(r,c_1,c_2)$ consists only of stable sheaves.
\begin{enumerate}
\item 
$M_H^X(r,c_1,c_2)$ is an irreducible symplectic manifold
which is deformation equivalent to 
$X^{[n]}$, where $n=\Delta/2-(r^2-1)$.
\item
If $\Delta/2-(r^2-1)>1$, then
$\theta_w$ is an isomorphism such that
$\theta_w$ preserves the Hodge structure
and the Beauville quadratic form 
$q_{M_H^X(r,c_1,c_2)}$ coincides with the quadratic form
associated to the Mukai lattice:
$\langle x^2 \rangle=q_{M_H^X(r,c_1,c_2)}(\theta_w(x))$.
If $\Delta/2-(r^2-1)=1$, then
$\theta_w$ is surjective with the kernel ${\Bbb Z}w$ 
and similar properties hold.

\end{enumerate}

\end{Theorem} 
 For the Euler characteristic of an irreducible symplectic 
manifold, we can use the following result due to Fujiki 
(cf. \cite[Corollary 23.18]{GHJ}).
\begin{Theorem}
For an irreducible symplectic manifold $M$,
there is a polynomial $f(x) \in {\Bbb Q}[x]$
such that for all $D \in H^2(M,{\Bbb Z})$,
\begin{equation*}
\int_M e^{D}\Todd(M)=f(q_M(D)),
\end{equation*} 
where $q_M$ is the Beauville quadratic form on $H^2(M,{\Bbb Z})$.
Obviously $f(x)$ is deformation invariant.
\end{Theorem}
Thus it is sufficient to compute the Euler characteristic of $\lambda(v)$ 
for the Hilbert scheme $X^{[n]}$ of $n$ points on a $K3$ surface $X$.
In this case,
the Euler characteristic is determined by \cite{EGL} 
(cf. \cite[Example 23.19]{GHJ}).
\begin{equation}\label{eq:rank1}
\chi(X^{[n]},\lambda(v))=
\binom{\frac{q_{X^{[n]}}(\lambda(v))}{2}+2+n-1}{n}.
\end{equation}
For $\phi(v)=(\rk v,-c_1(v),a)$ with 
$a=-\frac{1}{r}(c_1,c_1(v))-\frac{\rk v}{r}(\frac{(c_1^2)}{2}-c_2+r)=
-\frac{1}{r}(c_1,c_1(v))-
\frac{\rk v}{2r}(-\frac{\Delta}{r}+\frac{c_1^2}{r}+2r)$,
we have
$$
q_{M_H^X(r,c_1,c_2)}(\lambda(v))=
\langle \phi(v)^2 \rangle=\left(c_1(v)+\frac{\rk v}{r}c_1 \right)^2-
2(\rk v)^2\frac{\Delta/2-r^2}{r^2}.
$$
Therefore Proposition \ref{prop:K3} follows from \eqref{eq:rank1}.
\end{NB}
\begin{NB}
If $\Delta/2-(r^2-1)=0$, then
$M_H^X(r,c_1,c_2)$ is a point and the claim holds.
\end{NB}

\subsection{Nekrasov partition function}\label{subsec:Nek}
We briefly review the $K$-theoretic Nekrasov partition function in the case of rank $2$.
For more details see \cite[section 1]{NY3}.
Let $\ell_\infty$ be the line at infinity in $\P^2$.
Let $M(n)$ be the moduli space of
pairs $(E,\Phi)$, where $E$ is a rank $2$ torsion-free sheaf on $\P^2$ with $c_2(E)=n$,
which is locally free in a neighbourhood of $\ell_\infty$ and 
$\Phi:E|_{\ell_\infty}\to \oo_{\ell_\infty}^{\oplus 2}$ is an isomorphism.
$M(n)$ is a nonsingular quasiprojective variety of dimension $4n$.
The tangent space to $M(n)$ at $(E,\Phi)$ is $\Ext^1(E,E(-l_{\infty}))$.

Let $\Gamma:=\C^*\times \C^*$ and $\widetilde T:=\Gamma\times \C^*$.
$\widetilde T$ acts on $M(n)$ as follows:
For $(t_1,t_2)\in \Gamma$, let $F_{t_1,t_2}$ be the automorphism of 
$\P^2$ defined by 
$F_{t_1,t_2}([z_0,z_1,z_2])\mapsto[z_0,t_1z_1,t_2z_2]$, and for $e_2\in \C^*$ let
$G_{e_2}$ be the automorphism of $\oo_{\ell_\infty}^{\oplus 2}$ given by 
$(s_1,s_2)\mapsto (e_2^{-1}s_1,e_2s_2)$.
Then for $(E,\Phi)\in M(n)$ we put
$(t_1,t_2,e_2)\cdot(E,\Phi):=\big((F_{t_1,t_2}^{-1})^*E,\Phi'\big)$, where
$\Phi'$ is the composition 
\begin{diagram}[height=.8em,width=2em,scriptlabels]
(F^{-1}_{t_1,t_2})^*(E)|_{\ell_\infty}&\rTo^{(F^{-1}_{t_1,t_2})^*\Phi}& (F^{-1}_{t_1,t_2})^*\oo_{\ell_\infty}^{\oplus 2}&
\rTo{} &\oo_{\ell_\infty}^{\oplus 2}&\rTo^{G_{e_2}}& \oo_{\ell_\infty}^{\oplus 2}
\end{diagram}
where the middle arrow is the homomorphism given by the action.
%

\begin{Notation}
We denote $e_2$ the one-dimensional 
$\widetilde T$-module given by $(t_1,t_2,e_2)\mapsto e_2$.
and similar we write $t_i$  ($i=1,2$) for the  $1$-dimensional
$\widetilde T$ modules given by $(t_1,t_2,e_2)\mapsto t_i$.
We also write $e_1:=e_2^{-1}$.

Let $\ve_1,\ve_2,a$ be the coordinates on the Lie algebra of $\widetilde T$
corresponding to $t_1,t_2,e_2$. Then $\ve_1,\ve_2,a$ are generators of the equivariant
cohomology $H_{\widetilde T}^*(pt)$ of a point. We relate the two sets of variables
by 
$t_1=e^{\bbeta \ve_1}, t_2=e^{\bbeta \ve_2}, e_2=e^{\bbeta a}$, where $\bbeta\in\C$ is
a parameter. We write $a_1:=-a$, $a_2:=a$.
\end{Notation}

The instanton part of the $K$-theoretic partition function is defined as
\begin{equation}\label{eq:ZInstK}
\begin{split}
\Zin_K(\ve_1,\ve_2,a;\Lambda,\bbeta)&:=\sum_{n=0}^\infty
((\bbeta\Lambda)^4e^{-\bbeta(\ve_1+\ve_2)})^n\sum_{i} (-1)^i \ch H^i(M(n),\oo).
\end{split} 
\end{equation}
Here the character $\ch$ is a formal sum of weight spaces, which are all finite-dimensional by \cite[section 4]{NY2}.

Let $x,y$ be the coordinates on $\A^2=\P^2\setminus \ell_\infty$.
The fixpoint set $M(n)^{\widetilde T}$ is the set of 
$(\I_{Z_1},\Phi_1)\oplus (\I_{Z_2},\Phi_2)$, where the $\I_{Z_\alpha}$ are
 ideal sheaves
of zero dimensional schemes $Z_\alpha$ with support in the origin of $\A^2$ 
with $\len(Z_1)+\len(Z_2)=n$, and $\Phi_\alpha$ ($\alpha=1,2$) are isomorphisms of $\I_{Z_\alpha}|_{\ell_\infty}$ with the $\alpha$-th factor of 
$\oo_{\ell_\infty}^{\oplus 2}$.
Write $I_\alpha$  for the ideal of $Z_\alpha$ in $\C[x,y]$. Then the above is a fixpoint
if and only if $I_1$ and $I_2$ are generated by monomials in $x,y$.
The fixed point set $M(n)^{\widetilde T}$ is parametrized by the pairs of Young
diagrams  $\vec Y=(Y_1,Y_2)$ so that the ideal $I_\alpha$ is generated by the
$x^iy^i$ with $(i-1,j-1)$ outside $Y_i$. The total number of boxes  is 
$|\vec Y|:=|Y_1|+|Y_2|=n$.

We use the following notations:
For a Young diagram $Y$ let $\lambda_i$ the length of the $i^{\mathrm{th}}$ column.
Let $Y'$ be the transpose of $Y$ and let
$\lambda_j'$ be the length of the $j^{\mathrm{th}}$ column of $Y'$ (equal to the length of the 
$j^{\mathrm{th}}$ row of $Y$). 
For $s=(i,j)\in \Z_{\ge0}\times \Z_{\ge0}$ let 
$$a_Y(s):=\lambda_i-j, \quad l_Y(s)=\lambda_j'-i, \quad a'(s)=j-1,\quad l'(s)=i-1.
$$

Following \cite{NY3} let, for  $\alpha,\beta\in \{1,2\}$,  
\begin{equation}
\label{eq:nab}
\begin{split}
n_{\alpha,\beta}^{\vec{Y}}(\ve_1,\ve_2,a;\bbeta):=\prod_{s\in Y_\alpha}&
\Big(1-e^{-\bbeta(-l_{Y_\beta}(s)\ve_1+(a_{Y_{\alpha}}(s)+1)\ve_2+a_\beta-a_\alpha)}
\Big)\\
&\times \prod_{s\in Y_\beta}\Big(
(1-e^{-\bbeta((l_{Y_\alpha}(s)+1)\ve_1-a_{Y_{\beta}}(s)\ve_2+a_\beta-a_\alpha)}
\Big)
\end{split}
\end{equation}
be the $\widetilde  T$-equivariant character of 
$\bigl(\Ext^1(\I_{Z_{\alpha}},\I_{Z_{\beta}}(-\ell_{\infty}))\bigr)^\vee$. 
Then by the Atiyah-Bott Lefschetz fixed point formula we have
\begin{equation}\label{eq:ZinstKF}
\Zin_K(\ve_1,\ve_2,a;\Lambda,\bbeta)=
\sum_{\vec Y} \frac{((\bbeta\Lambda)^4e^{-\bbeta(\ve_1+\ve_2)})^{|\vec Y|}}
{\bigwedge_{-1}T^*_{\vec Y} M(n)}=
\sum_{\vec Y} \frac{((\bbeta\Lambda)^4e^{-\bbeta(\ve_1+\ve_2)})^{|\vec Y|}}
{\prod_{\alpha,\beta=1,2}n_{\alpha,\beta}^{\vec{Y}}(\ve_1,\ve_2,a;\bbeta)},
\end{equation}
where $\bigwedge_{-1}$ is the alternating sum of the exterior powers.

More generally we will consider the {\it partition function with 5D Chern-Simons 
term} (see \cite{Ta}):
Let $\E$ be the universal sheaf on  $\P^2\times M(n)$.
Consider the line bundle 
$${\cal L}:=\lambda_\E(\oo_{\P^2}(-\ell_\infty))^{-1}=(\det p_{2!}(\E\otimes p_1^*\oo_{\P^2}(-\ell_\infty)))^{-1}=
\det R^1p_{2*}(\E\otimes p_1^*\oo_{\P^2}(-\ell_\infty)).$$
For an integer $m$ consider the generating function
\begin{equation}\label{eq:Zinstm}
\begin{split}
\Zin_m(\ve_1,\ve_2,a;\Lambda,\bbeta,\tau):=\sum_{n=0}^\infty &
((\bbeta\Lambda)^4e^{-\bbeta(1+\frac{m}2)(\ve_1+\ve_2)})^n \sum_{i}
(-1)^i \ch H^i(M(n),{\cal L}^{\otimes m})
\\
 & \times \exp\Big(\tau\Big(-n+\frac{a^2}{\ve_1\ve_2}\Big)\Big).
\end{split} 
\end{equation}
\begin{NB}
The term $(1+\frac{m}2)$ is added. Oct. 18, H.N.   
\end{NB}
We denote $\Zin_m(\ve_1,\ve_2,a;\Lambda,\bbeta,0)$ simply by
$\Zin_m(\ve_1,\ve_2,a;\Lambda,\bbeta)$, in particular
$\Zin_0(\ve_1,\ve_2,a;\Lambda,\bbeta)=\Zin_K(\ve_1,\ve_2,a;\Lambda,\bbeta).$

We put 
\begin{equation*}
\begin{split}C^{\vec Y}_m(\ve_1,\ve_2,a;\bbeta,\tau)&:=\exp\Big(m\bbeta\sum_{\alpha=1}^2\sum_{s\in Y_\alpha}
(a_\alpha -l'(s)\ve_1-a'(s)\ve_2)\Big)\exp\Big(\tau\Big(-|\vec
Y|+\frac{a^2}{\ve_1\ve_2}\Big)\Big)
.
\end{split}
\end{equation*}
Then we get by localization
\begin{equation}\label{eq:ZinstmF}
  \Zin_m(\ve_1,\ve_2,a;\Lambda,\bbeta,\tau)
  = \sum_{\vec Y} \frac{((\bbeta\Lambda)^4e^{-\bbeta(1+\frac{m}2)(\ve_1+\ve_2)})^{|\vec Y|}C^{\vec Y}_m(\ve_1,\ve_2,a;\bbeta,\tau)}
{\prod_{\alpha,\beta=1,2}n_{\alpha,\beta}^{\vec{Y}}(\ve_1,\ve_2,a;\bbeta)}.
\end{equation}
(see also \cite{Ta}).%
\begin{NB}
$1+\frac{m}2$ added. The label \verb+eq:ZinstmF+ is moved, Oct. 18, H.N.  
\end{NB}
We briefly sketch the argument:
Let $\vec Y=(Y_1,Y_2)$ correspond to a fixpoint $(\I_{Z_1},\Phi_1)\oplus 
(\I_{Z_2},\Phi_2)$ of $M(n)$. By localization we have to show that
$$H^1(\P^2,(\I_{Z_1}\oplus \I_{Z_2})\otimes \oo(-\ell_\infty))=
\sum_{\alpha=1}^2 \sum_{s\in Y_\alpha}e_\alpha t_1^{-l'(s)}t_2^{-a'(s)},$$
as $\widetilde T$ modules.
The exact sequence $0\to (\I_{Z_1}\oplus \I_{Z_2})\otimes \I(-\ell_{\infty})\to
\oo(-\ell_{\infty})^{\oplus 2}\to \oo_{Z_1}\oplus \oo_{Z_2}\to 0$
induces an isomorphism 
$H^1(\P^2,(\I_{Z_1}\oplus \I_{Z_2})\otimes \oo(-\ell_\infty))\simeq
H^0(\oo_{Z_1})\oplus H^0(\oo_{Z_2})$.
We have seen that an equivariant basis of $H^0(\oo_{Z_\alpha})$
is the set $\bigl\{x^{l'(s)}y^{a'(s)}\bigm| s\in Y_\alpha\bigr\}$.
By definition $(t_1,t_2)\in\Gamma$ acts by multiplying $x$ by $t_1^{-1}$ and  $y$ by $t_2^{-1}$.
Finally by definition $e_\alpha$ acts $H^0(\oo_{Z_\alpha})$ by multiplying with $e_\alpha$. The claim follows.

For a variables $\tau_0,\tau_1$ let
\begin{equation}
\label{eq:FY}
\begin{split}
& E^{\vec Y}(\ve_1,\ve_2,a;\bbeta,\tau_0,\tau_1)
\\
:=& 
\exp\Big(
\sum_{\alpha=1}^2\sum_{\rho=0}^1\tau_\rho
\Big[\frac{e^{\bbeta a_\alpha}}{\bbeta^2 \ve_1\ve_2}
\Big(1-(1-e^{-\bbeta\ve_1})(1-e^{-\bbeta\ve_2})
\sum_{s\in Y_\alpha} e^{-\bbeta(l'(s)\ve_1+a'(s)\ve_2)}\Big)\Big]_{\rho}\Big).
\end{split}
\end{equation}
Here $[\cdot]_{\rho}$ means the part of degree $\rho$, where
$a,\ve_1,\ve_2$ have degree $1$.
This is $\exp(\sum_{\rho=0}^1 \tau_\rho \ch_{\rho+2}(\mathcal E)/[\C^2])$.
(See \cite[p.59]{NY2}.)
Then an easy computation gives that
\begin{equation}\label{eq:EKFE}
E^{\vec Y}(\ve_1,\ve_2,a;\bbeta,\tau,m)=C^{\vec Y}_{-m}(\ve_1,\ve_2,a;\bbeta,\tau)\exp\Big(m\bbeta\Big(|\vec Y|\frac{\ve_1+\ve_2}{2}+(|Y_2|-|Y_1|)
\frac{a^3}{6\ve_1\ve_2}\Big)\Big).
\end{equation}
\begin{NB}
Write $E^{\vec Y}(\ve_1,\ve_2,a;\bbeta,\tau_0,\tau_1)=\exp(\tau_0 G_0+\tau_1
G_1)$.
Then we get 
$G_0=-|\vec Y|+\frac{a^2}{\ve_1\ve_2}$,
\begin{align*}
G_1&=\sum_{\alpha=1}^2\sum_{s\in Y_\alpha}\bbeta\Big(l'(s)\ve_1+a'(s)\ve_2-a_\alpha+\frac{a_\alpha^3}{6\ve_1\ve_2}+(\ve_1+\ve_2)/2\Big)\\&=
\sum_{\alpha=1}^2\sum_{s\in Y_\alpha}\bbeta(l'(s)\ve_1+a'(s)\ve_2-a_\alpha)+\bbeta(|Y_2|-|Y_1|)\frac{a^3}{6\ve_1\ve_2}+\bbeta|\vec Y|(\ve_1+\ve_2)/2.
\end{align*}
Corrected 4.10.LG\end{NB}

As a power series in $\Lambda$, 
$\Zin_m(\ve_1,\ve_2,a;\Lambda,\bbeta,\tau)$
starts with $1$. 
Thus 
$$\Fin_m(\ve_1,\ve_2,a;\Lambda,\bbeta,\tau):=\log\Zin_m(\ve_1,\ve_2,a;\Lambda,\bbeta,\tau)
$$
is well-defined and we put $\Fin_m(\ve_1,\ve_2,a;\Lambda,\bbeta):=
\Fin_m(\ve_1,\ve_2,a;\Lambda,\bbeta,0)$, $\Fin_K(\ve_1,\ve_2,a;\Lambda,\bbeta):=
\Fin_0(\ve_1,\ve_2,a;\Lambda,\bbeta)$.



We define the perturbation part, see \cite[section 4.2]{NY3} for more
details. 
We set
\begin{equation}\label{gammati}\begin{split}
\gamma_{\ve_1,\ve_2}(x|\bbeta;\Lambda)&:=
\begin{aligned}[t]
\frac{1}{2\ve_1\ve_2}& \left(
-\frac{\bbeta}{6}\left(x+\frac{1}{2}(\ve_1+\ve_2)\right)^3 
 +x^2\log(\bbeta\Lambda)\right)
\\
&+\sum_{n \geq 1}\frac{1}{n}
\frac{e^{-\bbeta nx}}{(e^{\bbeta n\ve_1}-1)(e^{\bbeta n\ve_2}-1)},
\end{aligned}
\\
\widetilde{\gamma}_{\ve_1,\ve_2}(x|\bbeta;\Lambda)
&:=\gamma_{\ve_1,\ve_2}(x|\bbeta;\Lambda)+
\frac{1}{\ve_1 \ve_2} \left(\frac{\pi^2 x}{6 \bbeta}-\frac{\zeta(3)}{\bbeta^2} \right)\\
&\qquad \qquad +\frac{\ve_1+\ve_2}{2\ve_1 \ve_2} 
\left( x \log (\bbeta \Lambda)+\frac{\pi^2}{6\bbeta} \right)+
\frac{\ve_1^2+\ve_2^2+3\ve_1 \ve_2}{12 \ve_1 \ve_2} \log(\bbeta\Lambda)
\end{split}
\end{equation}
for $(x,\bbeta,\Lambda)$ in a neighbourhood  of $\sqrt{-1}\R_{>0}\times\sqrt{-1}\R_{<0}
\times\sqrt{-1}\R_{>0}$.
We formally expand $\ve_1 \ve_2\widetilde{\gamma}_{\ve_1,\ve_2}(x|\bbeta;\Lambda)$
as a power series of $\ve_1,\ve_2$ (around $\ve_1=\ve_2=0$).
Expanding
\begin{equation}\label{cm}
\frac{1}{(e^{\ve_1 t}-1)(e^{\ve_2t}-1)}=
\sum_{n \geq 0} \frac{c_n}{n!}t^{n-2},
\end{equation}
we obtain
$$
\sum_{n \geq 1}\frac{1}{n}
\frac{e^{-\bbeta nx}}{(e^{\bbeta n\ve_1}-1)(e^{\bbeta n\ve_2}-1)}
=\sum_{m \geq 0} \frac{c_m}{m!}\bbeta^{m-2} \mathrm{Li}_{3-m}(e^{-\bbeta x}),
$$
where $\Li_{3-m}$ is the polylogarithm (see \cite[Appendix B]{NY3}  for details).
Here we choose the branch of $\log$ by $\log(r\cdot e^{i\phi})=\log(r)+i\phi$
with $\log(r)\in \R$
for $\phi\in (-\pi/2,3\pi/2)$ and $r\in \R$.
We define  $\gamma_{\ve_1,\ve_2}(-x|\bbeta;\Lambda)$  by analytic continuation along circles in a 
counter-clockwise way. 
We then define the perturbation part of the partition function by
\begin{equation}\label{Fper}
\begin{split}
\Fper_K(\ve_1,\ve_2,x;\Lambda,\bbeta)&:=-\widetilde\gamma_{\ve_1,\ve_2}(2x|\bbeta;\Lambda)
-\widetilde\gamma_{\ve_1,\ve_2}(-2x|\bbeta;\Lambda),
\end{split}
\end{equation}
\begin{NB} The perturbation part for the CS term is deleted. Oct. 18,
  H.B.
\end{NB}
Then $\Fper_K(\ve_1,\ve_2,x;\Lambda,\bbeta)$ is  a formal power series in $\ve_1,\ve_2$ whose coefficients are 
holomorphic functions in $\Lambda\in \C\setminus \sqrt{-1}\R_{\le 0}$, $x\in \C\setminus \sqrt{-1}\R_{\le 0}$, $\bbeta\in \C$ with $|\bbeta|<\frac{\pi}{|x|}$.

Finally we define 
\begin{align*}
F_m(\ve_1,\ve_2,a;\Lambda,\bbeta,\tau)&:=\Fper_K(\ve_1,\ve_2,a;\Lambda,\bbeta)
+ \log\Zin_m(\ve_1,\ve_2,a;\Lambda,\bbeta,\tau), 
\\
F_m(\ve_1,\ve_2,a;\Lambda,\bbeta)&:=F_m(\ve_1,\ve_2,a;\Lambda,\bbeta,0),\
F_K(\ve_1,\ve_2,a;\Lambda,\bbeta):=F_0(\ve_1,\ve_2,a;\Lambda,\bbeta).
\end{align*}
Formally one defines $Z_m(\ve_1,\ve_2,a;\Lambda,\bbeta,\tau):=\exp(\Fper_m(\ve_1,\ve_2,a;\Lambda,\bbeta))\Zin_K(\ve_1,\ve_2,a;\Lambda,\bbeta,\tau)$,
and similarly for $Z_m(\ve_1,\ve_2,a;\Lambda,\bbeta)$, $Z_K(\ve_1,\ve_2,a;\Lambda,\bbeta)$.

\begin{NB}
It would be in some sense more natural to change the  definition of the "extended" partition function, so that it also is in $K$-theory (instead of just using the
extra factor from the topological version. I have attempted to do so, but was unable to do it. Thus I stick to my old definition. It would be easy if we wanted to compute holomorphic Euler characteristics
of bundles on  the moduli space of the form $p_!(q^*(v)\otimes \E)$. But we want instead to take the determinant bundle. The local contributions at the fixpoints of $M(n)$ 
before taking the determinant are equivariant vector spaces of infinite rank, so one would need some kind of renormalization to define the determinant, if this even can make sense. 

I sketch below what would be the correct definition if we just wanted to compute 
$p_!(q^*(v)\otimes \E)$. 

More generally we will consider the following:
Let
\begin{equation}
\label{eq:FY}
E^{\vec Y}(\ve_1,\ve_2,a;\bbeta,\tau):=
\tau \sum_{\alpha=1}^2 \
 \frac{e_\alpha}{(1-t_1)(1-t_2)}
\Big(1-(1-t_1)(1-t_2)
\sum_{s\in Y_\alpha} t_1^{l'(s)}t_2^{a'(s)}\Big).
\end{equation}
Let $p:\A^2\times M(n)\to M(n)$ be the projection, and let 
$\E$ be the universal sheaf on $X\times M(n)$. Analogously to \cite[(4.1)]{NY2} we can formally identify
$$E^{\vec Y}(\ve_1,\ve_2,a;\bbeta,\tau)=\tau \iota_{\vec Y}^*(p_!\E).$$
Then the instanton part of the partition function is defined as
\begin{equation}\label{eq:ZinstKFE}
\Zin_K(\ve_1,\ve_2,a;\Lambda,\bbeta,\tau):=
\sum_{\vec Y} \frac{((\bbeta\Lambda)^4e^{-\bbeta(\ve_1+\ve_2)})^{|\vec Y|}
E^{\vec Y}(\ve_1,\ve_2,a;\bbeta,\tau)}
{\prod_{\alpha,\beta=1,2}n_{\alpha,\beta}^{\vec{Y}}(\ve_1,\ve_2,a;\bbeta)}
\in \Q(\ve_1,\ve_2,a,\bbeta)[[\bbeta\Lambda,\tau]].
\end{equation}
In particular $\Zin_K(\ve_1,\ve_2,a;\Lambda,\bbeta,0)
=\Zin_K(\ve_1,\ve_2,a;\Lambda,\bbeta)$.
As a power series in $\Lambda$, 
$\Zin_K(\ve_1,\ve_2,a;\Lambda,\bbeta,\tau)$
starts with $1$. 
Thus 
$$\Fin_K(\ve_1,\ve_2,a;\Lambda,\bbeta,\tau):=\log\Zin_K(\ve_1,\ve_2,a;\Lambda,\bbeta,\tau)
\in \Q(\ve_1,\ve_2,a,\bbeta)[[\bbeta\Lambda,\tau]]$$
is well-defined and we put $\Fin_K(\ve_1,\ve_2,a;\Lambda,\bbeta):=\Fin(\ve_1,\ve_2,a;\Lambda,\bbeta).$

\end{NB}

\subsection{More on the partition function with 5D Chern-Simons term}
\label{subsec:MoreOn}

We explain how the known properties of the $K$-theoretic Nekrasov
partition function, obtained in \cite{NY3}, can be generalized to the
partition function with 5D Chern-Simons term, at least conjecturally.
Our explanation is mathematical, so a physical motivation can be found
in \cite{IMS,Ta} and the references therein.

This subsection is independent of the rest of this paper, and can be
safely skipped. We also keep the notation in \cite{NY3} except we set
$\q = \Lambda^{1/2r}$.

We consider the  general case $r\ge 2$, although we only consider the case $r=2$ in the main part of the paper.
Let $M(r,n)$ be the framed moduli space of rank $r$ torsion free
sheaves $E$ on $\P^2$ with $c_2(E) = n$.
Let $\E$ be the universal sheaf
on $\P^2\times M(r,n)$.  Consider the line bundle
$${\cal L}:=\lambda_\E(\oo_{\P^2}(-\ell_\infty))^{-1}=(\det p_{2!}(\E\otimes p_1^*\oo_{\P^2}(-\ell_\infty)))^{-1}=
\det R^1p_{2*}(\E\otimes p_1^*\oo_{\P^2}(-\ell_\infty)).$$
For an integer $m$ consider the generating function
\begin{equation}\label{eq:Zinstm'}
\begin{split}
\Zin_m(\ve_1,\ve_2,a;\Lambda,\bbeta)&:=\sum_{n=0}^\infty
((\bbeta\Lambda)^{2r}e^{-\bbeta(r+m)(\ve_1+\ve_2)/2})^n \sum_{i} (-1)^i \ch H^i(M(r,n),{\cal L}^{\otimes m}). 
\end{split} 
\end{equation}
\begin{NB}
I changed the definition of the partition function. I multiply
$e^{-\bbeta(r+m)(\ve_1+\ve_2)/2}$. The previous definition is
just $e^{-\bbeta r(\ve_1+\ve_2)/2}$. But this seems not compatible
with the blowup formula.
Oct. 14, H.N
\end{NB}
By the localization formula we have
\begin{equation}\label{eq:Zsum'}
\begin{split}
& \Zin_m(\ve_1,\ve_2,\vec{a};\Lambda,\bbeta)
\\
   =\; & \sum_{\vec{Y}} 
   \frac{((\bbeta\Lambda)^{2r} e^{-\bbeta(r+m)(\ve_1+\ve_2)/2})^{|\vec{Y}|}}
    {\displaystyle\prod_{\alpha,\beta}
    n^{\vec{Y}}_{\alpha,\beta}(\ve_1,\ve_2,\vec{a};\bbeta)}
   \exp\Big(m\bbeta\sum_\alpha \sum_{s \in Y_\alpha}
   (a_\alpha-l'(s)\ve_1-a'(s)\ve_2)\Big),
\end{split}
\end{equation}
where $\vec{Y}$ is an $r$-tuple of Young diagrams. The argument is the
same as in the rank $2$ case.

We have
\begin{equation}\label{eq:SerreDuality}
  \Zin_{-m} (-\ve_1,-\ve_2,-\vec{a};\Lambda,\bbeta)
  = \Zin_m(\ve_1,\ve_2,\vec{a};\Lambda,\bbeta).
\end{equation}
This is a consequence of Serre duality and the equality
$K_{M(r,n)} = e^{-r\bbeta(\ve_1+\ve_2)n}$ (\cite[Lemma~3.6]{NY3}).
But it also follows directly from
\begin{equation*}
   \prod_{\alpha,\beta} n_{\alpha,\beta}^{\vec{Y}}(-\ve_1,-\ve_2,-\vec{a};\bbeta)
   =
   e^{\bbeta r(\ve_1+\ve_2)|\vec{Y}|}
      \times
   \prod_{\alpha,\beta} n_{\alpha,\beta}^{\vec{Y}}(\ve_1,\ve_2,\vec{a};\bbeta).
\end{equation*}
\begin{NB}
We need to consider the pair $n_{\alpha,\beta}^{\vec{Y}}$
and $n_{\beta,\alpha}^{\vec{Y}}$. Oct. 13, H.N.
\end{NB}

\subsubsection{Correlation function on blow-up}

Let $X$ be the blow-up of ${\Bbb P}^2$ at the origin of 
${\Bbb C}^2$. Let $\widehat M(r,k,\widehat n)$ be the framed moduli
space on $X$.
We define the similar partition function
$\bZin_{m,k,d}(\ve_1,\ve_2,\vec{a};\Lambda,\bbeta)$ on $X$ by
considering
\begin{equation*}
    \sum_i (-1)^i \ch H^i(\widehat M(r,k,\widehat n),
    \widehat{\mathcal L}^{\otimes m}
    \otimes \mu(C)^{\otimes d}),
\end{equation*}
where $\widehat{\mathcal L}$ is defined as in the case of $\mathbb
P^2$ by taking the universal bundle $\widehat E$ over $X\times
\widehat M(r,k,\widehat n)$. (See \cite[\S2.1]{NY3}.)
As in \cite[\S2.2]{NY3}, we can write down
$\bZin_{m,k,d}(\ve_1,\ve_2,\vec{a};\Lambda,\bbeta)$ 
in terms of $\Zin_{m}(\ve_1,\ve_2,\vec{a};\Lambda,\bbeta)$.
The new factor comes from
\begin{equation*}
\begin{split}
   c_1\Big( \bigoplus_\alpha H^1({\cal O}_X(k_\alpha C))e_\alpha\Big)&=
   \sum_\alpha \frac{k_\alpha^3-k_\alpha}{6}(\ve_1+\ve_2)+ 
   \sum_\alpha \frac{k_\alpha(k_\alpha-1)}{2}a_\alpha
\\
   &= \sum_\alpha \frac{k_\alpha^3}{6}(\ve_1+\ve_2)+ 
   \sum_\alpha \frac{k_\alpha^2}{2}a_\alpha
   -\frac{1}2(\vec{k},\vec{a}).
\end{split}
\end{equation*}
Then the blowup formula is a slight modification of \cite[(2.2)]{NY3}:
\begin{multline}\label{eq:blow-up1}
  \bZin_{m,k,d}(\ve_1,\ve_2,\vec{a};\Lambda,\bbeta) =
  \sum_{\substack{\vec{k}\in\Z^r\\ \sum k_\alpha = k}}
   \frac{(e^{\bbeta(\ve_1+ \ve_2) (d - (r+m)/2)}
   (\bbeta\Lambda)^{2r})^{(\vec{k},\vec{k})/2}e^{\bbeta(\vec{k},\vec{a})(d-m/2)}}
   {\prod_{\vec\alpha\in\Delta}
     l^{\vec{k}}_{\vec{\alpha}}(\ve_1,\ve_2,\vec{a})}
\\
  \times
  \exp\left[m\bbeta \left(\frac16(\ve_1+\ve_2)\sum_\alpha k_\alpha^3
      + \frac12\sum_\alpha k_\alpha^2 a_\alpha\right)\right]
  \\
   \times
   \Zin_m(\ve_1,\ve_2-\ve_1,\vec{a}+\ve_1\vec{k};
    e^{\bbeta\ve_1 (d-(r+m)/2)/2r)} \Lambda,\bbeta)
\\ 
  \times
   \Zin_m(\ve_1-\ve_2,\ve_2,\vec{a}+\ve_2\vec{k};
    e^{\bbeta\ve_2 (d-(r+m)/2)/2r}\Lambda,\bbeta),
\end{multline}
where $(\vec{k},\vec{a}) = \frac1{2r}\sum_{\alpha,\beta}
(k_\alpha-k_\beta)(a_\alpha - a_\beta)$, and similarly for $(\vec{k},\vec{k})$.
Note that we need to normalize a vector $\vec{k} =
(k_\alpha)_{\alpha=1}^r$ with $\sum k_\alpha = k$ into $\vec{l} =
(k_1-\frac{k}r,\dots, k_r-\frac{k}r)$, as we assume $\sum a_\alpha = 0$.
(We took this normalization in \cite{NY3} without an explanation. It
was explained in \cite[\S6]{NY1}.)
Under this normalization we have
$(\vec{k},\vec{a}) = (\vec{l},\vec{a})$, 
$l^{\vec{k}}_{\vec{\alpha}}(\ve_1,\ve_2,\vec{a})
=l^{\vec{l}}_{\vec{\alpha}}(\ve_1,\ve_2,\vec{a})$,
$n^{\vec{Y}}_{\alpha,\beta}(\ve_1,\ve_2-\ve_1,\vec{a}+\ve_1\vec{k})
= n^{\vec{Y}}_{\alpha,\beta}(\ve_1,\ve_2-\ve_1,\vec{a}+\ve_1\vec{l})$,
etc. In particular, we simply replace $\vec{k}$ by $\vec{l}$
in the original partition function with $m=0$.%
\begin{NB}
For the original case, any expression is written in terms of $k_\alpha -
k_\beta$.
Oct. 24, H.N.  
\end{NB}
However the Chern-Simons term requires some care:
{\allowdisplaybreaks
\begin{equation*}
\begin{split}
  & 
  \begin{aligned}[t]
  & \exp \left[m\bbeta \left(\sum_\alpha \sum_{s\in Y_\alpha}
      \left( (a_\alpha+\ve_1 k_\alpha - l'(s)\ve_1 -
        a'(s)(\ve_2-\ve_1)
        \right)\right) \right]
    \\
  & \qquad
  = \exp\left[ \frac{\bbeta mk}r \ve_1 |\vec{Y}|\right]
  \exp \left[m\bbeta \left(\sum_\alpha \sum_{s\in Y_\alpha}
      \left( (a_\alpha+\ve_1 l_\alpha - l'(s)\ve_1 -
        a'(s)(\ve_2-\ve_1)
        \right)\right) \right],
  \end{aligned}
\\
  &
\begin{aligned}[t]
  & \sum_\alpha \left(
  \frac{k_\alpha^3}{6}(\ve_1+\ve_2)+\frac{k_\alpha^2}{2}a_\alpha\right)
\\
  &\qquad
   = \sum_\alpha \left(
   \frac{l_\alpha^3}{6}(\ve_1+\ve_2) + \frac{l_\alpha^2}{2}a_\alpha
   \right)
   + \left(\frac{k}{2r} (\vec{l},\vec{l}) + \frac{k^3}{6r^2}\right)
   (\ve_1+\ve_2)
   + \frac{k}r (\vec{l},\vec{a}).
  \end{aligned}
\end{split}
\end{equation*}
We} rewrite \eqref{eq:blow-up1} in terms of $\vec{l}$:
{\allowdisplaybreaks
\begin{equation}\label{eq:blow-up2}
\begin{split}
  \bZin_{m,k,d}(&\ve_1,\ve_2,\vec{a};\Lambda,\bbeta) 
    =
  \exp(\frac{k^3m\bbeta}{6r^2}(\ve_1+\ve_2))
\\
  & \times
  \sum_{\{\vec{l}\} = -k/r}
   \frac{(\exp\left[{\bbeta(\ve_1+ \ve_2) 
       (d+m\left(-\frac12+\frac{k}r\right)-\frac{r}2)}\right]
   (\bbeta\Lambda)^{2r})^{(\vec{l},\vec{l})/2}}
   {\prod_{\vec\alpha\in\Delta}
     l^{\vec{l}}_{\vec{\alpha}}(\ve_1,\ve_2,\vec{a})}
\\
  & \times
  \exp\left[\bbeta(\vec{l},\vec{a})(d+m(-\frac12+\frac{k}r))\right]
\\
  & \times
  \exp\left[m\bbeta \left(\frac16(\ve_1+\ve_2)\sum_\alpha l_\alpha^3
      + \frac12\sum_\alpha l_\alpha^2 a_\alpha\right)\right]
\\
 &
   \times
   \Zin_m(\ve_1,\ve_2-\ve_1,\vec{a}+\ve_1\vec{l};
    \exp\left[{\frac{\bbeta\ve_1}{2r}
        \left\{d+m\left(-\frac12+\frac{k}r\right)-\frac{r}2\right\}}\right]
    \Lambda,\bbeta)
\\ 
  &
  \times
   \Zin_m(\ve_1-\ve_2,\ve_2,\vec{a}+\ve_2\vec{l};
    \exp\left[{\frac{\bbeta\ve_2}{2r}
        \left\{d+m\left(-\frac12+\frac{k}r\right)-\frac{r}2\right\}}\right]
    \Lambda,\bbeta).
\end{split}
\end{equation}
Here} $\{\vec{l}\} = -{k}/r$ means that the fractional part of
$l_\alpha$ is independent of $\alpha$ and equal to $-{k}/r$.%
\begin{NB}
If we use
\begin{equation}\label{eq:original}
  \Zin_m(\ve_1,\ve_2,\vec{a};\Lambda,\bbeta)
  = \sum_n \chi(M(r,n),(\det V)^{\otimes m}) 
  (e^{-r(\ve_1+\ve_2)/2} (\bbeta\Lambda)^{2r})^n,
\end{equation}
\begin{multline*}
  \bZin_{m,d}(\ve_1,\ve_2,\vec{a};\Lambda,\bbeta) =
  \sum_{\vec{k}}
   \frac{(e^{\bbeta(\ve_1+ \ve_2) (d - r/2)} 
   (\bbeta\Lambda)^{2r})^{(\vec{k},\vec{k})/2}e^{\bbeta(\vec{k},\vec{a})(d-m/2)}}
   {\prod_{\vec\alpha\in\Delta} l^{\vec{k}}_{\vec{\alpha}}(\ve_1,\ve_2,\vec{a})}
\\
  \times
  \exp\left[m\left(\frac16(\ve_1+\ve_2)\sum_\alpha k_\alpha^3
      + \frac12\sum_\alpha k_\alpha^2 a_\alpha\right)\right]
\\
  \times
   \Zin_m(\ve_1,\ve_2-\ve_1,\vec{a}+\ve_1\vec{k};
    e^{\bbeta\ve_1 (d-r/2)/2r} \Lambda,\bbeta)
   \Zin_m(\ve_1-\ve_2,\ve_2,\vec{a}+\ve_2\vec{k};
    e^{\bbeta\ve_2 (d-r/2)/2r}\Lambda,\bbeta).
\end{multline*}
\end{NB}

By Serre duality we have
\begin{equation*}
   \bZin_{m,k,d}(\ve_1,\ve_2,\vec{a};\Lambda,\bbeta)
   = \bZin_{-m,k,r-d}(-\ve_1,-\ve_2,-\vec{a};\Lambda,\bbeta)
\end{equation*}
thanks to \cite[Lemma~3.6]{NY3}. It also follows from
\eqref{eq:blow-up1} and \eqref{eq:SerreDuality} together with
\begin{equation*}
   \prod_{\vec\alpha\in\Delta}
   l^{\vec{k}}_{\vec{\alpha}}(-\ve_1,-\ve_2,-\vec{a})
   = e^{-\bbeta r(\vec{k},\vec{a})}
   \prod_{\vec\alpha\in\Delta} l^{\vec{k}}_{\vec{\alpha}}(\ve_1,\ve_2,\vec{a}).
\end{equation*}

\subsubsection{The perturbation part}\label{subsubsec:pert}

In \cite[Sect.~4.2]{NY3} one of the reason of the introduction of the
perturbation part was to simplify the the blowup formula. As we have
an extra factor
\begin{equation*}
  \exp\left[m\bbeta \left(\frac16(\ve_1+\ve_2)\sum_\alpha l_\alpha^3
      + \frac12\sum_\alpha l_\alpha^2 a_\alpha\right)\right],
\end{equation*}
we need to modify the perturbation part so that it is absorbed in the
full partition function. The answer is the cubic term:
\begin{equation*}
   \exp\left[-m\bbeta\sum_{\alpha=1}^r \frac{a_\alpha^3}{6\ve_1\ve_2}\right].
\end{equation*}
We have the difference equation 
\begin{equation*}
      \left.\frac{x^3}{6 \ve_1\ve_2}\right|_{
   \substack{x\to x+l\ve_1\\ \ve_1\to\ve_1\\ \ve_2 \to \ve_2-\ve_1}}
   +
   \left.\frac{x^3}{6 \ve_1\ve_2}\right|_{
   \substack{x\to x+l\ve_2\\ \ve_1\to\ve_1-\ve_2\\ \ve_2 \to \ve_2}}
   -
   \frac{x^3}{6 \ve_1\ve_2}
   = -\frac{l^2 x}2 - \frac{(\ve_1+\ve_2)l^3}6,
\end{equation*}
We thus define
\begin{equation*}
\begin{split}
   F_m(\ve_1,\ve_2,\vec{a};\Lambda,\bbeta) \defeq & 
   \sum_{\vec{\alpha} \in \Delta} 
   -\widetilde{\gamma}_{\ve_1,\ve_2}(\langle\vec{a},\vec{\alpha}\rangle
   |\bbeta;\Lambda) 
   - m\bbeta \sum_{\alpha=1}^r \frac{a_\alpha^3}{6\ve_1\ve_2}
   + \log \Zin_m(\ve_1,\ve_2,\vec{a};\Lambda,\bbeta),\\
   \widehat{F}_{m,k,d}(\ve_1,\ve_2,\vec{a};\Lambda,\bbeta) \defeq & 
   \sum_{\vec{\alpha} \in \Delta} 
   -\widetilde{\gamma}_{\ve_1,\ve_2}(\langle\vec{a},\vec{\alpha}\rangle
   |\bbeta;\Lambda) 
   - m\bbeta \sum_{\alpha=1}^r \frac{a_\alpha^3}{6\ve_1\ve_1}
   + \log \Zin_{m,k,d}(\ve_1,\ve_2,\vec{a};\Lambda,\bbeta),
\end{split}
\end{equation*}
where $\widetilde{\gamma}_{\ve_1,\ve_2}$ is as in \eqref{gammati}.
Note that the term
$\sum_{\alpha=1}^r \frac{a_\alpha^3}{6\ve_1\ve_2}$
disappears%
\begin{NB} changed 30.10. LG\end{NB}
when $r=2$ thanks to the condition $a_1 + a_2 = 0$.%
\begin{NB}
As we change the definition of the partition function, the
perturbative part is changed. We don't have
$m\bbeta\frac{(\ve_1+\ve_2)x^2}{2\ve_1\ve_2}$ anymore !
Oct. 14, H.N. 
\end{NB}
We formally define
\begin{equation*}
\begin{split}
Z_m(\ve_1,\ve_2,\vec{a};\Lambda,\bbeta) \defeq & \exp( 
F_m(\ve_1,\ve_2,\vec{a};\Lambda,\bbeta)),\\
\bZ_{m,k,d}(\ve_1,\ve_2,\vec{a};\Lambda,\bbeta) \defeq & \exp(
\widehat{F}_{m,k,d}(\ve_1,\ve_2,\vec{a};\Lambda,\bbeta)).
\end{split}
\end{equation*}
%
The blowup formula is
\begin{equation}\label{eq:blow-up3}
\begin{split}
  \bZ_{m,k,d}(& \ve_1,\ve_2,\vec{a};\Lambda,\bbeta)
  = 
  \exp\left[\left\{
    -\frac{\left(4\left(d + m\left(-\frac12+\frac{k}r\right)\right) - r\right)
      (r-1)}{48}
   + \frac{k^3 m}{6r^2}
  \right\}
    \bbeta(\ve_1+\ve_2)\right]
\\
&
   \times
   \sum_{\{\vec{l}\} = -k/r}\!\!
   \begin{aligned}[t]
   & Z_m(\ve_1,\ve_2-\ve_1,\vec{a}+\ve_1\vec{l};
    \exp\left[{\frac{\bbeta\ve_1}{2r}
        \left\{d+m\left(-\frac12+\frac{k}r\right)-\frac{r}2\right\}}\right]
    \Lambda,\bbeta)
\\
  & \times 
   Z_m(\ve_1-\ve_2,\ve_2,\vec{a}+\ve_2\vec{l};
    \exp\left[{\frac{\bbeta\ve_2}{2r}
        \left\{d+m\left(-\frac12+\frac{k}r\right)-\frac{r}2\right\}}\right]
    \Lambda,\bbeta).
   \end{aligned}
\end{split}
\end{equation}
This is exactly the same as \cite[(4.9)]{NY3} with the replacement
$d\to d+m(-1/2+k/r)$.

\subsubsection{A conjectural blowup equation}

For the original $K$-theoretic Nekrasov partition function we have a
blowup equation \cite[Th.~2.4 and (4.9)]{NY3}, which determines the
partition function from its perturbative part. It was derived from 
vanishing of higher direct image sheaves of a determinant line bundle%
\begin{NB} changed 15.11 LG\end{NB}
$\mu(C)$ with respect to the projection $\widehat\pi\colon \widehat
M(r,0,n)\to N(r,n)$, where $N(r,n)$ is the Uhlenbeck compactification
of the framed moduli space of locally free sheaves on $\mathbb P^2$,
denote by $M_0(r,n)$ in \cite{NY3}.

The proof of the vanishing theorem cannot be carried over to the
partition functions with Chern-Simons terms. But a numerical
computation suggests
\begin{equation}\label{eq:conjvanish}
\bZin_{m,0,d}(\ve_1,\ve_2,\vec{a};\Lambda,\bbeta) =
\Zin_m(\ve_1,\ve_2,\vec{a};\Lambda,\bbeta)
\qquad\text{for $0\le d\le r$, $|m|\le r$.}
\end{equation}
This is exactly the same what we have proved
for the original K-theoretic partition function, i.e.\  $m=0$ in \cite{NY3}.
It seems likely that the left hand side can be always written in terms
of the correlation function, which is the holomorphic Euler
characteristic of certain (virtual) bundles on $M(r,n)$. But
it can be written as above only in the limited range of $d$ and $m$.
In fact, we check the above equation holds in a slightly wider
situation when $r=2$, $m=1$: it seems to hold for $0\le d\le
3=r+m$. But we also check that when $r=2$, $m=2$, the above is not
true for $d=4$.

We have the following analogs of \cite[Lemma~4.3, Theorem~4.4]{NY3}:
\begin{Proposition}\label{prop:regular'}
Suppose \eqref{eq:conjvanish} holds and assume $|m| < r$. Then

\textup{(1)}
\(
  \Zin_m(\ve_1,-2\ve_1,\vec{a};\Lambda,\bbeta)
  = \Zin_m(2\ve_1,-\ve_1,\vec{a};\Lambda,\bbeta).
\)

\textup{(2)} 
\(
   \ve_1\ve_2 \log \Zin_m(\ve_1,\ve_2,\vec{a};\Lambda,\bbeta)
\)
is regular at $(\ve_1,\ve_2) = (0,0)$.
\end{Proposition}

We only give the proof of (1), as the proof of (2) is exactly the same
as the original.

\begin{proof}
By \eqref{eq:SerreDuality} we may assume $m\le 0$. By the assumption
we have \eqref{eq:conjvanish} for $d=0$ and $d=r+m$. Note that 
we have $0\neq r+m$ as $m\neq -r$.

Let us put $\bbeta=1$ for brevity.
We take the difference of both sides of \eqref{eq:blow-up1} with
$d=r+m$, $0$ after setting $\ve_2 = -\ve_1$. We have
\begin{equation*}
\begin{split}
  & \left(Z_n(\ve_1,-2\ve_1,\vec{a}) - Z_n(2\ve_1,-\ve_1,\vec{a})\right)
  \left(e^{(r+m)n \ve_1/2} - e^{-(r+m)n \ve_1/2}\right)
\\
  =\; &
  -\sum_{\substack{(\vec{k},\vec{k})/2+l+l' = n\\
      l\neq n, l'\neq n}}
\begin{aligned}[t]
  & \frac{e^{r{(\vec{k},\vec{a})}/{2}}
   Z_{l'}(\ve_1,-2\ve_1,\vec{a}+\ve_1\vec{k})
   Z_l(2\ve_1,-\ve_1,\vec{a}-\ve_1\vec{k})}
  {\prod_{\vec\alpha\in\Delta} l^{\vec{k}}_{\vec{\alpha}}(\ve_1,-\ve_1,\vec{a})}  \\
  & \qquad
  \times\exp\left[ \frac{m\bbeta}2 \sum_\alpha k_\alpha^2
    a_\alpha\right] \\
  & \qquad
  \times\left(e^{(r+m){(\vec{k},\vec{a})}/{2}}
    e^{{(r+m)(l'-l)}\ve_1/{2}}-
    e^{{-(r+m)(\vec{k},\vec{a})}/{2}}e^{{-(r+m)(l'-l)}\ve_1/{2}}\right),
\end{aligned}
 \end{split}
\end{equation*}
where we expand $\Zin_m$ as
\begin{equation*}
    \Zin_m(\ve_1,\ve_2,\vec{a};\Lambda,\bbeta=1)
    = \sum_n Z_n(\ve_1,\ve_2,\vec{a}) \Lambda^{2rn}.
\end{equation*}

\begin{NB}
{\bf Correction}: In the proof of \cite[Lemma~4.3]{NY3} we expand
$\Zin = \sum_n \left(
  \q\bbeta^{2r}e^{-r\bbeta(\ve_1+\ve_2)/2}\right)^n
  Z_n$.
But this was wrong and we expand it as above.
\end{NB}

Let us show that
\(
Z_n(\ve_1,-2\ve_1,\vec{a})
   = Z_n(2\ve_1,-\ve_1,\vec{a})
\)
by using the induction on $n$. It holds for $n = 0$ as $Z_0 = 1$.
Suppose that it is true for $l,m< n$. Then the right hand side of the
above equation vanishes, as terms with $(\vec{k},l,l')$ and $(-\vec{k},
l', l)$ cancel thanks to \cite[Lemma~4.1(1) and (4.2)]{NY3}, and the
term $(0,l,l)$ is $0$. Therefore it is also true for $n$.
\end{proof}

\begin{NB}
This argument cannot be applied to the original
definition~\eqref{eq:original} with $d=r$, $0$: We need to have
$e^{r{(\vec{k},\vec{a})}/{2}}$ to use \cite[(4.2)]{NY3}. But then
$e^{(\vec{k},\vec{a})(r-m/2)} = e^{(\vec{k},\vec{a})r}
e^{-(\vec{k},\vec{a})m/2}$ and
$e^{(\vec{k},\vec{a})(0-m/2)} = e^{(\vec{k},\vec{a})r}
e^{-(\vec{k},\vec{a})(r+m/2)}$.
\end{NB}


We expand 
$\ve_1\ve_2 \log Z(\ve_1,\ve_2,\vec{a};\Lambda,\bbeta)$
as in \eqref{eq:genus_expansion}.
The following can be proved exactly as in \cite[(4.11)]{NY3}:
\begin{Proposition}\label{prop:contact}
Suppose \eqref{eq:conjvanish} holds and assume $|m| < r$. Then
\begin{equation*}
\begin{split}
  & \exp\left[- \frac{\bbeta^2}{8r^2}\left(d-\frac{r+m}{2} \right)^2
    \frac{\partial^2 \F_0}{(\partial \log \Lambda)^2}
  \right]
  \Theta_E \left(\left. - \frac1{2\pi\sqrt{-1}}
      \frac{\bbeta}{2r} \left(d-\frac{r+m}{2} \right)
      \frac{\partial^2 \F_0}
      {\partial\log\Lambda\partial\vec{a}}
  \right| \tau(\bbeta) \right)
\end{split}
\end{equation*}
is independent of $d=0,\dots,r$. Here
\begin{equation*}
  \tau(\bbeta) = - \frac1{2\pi\sqrt{-1}} \frac{\partial^2 \F_0}
 {(\partial \vec{a})^2}
\end{equation*}
and $\Theta_E$ is the Riemann theta function with the characteristic
${}^t\left(\frac{1}{2},\frac{1}{2},\dots,\frac{1}{2} \right)$. 
\textup(See \cite[Appendix~B]{NY2} for convention.\textup)
\end{Proposition}

We call this the {\it contact term equation\/}. 

\begin{NB}
The range of $d$ is wider for the SW prepotential. But it is probably
easy to show that the range can be extended automatically from the
above range to the wider range, thanks to the symmetry of the theta
function and \eqref{eq:SerreDuality}. But it does not {\it a priori\/}
imply that the blowup equation holds in the wider range, i.e.\ we may
not have the equation before taking the limit. I still need to study
the symmetry of the partition function on the blowup.
Oct. 15, H.N.
\end{NB}

As $\Theta_E$ is an even function, the above holds for $d$ if and only
if it holds for $r+m-d$. In particular, the above expression is
independent of $0 \le d \le r+m$ for $m\ge 0$, and
$m\le d\le r$ for $m\le 0$.%
\begin{NB}
  Added on Oct. 16, H.N. It seems that this is not related to the
  Serre duality directly. Anyway we do not know this holds for the
  original blowup equation. Kota's numerical calculation check this
  for $r=2$, $m=\pm 1$.
\end{NB}

We will prove that the Seiberg-Witten prepotential defined via the
periods of hyperelliptic curves satisfies the same equation and has
the same perturbation part in \secref{sec:SWcurve}.
As the contact term equation determines
the instanton part of the prepotential recursively from the
perturbation part, we get
\begin{Theorem}
  Suppose \eqref{eq:conjvanish} holds and assume $|m| < r$. Then $\F_0$
  coincides with the Seiberg-Witten prepotential defined in
  \eqref{eq:SWprep}.
\end{Theorem}

As we have \eqref{eq:conjvanish} for the case $m=0$, we have the
assertion without the condition in this case. This is the proof of
Nekrasov's conjecture for the $K$-theoretic partition function
\cite{Nek}. See \cite{NO} for another proof.

By \propref{prop:regular'}(1) the next coefficient
$H(\vec{a};\Lambda,\bbeta)$ of the
expansion~\eqref{eq:genus_expansion} comes from the perturbation part:

\begin{Proposition}
Suppose \eqref{eq:conjvanish} holds and assume $|m| < r$.
\begin{equation*}
   H(\vec{a};\Lambda,\bbeta) = -\pi\sqrt{-1}\langle\vec{a},\rho\rangle.
\end{equation*}
\end{Proposition}
\begin{NB}
The term
\begin{equation*}
   - \frac{m \bbeta}4 \sum_{\alpha=1}^r a_\alpha^2.
\end{equation*}
disappears as we change the definition of the partition function and
the perturbation part.
Oct. 15, H.N.
\end{NB}

\subsubsection{Genus $1$ parts}

Next we turn to the genus $1$ parts of the
expansion~\eqref{eq:genus_expansion}. When $r=2$, $m=0$, we determined
$A$, $B$ explicitly as theta constants in \cite{NY3}. So we assume
$r=2$, $m=1$.
Let $F_1 = A - \frac23 B$, $G = \frac13 B$.

We have \cite[(4.11)]{NY3} if we replace $d$ by $d-\frac{m}2$. (Note
that this is $k=0$ case.)
Taking
$d=1$ (and $r=2$, $m=1$) we have
\begin{equation}\label{eq:genus1}
\exp(G-F_1)
=
\begin{NB}
\exp(B - A) =
\end{NB}
\exp\left[  -\frac{\bbeta^2}{128}
  \frac{\partial^2\F_0}{(\partial\log\Lambda)^2}\right]
\theta_{01}\left(\left.
  \frac{\bbeta}{16\pi\sqrt{-1}}
  \frac{\partial\F_0}{\partial\log\Lambda\partial a}\right|
  \tau
\right).
\end{equation}

We assume
\begin{equation}\label{eq:conjvanish2}
   \bZin_{m,k,d}(\ve_1,\ve_2,\vec{a};\Lambda,\bbeta) = 0
\end{equation}
for $0 < k < r$, $0 < d < r$.%
\begin{NB}
  Then we have \cite[the last displayed equation in p.514]{NY3} if we
  replace $d$ by $d+m(-1/2+k/r)$. We take $k=1$, $d=1$ (and $r=2$,
  $m=1$). Then we have
\begin{equation*}
   0 = \theta_{11}(0|\tau).
\end{equation*}
This is of course a trivial identity.
\end{NB}
Then we have \cite[the first displayed equation in p.515]{NY3} if we
replace $d$ by $d+m(-1/2+k/r)$. We take $k=1$, $d=1$ (and $r=2$,
$m=1$). Then we have exactly the same equation as in
\cite{NY3}. Therefore we get
\begin{equation*}
  G + F_1 = -\frac13 \log\left(
    -2\pi q^{\frac{1}{8}}\prod_{d=1}^{\infty}(1-q^{d})^3 \right)+C
\end{equation*}
where $C$ is a function on $\Lambda$. Here $q =
\exp(2\pi\sqrt{-1}\tau) = \exp(-d^2 F/da^2)$ and the convention is
different from that in \cite{NY3}.
Combining with \eqref{eq:genus1}, we get
\begin{equation*}
   \exp F_1 = C' q^{-1/48} \prod_{d=1}^\infty (1 - q^d)^{-1/2}
   \exp\left[ \frac{\bbeta^2}{256}
  \frac{\partial^2\F_0}{(\partial\log\Lambda)^2}\right]
  \theta_{01}\left(\left.
      \frac{\bbeta}{16\pi\sqrt{-1}}
      \frac{\partial\F_0}{\partial\log\Lambda\partial a}\right|
    \tau
  \right)^{-1/2}
\end{equation*}
for $C' = C'(\Lambda)$. By the same argument in \cite[p.515]{NY3} we
have $C' \equiv 1$. 
Let us briefly recall the argument and explain how it is modified in
our case. The proof is based on the observation that $\eta(\tau/2)\exp
F_1$ depends on $\Lambda$ in the form $\C[[\zeta_{1,2}\Lambda^{4}]]$,
where $\zeta_{1,2} = \frac{\bbeta}{1-e^{2\bbeta a}}$ (see
\subsecref{subsec:pm1}).%
\begin{NB}
({\bf Correction} to \cite{NY3}: Replace
$\prod_{\alpha<\beta}\zeta_{\alpha,\beta}\Lambda$ by
$(\prod_{\alpha<\beta}\zeta_{\alpha,\beta}) \Lambda^{2r}$ in two line
above of (5.3), (5.3), and (5.5). Note also that $\Fin_K$ in
\cite{NY3} is $\ve_1\ve_2\Fin_K$ here.)
\end{NB}
There is an extra factor $\exp(m \bbeta a(|Y^2|-|Y^1|))$ coming from
the Chern-Simons terms. Hence the coefficient of $\Lambda^{4n}$ is
divisible by $\exp({mn}\bbeta a) \zeta_{1,2}^n$.  Under our assumption
$m=1$, we cannot get a term which is constant with respect to $a$.
Therefore
\begin{equation}\label{eq:AB}
\begin{split}
   \exp F_1 & = q^{-1/48} \prod_{d=1}^\infty (1 - q^d)^{-1/2}
   \exp\left[ \frac{\bbeta^2}{256}
  \frac{\partial^2\F_0}{(\partial\log\Lambda)^2}\right]
  \theta_{01}\left(\left.
      \frac{\bbeta}{16\pi\sqrt{-1}}
      \frac{\partial\F_0}{\partial\log\Lambda\partial a}\right|
    \tau
  \right)^{-1/2},
\\
   \exp G &= q^{-1/48} \prod_{d=1}^\infty (1 - q^d)^{-1/2}
   \exp\left[ -\frac{\bbeta^2}{256}
  \frac{\partial^2\F_0}{(\partial\log\Lambda)^2}\right]
  \theta_{01}\left(\left.
      \frac{\bbeta}{16\pi\sqrt{-1}}
      \frac{\partial\F_0}{\partial\log\Lambda\partial a}\right|
    \tau
  \right)^{1/2}.
\end{split}
\end{equation}

\section{Computation of the wallcrossing in terms of Hilbert schemes}\label{sec:Hilb}

Let $X$ be a simply connected 
smooth projective  surface with $p_g=0$. 
In this section we will compute the wallcrossing of the $K$-theoretic Donaldson invariants of $X$
in terms of the holomorphic Euler characteristic of certain sheaves on Hilbert schemes of points on $X$. 
Later we will specialize to the case that $X$ is a
smooth toric surface and relate this result to the $K$-theoretic Nekrasov partition
function.

\begin{Notation}
Let $t$ be a variable. If $Y$ is a variety and 
$b\in H^*(Y)[t]$, we denote
by
$[b]_d$ its part of degree $d$, where elements in 
$H^{2n}(Y)$ have degree $n$ 
and $t$ has degree $1$.

If $R$ is a ring, $t$ a variable and $b\in R((t))$, we will denote for 
$i\in \Z$ by
$[b]_{t^i}$ the coefficient of $t^i$ of $b$.

If $E$ is a vector bundle of rank $r$ on $Y$, let $\bigwedge_{-t}E:=
\sum_{i}(-1)^i \Lambda^i (E) t^i\in K(Y)[t]$, and let 
$S_t(E):=\sum_{i}S^i(E)t^i$, where $S^i(E)$ is the $i^{\mathrm{th}}$ symmetric power of $E$.
Note that $S_t(E)=\frac{1}{\bigwedge_{-t}(E)}$.
\end{Notation}

\subsection{The wallcrossing term}\label{wallmod}
Denote by $\cc$ the ample cone of $X$. Then $\cc$ has a chamber structure:
For a class $\xi\in H^2(X,\Z)\setminus \{0\}$ let $W^\xi:=\big\{ x\in \cc\bigm| \< x,\xi\>=0\big\}$.
Assume $W^\xi\ne \emptyset $. Then we call $\xi$ a {\it class of type} $(c_1,d)$ and call 
$W^\xi$ a {\it wall of type} $(c_1,d)$ if the following conditions hold
\begin{enumerate}
\item 
$\xi+c_1$ is divisible by $2$ in $H^2(X,\Z)$,
\item $d+3+\xi^2\ge 0$.
\end{enumerate}
We call $\xi$ a {\it class of type} $c_1$, if $\xi+c_1$ is divisible by $2$ in $H^2(X,\Z)$.
The {\it chambers of type} $(c_1,d)$ are the connected components of the complement
of the walls of type $(c_1,d)$ in $\cc$.
Then $M_H^X(c_1,d)$ depends only on the chamber of type $(c_1,d)$ of $H$.

Let $\xi\in  H^2(X,\Z)$ by a class of type $c_1$. 
We say that $\xi$  {\it good\/} and $W^\xi$ is a {\it good wall\/} if
 $D+K_X$ is not effective for any divisor $D$ with $W^{c_1(D)}=W^\xi$.
A sufficient condition for $\xi$ to be good  is that $W^\xi$ contains an 
ample divisor $H$ with $H\cdot K_X<0$. 
One can show that an ample divisor $H$ is general with respect to 
$(2,c_1,c_2)$ if and only if $H$ lies in a chamber of type  $(c_1,4c_2-c_1^2-3)$.

Let $\xi$ be a  class of type $c_1$.
Let $X^{[n]}$ be the Hilbert scheme of subschemes of length $n$ on $X$.
Let $Z_n(X)\subset X\times X^{[n]}$ be the universal subscheme. 
Let  $\I_1$ (resp.~$\I_2$) be the sheaf $p_{1,2}^*(\I_{Z_n(X)})$ (resp.~
$p_{1,3}^*(\I_{Z_m(X)})$ on $X\times X^{[n]}\times X^{[m]}$.
We also denote $\F_1:=\I_1(\frac{c_1+\xi}{2})$ and
$\F_2:=\I_{Z_2}(\frac{c_1-\xi}2,)$. Note that $X^{[n]}=M^X_H(1,\frac{c_1+\xi}2,n)$
and $X^{[m]}=M^X_H(1,\frac{c_1-\xi}2,m)$ and  $\F_1$, $\F_2$ are the corresponding universal sheaves. Let $f_1,f_2\in K(X)$ be the classes of elements of 
$M^X_H(1,\frac{c_1+\xi}2,n)$ and $M^X_H(1,\frac{c_1-\xi}2,m)$ respectively.

Let $p:X\times X^{[n]}\times X^{[m]}\to X^{[n]}\times X^{[m]}$,  
$q:X\times X^{[n]}\times X^{[m]}\to X$ be the projections. Let
$
\AA_{\xi,-}:=-p_!(\I_2^\vee\otimes \I_1\otimes q^!\xi),\
\AA_{\xi,+}:=-p_!(\I_1^\vee\otimes \I_2\otimes q^!\xi^\vee)\in K(X^{[n]}\times X^{[m]})
$.
We also just write 
$\AA_-$ and $\AA_+$ instead
of  $\AA_{\xi,-}$, $\AA_{\xi,+}$. 
\begin{NB} I redefined $\AA_+,\AA_-$ in this way so that it also works
if the wall is not good. LG 27.8.06\end{NB}

Now assume  $\xi$ is good. Then 
$\Ext_p^0(\I_2,\I_1(\xi))=
\Ext_p^2(\I_2,\I_1(\xi))=0$ and we will write $\AA_{\xi,-}$ for its
 representative $\Ext^1_{p}(\I_{2},\I_{1}(\xi))$, which is a  locally free sheaf
 on $X^{[n]}\times X^{[m]}$.
 Similarly we write  $\AA_{\xi,+}$ for the locally free sheaf $\Ext^1_{p}(\I_{1},\I_{2}(-\xi))$, and we put   $\P_-:=\P(\AA_-^\vee)$ and $\P_+:=\P(\AA_+^\vee)$
 (we use the Grothendieck notation, i.e. this is the bundle of $1$-dimensional quotients). Let $\pi_{\pm}:\P_{\pm}
\to X^{[n]}\times X^{[m]}$ be the projection. 

\begin{Definition}
Fix $c_1\in H^2(X,\Z)$, and let $v\in K(X)$. 
Let $\xi\in H^2(X,\Z)$ be a  class of type $c_1$. 
We denote $\chi(f_1\otimes v)=\chi(\I_{Z_1}(\frac{c_1+\xi}2\otimes v))$, 
$\chi(f_2\otimes v)=\chi(\I_{Z_2}(\frac{c_1-\xi}2\otimes v))$
for $(Z_1,Z_2)\in X^{[n]}\times X^{[m]}$.
By the Riemann-Roch-Theorem we see that
\begin{equation}\label{eq:f2ud}
\frac{1}{2}(\chi(f_2\otimes v)-\chi(f_1\otimes v))=
-\Bigl\<\frac{\xi}{2},v^{(1)}\Bigr\>+\frac{\rk(v)}{2}(n-m).
\end{equation} 
In particular it  only depends on $\rk(v)$ and $c_1(v)$, and it is independent of $n,m$ if
$\rk(v)=0$.

\begin{NB} Strictly speaking we should write $f_1(n,m)$, $f_2(n,m)$
because $f_1,f_2$ depend on $n,m$, but I hope that this will not lead to confusion. What do you think?
(LG 21.8.06)\end{NB}

%
The {\em wallcrossing terms} are
\begin{equation}\label{wallct}
\begin{split}&\Delta^X_{\xi,T}(v;\Lambda):=\sum_{{n,m\ge 0}\atop {d=4(n+m)+\xi^2-3}}\frac{\Lambda^d}{T^{\frac{1}{2}(\chi(f_2\otimes v)-\chi(f_1\otimes v))}}
\chi\Bigl(X^{[n]}\times X^{[m]}, 
\frac{\lambda_{\F_1}(v)\otimes \lambda_{\F_2}(v)}
{\bigwedge_{-T}(\AA^\vee_{\xi,+})
\bigwedge_{-T^{-1}}(\AA^\vee_{\xi,-})}\Bigr),\\
&\Delta^X_{\xi}(v;\Lambda):=[\Delta^X_{\xi,T}(v;\Lambda)]_{T^0}-[\Delta^X_{\xi,T}(v;\Lambda)]_{(T^{-1})^0}.
\end{split}
\end{equation}
Here the right hand side of the first equation is understood as a
rational function in $T^{1/2}$ as follows, and $[\bullet]_{T^0}$,
$[\bullet]_{(T^{-1})^0}$ denote the constant terms of the expansions
at $T^{1/2}=0$, $T^{1/2}=\infty$ respectively. We formally apply the
Hirzebruch-Riemann-Roch theorem to get
\begin{multline*}
\chi\Bigl(X^{[n]}\times X^{[m]}, 
\frac{\lambda_{\F_1}(v)\otimes \lambda_{\F_2}(v)}
{\bigwedge_{-T}(\AA^\vee_{\xi,+})
\bigwedge_{-T^{-1}}(\AA^\vee_{\xi,-})}\Bigr)
\\
= \int_{X^{[n]}\times X^{[m]}}
   \frac{\ch(\lambda_{\F_1}(v)) \ch(\lambda_{\F_2}(v))}
   {\ch \bigwedge_{-T}(\AA^\vee_{\xi,+}) \ch\bigwedge_{-T^{-1}}(\AA^\vee_{\xi,-})}
   \Todd (X^{[n]}\times X^{[m]}).
\end{multline*}
Then we consider 
$\ch \bigwedge_{-T}(\AA^\vee_{\xi,+})$,
$\ch\bigwedge_{-T^{-1}}(\AA^\vee_{\xi,-})$ as
$\End(H^*(X^{[n]}\times X^{[m]}))$-valued Laurent polynomials. Their
inverses are defined as their cofactor matrices divided by their
determinants (which are equal to $(1-T)^{\rk(\AA^\vee_{\xi,+})}$,
$(1-T^{-1})^{\rk(\AA^\vee_{\xi,-})}$) respectively. Then their inverse
are in $\End(H^*(X^{[n]}\times X^{[m]}))\otimes_{\Q} \Q(T)$. Thus the
integral  is an element of $\Q(T)$. This way of  understanding the formula will become
more apparent  when we will consider the equivariant wallcrossing term
in \secref{sec:comparison}. In that case we can interpret the formula so that the computation
is done in the localized equivariant $K$-theory.

\begin{NB}
Differently from the topological case, the wallcrossing terms will also depend 
on $c_1$ and not only on $\xi$, although I expect to prove later that they are
independent of $c_1$ in case $\rk(v)=0$, which is the only  case where one 
can complete the computation. Thus it might be necessary to put $c_1$ also in the
notation (what do you think)?

There is one alternative here: One could always change $v$ formally into 
$v\otimes c_1/2$. This means that we actually do not change anything, but
notationally we pretend that we have tensorized the bundle by $-c_1/2$, and 
$v$ by $c_1/2$. This makes the dependence on $c_1$ clear, and later we will 
see it is also the way that $c_1$ figures in the computations, because in
the Nekrasov partition function $c_1$ is always normalized to be $0$.
\end{NB} 

The expansions at $T=0$, $T=\infty$ can be also understood
differently. Note that for a vector bundle $E$ of rank $r$ we have
\begin{NB}
Slightly changed according to the above paragraph. Nov. 6, HN. Here is
the original:

$\Delta^X_{\xi,T}(v;\Lambda)$ is {\it a priori\/} just a formal expression, however
note that for a vector bundle $E$ of rank $r$ we have
\end{NB}
\begin{equation}
\frac{1}{\bigwedge_{-T}(E)}=S_T(E), \qquad
\frac{1}{\bigwedge_{-T^{-1}}(E)}=\frac{(-T)^r}{ \det(E) \otimes \bigwedge_{-T} E^\vee}=
(-T)^r\det(E^\vee) \otimes S_T(E^\vee).
\end{equation}
Thus $\Delta^X_{\xi,T}(v;\Lambda)$
can be developed as Laurent series both in $T^{1/2}$ and  in $\frac{1}{T^{1/2}}$
\begin{equation}
\label{delT}
\begin{split}
\Delta^X_{\xi,T}(v;\Lambda)=\sum_{{n,m\ge 0}\atop {d=4(n+m)-\xi^2-3}} &
\frac{\Lambda^d (-T)^{\rk(\AA_-)}}{T^{\frac{1}{2}(\chi(f_2\otimes v)-\chi(f_1\otimes v))}}
\chi\Bigl(X^{[n]}\times X^{[m]},\lambda_{\F_1}(v)\otimes \lambda_{\F_2}(v)\otimes \\
&\det(\AA_{\xi,-})\otimes S_T(\AA_{\xi,+}^\vee)\otimes
S_T(\AA_{\xi,-})\Bigr)\in \Z((T^{\frac{1}{2}}))[[\Lambda]],\\
\end{split}
\end{equation}
\begin{equation}
\label{delT1}
\begin{split}
\Delta^X_{\xi,T}(v;\Lambda)=\sum_{{n,m\ge 0}\atop {d=4(n+m)-\xi^2-3}} &
\frac{\Lambda^d (-T^{-1})^{\rk(\AA_+)}}{(T^{-1})^{\frac{1}{2}(\chi(f_1\otimes v)-\chi(f_2\otimes v))}}
\chi\Bigl(X^{[n]}\times X^{[m]},\lambda_{\F_1}(v)\otimes \lambda_{\F_2}(v)\\
&\otimes\det(\AA_{\xi,+}) \otimes S_{T^{-1}}(\AA_{\xi,+})\otimes
S_{T^{-1}}(\AA_{\xi,-}^\vee)\Bigr) \in\Z((T^{-\frac{1}{2}}))[[\Lambda]].
\end{split}
\end{equation}
Then
$[\Delta^X_{\xi,T}(v;\Lambda)]_{T^0}$ is equal to the coefficient of $T^0$ of \eqref{delT} and
$[\Delta^X_{\xi,T}(v;\Lambda)]_{(T^{-1})^0}$ is equal to the
coefficient of $(\frac{1}{T})^0$ of \eqref{delT1}.
However note that it was not clear the expressions are
in $\Z((T))[[\Lambda]]$ or $\Z((T^{-1}))[[\Lambda]]$ in the original
formulation in terms of the Hirzebruch-Riemann-Roch theorem.
\end{Definition}

\begin{NB}
As you see I removed the conditions on $c_1(v),\rk(v)$. I add a remark
in which I explain that in case $v\in K_c$ I get the same as before and 
that otherwise I usually get $0$.

I also made the definition in case the wall is not good
 (21.8.06 LG)
\end{NB}

\begin{Remark}\label{integral}
Fix $c_1,d$ and let $c\in K(X)$ be the class of an element in $M_H^X(c_1,d)$.
 Let $v\in K(X)$.
Then either $\frac{1}{2}(\chi(f_2\otimes v)-\chi(f_1\otimes v))\in \Z$ for all
$n,m\in \Z_{\ge 0}$ with $4(n+m)-\xi^2-3=d$, or $\frac{1}{2}(\chi(f_2\otimes v)-\chi(f_1\otimes 0))\in \Z+\frac{1}{2}$ for all such
$n,m$.
In the second case the coefficients of $\Lambda^d$ of 
$[\Delta_{\xi,T}(v,\Lambda)]_{T^0}$ and $[\Delta_{\xi,T}(v,\Lambda)]_{(T^{-1})^0}$ are trivially $0$.

On the other hand, if $v\in K_c$, then $\chi(f_2\otimes v)=-\chi(f_1\otimes v)$
and thus $\frac{1}{2}(\chi(f_2\otimes v)-\chi(f_1\otimes v))
=\chi(f_2\otimes v)\in \Z$.
\end{Remark}


\begin{Remark}\label{rem:polynom}
Let $v\in K(X)$ be a class of rank $0$.
Let $\xi$ be a  wall of type $(c_1,d)$. Let $l=\frac{d+3+\xi^2}{4}$.
Fix $l\ge 0$. Write $d:=4l-\xi^2-3$. Note that by \cite[Lemma 4.3]{EG1}
\begin{equation}\label{eq:r} \rk(\AA_-)=-\frac{\xi(\xi-K_X)}{2}+l-1,\quad 
\rk(\AA_+)=-\frac{\xi(\xi+K_X)}{2}+l-1.
\end{equation}
Note that by definition the coefficient of $\Lambda^d$ of $\Delta^X_{\xi}(v;\Lambda)$
is zero if $-\rk(\AA_+)<-\<\xi/2,c_1(v)\> <\rk(\AA_-)$. By \eqref{eq:r}
it thus follows that 
that the coefficient of $\Lambda^d$ of $\Delta^X_{\xi}(L;\Lambda)$ is $0$ unless
$0\le l\le |\< \frac{\xi}{2},c_1(v)+K_X\>|+1+\frac{\xi^2}{2}$, which is equivalent to
$-\xi^2-3\le d\le \xi^2+|\<2\xi,c_1(v)+K_X\>|+1$.  In particular  $\Delta^X_{\xi}(v;\Lambda)\in \C[\Lambda]$.
\end{Remark}

The aim of this section is to prove that the wallcrossing for the 
$K$-theoretic Donaldson invariants can be expressed as a sum over $\Delta^X_{\xi}(v;\Lambda)$.

\begin{Proposition}\label{wallcr} Fix $c_1,d$, let $c\in K(X)$ be the class of 
an element of $M^X_H(c_1,d)$.
Let $v\in K_c$.
Let $H_-$, $H_+$ be ample divisors on $X$, which do not lie on a wall 
of type $(c_1,d)$. Let $B_+$ be the set of classes $\xi$ of type $(c_1,d)$
with $\<\xi\cdot H_+\>>0 >\<\xi\cdot H_-\>$. Assume that all classes in $B_+$ 
are good. 
Then  
\begin{align*}
\chi(M^X_{H_+}(c_1,d),\lambda(v))-\chi(M^X_{H_-}(c_1,d),\lambda(v))&=\sum_{\xi\in B_+}\bigl[\Delta^X_{\xi}(v;\Lambda)\bigr]_{\Lambda^d}.
\end{align*}
\end{Proposition}

\begin{NB} Below is the old version with $u_d$.

\begin{Proposition}\label{wallcr}
Let $v\in K(X)$.
Let $H_-$, $H_+$ be ample divisors on $X$, which do not lie on a wall 
of type $c_1$. Let $B_+$ be the set of all 
classes $\xi$ of type $c_1$
with $\<\xi\cdot H_+\>>0 >\<\xi\cdot H_-\>$. Assume that all classes in $B_+$ 
are good. 
Then  
\begin{align*}
\chi_{c_1}^{H_+}(v; \Lambda)-\chi_{c_1}^{H_-}(v; \Lambda)&=\sum_{\xi\in B_+}\Delta^X_{\xi}(v;\Lambda).
\end{align*}
\end{Proposition}
\end{NB}


In the rest of this section we will show Prop.~\ref{wallcr}. 

$M^X_H(c_1,d)$ and thus $\chi(M_H^X(c_1,d),\lambda(v))$ is constant as long as $H$ stays in the same  chamber of type $(c_1,d)$ 
and only changes when
$H$ crosses a wall of type $(c_1,d)$. 
By \cite{EG1}, \cite{FQ} the change of the moduli spaces can be described as follows.
Let $B_d$ be the set of all $\xi\in B_+$ which define a wall of type $(c_1,d)$.
For the moment assume for simplicity that $B_d$ consists of a single element
$\xi$.
Let $l:=(d+3+\xi^2)/4\in \Z_{\ge 0}$.
Write $M_{0,l}:=M^X_{H_-}(c_1,d)$.  Then successively for all $n=0,\ldots,l$ write 
$m:=l-n$. Then one has the following: 
$M_{n,m}$ contains a closed subscheme $E_-^{n,m}$ isomorphic to $\P_{-}^{n,m}$
and $M_{n,m}$ is nonsingular in a neighbourhood of $E_-^{n,m}$.  
Let $\widehat M_{n,m}$ be the blow up of $M_{n,m}$  along $E_-^{n,m}$.
The exceptional divisor is isomorphic to the fibre product  
$D^{n,m}:=\P_-^{n,m}\times_{X^{[n]}\times X^{[m]}} 
\P_{+}^{n,m}$. We can blow down $\widehat M_{n,m}$  in $D^{n,m}$ 
in the other fibre direction to obtain a new variety  $M_{n+1,m-1}$. 
The image of $D^{n,m}$ is a closed subset $E_+^{n,m}$ isomorphic to 
$\P_+^{n,m}$ and $M_{n+1,m-1}$ is smooth in a neighbourhood of $E_+^{n,m}$.

The transformation from $M_{n,m}$ to $M_{n+1,m-1}$  does not have to be  birational. 
It is possible that $E_+^{n,m}=\emptyset$, i.e. $\AA_+=0$. 
As $\rk(\AA_-)+\rk(\AA_+)+2l=d+1$, this happens if and only if $E_-^{n,m}$ has dimension $d$
and thus by the smoothness of  $M_{n,m}$ near $E_-^{n,m}$, we get that $E_-^{n,m}$
is a connected component of  $M_{n,m}$. Then blowing up along $E_-^{n,m}$ just means
deleting $E_-^{n,m}$. Thus in this case $M_{n+1,m-1}=M_{n,m}\setminus E_-^{n,m}$.
Similarly  we have $E_-^{n,m}=\emptyset$, i.e. $\AA_-=0$, if and only if 
$E_+^{n,m}$ is a 
connected component of $M_{n+1,m-1}$ and  $M_{n+1,m-1}=M_{n,m}\sqcup E_+^{n,m}$.
Below, if the transformation from $M_{n,m}$ to $M_{n+1,m-1}$ is birational, we say 
we are in case (1), otherwise in case (2).
Finally we have $M_{l+1,-1}=M_{H_+}^X(c_1,d)$. 
If $B_d$ consists of more than one element, one obtains $M_{H_+}(c_1,d)$
from $M_{H_-}(c_1,d)$  by iterating  this procedure in a suitable order over  all
$\xi\in B_+$.

Fix $\xi$ in $B_d$. 
Fix $n,m\in \Z_{\ge 0}$ with $n+m=l:=(d+3+\xi^2)/4$. 
We write $M_-:=M_{n,m}$, $M_+:=M_{n+1,m-1}$. Let $\overline \E_\pm$ be 
universal sheaves on $X\times M_\pm$ respectively.
Let $E_-:=E_-^{n,m}$, $E_+=E_+^{n,m}$.
Let $\widetilde M$ be the blowup of $M_-$ along $E_-$, and denote
by $D$ the exceptional divisor (which is also the exceptional divisor 
of the blowup of $M_+$ along $E_+$). Write $D':=X\times D$ and let $j:D\to \widetilde M$,
$j':X\times D\to X\times \widetilde M$ be the embeddings.
Let $\E_-$, $\E_+$ be the pullbacks of $\overline \E_-$, $\overline 
\E_+$ to $X\times \widetilde M$.

\begin{Notation} 
We denote by $T_-$ (resp.\ $T_+$) the universal quotient line bundle on 
$\P_-=\P(\AA_-^\vee)$ (resp. on $\P_+=\P(\AA_+^\vee)$).
For a class $a\in H^*(X)$ we also denote by $a$ its pullback to $X\times Y$
for a variety $Y$.
We write $\I_1,\I_2$ also for the pullback of $\I_1$, $\I_2$ to $D'$ and we write
$T_+,T_-$ also for their pullbacks to $D$ and $D'$.
\end{Notation}

By the condition $\chi(c\otimes v)=0$, we can 
replace $\frac{1}{2}(\chi(f_2\otimes v)-\chi(f_1\otimes v))$ by 
$\chi(f_2\otimes v)$.
We will show
\begin{equation}
\label{restr}
\begin{split}
&\chi(M_+,\lambda_{\overline \E_+}(v))-\chi(M_-,\lambda_{\overline \E_-}(v))
=\chi\Biggl(X^{[n]}\times X^{[m]},\lambda_{\F_1}(v)\otimes \lambda_{\F_2}(v)\\
&\otimes\Big(\Big[\frac{(-t)^{\rk(\AA_-)}S_{t}(\AA_+^\vee)\otimes S_{t}(\AA_-)\otimes \det(\AA_-)}{t^{\chi(f_2\otimes v)}}
\Big]_{t^0}-\Big[\frac{(-t)^{\rk(\AA_+)}
S_{t}(\AA_-^\vee)\otimes S_t(\AA_+)\otimes \det(\AA_+)}{t^{-\chi(f_2\otimes v)}}
\Big]_{t^0}\Big)\Biggr).
\end{split}
\end{equation}
Formula (\ref{restr}) implies Proposition \ref{wallcr} by summing over all 
$\xi\in B_+$, 
and over all $n,m$ with $n+m=(d+\xi^2+3)/4$.

Assume first that we are in case (1). Note that this is equivalent to 
both $\rk(\AA_-)$ and $\rk(\AA_+)$ strictly positive, 
and then it is evident that the first (resp. second)
summand on the left hand side of \eqref{restr} vanishes if $\chi(f_2\otimes v)\le 0$
(resp. if $\chi(f_2\otimes v)\ge 0$).
Let $\pi_{\pm}:\widetilde M\to M_{\pm}$  be the blowup morphisms.
By \cite[Prop.~VI.4.1]{FuLa} and its proof,  $R^i\pi_{\pm*} \oo_{\widetilde M}=
0$ for $i>0$, and  $\pi_{\pm*}\oo_{\widetilde M}=\oo_{M_{\pm}}$.
Thus  the  projection formula  gives
$\chi(M_{\pm},L)=\chi(\widetilde M,\pi_{\pm}^*L)$ for any line bundle $L$ on  
$M_{\pm}$.
Therefore it is enough to prove (\ref{restr})
with the left-hand side replaced by 
$\chi(\widetilde M,\lambda_{\E_+}(v))-\chi(\widetilde M,\lambda_{\E_-}(v)).$

\begin{Lemma}\label{lamD} In $K(\widetilde M)$ we have
$$\lambda_{\E_+}(v)-\lambda_{\E_-}(v)=
j_*\Big(\Big(\frac{t^{\chi(f_2\otimes v)}-s^{-\chi(f_2\otimes v)}}{1-st}\pi^*(\lambda_{\F_1}(v)\otimes \lambda_{\F_2}(v)\Big)|_{{s=T_-}\atop{t=T_+}}\Big).$$
\end{Lemma}
\begin{proof}
By \cite[section 5]{EG1} there exists a line bundle $\mu$ on $D$ such that
\begin{equation}
\label{EF}
(j')^*\E_-=\F_1\otimes \mu+\F_2\otimes T_-^{-1}\otimes \mu\hbox{ in $K(D')$},
\quad 
\E_+=\E_--j'_*(\F_2\otimes T_-^{-1}\otimes \mu) \hbox{ in  $K(X\times \widetilde M)$}.
\end{equation}
Thus we get 
$\lambda_{\E_+}(v)=\lambda_{\E_-}(v)\otimes
\det\bigl(p_!(v\otimes j'_*(\F_2\otimes T_-^{-1}\otimes \mu))\bigr)^{-1}$.
Note that $p_!(v\otimes j'_*(\F_2\otimes T_-^{-1}\otimes \mu))$ is a coherent sheaf
of rank $\chi(v\otimes f_2)$ on $D$. As $D$ is a Cartier divisor, it follows
that 
$$\lambda_{\E_+}(v)=\lambda_{\E_-}(v)\otimes
\det(\chi(f_2\otimes v) [\oo_D])^{-1}=\lambda_{\E_-}(v)\otimes
\det(\oo_{\widetilde M}(-D))^{\chi(f_2\otimes v)}.$$
For the second equality we have used that $\oo_D=\oo_{\widetilde M}-\oo_{\widetilde M}(-D)$ in $K(\widetilde M)$ and thus $\det(\oo_D)=\det(\oo_{\widetilde M}(-D))^{-1}$.
Thus we get in $K(\widetilde M)$ that 
\begin{align*}
\lambda_{\E_+}(v)-\lambda_{\E_-}(v)&=(\oo_{\widetilde M}(-D)^{\chi(f_2\otimes v)}-1)\otimes 
\lambda_{\E_-}(v)
=j_*\Big(\Big(\frac{t^{\chi(f_2\otimes v)}-1}{1-t}j^*(\lambda_{\E_-}(v))\Big)|_{t=T_+\otimes T_-}\Big).\end{align*}
In the last step we have used that for a locally free sheaf $\G$ on 
$\widetilde M$ we have $(1-\oo_{\widetilde M}(-D))\otimes\G =j_*(j^*\G)$ in 
$K(\widetilde M)$.
As the determinant bundles are compatible with pullback 
we obtain by \eqref{EF} that 
\begin{equation}\label{eq:int}
\begin{split}
j^*(\lambda_{\E_-}(v))&=
\lambda_{\pi^*(\F_1)\otimes \mu}(v)\otimes \lambda_{\pi^*(\F_2)\otimes T_-^{-1}\otimes \mu}(v)\\
&=
\pi^*(\lambda_{\F_1}(v)\otimes \lambda_{\F_2}(v))\otimes \mu^{\chi(f_1\otimes v)+\chi(f_2\otimes v)} T_-^{-\chi(f_2\otimes v)}\\
&=
\pi^*(\lambda_{\F_1}(v)\otimes \lambda_{\F_2}(v))\otimes T_-^{-\chi(f_2\otimes v)}.
\end{split}\end{equation}
In the last step we use that $\chi(f_1\otimes v)+\chi(f_2\otimes v)=0$. 
The result follows.
\end{proof}

In case $\chi(f_2\otimes v)=0$ the left hand side of Lemma \ref{lamD}  is obviously $0$.
Thus we only need to show \eqref{restr} in the cases $\chi(f_2\otimes v)> 0$
and $\chi(f_2\otimes v)<0$.

{\bf (a)} $\chi(f_2\otimes v)< 0$:
We apply the formula
\begin{equation}\label{eq:rat}
\frac{y^{-c}-x^{c}}{1-xy}=y^{-1}\frac{x^{c}-y^{-c}}{x-y^{-1}}=\sum_{{a+b=c}\atop
{a\ge 0,\ b>0}} x^ay^{-b},\quad c\in \Z_{>0}
\end{equation} 
for $x=T_-$, $y=T_+$ to Lemma \ref{lamD} to obtain
\begin{align*}\lambda_{\E_+}(v)-\lambda_{\E_-}(v)&=
j_*\Big(\sum_{{a+b=-\chi(f_2\otimes v)}\atop {a\ge 0,\ b>0}}T_-^{ a}\otimes 
T_+^{-b}\otimes \pi^*(\lambda_{\F_1}(v)\otimes \lambda_{\F_2}(v))\Big)
\end{align*} 
in $K(\widetilde M)$.
Let $\E$ be a vector bundle of rank $e$ on a variety $Y$ and let 
$p:\P(\E^\vee)\to Y$ be the projection and $\oo(1)$ the universal quotient line 
bundle on $\P(\E^\vee)$. Then by \cite[Ex.~III.8.4]{Ha}
$$p_!(\oo(n))=\begin{cases} S^n(\E^\vee) & n\ge 0,\\(-1)^{e-1} S^{-n-e}(\E)\otimes \det(\E)&
n\le -e,\\ 0& \hbox{otherwise.}\end{cases}$$
Let $\pi:D=\P(\AA_+^\vee)\times_{X^{[n]}\times X^{[m]}}\P(\AA_+^\vee)\to X^{[n]}
\times X^{[m]}$  be the projection. Then we get using the projection formula
\begin{equation}\label{eq:SST}
\begin{split}
&\chi(\widetilde M,\lambda_{\E_+}(v))-\chi(\widetilde M,\lambda_{\E_-}(v))=
\sum_{{a+b=-\chi(f_2\otimes v)}\atop {a\ge 0,\ b>0}}\chi\big(D,T_-^{ a}\otimes 
T_+^{-b}\otimes \pi^*(\lambda_{\F_1}(v)\otimes \lambda_{\F_2}(v))\big)\\
&=-\sum_{{a+b=-\chi(f_2\otimes v)}\atop {a\ge 0,\ b>0}}
(-1)^{\rk(\AA_+)}\chi\Bigl(X^{[n]}\times X^{[m]}, \lambda_{\F_1}(v)\otimes \lambda_{\F_2}(v)\\
&\qquad \qquad \qquad\qquad\qquad\otimes S^a(\AA_-^\vee)\otimes S^{b-\rk(\AA_+)}(\AA_+) \otimes\det(\AA_+)
\Bigr)\\
&=\chi\Bigl(X^{[n]}\times X^{[m]},\lambda_{\F_1}(v)\otimes \lambda_{\F_2}(v)\otimes\Big[-\frac{(-t)^{\rk(\AA_+)}\otimes S_{t}(\AA_-^\vee)\otimes S_{t}(\AA_+)\otimes 
\det(\AA_+)}{t^{-\chi(f_2\otimes v)}}\Big]_{t^0}\Bigr).
\end{split}\end{equation}

{\bf (b)} $\chi(f_2\otimes v)> 0$:
The formula \eqref{eq:rat}
for $x=T_+$, $y=T_-$, gives
\begin{align*}\lambda_{\E_+}(v)-\lambda_{\E_-}(v)&=-
j_*\Big(\sum_{{a+b=\chi(f_2\otimes v)}\atop {a\ge 0,\ b>0}}T_+^{ a}\otimes 
T_-^{-b}\otimes \pi^*(\lambda_{\F_1}(v)\otimes \lambda_{\F_2}(v))\Big).
\end{align*} 
Then the same arguments as in the case $\chi(f_2\otimes v)< 0$
show that 
\begin{align*}
&-\sum_{{a+b=\chi(f_2\otimes v)}
\atop {a\ge 0,\ b >0}}
\chi\Big(D, \pi^*(\lambda_{\F_1}(v)\otimes \lambda_{\F_2}(v))
\otimes T_-^{ a}\otimes T_+^{ -b})\Big)\\
&=
\chi(X^{[n]}\times X^{[m]},\lambda_{\F_1}(v)\otimes \lambda_{\F_2}(v)
\otimes\Big[\frac{(-t)^{rk(\AA_-)}\otimes S_{t}(\AA_+^\vee)\otimes S_{t}(\AA_-)\otimes \det(\AA_-)\big)}
{t^{\chi(f_2\otimes v)}}
\Big]_{t^0}\Big).
\end{align*}

In case (2), we can assume by symmetry that $\P_+=\emptyset$, thus $\AA_{+}=0$ 
and $\AA_-$ has rank $d+1-2(n+m)$. 
Then we have 
\begin{align*}
\chi(M_+, \lambda_{\overline \E_+}(v))-\chi(M_-, \lambda_{\overline 
\E_-}(v))
=-\chi(\P_-,\lambda_{(j')^*(\overline \E_-)}(v))
\end{align*}
where  $j':X\times \P_-\to X\times M_-$ is the inclusion. 
The same argument as in the proof of \eqref{eq:int} shows that 
\begin{align*}
\lambda_{(j')^*(\overline \E_-)}(v)=
\pi_-^*(\lambda_{\F_1(v)}\otimes \lambda_{\F_2}(v))\otimes T_-^{ -\chi(f_2\otimes v)}.
\end{align*}
Note that differently from case (1) this  is not zero when $\chi(f_2\otimes v)=0$.
Now the same arguments as in the proof of \eqref{eq:SST}
show that 
$-\chi\bigl(\P_-,\lambda_{(j')^*(\overline \E_-)}(v)\bigr)$ is equal to 
$$-
\chi\bigl(X^{[n]}\times X^{[m]}, \lambda_{\F_1}(v)\otimes \lambda_{\F_2}(v)\otimes 
S^{-\chi(f_2\otimes v)}(\AA_-^\vee))
$$
in case  $\chi(f_2\otimes v)\le 0,$ and to 
$$(-1)^{\rk(\AA_-)}\chi(X^{[n]}\times X^{[m]},
\lambda_{\F_1}(v)\otimes \lambda_{\F_2}(v)\otimes S^{\chi(f_2\otimes v)-\rk(\AA_-)}(\AA_-)\otimes 
\det(\AA_-))$$ in case $\chi(f_2\otimes v)>0$.
As $S_t(\AA_+)=1$, this shows \eqref{restr} also  in  case (2) and thus finishes the 
proof of Proposition \ref{wallcr}.

\section{Comparison with the partition function}\label{sec:comparison}

For the next two sections (except in \subsecref{subsec:nontoric}) let $X$ be a smooth
projective toric surface over $\C$, in particular $X$ is simply
connected and $p_g(X)=0$.  $X$ carries an action of
$\Gamma:=\C^*\times \C^*$ with finitely many fixpoints, which we will
denote by $p_1,\ldots,p_\chi$, where $\chi$ is the Euler number of
$X$.  Let $w(x_i)$, $w(y_i)$ the weights of the $\Gamma$-action on
$T_{p_i}X$. Then there are local coordinates $x_i,y_i$ at $p_i$, so
that $(t_1,t_2) x_i=e^{-w(x_i)}x_i$ and $(t_1,t_2)
y_i=e^{-w(y_i)}y_i$.%
\begin{NB} changed 21.8. LG, please check.
\end{NB}
By definition $w(x_i)$ and $w(y_i)$
are linear forms in $\ve_1$ and $\ve_2$. For $\beta\in H^*_\Gamma(X)$ or $\beta\in H_*^\Gamma(X)$, we denote by
$\iota_{p_i}^*\beta$ its  pullback to the fixpoint $p_i$.
More generally, if $\Gamma$ acts on a nonsingular variety $Y$ and $W\subset Y$ is 
invariant under the $\Gamma$-action, we denote by 
$\iota_W^*:H^*_\Gamma(Y)\to H^*_\Gamma(W)$ the pullback homomorphism.

Note that $T_X$ and the canonical bundle  are canonically equivariant.
Thus any polynomial in the Chern classes $c_i(X)$ and $K_X$ is canonically an element
of $H^*_\Gamma(X)$.  

\subsection{Equivariant $K$-theoretic Donaldson invariants and equivariant wallcrossing}
For $t\in \Gamma$ denote by $F_t$ the automorphism $X\to X; x\mapsto t\cdot x$.
Then $\Gamma$ acts on $X^{[n]}\times X^{[m]}$ by $t\cdot(\I_{Y_1},\I_{Y_2})=((F_{t}^{-1})^*\I_{Y_1},(F_{t}^{-1})^*\I_{Y_2})$ and on 
$X\times X^{[n]}\times X^{[m]}$%
\begin{NB} corrected 2.10.06, 30.10 LG\end{NB}
 by $t\cdot(x,\I_{Y_1},\I_{Y_2})=(F_t(x),(F_{t}^{-1})^*\I_{Y_1},
(F_{t}^{-1})^*\I_{Y_2})$   and the sheaves $\I_{1}$,
$\I_{2}$ are $\Gamma$-equivariant. If we choose an equivariant lifting of $c_1$ and $\xi$, then also $\F_1$, $\F_2$ are $\Gamma$-equivariant sheaves.

We write $X_2:=X\sqcup X$ and $X_2^{[l]}:=\coprod_{n+m=l} X^{[n]}\times X^{[m]}$.
The fixpoints of the $\Gamma$-action on $X_2^{[l]}$ are the pairs 
$(Z_1,Z_2)$ of zero-dimensional subschemes with support in 
$\{p_1,\ldots,p_\chi\}$ with $\len(Z_1)+\len(Z_2)=l$ and such that each
$I_{Z_\alpha,p_i}$ is generated by monomials in $x_i,y_i$.
We associate to 
$(Z_1,Z_2)$ the $\chi$-tuple $(\vec Y^1,\ldots,\vec Y^\chi)$ with 
$\vec Y^i=(Y^i_1,Y^i_2)$, where
$$Y_\alpha^i=\big\{ (n,m)\in \Z_{>0}\times \Z_{>0} \bigm|
x_i^{n-1}y_i^{m-1}\not\in I_{Z_\alpha,p_i}\big\}.$$ 
We write $|Y^i_\alpha|$ for the number of elements of $Y^i_\alpha$ and $|\vec Y^i|:=|Y^i_1|+|Y^i_2|$.
This gives a bijection from the fixpoint set $(X_2^{[l]})^\Gamma$ to the 
set of the $\chi$-tuples of pairs of Young diagrams
$(\vec Y^1,\ldots,\vec Y^\chi)$, with $\sum_{i} |\vec Y^i|=l$.

Similarly $\Gamma$ acts on 
$X\times M^H_X(c_1,c_2)$ by $t\cdot(x,E)=(F_t(x),(F_{t}^{-1})^*E)$. Assume for the moment that there exist a universal sheaf  $\E$ over $X\times M^H_X(c_1,d)$,
then one can show that $\E$ has a lifting to a $\Gamma$-equivariant sheaf, unique up
to twist by a character. 

The definition of the determinant bundles and the $K$-theoretic Donaldson invariants
is easily generalized to the equivariant case. If $Y$ is a variety with an action of $\Gamma$, we denote by $K^\Gamma(Y)$, $K^{0\Gamma}(X)$ the Grothendieck groups of $\Gamma$-equivariant coherent sheaves and $\Gamma$-equivariant locally free sheaves respectively. $\chi(u\otimes v):K^\Gamma(X)^2\to \Z$%
\begin{NB} corrected 2.10. LG \end{NB} 
is still a quadratic form. 
The formula  
\eqref{eq:lambdaE} defines a homomorphism $K^\Gamma(X)\to \Pic^\Gamma(S)$, where now $S$ is a scheme with a $\Gamma$-action, and $\E$ a flat family of 
$\Gamma$-equivariant coherent sheaves of class $c\in K(X)_{{\rm num}}$ on $X$, flat over $S$. For $c\in K(X)_{{\rm num}}$ we define $K_c^\Gamma,K_{c,H}^\Gamma\subset 
K^\Gamma(X)$
by the same formula as in section \ref{sec:detbun}.
In the same way as in \ref{sec:detbun}, there are homomorphisms
$\lambda\colon K_c^\Gamma\to \Pic^\Gamma(M_H^X(r,c_1,c_2)_s)$, $\lambda\colon K_{c,H}^\Gamma
\to \Pic^\Gamma(M^X_H(r,c_1,c_2))$, which commute with the inclusions 
$K^\Gamma_{c,H}\subset K^\Gamma_{c}$ and 
$\Pic^{\Gamma}(M^X_H(r,c_1,c_2)_s)\subset \Pic^{\Gamma}(M^X_H(r,c_1,c_2))$.
If $H$ is general with respect to $(r,c_1,c_2)$, then 
$\lambda\colon K_{c,H}^\Gamma
\to \Pic^\Gamma(M^X_H(r,c_1,c_2))$ can be extended to $K_{c}$.
For a flat family $\E$ of equivariant stable sheaves on $X$ parametrized by  $S$, $\lambda$ and $\lambda_\E$ commute with the pullback 
$\phi^*_\E:\Pic^\Gamma(M_H^X(r,c_1,c_2))\to \Pic^\Gamma(S)$ by the  classifying morphism. 

Let $v\in K^\Gamma_c$, where $c$ is the class of an element of 
$M_H^X(c_1,d)$, where $d=4c_2-c_1^2-3$.  
Assume that $H$ is general with respect to $(2,c_1,c_2)$.
 If $Y$ is a variety with a $\Gamma$-action and $w\in K^{0\Gamma}(Y)$,
we denote 
\begin{equation}
\label{chitil}
\widetilde \chi(Y,w):=\pi_!(w)\in \C[t_1^{\pm 1},t_2^{\pm 1}],
\end{equation}
where  $\pi:Y\to pt$ is the projection to a point. 
The {\em equivariant $K$-theoretic Donaldson invariant} of $X$ with respect to 
$v,c_1,d,H$ is 
$\widetilde \chi(M^X_H(c_1,d),\lambda(v))$.%
\begin{NB} cut, 2.10.06 LG \end{NB}
If $L$ is a $\Gamma$ equivariant line bundle on $X$ with $\<c_1(L),c_1\>$ even, let  $v(L)\in K_c^\Gamma$ be an equivariant lift of the class defined by  
\eqref{eq:uL} and $\mu(L):=\lambda(v(L))\in \Pic^\Gamma(M^X_H(c_1,d))$.
We put
$\widetilde \chi(M^X_H(c_1,d),\mu(L))$
and 
$\widetilde \chi_{c_1}^H(L;\Lambda):=\sum_{d\ge 0}\Lambda^d  \widetilde \chi(M^X_H(c_1,d),\mu(L)).$

\begin{NB}
I wanted not to define $u_d$ in the equivariant case. I will only deal with the case
$v\in K_c$, either if $\rk(v)=0$, or fixing $d$.
I keep here the old version in case we need it.

Fix a fixpoint $x\in X^\Gamma$. Assume that $u$ satisfies condition \eqref{even}.
Let  $d\in \Z$ with $d\equiv -c_1^2-3$ modulo $4$. 
Let $c\in K(X)_{{\rm num}}$ be the class of an element of $M^X_H(c_1,d)$.
Then we define 
$u_d:=u-\frac{\chi(v\otimes c)}{2}[\oo_x]\in K^\Gamma_c$, and 
$
\widetilde \chi_{c_1}^H(v;\Lambda):=\sum_{d\ge 0}\Lambda^d  \widetilde 
\chi(M^X_H(c_1,d),\lambda(v_d)).
$
Finally, if $L$ is a $\Gamma$ equivariant line bundle on $X$ with $\<c_1(L),c_1\>$ even, we define $u(L)\in K_c^\Gamma$ by the formula 
\eqref{eq:uL} and $\mu(L):=\lambda(v(L))\in \Pic^\Gamma(M^X_H(c_1,d))$.

\end{NB}

\begin{NB}
This remark is mostly for later use when we want to study things like an equivariant 
version of the strange duality conjecture. 
For this we will be interested in the $\widetilde \chi_{c_1}^H(L;\Lambda)$ and thus in 
$v(L)$. Note that $v(L)$ does not depend on $d$, but it is still not canonical,
because it depends an equivariant lift of the class of a point. 

\begin{Remark}
We will later see that the nonequivariant wallcrossing for $\widetilde \chi_{c_1}^H(L;\Lambda)$ (and more generally for any class $v\in K(X)$ of rank $0$) does not change
when one changes $c_1$ but leaves $\xi$ fixed (i.e. one twists $\E$ with an equivariant line bundle). 
Let me put $v(L)=-(1-L^{-1})$ (i.e. without any class of the  structure sheaf of a point).
In many cases one will be able to do this in such a way that $\chi(v(L)\otimes c)=0$, i.e. so that $v(L)\in K_c$ or in $K_{c,H}$.
Lets look at this in the  case of $\P^2$. 
\begin{enumerate}
\item 
Look at the case $c_1$ odd. Only $v(\oo(2n))$ with $n\in \Z$ give line bundles 
on the moduli space. When we write the moduli space as 
$M_H((2n-3)H,d)$, then we get indeed $\chi(v(L)\otimes c)=0$ for $c$ the class of an
element in $M_H((2n-3)H,d)$. This can be seen as follows. We know that the  condition 
is independent of $d$. Thus we can assume that $c$ is the class of $\oo(n-2)\oplus 
\oo(n-1)$. Then $$\chi(v(L)\otimes (\oo(n-2)\oplus 
\oo(n-1))=-\chi(\oo(n-2))-\chi(\oo(n-1))+\chi(\oo(-n-1))+\chi(\oo(-n-2))=0$$
by Serre Duality.  Thus we get $v(L)\in K_c$. This means that we get a canonical
equivariant lift of the determinant bundle.

Alternatively we can also just look at $M_H(H,d)$ and replace $v(\oo(2n))$ by 
$\oo(n-2)-\oo(-n-2)$.
\item 
Now look at the case $c_1$ is even. 
Let us deal with $v(\oo(2n+1))$.
Then we can write the moduli space as $M_H((2n-2)H,d)$, and we again only have to 
check for the class of $\oo(n-1)\oplus \oo(n-1)$. Then we get
$$\chi(v(L)\otimes (\oo(n-1)\oplus 
\oo(n-1))=-\chi(\oo(n-1))-\chi(\oo(n-1))+\chi(\oo(-n-2))+\chi(\oo(-n-2))=0$$
again by Serre duality. It does not seem that something similar works for 
$\oo(2n)$.

Alternatively we can also look at $M_H(0,d)$ and replace $v(\oo(2n+1))$
by $\oo(n-1)- \oo(-n-2)$.
\end{enumerate}

One  reason why this should be interesting is the following.
Assume we are able to prove the Strange duality conjecture (maybe just
for some fixed small n=1,2). This means that we compute the holomorphic Euler 
characteristic of the corresponding line bundle to be equal to that of a a certain line
bundle on some projective space. On the other hand at least for small $n$ it is
know that there is a canonical injective homomorphism from one space of 
sections to the other (and the higher cohomology vanishes). This means that 
the space of sections of suitable Donaldson line bundle on a moduli space $\P^2$ is canonically
isomorphic to the space of sections of a certain line bundle on $\P^m$.
If the  line bundle on the moduli space on $\P^2$ is given a canonical lifting, we should
get that the isomorphism is equivariant. Thus we will also get a very simple formula
for the character of the space of sections of the Donaldson line bundle.
On the other  hand this result will give very nontrivial relations for the  
$K$-theoretic Nekrasov partition function (for all genera), which might
actually help to compute the higher genus parts.

\end{Remark}
\end{NB}

\begin{Definition}\label{eqwall}
Let $v\in K(X)$. 
Let $\xi\in H^2(X,\Z)$ be an equivariant lifting of a class of type $c_1$. 
Then $\I_1$, $\I_2$, $\F_1$, $\F_2$. $\AA_{\xi,+}$ and $\AA_{\xi,-}$ are in a natural way equivariant sheaves on $X^{[n]}\times X^{[m]}$ (resp. elements in 
$K^\Gamma(X^{[n]}\times X^{[m]})$), and the equivariant wallcrossing terms $\widetilde\Delta_{\xi,T}^X(v;\Lambda)$,
$\widetilde \Delta_{\xi}^X(v;\Lambda)$  are defined by the right-hand side of 
formulas (\ref{wallct}), with the holomorphic  Euler characteristic $\chi$ replaced by 
the equivariant pushforward $\widetilde \chi$ to a point.
Now $\widetilde\Delta_{\xi,T}^X(v;\Lambda)$ can be understood by
localization in equivariant $K$-theory on $X^{[n]}\times X^{[m]}$. 
Then $\widetilde\Delta_{\xi,T}^X(v;\Lambda)\in
\Lambda^{-\xi^2-3}\Q(t_1,t_2,T^{\frac{1}{2}})[[\Lambda]]$.%
\begin{NB}
Changed. 
Nov. 6, HN
Original:

Now $\widetilde\Delta_{\xi,T}^X(v;\Lambda)$ is not only a formal expression:
application of localization in $K$-theory on $X^{[n]}\times X^{[m]}$ 
gives $\widetilde\Delta_{\xi,T}^X(v;\Lambda)\in
\Lambda^{-\xi^2-3}\Q(t_1,t_2,T^{\frac{1}{2}})[[\Lambda]]$.
\end{NB}
Then using  \eqref{delT} we can view $\widetilde\Delta_{\xi,T}^X(v;\Lambda)$ as an element
of $\Lambda^{-\xi^2-3}\Q[t_1^{\pm 1},t_2^{\pm 1}]((T^{\frac{1}{2}}))[[\Lambda]]$, and 
$[\widetilde\Delta_{\xi,T}^X(v;\Lambda)]_{T^0}$ is its coefficient of $T^0$.
Similarly using 
\eqref{delT1}, $\widetilde\Delta_{\xi,T}^X(v;\Lambda)$ is an element of 
$\Lambda^{-\xi^2-3}\Q[t_1^{\pm 1}, t_2^{\pm 1}]((T^{-\frac{1}{2}}))[[\Lambda]]$, and 
$[\widetilde\Delta_{\xi,T}^X(v;\Lambda)]_{(T^{-1})^0}$ is its coefficient of $(T^{-1})^0$.
In particular $\widetilde\Delta_{\xi}^X(v;\Lambda)\in \Q[t_1^{\pm 1},t_2^{\pm 1}][[\Lambda]]$, 
and $\widetilde\Delta_{\xi}^X(v;\Lambda)|_{t_1=t_2=1}=\Delta_{\xi}^X(v;\Lambda)$.
\end{Definition}

Let $c\in K(X)$ be the class of an element of $M^X_H(c_1,d)$.
In the same way as in Remark \ref{integral}, we see that 
the coefficient of $\Lambda^d$ of $\widetilde\Delta_{\xi,T}^X(v;\Lambda)$ is 
either in $T^{\frac{1}{2}}\Q(t_1,t_2,T)$ (and the coefficient of $\Lambda^d$ of  $\widetilde\Delta_{\xi}^X(v;\Lambda)$ is $0$)
or in $\Q(t_1,t_2,T)$. If $v\in K_c^\Gamma$, then the coefficient is in $\Q(t_1,t_2,T)$.

Let $v\in K_c^\Gamma$. 
Under the assumptions of Proposition \ref{wallcr}
let $\widetilde B_+$ be a set consisting  of one equivariant lift
 $\xi$  for each class of type $(c_1,d)$
with $\<\xi\cdot H_+\>>0 >\<\xi\cdot H_-\>$. 
Then the same proof as before (with all sheaves and classes replaced by their
equivariant versions) shows that 
\begin{align*}
\widetilde \chi(M^X_{H_+}(c_1,d),\lambda(v))-\widetilde \chi(M^X_{H_-}(c_1,d),\lambda(v))&=\sum_{\xi\in B_+}\bigl[\widetilde\Delta^X_{\xi}(v;\Lambda)\bigr]_{\Lambda^d}.
\end{align*}

\begin{NB}
Version with $u_d$.

Then the same proof as that of Proposition \ref{wallcr}
(with all sheaves and classes replaced by the equivariant versions)  shows
the following.  

\begin{Lemma}\label{ewallcr} Fix $c_1,d$, let $c\in K(X)$ be the class
of an element of $M_H^X(c_1,d)$.
Let $v\in K_c^\Gamma$.
Let $H_-$, $H_+$ be ample divisors on $X$, which do not lie on a wall 
of type $(c_1,d)$. Let $B_+$ be a set consisting  of one equivariant lift
 $\xi$  for each class of type $(c_1,d)$
with $\<\xi\cdot H_+\>>0 >\<\xi\cdot H_-\>$. Assume that all classes in $B_+$ 
are good. 
Then  
\begin{align*}
\widetilde \chi(M^X_{H_+},\lambda(v))-\widetilde \chi(M^X_{H_-},\lambda(v))&=\sum_{\xi\in B_+}\bigl[\widetilde \Delta^X_{\xi}(v;\Lambda)\bigr]_{\Lambda^d}.
\end{align*}
\end{Lemma}
\end{NB}


Now we want to give a formula expressing $\widetilde \Delta^X_{\xi,T}(v;\Lambda)$ in terms 
of the K-theoretic Nekrasov partition function $Z_K$.
 For the rest of this section let $\xi$ be an equivariant lift of a class of type
$c_1$, and let $v\in K^\Gamma(X)$. 
We first give, up to a correction term,  an 
expression for  $\widetilde \Delta^X_{\xi,T}(v,\Lambda)$ in terms of the instanton
part. Then we show
that this correction term is given  by the perturbation part.

\begin{Theorem}\label{main1} Let $v\in K^\Gamma(X)$. 
\begin{align*}&\widetilde \Delta_{\xi,e^{-\bbeta t}}^X(v;\bbeta\Lambda)|_{{t_1\to e^{\bbeta\ve_1}}\atop{t_2\to e^{\bbeta\ve_2}}}\\
&=
\frac{1}{\bbeta\Lambda}\exp\Bigg(\bbeta\Big(\frac{\<K_X^3\>}{48}-\frac{\<\Todd_2(X) K_X\>}2+\<[2\ch(v)\exp(c_1/2)\Todd(X)]_3\>\Big)\\
&\quad +\sum_{i=1}^\chi F_{-\rk(v)}\big(w(x_i),w(y_i),
\hbox{$\frac{t-\iota_{p_i}^*\xi}{2}$};\Lambda e^{-\bbeta\iota_{p_i}^*K_X/4},\bbeta\iota_{p_i}^*(c_1(v)+\hbox{$\frac{\rk(v)}{2}$}(c_1-K_X)))\Bigg).
\end{align*}
\end{Theorem}

\begin{NB}
The term $(1+\frac{\rk(v)}{2})$ was deleted. (twice) Oct. 18, H.N.
Original was
\begin{align*}&\widetilde \Delta_{\xi,e^{-\bbeta t}}^X(v;\bbeta\Lambda)|_{{t_1\to e^{\bbeta\ve_1}}\atop{t_2\to e^{\bbeta\ve_2}}}\\
&=
\frac{1}{\bbeta\Lambda}\exp\Bigg(\bbeta\Big(\frac{\<K_X^3\>}{48}-(1+\frac{\rk(v)}{2})\frac{\<\Todd_2(X)
   K_X\>}2+\<[2\ch(v)\exp(c_1/2)\Todd(X)]_3\>\Big)\\
&\quad +\sum_{i=1}^\chi F_{-\rk(v)}\big(w(x_i),w(y_i),
\hbox{$\frac{t-\iota_{p_i}^*\xi}{2}$};\Lambda e^{-\bbeta(1+\frac{\rk(v)}{2})\iota_{p_i}^*K_X/4},\bbeta\iota_{p_i}^*(c_1(v)+\hbox{$\frac{\rk(v)}{2}$}(c_1-K_X)))\Bigg).
\end{align*}
\end{NB}

Note that the left-hand side 
lies in $(\bbeta\Lambda)^{-\xi^2-3}\Q(e^{\ve_1},e^{\ve_2},e^{t\bbeta})
[[\bbeta\Lambda]]$.%
\begin{NB} corrected 2.10.LG \end{NB}
In the course
of the proof we will also have to show how one can interpret the right-hand
side, so that both sides lie in the same ring. 
\begin{NB}
This is certainly incorrect. We have to be much more careful with the coefficients
and use the calculations of Kota.

Now I instead think it is correct 2.10.06 LG
\end{NB}


\begin{Lemma}\label{instan} 
Let $M$ be a $\Gamma$-equivariant line bundle on $X$, with 
$c_1(M)=\xi$. 
Then in 
$(\bbeta\Lambda)^{-\xi^2-3}\Q(e^{\bbeta\ve_1},e^{\bbeta\ve_2},e^{\bbeta t})[[\bbeta\Lambda]]$ we have
\begin{align*}
&\widetilde \Delta_{\xi,e^{-\bbeta t}}^X(v,\bbeta \Lambda)|_{{t_1\to e^{\bbeta\ve_1}}\atop{t_2\to e^{\bbeta\ve_2}}}=\exp(2\bbeta \<v^{(3)}\>)\\
&\cdot \frac{\prod_{i=1}^\chi \Zin_{-\rk(v)}\big(w(x_i),w(y_i),
\frac{t-\iota_{p_i}^*\xi}{2};\Lambda e^{-\bbeta\iota_{p_i}^*(K_X)/4},\bbeta,
\bbeta\iota_{p_i}^*(c_1(v)+\frac{\rk(v)}{2}(c_1-K_X)\big)
} {(\bbeta\Lambda)^{\xi^2+3}
\bigwedge_{-e^{-\bbeta t}}(-\widetilde\chi(X,M^\vee)^\vee)\bigwedge_{-e^{\bbeta t}}(-\widetilde \chi(X,M)^\vee)
}.
\end{align*}
\end{Lemma}

\begin{NB}
The term $(1+\frac{\rk(v)}{2})$ is deleted. Oct. 18, H.N.
\end{NB}

\begin{proof}
Following \cite{GNY}, we denote $C(0):=\ch(\I_1)e^{\xi/2}+
\ch(\I_2)e^{-\xi/2}$ on $X\times X^{[n]}\times X^{[m]}$, and 
$C_i(0):=[C(0)]_i$.
The Grothendieck-Riemann-Roch theorem implies that 
\begin{equation}\label{FC0}
\begin{split}
\ch(\lambda_{\F_1}(v)\otimes \lambda_{\F_2}(v))&=
\exp([p_*(q^*(\ch(v)) \ch(\F_1\oplus \F_2)\Todd(X)]_1)\\
&=\exp([p_*(q^*(\ch(v)) C(0) e^{c_1/2}\Todd(X))]_1)\\
&=\exp(C_3(0)/\rk(v)+C_2(0)/v^{(1)}+2/v^{(3)}).
\end{split}
\end{equation}
Let $(Z_1,Z_2)\in (X^{[n]}\times X^{[m]})^\Gamma$ correspond to 
$(\vec Y^1,\ldots,\vec Y^\chi)$. 
By \cite[Lemma 3.4]{GNY} the cotangent space and the fibres
of $\AA_{\pm}^\vee$ at $(Z_1,Z_2)$ are 
\begin{equation}
\label{eq:denom}
\begin{gathered}
\bigwedge\nolimits_{-1}T^*_{(Z_1,Z_2)}X^{[n]}\times X^{[m]}|_{{t_1\to
    e^{\bbeta\ve_1}}\atop{t_2\to e^{\bbeta\ve_2}}}
=\prod_{i=1}^\chi\prod_{\gamma=1}^2
n_{\gamma,\gamma}^{\vec Y_i}(w(x_i),w(y_i),\hbox{$\frac{t-\iota_{p_i}\xi}{2}$};\bbeta),
\\
\bigwedge\nolimits_{-e^{-t\bbeta}}\AA_{+}^\vee(Z_1,Z_2)
|_{{t_1\to e^{\bbeta\ve_1}}\atop{t_2\to e^{\bbeta\ve_2}}}
=\bigwedge\nolimits_{-e^{-t\bbeta}}-\widetilde \chi(X,M^\vee)^\vee |_{{t_1\to e^{\bbeta\ve_1}}\atop{t_2\to e^{\bbeta\ve_2}}}\prod_{i=1}^\chi
n_{1,2}^{\vec Y_i}\big(w(x_i),w(y_i),\hbox{$\frac{t-\iota_{p_i}\xi}{2}$};\bbeta),
\\
\bigwedge\nolimits_{-e^{ t\bbeta}}\AA_-^\vee(Z_1,Z_2)
|_{{t_1\to e^{\bbeta\ve_1}}\atop{t_2\to e^{\bbeta\ve_2}}}
=\bigwedge\nolimits_{-e^{ t\bbeta}}-\widetilde \chi(X,M)^\vee|_{{t_1\to e^{\bbeta\ve_1}}\atop{t_2\to e^{\bbeta\ve_2}}}\prod_{i=1}^\chi
n_{2,1}^{\vec Y_i}\big(w(x_i),w(y_i),\hbox{$\frac{t-\iota_{p_i}\xi}{2}$};\bbeta).
\end{gathered}
\end{equation}
\begin{NB}
The first formula of \eqref{eq:denom} is obvious and the other two 
are obvious with $t$ replaced by $0$ on both sides. To get the formula with $t$
note that replacing $a$ by $t/2+a$ will multiply every factor of 
$n_{1,2}(\ve_1,\ve_2,a;\bbeta)$ in 
\eqref{eq:nab} by $e^{-\bbeta t}$ and every factor of $n_{2,1}(\ve_1,\ve_2,a;\bbeta)$
by $e^{\bbeta t}$.
\end{NB}

By \eqref{eq:EKFE} and \cite[(3.11)]{GNY}  we get 
{\allowdisplaybreaks
\begin{equation}
\label{eq:EY}
\begin{split}
&\prod_{i=1}^\chi \exp\Big(\bbeta\rk(v)|\vec Y_i|\frac{w(x_i)+w(y_i)}{2}
\Big)C_{-\rk(v)}^{\vec Y}(w(x_i),w(y_i),\hbox{$\frac{t-\iota_{p_i}^*\xi}2$};\bbeta,
\bbeta\iota_{p_i}^*v^{(1)})\\
&=\prod_{i=1}^\chi E^{\vec Y}(w(x_i),w(y_i),\hbox{$\frac{t-\iota_{p_i}^*\xi}2$};\bbeta,
\bbeta\iota_{p_i}^*v^{(1)},\rk(v))\\
&=\iota_{(Z_1,Z_2)}^*\exp\Bigl(
\begin{aligned}[t]
\bigl[\ch(\I_1)
e^{\frac{\xi-\bbeta t}2}\oplus &
\ch(\I_2)e^{\frac{\bbeta t-\xi}2}\bigr]_{2}/v^{(1)}
\\
&+\bigl[\ch(\I_1)
e^{\frac{\xi-\bbeta t}2}\oplus 
\ch(\I_2)e^{\frac{\bbeta t-\xi}2}\bigr]_{3}/\rk(v)\Bigr)|_{{\ve_1\to \bbeta \ve_1}\atop {\ve_2\to
\bbeta \ve_2}}
\end{aligned}
\\
&=\iota_{(Z_1,Z_2)}^*\Bigl(\exp\bigl( C_2(0)/v^{(1)}+C_3(0)/\rk(v)\bigr)\Bigr)_{{\ve_1\to \bbeta\ve_1}\atop {\ve_2\to \bbeta\ve_2}}(e^{\bbeta t})^{\frac{1}{2}(\chi(f_2\otimes v)-\chi(f_1\otimes v))}\\
&=\exp(-2\bbeta \<v^{(3)}\>) 
\iota_{(Z_1,Z_2)}^*\Bigl(\ch(\lambda_{\F_1}(v)\otimes \lambda_{\F_2}(v)))\Bigr)_{{\ve_1\to\bbeta\ve_1}\atop {\ve_2\to \bbeta \ve_2}}(e^{\bbeta t})^{\frac{1}{2}(\chi(f_2\otimes v)-\chi(f_1\otimes v))}.
\end{split}
\end{equation}
In} the fourth line  we use that $\ch_0(\I_{\alpha})=1$, $\ch_1(\I_\alpha)=0$ for $\alpha=1,2$, and that $\ch_2(\I_1)/1=-n$, $\ch_2(\I_2)/1=-m$
and thus
\begin{align*} \iota_{(Z_1,Z_2)}^*\Bigl(\bigl[\ch(\I_1)
e^{\frac{\xi-\bbeta t}2}\oplus 
\ch(\I_2)e^{\frac{\bbeta t-\xi}2}\bigr]_{2}/v^{(1)}\Bigr)
&=\iota_{(Z_1,Z_2)}^* \bigl(C_2(0)/v^{(1)}\bigr)
-\<\xi/2,v^{(1)}\> \bbeta t,\\
\iota_{(Z_1,Z_2)}^*\Bigl(\bigl[\ch(\I_1)
e^{\frac{\xi-\bbeta t}2}\oplus 
\ch(\I_2)e^{\frac{\bbeta t-\xi}2}\bigr]_{3}/\rk(v)\Bigr)&=\iota_{(Z_1,Z_2)}^* 
\bigl(C_3(0)/\rk(v)\bigr)+(n-m)\frac{\rk(v)}{2}\bbeta t, 
\end{align*}
and formula \eqref{eq:f2ud}. In the last line of \eqref{eq:EY} we use  \eqref{FC0}.

Write $|Y|:=|\vec Y_1|+\ldots+|\vec Y_\chi|$, and 
write $(Z^Y_1,Z^Y_2)$ for the point of  $X^{[n]}\times X^{[m]}$ with $n+m=|Y|$ 
determined by an $\chi$-tuple $Y= (\vec Y_1,\ldots,\vec Y_\chi)$ of pairs of
Young diagrams. Using that $\iota_{p_i}^*K_X=-w(x_i)-w(y_i)$, we get by localization, \eqref{eq:EY} and \eqref{eq:ZinstmF} 
{\allowdisplaybreaks
\begin{align*}
&(\bbeta\Lambda)^{-\xi^2-3}\frac{\prod_{i=1}^\chi
\Zin_{-\rk(v)}\big(w(x_i),w(y_i),\frac{t-\iota_{p_i}^*\xi}{2};\Lambda 
e^{-\bbeta\iota^*_{p_i} K_X/4},\bbeta,
\bbeta\iota_{p_i}^*v^{(1)}\big)}{\bigwedge_{e^{-\bbeta t}}-\widetilde \chi(X,M^\vee)^\vee 
\bigwedge_{-e^{\bbeta t}}
-\widetilde \chi(X,M)^\vee})|_{{t_1\to e^{\bbeta\ve_1}}\atop{t_2\to e^{\bbeta\ve_2}}}\\
&=
\begin{aligned}[t]
& \sum_{Y=(\vec Y_1,\ldots,\vec Y_\chi)}
(\bbeta\Lambda)^{4|Y|-\xi^2-3}
\\
&\quad \times
\frac{\prod_{i=1}^\chi
 \exp\Big(\bbeta\rk(v)|\vec Y_i|\frac{w(x_i)+w(y_i)}{2}
\Big)C_{-\rk(v)}^{\vec Y}(w(x_i),w(y_i),\hbox{$\frac{t-\iota_{p_i}^*\xi}2$};\bbeta,
\bbeta\iota_{p_i}^*v^{(1)})}
{\big(\bigwedge_{-1}(T^*_{(Z_1^Y,Z^Y_2)}X_2^{[|Y|]})
\bigwedge_{-e^{-\bbeta t}}\AA_+^\vee(Z_1^Y,Z_2^Y)
\bigwedge_{-e^{\bbeta t}}
\AA_-^\vee(Z_2^Y,Z_1^Y))\big)|_{{t_1\to e^{\bbeta\ve_1}}\atop{t_2\to
    e^{\bbeta\ve_2}}}}
\end{aligned}
\\
&
\begin{aligned}[t]
=\exp(-2\bbeta \<v^{(3)}\>)
\!\!\! & \sum_{{n,m\ge 0}\atop {d=4(n+m)-\xi^2-3}}\!\!\! (\bbeta\Lambda)^{d}
\\
& \quad \times\widetilde\chi\Bigl(X^{[n]}\times X^{[m]},
\frac{\lambda_{\F_1}(v)\otimes \Lambda_{\F_2}(v) }{e^{-\bbeta t(\frac{1}{2}
\chi(f_2\otimes v)-\chi(f_1\otimes v)))} \bigwedge_{-e^{-\bbeta t}} \AA_+^\vee\bigwedge_{-e^{\bbeta t}} \AA_-^\vee
}\Bigr)|_{{t_1\to e^{\bbeta\ve_1}}\atop {t_2\to e^{\bbeta\ve_2}}}
\end{aligned}
\\
&=
\widetilde\Delta_{\xi,e^{-\bbeta t}}^X(v,\bbeta\Lambda)|_{{t_1\to
e^{\bbeta\ve_1}}\atop{t_2\to e^{\bbeta\ve_2}}} \exp(-2\bbeta \<v^{(3)}\>).
\end{align*}
In} the third line we use \eqref{eq:EY} and equivariant localization.
\begin{NB}
The term $(1+\frac{\rk(v)}{2})$ was deleted in the first line. 
I also add several line breaks.
Oct.18, H.N.
\end{NB}
\end{proof}

Now we identify the contribution of the perturbation part.
Let  ${\widetilde \oo}$ be the ring of holomorphic functions in  $(x,\bbeta,t)$ 
in a neighborhood of $\sqrt{-1}\R_{>0}\times \sqrt{-1}\R_{< 0}\times \sqrt{-1}\R_{> 0}$.

\begin{Lemma}\label{pertur}
\begin{align*}\sum_{i=1}^\chi& \Fper(w(x_i),w(y_i),\hbox{$\frac{t-i_{p_i}^*\xi}{2}$};\Lambda e^{-\bbeta\iota^*_{p_i}K_X/4},\bbeta)
=\big(-(\chi(M)+\chi(M^\vee)\big)\log(\bbeta\Lambda)\\&
-\frac{\bbeta\<K_X^3\>}{48}+\frac{\bbeta}{2}\<\Todd_2(X) K_X\>
+\log\Big(\frac{1}{\bigwedge_{-e^{\bbeta t}} -\widetilde \chi(X,M)^\vee
\bigwedge_{-e^{-\bbeta t} }-\widetilde \chi(X,M^\vee)^\vee}\Big)|_{{t_1\to e^{\bbeta\ve_1}}\atop
{t_2\to e^{\bbeta\ve_2}}}.
\end{align*}
holds in $\widetilde \oo[[\ve_1,\ve_2]][\prod_i (w(x_i)w(y_i))^{-1}]$.
\end{Lemma}

\begin{NB}
I delete $(1+\frac{\rk(v)}{2})$ twice, and delete the subscript
${}_{-\rk(v)}$ from $\Fper$. Oct. 18, H.N.
Original:
\begin{align*}\sum_{i=1}^\chi& \Fper_{-\rk(v)}(w(x_i),w(y_i),\hbox{$\frac{t-i_{p_i}^*\xi}{2}$};\Lambda e^{-\bbeta(1+\frac{\rk(v)}{2})\iota^*_{p_i}K_X/4},\bbeta)
=\big(-(\chi(M)+\chi(M^\vee)\big)\log(\bbeta\Lambda)\\&
-\frac{\bbeta\<K_X^3\>}{48}+\frac{\bbeta}{2}(1+\frac{\rk(v)}{2})\<\Todd_2(X) K_X\>
+\log\Big(\frac{1}{\bigwedge_{-e^{\bbeta t}} -\widetilde \chi(X,M)^\vee
\bigwedge_{-e^{-\bbeta t} }-\widetilde \chi(X,M^\vee)^\vee}\Big)|_{{t_1\to e^{\bbeta\ve_1}}\atop
{t_2\to e^{\bbeta\ve_2}}}.
\end{align*}
\end{NB}


\begin{NB} I have not yet thought carefully about where precisely these equalities
hold (mostly how to improve  where they hold), and where precisely we want
Theorem \ref{main1} to hold. What I wrote comes precisely out of the
proof of \cite[Lem.~3.14]{GNY}.
\end{NB}
\begin{proof}
By \cite[(3.17)]{GNY} we get that $\sum_{i=1}^\chi \Fper_K(w(x_i),w(y_i),\hbox{$\frac{t-i_{p_i}^*\xi}{2}$};\Lambda,\bbeta)$ is given by the same formula with
 $\frac{\bbeta}{2}\<\Todd_2(X) K_X\>$ replaced by 
$-\frac{\bbeta}{4}\<\xi^2 K_X\>+\frac{\bbeta}{2}\<\xi K_X\> t$.%
\begin{NB}
$(1+\frac{\rk(v)}{2})$ is deleted. Oct. 18, H.N.
\end{NB}
Note that by \eqref{gammati}, when changing $\Lambda$ to $\Lambda
e^{-\bbeta K_X/4}$,%
\begin{NB}
$(1+\frac{\rk(v)}{2})$ is deleted. Oct. 18, H.N.  
\end{NB}
the result changes by adding
\begin{align*}
-\bbeta\sum_{i=1}^\chi & \frac{ w(x_i)+w(y_i)}{4w(x_i)w(y_i)}
\Big((t-\iota_{p_i}^*\xi)^2+\frac{w(x_i)^2+w(y_i)^2+3w(x_i)w(y_i)}{6}\Big)\\
&=
\Big(\frac{\bbeta}{4}\<\xi^2 K_X\>-\frac{\bbeta}{2}\<\xi K_X\> t +\frac{\bbeta}{2}
\<\Todd_2(X) K_X\>\Big).
\end{align*}
\begin{NB}
$(1+\frac{\rk(v)}{2})$ is deleted in both sides. Oct. 18, H.N.
\end{NB}
\begin{NB}
I comment out the following, Oct. 18, H.N:

By formula \eqref{Fper} and localization, we get that 
\begin{align*}
\sum_{i=1}^\chi &\Fper_{-\rk(v)}(w(x_i),w(y_i),\hbox{$\frac{t-i_{p_i}^*\xi}{2}$};\Lambda,\bbeta)-\sum_{i=1}^\chi \Fper_{K}(w(x_i),w(y_i),\hbox{$\frac{t-i_{p_i}^*\xi}{2}$};\Lambda,\bbeta)
\\&=\bbeta\frac{\rk(v)}{2}\Big(-\frac{\<K_X\xi^2\>}{4}
+\frac{\<K_X\xi\>t}{2}\Big).
\end{align*}
\end{NB}
The result follows.
\end{proof}

\begin{NB} I added the following remark in the spirit of our previous paper to show why the theorem follows from these two lemmas.
\end{NB}

Writing $\ch(-\widetilde \chi(X,M))=\sum_{i=0}^\ell e^{\alpha_j}$, 
$\ch(-\widetilde \chi(X,M^\vee))=\sum_{i=0}^{\ell'} e^{\alpha'_k}$, we see that
\begin{equation}
\label{reg}\begin{split}\log\Big(&\frac{\bbeta^{-(\chi(M)+\chi(M^\vee))}}{\bigwedge_{-e^{\bbeta t}} -\widetilde \chi(X,M)^\vee
\bigwedge_{-e^{-\bbeta t} }-\widetilde \chi(X,M^\vee)^\vee}\Big)|_{{t_1\to e^{\bbeta\ve_1}}\atop
{t_2\to e^{\bbeta\ve_2}}}\\&=\sum_{j=1}^\ell 
\log\Big(\frac{\bbeta}{1-e^{-(\alpha_j-t)\bbeta}}\Big)+\sum_{k=1}^{\ell'}
\log\Big(\frac{\bbeta}{1-e^{-(\alpha'_k+t)\bbeta}}\Big).
\end{split}\end{equation}
Let $\widetilde \oo$ denote the ring of holomorphic functions in $(t,\Lambda,\bbeta)$ in an open subset of $\C^3$ which contains for any 
$(t,\Lambda)\in (\C\setminus \R_{\le 0})^2$ an open neighbourhood of $\bbeta=0$. Then \eqref{reg} shows that the left hand side of Lemma 
\ref{pertur} lies in $\widetilde \oo[[\ve_1,\ve_2]]$.
Thus we can view also the left hand side of Lemma 
\ref{pertur} to lie in $\widetilde  \oo[[\ve_1,\ve_2]]$, and we can take the
exponential of both sides of the equation. Note that the exponential of the 
right hand side lies in $\Q(e^{\bbeta\ve_1},e^{\bbeta\ve_2},e^{\bbeta t})[[\bbeta\Lambda]]$. With this remark
Theorem \ref{main1} follows from Lemma \ref{instan} and Lemma \ref{pertur}.

Now we  express $\widetilde\Delta_{\xi,T}^X(v,\Lambda)$ in terms of the  
$Z_{-\rk(v)}(\ve_1,\ve_2,a;\Lambda,\bbeta)$. 
The Nekrasov conjecture determines the lowest order terms 
in $\ve_1,\ve_2$ of 
$F_{-\rk(v)}(\ve_1,\ve_2,a;\Lambda,\bbeta)$, but not of $F_{-\rk(v)}(\ve_1,\ve_2,a;\Lambda,\bbeta,\tau)$.

\begin{Corollary}\label{deltaq}
Let $v\in K^\Gamma(X)$.%
\begin{NB} Corrected 4.10.06 LG\end{NB}
Then 
\begin{align*}
\widetilde\Delta_{\xi,e^{-\bbeta t}}^X(v,\bbeta\Lambda)&=
\frac{1}{\bbeta\Lambda}\exp\Big(
\begin{aligned}[t]
\bbeta\Big(\frac{\<K_X^3\>}{48}-&\frac{1}{2}\big\<\Todd_2(X)(K_X+c_1(v)+\frac{\rk(v)}{2}(c_1-K_X))\big\>
\\
&+2\<[\ch(v)e^{c_1/2}\Todd(X)]_3\>\Big)\Big)
\end{aligned}
\\
&\times \Big(\sum_{i=1}^\chi
F_{-\rk(v)}\big(w(x_i),w(y_i),\hbox{$\frac{t-\iota_{p_i}^*\xi}{2}$};\Lambda e^{-\frac{\bbeta}{4}\iota_{p_i}^*(K_X+c_1(v)+\frac{\rk(v)}{2}(c_1-K_X))}\big)\Big).
\end{align*}
\end{Corollary}

\begin{NB}
I add $-\frac{\rk(v)}2 K_X$ according to my change of the definition
of the partition function. Oct. 18, H.N.

Original:
\begin{align*}
\widetilde\Delta_{\xi,e^{-\bbeta t}}^X(v,\bbeta\Lambda)&=
\frac{1}{\bbeta\Lambda}\exp\Big(\bbeta\Big(\frac{\<K_X^3\>}{48}-\frac{1}{2}\big\<\Todd_2(X)(c_1(v)+K_X+\frac{\rk(v)}{2}c_1)\big\> +2\<[\ch(v)e^{c_1/2}\Todd(X)]_3\>\Big)\Big)\\
&\times \Big(\sum_{i=1}^\chi
F_{-\rk(v)}\big(w(x_i),w(y_i),\hbox{$\frac{t-\iota_{p_i}^*\xi}{2}$};\Lambda e^{-\frac{\bbeta}{4}\iota_{p_i}^*(c_1(v)+K_X+\frac{\rk(v)}{2}c_1)}\big)\Big).
\end{align*}
\end{NB}

\begin{proof} Let $\tau,\sigma$ be variables.
In the same way as in  \cite[section 4.5]{NY2}, we see that 
$$\Zin_m(\ve_1,\ve_2,a;\Lambda e^{-\sigma/4},
\bbeta, \tau)=
\exp\Big(\frac{\tau a^2}{\ve_1\ve_2}\Big)\Zin_m(\ve_1,\ve_2,a;\Lambda 
e^{-(\tau+\sigma)/4},\bbeta).$$
On the other hand, by \cite[formula after (4.12)]{NY3},
we get that 
\begin{align*}
\Fper_m(\ve_1,\ve_2,a;\Lambda e^{-\sigma/4},\bbeta ,\tau)=
\Fper_m(\ve_1,\ve_2,a;\Lambda e^{-(\tau+\sigma)/4},\bbeta)
-\frac{\tau a^2}{\ve_1\ve_2}-\frac{\tau (\ve_1^2+\ve_2^2+3\ve_1\ve_2)}
{24\ve_1\ve_2}.
\end{align*}
The result follows by 
%
by localization and Theorem \ref{main1}. \end{proof}
\begin{NB}
We get
\begin{align*}
&\sum_{i=1}^\chi F_{-\rk(v)}\big(w(x_i),w(y_i),
\hbox{$\frac{t-\iota_{p_i}^*\xi}{2}$};\Lambda e^{-\bbeta(1+\frac{\rk(v)}{2})\iota_{p_i}^*K_X/4},\bbeta\iota_{p_i}^*(c_1(v)+\hbox{$\frac{\rk(v)}{2}$}(c_1-K_X)))\\
&=-\frac{\bbeta}{2}\big\<\Todd_2(X)(c_1(v)+\frac{\rk(v)}{2}(c_1-K_X)\big\> 
+\sum_{i=1}^\chi F_{-\rk(v)}\big(w(x_i),w(y_i),
\hbox{$\frac{t-\iota_{p_i}^*\xi}{2}$};\Lambda e^{-\bbeta(1+\frac{\rk(v)}{2})\iota_{p_i}^*K_X/4-\bbeta(c_1(v)+\hbox{$\frac{\rk(v)}{2}$}(c_1-K_X))/4}\big).
\end{align*}
and 
$$(1+\frac{\rk(v)}{2})\iota_{p_i}^*K_X+\iota_{p_i}^*(c_1(v)+\hbox{$\frac{\rk(v)}{2}$}(c_1-K_X))=\iota_{p_i}^*(c_1(v)+K_X+\frac{\rk(v)}{2}).$$
\end{NB}

\section{Explicit formulas in terms of modular forms}\label{sec:modular}

The result of \cite{NY3} together with \secref{sec:SWcurve} implies
that the following solution of Nekrasov's conjecture and its
refinement are true for the K-theoretic partition function when $m=0$:
\begin{enumerate}
\item $\ve_1\ve_2 F_{m}(\ve_1,\ve_2,a;\Lambda)$ is regular at
$\ve_1$, $\ve_2=0$, 
\item $\mathcal F_0(a;\Lambda)$ is the Seiberg-Witten prepotential 
associated with the Seiberg-Witten curve 
$Y^2 = P(X)^2 - 4(-X)^{2+m}(\bbeta\Lambda)^4$,
\item $H$ comes only from the perturbation part, i.e.\
\(
   H(a,\Lambda) = \pi\sqrt{-1}a,
\)%
\begin{NB}
$H$ was originally
\begin{equation*}
   H(a,\Lambda) = \pi\sqrt{-1}a - m\bbeta \frac{a^2}2.
\end{equation*}
But the second term disappear as we changed the definition of the
partition function.
Oct. 19, H.N.
\end{NB}
\item $\exp A = \left(\frac{2}{\theta_{00}\theta_{10}}\right)^{1/2}$,
$\exp B = \theta_{01}\exp A$, where the $\theta_{**}$ are theta functions
with variable $q=e^{2\pi\sqrt{-1}\tau}$, where $\tau$ is the period of
the above Seiberg-Witten curve, i.e.\
$\tau=-\frac1{2\pi\sqrt{-1}}\frac{\partial^2 \F_0}{\partial a^2}$.
\end{enumerate}
Here $\mathcal F_0$, $H$, $A$, $B$ are given by the expansion
\begin{multline}\label{eq:genus_expansion}
    \ve_1\ve_2 F_m(\ve_1,\ve_2,a;\Lambda,\bbeta)
\\
    = \mathcal F_0(a;\Lambda,\bbeta) + (\ve_1+\ve_2) H(a;\Lambda)
    + \ve_1\ve_2 A(a;\Lambda,\bbeta)
    + \frac{\ve_1^2+\ve_2^2}3 B(a;\Lambda,\bbeta) + \cdots.
\end{multline}
When $|m| < 2$, the above (1)$\--$(3) follow from a conjectural blowup
equation~(\ref{eq:conjvanish}) as we explained in
\subsecref{subsec:MoreOn}.
The analogue of the statement (4) is \eqref{eq:AB} which follows from
the conjecture (\ref{eq:conjvanish2}).
In the above we implicitly assume $|m|\le 2$ as the Seiberg-Witten
curve changes the genus otherwise.
According to a physical argument \cite{IMS,Ta}, the remaining case
$m=\pm 2$ is similar to the case $|m| < 2$, in particular (1),(2)
should be true. (These probably follow from the approach in
\cite{NO}.) But we believe that the blowup equation must be modified,
and (3) is probably {\it not\/} true.
\begin{NB}
The explanation is modified and expanded. Oct. 20, H.N.
\end{NB}

In the following we {\it assume\/} the above (1)$\--$(3)%
\begin{NB} changed 30.10 LG\end{NB}
and
\eqref{eq:AB} are also true for $m=\pm 1$.

Once we have the above (1)$\--$(3),%
\begin{NB} changed 30.10 LG\end{NB} 
then the same argument as in
\cite[proof of Thm.~4.2, in particular of (4.12)]{GNY} gives

\begin{Corollary}\label{cor:formula1}
\begin{equation*}
\begin{split}
\left.\Delta_{\xi,e^{-\bbeta t}}^X(v,\bbeta\Lambda)\right|_{t=2a}
&= \frac1{\bbeta\Lambda}
\sqrt{-1}^{\<\xi,K_X\>} q^{-\frac12 \left(\frac\xi2\right)^2}
\\
& \quad\times \exp\Bigg[
\frac{\bbeta}8 \frac{\partial^2\F_0}{\partial a\partial\log\Lambda}
\left\<\xi(K_X+c_1(v)+\frac{\rk(v)}2 (c_1-K_X))\right\>
\\
&\qquad\quad
+\frac{\bbeta^2}{32} \frac{\partial^2\F_0}{(\partial\log\Lambda)^2}
\left\<(K_X+c_1(v)+\frac{\rk(v)}2 (c_1-K_X))^2\right\>
+ \chi A + \sigma B
\Bigg].
\end{split}
\end{equation*}
\end{Corollary}

\begin{NB}
Corrected according to the change of the partition function by H.N,
Oct. 18. 
\end{NB}

We have expressed the wallcrossing $\Delta^X_{\xi,e^{-\bbeta
    t}}(v,\bbeta\Lambda)$ in terms of the partition function with 5D
  Chern-Simons term. As in \cite[\S4]{GNY} we use the Nekrasov
  conjecture to give an explicit formula in terms of $q$-development
  of modular forms.

We identify $t/2$ with $a$ hereafter.

\begin{Theorem}\label{thm:q-expansion}

\textup{(1)} 
Let 
\(
   \Delta_{\xi,e^{-2\bbeta a}}^X(v,\bbeta\Lambda)
   = \sum_{n\ge 0} \Delta_n \Lambda^{4n-\xi^2-3}.
\)%
\begin{NB}
Corrected. $n\to 4n$. Nov. 9, HN  
\end{NB}
Then $\Delta_n$ is equal to $0$ 
if $\langle\xi, c_1(v)+\frac{\rk(v)}2 (c_1-K_X)\rangle+{\rk(v)} n$ is
odd, and equal to the coefficient in
\begin{equation*}
  2 \Coeff_{(q^{1/8})^0}
  \left[ \left.\Delta^X_{\xi,e^{-2\bbeta a}}(v,\Lambda)
      \frac{a^2}{\Lambda}
    \right|_{a=a(q^{1/8},\Lambda)}
    q^{1/8}\frac{\partial\left(\frac{\Lambda}a\right)}{\partial (q^{1/8})}
    \right]
\end{equation*}
otherwise.
\begin{NB}
The change of variable is made as $\frac\Lambda{a}\leftrightarrow
q^{1/8}$ instead of ${a}\leftrightarrow q^{1/8}$.
Oct. 20, H.N.
\end{NB}

\textup{(2)} Suppose $\rk(v) = -m = 0$. Then the terms in $[\ \ ]$ above
are given in explicit modular forms in $\C((q^{1/8}))[[\Lambda]]$.
\end{Theorem}

Here the change of variable from $\frac\Lambda{a}$ to $q^{1/8}$ will
be explained later during the proof. It will be done in several steps
in \S\S\ref{subsec:pm1},\ref{subsec:a0toinfty},\ref{subsec:changevar}.
The explicit forms stated in (2) will be given in
\subsecref{subsec:explicit}.

\begin{NB}
If I am not mistaken, the proof shows actually more: one does not need any 
assumption on the parities, however the statement changes:
When either $\rk(v)$ or $\langle\xi, c_1(v)+\frac{\rk(v)}2
(c_1-K_X)\rangle$ is not even, the coefficient of
$\Lambda^{4n-\xi^2-3}$ in the left hand side is $0$ or
equal to the coefficient in the right hand side, depending on
the parity of $(-1)^{\langle\xi, c_1(v)+\frac{\rk(v)}2
(c_1-K_X)\rangle+{\rk(v)} n}$. I do not know whether one should explicitly state
that, because it is a stronger result, or whether you find it too clumsy.
5.9. LG 
\end{NB}

\begin{NB}
I have changed the statement. I hope my statement is not too
clumsy....
Oct. 19. H.N.
\end{NB}

For $\rk(v) = \pm 1$, the terms are written in terms of the
Seiberg-Witten prepotential $\F_0$, but we do not know how to write
them {\it explicitly\/} in terms of $q^{1/8}$ and $\Lambda$ at this
moment. This is a problem is about elliptic integrals and modular forms.
\begin{NB}
I add a comment on $m=1$ case. I am not sure I explained well....
Oct. 26, H.N. Changed slightly 15.11.LG
\end{NB}

\subsection{From the residues at $e^{\bbeta a} = 0,\infty$ to
the residue at $e^{\bbeta a} = 1$}\label{subsec:pm1}

Let 
\(
   \Delta_{\xi,e^{-2\bbeta a}}^X(v,\bbeta\Lambda)
   = \sum_{n\ge 0} \Delta_n \Lambda^{4n-\xi^2-3}.
\)
\begin{NB}
Corrected. $n\to 4n$. Nov. 9, HN  
\end{NB}

\begin{Proposition}\label{prop:regular}
\textup{(1)} 
The coefficient $\Delta_n$ is a rational function in $e^{\bbeta a}$%
\begin{NB}
$= e^{\bbeta t/2}$
\end{NB},
which is regular on $\proj^1\setminus \{0,\infty,1,-1\}$.

\textup{(2)} $\Delta_n$ is multiplied by $(-1)^{\langle\xi,
c_1(v)+\frac{\rk(v)}2 (c_1-K_X)\rangle+{\rk(v)} n}$ under the replacement
$e^{\bbeta a} \mapsto - e^{\bbeta a}$.
\end{Proposition}

\begin{NB}
Probably it is not true that 
$\Delta_{\xi,e^{-\bbeta t}}^X(v,\bbeta\Lambda)$ is a function in
$e^{-\bbeta t}$, contrary to the notation.
\end{NB}
\begin{NB} In fact the point is that with my definition 
$\Delta_{\xi,T}^X$ is a function of $T^{\frac{1}{2}}$ and not of $T$.
Maybe the correct thing would be to change $T$ to $T^2$ in my definition.
5.9. LG
\end{NB}

\begin{Corollary}\label{cor:formula2}
Assume $\rk(v)$ and $\langle\xi, c_1(v)+\frac{\rk(v)}2
(c_1-K_X)\rangle$ are even. Then
\begin{equation*}
\begin{split}
  \Delta^X_\xi(v;\bbeta\Lambda)
&=
\Res_{e^{-\bbeta a}=0} \Delta^X_{\xi,e^{-2\bbeta a}}(v;\bbeta\Lambda) 
\frac{d e^{-\bbeta a}}{e^{-\bbeta a}}
+
\Res_{e^{-\bbeta a}=\infty} \Delta^X_{\xi,e^{-2\bbeta a}}(v;\bbeta\Lambda) 
\frac{d e^{-\bbeta a}}{e^{-\bbeta a}}
\\
&=
- 2 \Res_{e^{-\bbeta a}=1} \Delta^X_{\xi,e^{-2\bbeta a}}(v;\Lambda) 
\frac{d e^{-\bbeta a}}{e^{-\bbeta a}}.
\end{split}
\end{equation*}
\end{Corollary}

The first equality follows from (1) (and $T = e^{-\bbeta t} =
e^{-2\bbeta a}$). The second equality follows from (1) and the residue
theorem, together with (2). This corollary means that we can move the
position taking residues from $0$, $\infty$ to $1$.

When either $\rk(v)$ or $\langle\xi, c_1(v)+\frac{\rk(v)}2
(c_1-K_X)\rangle$ is not even, the coefficient of
$\Lambda^{4n-\xi^2-3}$%
\begin{NB} $n\to 4n$ 15.11 LG\end{NB}
in the left hand side is $0$ or
equal to the coefficient in the right hand side, depending on
the parity of $(-1)^{\langle\xi, c_1(v)+\frac{\rk(v)}2
(c_1-K_X)\rangle+{\rk(v)} n}$.
We assume that both are even for brevity in the above corollary, but
it is clear that we have a statement like in
\thmref{thm:q-expansion}(1).
\begin{NB}
I have changed the explanation a little bit. Oct. 19, H.N.
\end{NB}

Before starting the proof of \propref{prop:regular} we give new
variables so that the partition function becomes homogeneous.

Recall we set $a_1 = -a$, $a_2 = a$. Following \cite[\S5]{NY3}, we set
\begin{equation*}
   \zeta_{\alpha,\beta}:=\frac{\bbeta}{1-e^{-(a_\alpha-a_\beta)\bbeta}}.
\end{equation*}
We first consider the case when the 5D Chern-Simons term is {\it
  not\/} included. 

By \cite[(5.3)]{NY3}%
\begin{NB}
({\bf Correction} to \cite{NY3}: Replace
$\prod_{\alpha<\beta}\zeta_{\alpha,\beta}\Lambda$ by
$(\prod_{\alpha<\beta}\zeta_{\alpha,\beta}) \Lambda^{2r}$ in two line
above of (5.3), (5.3), and (5.5). Note also that $\Fin_K$ in
\cite{NY3} is $\ve_1\ve_2\Fin_K$ here.)
\end{NB}
we have
\begin{equation*}
   \ve_1\ve_2 \Fin_K \in \C[\zeta_{1,2},\zeta_{2,1},\bbeta]
   [[\ve_1,\ve_2,\zeta_{1,2}\Lambda^4]].
\end{equation*}
We assign degrees as $\deg \ve_1=\deg \ve_2=\deg \Lambda=1$ and $\deg
\bbeta=\deg \zeta_{\alpha,\beta}=-1$. Then $\Zin_K$ is homogeneous of
degree $0$, and hence $\ve_1\ve_2\Fin_K$ is of degree
$2$. Let 
\[
   \Finz_0 := \left.\ve_1\ve_2 \Fin_K\right|_{\ve_1=\ve_2=0}
   = \sum_{n \ge 1} \Finz_n (\bbeta\Lambda)^{4n}.
\]
Then the coefficient $\Finz_n$ is a homogeneous polynomial of $\bbeta$
and $\zeta_{\alpha,\beta}$ of degree $2-4n$.
When we exchange $a_1$ and $a_2$, $\zeta_{2,1}$ and $\zeta_{1,2}$ are
exchanged accordingly. Since $\Finz_n$ is symmetric in $a_1$, $a_2$,
$\Finz_n$ is symmetric in $\zeta_{1,2}$ and $\zeta_{2,1}$.  By the
equality $\zeta_{2,1}=\bbeta-\zeta_{1,2}$, we see that there is a
weighted homogeneous polynomial $A_{4n-2}(x,y) \in {\mathbb C}[x,y]$
of degree $4n-2$ with $\deg x=1$ and $\deg y=2$ such that
\begin{equation*}
\Finz_n =A_{4n-2}(\bbeta,\zeta_{1,2}\zeta_{2,1}).
\end{equation*}
Moreover, as $\Fin_K$ is a formal power series in
$\zeta_{1,2}\Lambda^4$ by \cite[(5.3)]{NY3}, $\Finz_n$ is divisible by
$\left(\zeta_{1,2}\zeta_{2,1}\right)^n$.

We further introduce
\begin{equation*}
   z :=\frac{-\sqrt{-1}\bbeta \Lambda}{e^{\bbeta a_1} - e^{\bbeta a_2}}.
\end{equation*}%
\begin{NB}
My $z$ is Kota's $z\Lambda$. Aug.29 H.N.
\end{NB}
We have $z^2=\zeta_{1,2}\zeta_{2,1}\Lambda^2$. From the above
consideration we have
\begin{equation}\label{eq:zeta'}
   \Finz_0 \in z^2\Lambda^2\C[\bbeta,\Lambda][[z^2]].
\end{equation}
\begin{NB}
\begin{equation*}
\begin{split}
   & A_{4n-2}(\bbeta,\zeta_{1,2}\zeta_{2,1})(\bbeta\Lambda)^{4n}
   = \left(a (\zeta_{1,2}\zeta_{2,1})^{2n-1} +
   b(\zeta_{1,2}\zeta_{2,1})^{2n-2}\bbeta^2 + \cdots
   + c(\zeta_{1,2}\zeta_{2,1})^n \bbeta^{2n-2}\right)
   \times (\bbeta\Lambda)^{4n}
\\
   = \; & a z^{4n-2} \bbeta^{4n}\Lambda^2 + b z^{4n-4}\bbeta^{4n}\Lambda^4
   + \cdots + c z^{2n} \bbeta^{4n}\Lambda^{2n}.
\end{split}
\end{equation*}
Sep. 1, H.N. 
\end{NB}
As
\(
   \pd{}{a}z = -z (\zeta_{1,2}-\zeta_{2,1}),
\)
and
\(
   \pd{}{a}(\zeta_{1,2}-\zeta_{2,1})
   = 4 (z/\Lambda)^2,
\)
we have
\begin{equation}\label{eq:zeta''}
\begin{split}
   \pd{\Finz_0}{a} &\in (\zeta_{1,2}-\zeta_{2,1}) z^2 \Lambda^2
   \C[\bbeta,\Lambda][[z^2]],
\\
   \pd{\Finz_0}{a^2} &\in z^2 \C[\bbeta,\Lambda][[z^2]].
\end{split}
\end{equation}
Even when we include the 5d Chern-Simons term \eqref{eq:ZinstmF}, we
can repeat the above proof. We only need to change
$\C[\bbeta,\Lambda]$ by $\C[\bbeta,\Lambda,e^{\pm \rk(v)\bbeta a}]$.

\begin{proof}[Proof of \propref{prop:regular}]
Let us look at the expression of $\Delta_{\xi,e^{-2\bbeta
a}}^X(v,\bbeta\Lambda)$ given in \corref{cor:formula1}. We will write
it as a multiple of an explicit rational function in $e^{\bbeta a}$
and a formal power series in $z$. The explicit function comes from the
perturbation part of the partition function.

First note that
\(
  \frac{\partial^2\F_0}{(\partial \log\Lambda)^2}
\)
consists only of the instanton part. Therefore \eqref{eq:zeta'}
implies 
\begin{equation*}
   \frac{\partial^2\F_0}{(\partial \log\Lambda)^2}
   \in z^2 \Lambda^2\C[\bbeta,\Lambda,e^{\pm \rk(v)\bbeta a}][[z^2]].
\end{equation*}
Next we have
\begin{equation}\label{eq:q_expansion}
\begin{split}
q^{{1}/{8}}&= \left(\frac{-\sqrt{-1}\bbeta \Lambda}
{e^{-\bbeta a}-e^{\bbeta a}}\right)
\exp\left(-\frac{1}{8}\frac{\partial^2 \Finz}{\partial a^2}\right)
\in z\left(1+z^2\C[\bbeta,\Lambda,e^{\pm \rk(v)\bbeta a}][[z^2]]\right).
\end{split}
\end{equation}
from \eqref{eq:zeta''}.%
\begin{NB}
Corrected $a_1 \to -a$, $a_2 \to a$  by H.N. Oct. 19
\end{NB}
\begin{NB} Maybe that is too trivial, but would it be good to give a reference
to a formula or a sketch, why  $\left(\frac{-\sqrt{-1}\bbeta \Lambda}
{e^{\bbeta a_1}-e^{\bbeta a_2}}\right)$ is the contribution of the perturbation term.
5.9. LG
\end{NB}
Here we have used
\begin{equation*}
   \overline\gamma_0''(x|\bbeta;\Lambda)=2 \log\left(
   \frac{-\sqrt{-1}\bbeta\Lambda}{
     e^{\bbeta x/2} - e^{-\bbeta x/2}}\right)
\end{equation*}
(cf.\ \eqref{eq:gamma'}) to calculate the first term coming from the
perturbation part.

\begin{NB}
I add an explanation according to Lothar's suggestion. Oct. 19 H.N.
\end{NB}

Next consider the genus $1$ parts. When $\rk(v)=0$, we have
\begin{equation*}
  \exp (\chi A + \sigma B) = 
  \left(\frac2{\theta_{00}\theta_{10}}\right)^2
  \theta_{01}^\sigma \in 4z^{-2}\left(1+z^2\C[\bbeta,\Lambda,e^{\pm \rk(v)\bbeta
a}][[z^2]]\right).
\end{equation*}
The case $\rk(v)=\pm 1$ is similar thanks to \eqref{eq:AB}.

Finally again by \eqref{eq:zeta''} we have
\begin{equation*}
\begin{split}
& \exp
\left(
  \frac\bbeta8 \frac{\partial^2 {\mathcal F_0}}
  {\partial a \partial \log\Lambda} 
\langle\xi, K_X+c_1(v)+\frac{\rk(v)}2 (c_1-K_X)\rangle
\right)
=
\left(e^{-\bbeta a}\right)^{
N
}
\exp 
\left( N \frac\bbeta8 \frac{\partial^2 \Finz_0}
  {\partial a \partial \log\Lambda} 
\right)
\\
\in\; &
\left(e^{-\bbeta a}\right)^{
N
}
\C[\bbeta,\Lambda,e^{\pm \rk(v)\bbeta a}][[(\zeta_{1,2}-\zeta_{2,1})z^2\Lambda^2,z^2]],
\end{split}
\end{equation*}
with $N = \langle\xi, K_X+c_1(v)+\frac{\rk(v)}2 (c_1-K_X) 
\rangle$.%
\begin{NB}
The formula was corrected according to Lothar's message.
Oct.30, H.N.  
\end{NB}
Note that
$\langle\xi, c_1-K_X\rangle \equiv
\langle\xi,\xi-K_X\rangle \equiv 0 \mod 2$, where the first equality
follows from the assumption (\subsecref{wallmod}(2)), and the second
from the Riemann-Roch theorem. Therefore $N$ is an integer. 

As $z = \frac{-\sqrt{-1}\bbeta\Lambda}{e^{-\bbeta a} - e^{\bbeta a}}$,
$\zeta_{1,2}-\zeta_{2,1}=-\bbeta \frac{e^{2\bbeta a}+1}{e^{2\bbeta a}-1}$, the
statement~(1) becomes clear now.

Let us check the statement~(2). We substitute $e^{\bbeta a}$ by
$-e^{\bbeta a}$. Then $z$ changes the sign and $\zeta_{1,2}$,
$\zeta_{2,1}$ are invariant. Therefore the change of the instanton
part of $\Delta_{\xi,e^{-2\bbeta a}}^X(v,\bbeta\Lambda)$ comes only
from $C^{\vec{Y}}_{-\rk(v)}(\ve_1,\ve_2,a;\bbeta,\tau)$ in
\eqref{eq:ZinstmF}. It is multiplied by
\(
   (-1)^{\rk(v)(|Y^1|+|Y^2|)}.
\)
The perturbation part of $\Delta_{\xi,e^{-2\bbeta
a}}^X(v,\bbeta\Lambda)$ is multiplied by
\begin{equation*}
   (-1)^{N + \langle\xi^2\rangle}
   = (-1)^{\langle\xi, c_1(v)+\frac{\rk(v)}2 (c_1-K_X)\rangle}.
\end{equation*}
Altogether the coefficient of $\Lambda^{4n-\xi^2-3}$%
\begin{NB} $n\to 4n$ 15.11. LG\end{NB}
 in $\Delta_{\xi,e^{-2\bbeta
a}}^X(v,\bbeta\Lambda)$ is multiplied by
\begin{equation*}
  (-1)^{\langle\xi, c_1(v)+\frac{\rk(v)}2 (c_1-K_X)\rangle+{\rk(v)} n}.
\end{equation*}
\end{proof}

\subsection{From the expansion at $a=0$ to $a=\infty$}\label{subsec:a0toinfty} 

We set $\bbeta = 1$ hereafter.

We expand $\Delta^X_{\xi,e^{-2 a}}(v,\Lambda)$ at $a = 0$:
\begin{equation}\label{eq:expansion}
   \Delta_{\xi,e^{-2 a}}^X(v;\Lambda)
   = \sum_{\substack{n\ge 0 \\ m\in\Z}} 
   \Delta_{m,n} a^m \Lambda^{4n-\xi^2-3}
   \in \Lambda^{-\xi^2-3}\C((a))[[\Lambda]].
\end{equation}
Then
\begin{equation*}
   \Res_{e^{a}=1} \Delta^X_{\xi,e^{-2a}}(v;\Lambda) 
   \frac{d e^{a}}{e^{a}}
   = \Coeff_{(a)^0}
   \left[
     \Delta^X_{\xi,e^{-2 a}}(v;\Lambda) \times a\right]
   = \sum_n \Delta_{-1,n} \Lambda^{4n-\xi^2-3}.
\end{equation*}
\begin{NB}
Corrected (twice). $n\to 4n$. Nov. 9, HN  
\end{NB}

\begin{Proposition}
 $\Delta^X_{\xi,e^{-2 a}}(v;\Lambda)$ is in
$\Lambda^{-\xi^2-3}\C[[\frac\Lambda{a},a]]$, i.e.\ $\Delta_{m,n}=0$
unless $m\ge -n$ in \eqref{eq:expansion}.    
\end{Proposition}

This is a consequence of the proof of \propref{prop:regular}. The key
observation is that $z$, $(\zeta_{1,2}-\zeta_{2,1}) \Lambda \in
\frac{\Lambda}a \C[[a]]$.

We rewrite the above expansion as
\begin{equation*}
  \Delta^X_{\xi,e^{-2 a}}(v;\Lambda) \times a
=  \sum_{\substack{n\ge 0 \\ m+n\ge 0}} 
   \Delta_{m,n} a^{m+1} \Lambda^{4n-\xi^2-3}
= \sum_{\substack{n\ge 0 \\ m+n\ge 0}} 
   \Delta_{m,n} \left(\frac{\Lambda}{a}\right)^{-m-1} \Lambda^{4n+m+1-\xi^2-3}.
\end{equation*}%
\begin{NB}
Corrected (twice). $n\to 4n$. Nov. 9, HN  
\end{NB}
The last expression is an element in 
\( 
  \Lambda^{-\xi^2-2} \C((\frac{\Lambda}{a}))[[\Lambda]],
\)
and
\(
  \sum_n \Delta_{-1,n} \Lambda^{4n-\xi^2-3}
\)%
\begin{NB}
Corrected. $n\to 4n$. Nov. 9, HN  
\end{NB}
is equal to its coefficient of
$\left(\frac{\Lambda}{a}\right)^0$. Thus we get

\begin{Corollary}\label{cor:step2}
\begin{equation*}
   2 \Res_{e^{a}=1} \Delta^X_{\xi,e^{-2a}}(v;\Lambda) 
   \frac{d e^{a}}{e^{a}}
   = 2 \Coeff_{(\frac{\Lambda}a)^0}
   \left[
     \Delta^X_{\xi,e^{-2 a}}(v;\Lambda) \times \frac{a}\Lambda \Lambda\right].
\end{equation*}
\end{Corollary}

\subsection{From $a=\infty$ to $q=0$}\label{subsec:changevar}

By \eqref{eq:q_expansion} we have the following expansion in 
$\C[[\frac\Lambda{a},a]]$:
\begin{equation*}
   q^{1/8} = \frac{\sqrt{-1}\Lambda}{2a}
   \left(1 + O(a,\frac\Lambda{a})\right).
\end{equation*}%
\begin{NB}
Kota also observed that $\frac{\Lambda}{q^{1/8}} =
\frac{2a}{\sqrt{-1}}(1+ O(a,\frac\Lambda{a}))$. 
H.N.  Aug. 29 
\end{NB}
As in the previous subsection, we consider this as an element in
$\C((\frac\Lambda{a}))[[\Lambda]]$. Then we have
\begin{equation*}
\begin{split}
   q^{1/8} &= q_0(\frac\Lambda{a}) + q_1(\frac\Lambda{a})\Lambda
   + \cdots,
\\
   q_0(\frac\Lambda{a}) &= \frac{\sqrt{-1}\Lambda}{2a}
   + a_2 \left(\frac\Lambda{a}\right)^2 
   + a_3 \left(\frac\Lambda{a}\right)^3  + \cdots .
\end{split}
\end{equation*}
From this we see that $\C((\frac\Lambda{a}))[[\Lambda]] \cong
\C((q_0))[[\Lambda]]\cong \C((q^{1/8}))[[\Lambda]]$.%
\begin{NB}
Kota also observed
$\C[[\frac\Lambda{a},a]]\cong\C[[q^{1/8},\frac{\Lambda}{q^{1/8}}]]$ from
the above expansion.
H.N. Aug.29
\end{NB}
We now change the variable from
$\frac{\Lambda}{a_2}$ to $q^{1/8}$ by the following lemma:
\begin{Lemma}\label{lem:change}
Let us consider the change of the variable from $x$ to $y$ given by
$y = y(x,\Lambda) = y_0(x) + y_1(x)\Lambda + \cdots
\in \C((x))[[\Lambda]]$. Assume $y_0(x) = x + a_2 x^2 + \cdots
\in x(1 + x \C[[x]])$.
Let $f(y,\Lambda)\in \C((y))[[\Lambda]]\cong \C((x))[[\Lambda]]$.
Then
\begin{equation*}
   \Coeff_{y^0} \left[ y f(y,\Lambda) \right]
   = \Coeff_{x^0} \left[ x f(y(x,\Lambda),\Lambda) \frac{dy}{dx} \right].
\end{equation*}
\end{Lemma}

This lemma just means the invariance of the residue under the change
of variables. As we have an extra parameter $\Lambda$ which does not
appear in the usual setting, we give a proof.

\begin{proof}
It is enough to check the case $f(y,\Lambda) = y^{m-1}$ for $m\in \Z$.
First suppose $m\neq 0$. Then the left hand side is equal to $0$. On
the other hand,
\begin{equation*}
      y(x,\Lambda)^{m-1} \frac{dy}{dx}
  =  \frac1m \frac{d}{dx} \left( y(x,\Lambda)^m \right)
\end{equation*}
does not contain the term $x^{-1}$, as it is a derivative of a formal
power series in $x$. Therefore the right hand side is also $0$.
\begin{NB}
We have used the Leibniz rule here (and the induction on $m$). To
check it, we cut $y$ modulo $\Lambda^N$, apply the usual Leibniz
rule for polynomial and make $N\to \infty$. Sep.3., H.N.
\end{NB}

Next suppose $m=0$. Then the left hand side is $1$. Let us consider 
\begin{equation*}
   \log \frac{y(x,\Lambda)}{y_0(x)}
   = \log\left(
     1 + \frac{y_1(x)}{y_0(x)}\Lambda + \cdots
     \right).
\end{equation*}
This is well-defined in $\C((x))[[\Lambda]]$. Then we have
\begin{equation*}
\begin{split}
   \frac1{y(x,\Lambda)} \frac{dy}{dx}
   & = \frac1{y_0(x)} \frac{d y_0(x)}{dx}
   + \frac{d}{dx} \left\{\log\left(
     1 + \frac{y_1(x)}{y_0(x)}\Lambda + \cdots
     \right)\right\}
\\
   &= \frac1x \left(1+a_2 x + \cdots\right)^{-1}
   \left(1 + 2 a_2 x + \cdots\right) 
   + \frac{d}{dx} \left\{\log\left(
     1 + \frac{y_1(x)}{y_0(x)}\Lambda + \cdots
     \right)\right\}.
\end{split}
\end{equation*}
The second term does not contain the term $x^{-1}$ by the same reason
as above. Therefore we get $x^{-1}$ only from the first
term. Hence we have found that the right hand side is also equal to $1$.
\end{proof}

Applying this to the right hand side of \corref{cor:step2} we get
\begin{equation*}
  2 \Coeff_{(\frac{\Lambda}a)^0}
   \left[
     \Delta^X_{\xi,e^{2 a}}(u;\Lambda) \times a\right]
   = 
  2 \Coeff_{(q^{1/8})^0}
   \left[\left.
     \Delta^X_{\xi,e^{2 a}}(u;\Lambda) 
     \times \left(\frac{a}\Lambda\right)^2 \Lambda
    \right|_{\frac\Lambda{a}=\frac\Lambda{a}(q^{1/8},\Lambda)}
    q^{1/8}
    \frac{d(\frac\Lambda{a})}{d(q^{1/8})}\right].
\end{equation*}
This completes the proof of \thmref{thm:q-expansion}(1).

\begin{NB}
Originally the following is attached at the end of the proof. But is
deleted. Oct. 20, H.N.

Finally substituting 
\begin{equation*}
   \frac{d \left(\frac\Lambda{a}\right)}{d(q^{1/8})}
   = - \frac\Lambda{a^2} \frac{da}{d(q^{1/8})}
\end{equation*}
into this equation, we get \thmref{thm:q-expansion}.
\verb+\begin{NB}+
Note   
\begin{equation*}
   a = \frac{a}{\Lambda}\Lambda \in \C((\frac\Lambda{a}))[[\Lambda]].
\end{equation*}
We apply the Leibnitz rule to $a\times \frac\Lambda{a} = \Lambda$
(constant) to get the above equality.

Kota used the smaller ring $\C[[q^{1/8},\frac{\Lambda}{q^{1/8}}]]$,
but I do not think that we need it.

H.N. Aug. 29
\verb+end{NB}+
\end{NB}

\begin{NB}
The following subsection was moved to the end of the section for the SW
curve.

\subsection{rank $2$ case}
We assume $r=2$ in this subsection.

.....  
\end{NB}

\subsection{Explicit expressions}\label{subsec:explicit}

Our remaining task is to express the terms in $[\ \ ]$ of the right hand
side of \thmref{thm:q-expansion} in explicit forms in
$\C((q^{1/8}))[[\Lambda]]$. We suppose $m = -\rk(v) = 0$ in this subsection.

By \cite[\S5]{NY3} $\exp(\chi A + \sigma B)$ can be written explicitly
in terms of $q^{1/8}$. So we only need to express
$q^{1/8}\frac{\partial{(\Lambda/a)}}{\partial (q^{1/8})}$,
$\frac{\partial^2\F_0}{\partial a\partial\log\Lambda}$,
and
$\frac{\partial^2\F_0}{(\partial\log\Lambda)^2}$.
The expressions will be given in \eqref{eq:dadq}, \eqref{eq:h},
\eqref{eq:T} respectively.

\subsubsection{The term $q^{1/8}\frac{\partial{(\Lambda/a)}}{\partial (q^{1/8})}$}
We consider $a$ defined as a period of the Seiberg-Witten curve as in
\secref{sec:SWcurve}. In particular, we are in the region $D^*$ such
that $\sqrt{-1}a$ has a large real part and $0 < |\Lambda| \ll 1$. We
will compute $q^{1/8}\frac{\partial{(\Lambda/a)}}{\partial (q^{1/8})}$ first in
this region and then see later that the computation holds in
$\C((q^{1/8}))[[\Lambda]]$.

For simplicity we introduce a variable $u$ by
\begin{equation*}
   u := - \frac{\theta_{00}^4 +
   \theta_{10}^4}{\theta_{00}^2\theta_{10}^2}\bbeta^2\Lambda^2
   \in \bbeta^2 \Lambda^2\C((q^{1/8}))
\end{equation*}
where $\theta$-functions are evaluated at $(0,\tau)$. This definition is
motivated by a fundamental variable in the homological version (see
\cite[(4.1)]{GNY}).%
\begin{NB}
I have multiplied $u$ by $\bbeta^2\Lambda^2$ from Kota's note.
Aug. 31, H.N. 
\end{NB}
By \eqref{eq:U_1} we have
\begin{equation*}
   U_1 = \pm 2 \sqrt{1 + u + \bbeta^4\Lambda^4}.
\end{equation*}
By a certain standard equality for $\theta$-functions (cf.\
\cite[p.29]{GNY}) we have
\begin{equation*}
   \frac{d u}{d\tau} 
   = -\frac{\bbeta^2\Lambda^2 \pi}{2\sqrt{-1}}
   \frac{\theta_{01}^8}{\theta_{00}^2\theta_{10}^2}.
\end{equation*}
Combining this with \eqref{eq:dadU}, we get
\begin{equation*}
   \frac{d a}{d \tau} = \frac{d a}{d U_1}
   \frac{d U_1}{d\tau} = \pm \frac{\pi\Lambda}{4}
   \frac{\theta_{01}^8}{\theta_{00}\theta_{10}}
   \frac1{\sqrt{1+ u + \bbeta^4\Lambda^4}}.
\end{equation*}
Therefore we have
\begin{equation}\label{eq:dtauda}
   \left(\frac{d\tau}{da}\right)^2
   = \frac{16}{\pi^2\Lambda^2}
   \frac{\theta_{00}^2\theta_{10}^2}{\theta_{01}^{16}}
   \left(1 + u + \bbeta^4\Lambda^4\right).
\end{equation}
This is {\it a priori\/} an equality on $D^*$. However both
sides extend%
\begin{NB} corrected 30.10 LG\end{NB}
to $\Lambda=0$: The right hand side is a function in
$q^{1/8}$ and we have 
\( 
   q^{1/8} \sim
   \frac{-\sqrt{-1}\bbeta\Lambda}{e^{\bbeta a_1}-e^{\bbeta a_2}}.
\)
Here $\sim$ means the equality up to the instanton part.
\begin{NB} did we introduce the notation $\sim$?
5.9 LG
\end{NB}
\begin{NB}
I have introduced $\sim$.
Oct. 19, H.N.  
\end{NB}

Therefore $\theta_{10}/\Lambda$ and $u$ are regular at $\Lambda=0$,
hence so is the right hand side.
The left hand side is a triple derivative of the prepotential with
respect to $a$, and hence has no perturbation part. Thus it is
regular at $\Lambda=0$.
Therefore \eqref{eq:dtauda} holds even at $\Lambda=0$.

Next we consider the coefficients of $\Lambda^k$ for both sides%
\begin{NB} I think it should be always "both sides" and not both
hand sides 5.9. LG \end{NB}
of \eqref{eq:dtauda}. The equality holds {\it a priori\/} for $a$ such
that $\sqrt{-1}a$ has a large real part. However  both sides
are rational functions in $e^{\bbeta a}$: This claim can be checked as
above. The left hand side has no perturbation part, so the claim was
proved during the proof of \propref{prop:regular}. The right hand side
is a function in $q^{1/8}$, hence the claim was again proved during
the proof of \propref{prop:regular}. Considering the expansion at
$a=0$, we conclude that \eqref{eq:dtauda} holds in $\C((a))[[\Lambda]]$.

From the discussion in \subsecref{subsec:a0toinfty} we see that
both sides of \eqref{eq:dtauda} are in
$\frac1{a^2}\C[[a,\frac\Lambda{a}]]$. Therefore \eqref{eq:dtauda}
holds in $\frac1{a^2}\C[[a,\frac\Lambda{a}]]$, and
hence in $\C((\frac{\Lambda}{a}))[[\Lambda]]$. We now change the
variable from $\Lambda/a$ to $q^{1/8}$ as in
\subsecref{subsec:changevar} and use the composition law%
\begin{NB}
The proof of the composition law (it is already used in the proof of
\lemref{lem:change}):

Under the setting in \lemref{lem:change} we want to prove
\begin{equation*}
  \frac{d f(y(x,\Lambda),\Lambda)}{dx}
  = \frac{d f(y,\Lambda)}{dy} \frac{dy}{dx}
\end{equation*}
for $x = \frac\Lambda{a}$, $y=q^{1/8}$, $f(y,\Lambda) =
\left(\frac\Lambda{a}\right)(q^{1/8})$. We may assume $f(y,\Lambda) =
y^m$ cutting $f$ modulo $\Lambda^N$ for $N\gg 0$ as before. Then make
$N\to \infty$. Finally the formula for $f(y,\Lambda) = y^m$ follows
from Leibniz rule and the induction on $m$. 

Sep. 3, H.N.

The following was what Kota added to my NB above. (The main comment is
already treated. I am not sure I understood what Kota said, but I
changed the change of variable $a\to \frac{\Lambda}a$ at the end.
Oct. 21, H.N.

By the way, it might be better to show
$d\tau/d(\Lambda/a)=(d\tau/da)(-\Lambda(a/\Lambda)^2)$
first (this is clear), and apply Lemma 1.10 to
$q$ and $\Lambda/a$. Thus we remove the last sentence in section 1.3.
2006/9/5, K.Y.  
\end{NB}
to get
\begin{equation*}
\left( \frac{a^2}{\Lambda}q^{1/8}
  \frac{d\left(\frac{\Lambda}a\right)}{d(q^{1/8})}\right)^2 
= \left(-\sqrt{-1}
  \frac{\theta_{01}^8\Lambda}{\theta_{00}\theta_{10}}\right)^2
  \frac1{1+u+\bbeta^4\Lambda^4}.
\end{equation*}%
\begin{NB}
More precise explanation: as \eqref{eq:dtauda} holds in
$\frac1{a^2}\C[[a,\frac{\Lambda}a]]$, we first observe
\begin{equation*}
  \frac{d \tau}{d \left(\frac{\Lambda}a\right)}
  = - \frac{d\tau}{da} \Lambda \left(\frac{a}\Lambda\right)^2.
\end{equation*}
Therefore
\begin{equation*}
    \left(\frac{d\tau}{d \left(\frac{\Lambda}a\right)}\right)^2
    = \frac{16}{\pi^2}
   \frac{\theta_{00}^2\theta_{10}^2}{\theta_{01}^{16}}
   \left(1 + u + \bbeta^4\Lambda^4\right)
   \left(\frac{a}\Lambda\right)^4.
\end{equation*}
We then use the composition law for $\frac{\Lambda}a \leftrightarrow
q^{1/8}$ to get the above.
Oct. 30, H.N.
\end{NB}
As $a \sim \frac{\sqrt{-1}}2 \frac{\Lambda}{q^{1/8}}$, we can
determine the branch of the square root to get
\begin{equation}\label{eq:dadq}
   \frac{a^2}{\Lambda} q^{1/8} \frac{d\left(\frac{\Lambda}a\right)}{d(q^{1/8})}
   = \sqrt{-1}\frac{\theta_{01}^8\Lambda}{\theta_{00}\theta_{10}}
   \sum_{n\ge 0}\binom{-\frac{1}{2}}{n}
   (u+\bbeta^4\Lambda^4)^n.
\end{equation}
This is an equality in $\C((q^{1/8}))[[\Lambda]]$.
\begin{NB}
Corrected by H.N, Oct. 20  
\end{NB}

\subsubsection{The term $\frac{\partial^2\F_0}{\partial a\partial\log\Lambda}$}
Let
\(
   h := - \frac14 \frac{\partial^2 \mathcal F_0}
                       {\partial a\partial\log\Lambda}
      = \frac{\pi\sqrt{-1}}2\pd{a^D}{\log\Lambda}.
\)
Let us rewrite \eqref{eq:sn} in terms of $\sn$ associated with the
elliptic curve with period $\tau$. (Be aware that we have used $\sn$
with period $-2/\tau$ before.) We get
\begin{equation*}
    - \frac{\theta_{10}}{\theta_{00}}
     \sn(\theta_{00}^2 \frac{\bbeta h}{2\sqrt{-1}},\kappa(\tau))
     = 
     \frac{\theta_{11}(\frac{\bbeta h}{2\pi\sqrt{-1}})}
     {\theta_{01}(\frac{\bbeta h}{2\pi\sqrt{-1}})} = - \bbeta\Lambda.
\end{equation*}
Therefore
\begin{equation*}
   \theta_{00}^2\frac{\bbeta h}{2\sqrt{-1}}
   = \int_0^{\frac{\theta_{00}}{\theta_{10}}\bbeta\Lambda}
   \frac{dx}{\sqrt{(1-x^2)(1 - \kappa^2 x^2)}}
   = \frac{\theta_{00}}{\theta_{10}} \int_0^{\bbeta\Lambda}
    \frac{dx}{\sqrt{1 + \frac{u}{\bbeta^2\Lambda^2} x^2 + x^4}}.
\end{equation*}
Therefore
\begin{equation*}
  h = 
  \frac{2\sqrt{-1}}{\bbeta \theta_{00} \theta_{10}}
  \int_0^{\bbeta\Lambda}
    \frac{dx}{\sqrt{1 + \frac{u}{\bbeta^2\Lambda^2} x^2 + x^4}}.
\end{equation*}
Using
\(
   \frac1{\sqrt{1 + \frac{u}{\bbeta^2\Lambda^2} x^2 + x^4}}
   = \sum_{n\ge 0, n\ge k\ge 0} \binom{-\frac{1}{2}}{n}
   \binom{n}{k}
  \left(\frac{u}{\bbeta^2\Lambda^2}\right)^k x^{4n-2k},
\)
we get 
\begin{equation}\label{eq:h}
    h = \frac{2\sqrt{-1}}{\bbeta \theta_{00} \theta_{10}}
    \sum_{\substack{n\ge 0 \\ n\ge k\ge 0}}
    \binom{-\frac{1}{2}}{n}\binom{n}{k}
    \frac{u^k (\bbeta\Lambda)^{4(n-k)+1}}{4n-2k+1}.
\end{equation}%
\begin{NB} 
A slight simplification is 
\begin{equation}\label{eq:h}
    h = \frac{2\sqrt{-1}}{\bbeta \theta_{00} \theta_{10}}
    \sum_{\substack{n\ge 0 \\ n\ge k\ge 0}}
    \binom{-\frac{1}{2}}{n}\binom{n}{k}
    \frac1{4n-2k+1} 
    u^k (\bbeta\Lambda)^{4(n-k)+1}.
\end{equation}
5.9. LG

I have changed the expression according to Lothar's comment. Nov. 6, HN
\end{NB}
This gives us an explicit expression in terms of $q^{1/8}$ as, e.g.\ 
\begin{equation*}
   \bbeta\theta_{00}\theta_{10} h
   = 2\sqrt{-1} \left(\bbeta\Lambda - \frac{u}6 \bbeta\Lambda +
     \cdots\right). 
\end{equation*}

\begin{NB}
It is better to study the behaviour of $h$ as in Kota's note (at
least in \verb+NB+). Aug. 31, H.N. 
$-\bbeta \Lambda$ in the integral seems to be $\bbeta \Lambda$
and the sign of $u/(\bbeta\Lambda)^2$ seems to be different from my note.. 
2006/9/1, K.Y.
I have corrected the two signs, and add an explicit formula of the expansion. Sep.
3, H.N.
\end{NB}

\begin{NB}
As $\bbeta\to 0$, we have
\begin{equation*}
   h\to \frac{du}{da} = \frac{2\sqrt{-1}}{\theta_{00}\theta_{10}}\Lambda.
\end{equation*}
We pick up the first term of the above expansion.
Sep.\ 3. H.N.
\end{NB}

\subsubsection{The term $\frac{\partial^2\F_0}{(\partial\log\Lambda)^2}$.}
We use
\begin{equation*}
    \frac{\theta_{11}(\frac{\bbeta h}{2\pi\sqrt{-1}},\tau)}
    {\theta_{01}(0,\tau)} 
    = 
        \frac{\theta_{11}(\frac{\bbeta h}{2\pi\sqrt{-1}},\tau)}
    {\theta_{01}(\frac{\bbeta h}{2\pi\sqrt{-1}},\tau)} 
    \frac{\theta_{01}(\frac{\bbeta h}{2\pi\sqrt{-1}},\tau)}
    {\theta_{01}(0,\tau)} 
    = -\bbeta\Lambda 
    \exp\left(\frac{\bbeta^2}{32} \frac{\partial^2\mathcal F_0}
      {(\partial\log\Lambda)^2}\right),
\end{equation*}
where the second equality follows from \eqref{eq:sn} and
\eqref{eq:contact}. We use the formula (see \cite[21$\cdot$43]{W-W}):
\begin{equation*}
    \frac{\theta_{11}(z,\tau)}
    {\theta_{11}'(0,\tau)} 
    = z \exp\left(-\sum_{k=1}^\infty \frac{G_{2k}(\tau)}{2k}z^{2k}\right),
\end{equation*}
where $G_{2k} = 2\zeta(2k)E_{2k}$ are Eisenstein series, and $E_{2k}$
are normalized Eisenstein series. Using Jacobi's derivative formula
(\cite[21$\cdot$41]{W-W}), we get
\begin{equation}\label{eq:T}
    \frac{\bbeta^2}{32} \frac{\partial^2\mathcal F_0}
      {(\partial\log\Lambda)^2}
      = \log\left[\frac{\theta_{00}\theta_{10}
          h}{2\sqrt{-1}\Lambda}\right]
      - \sum_{k=1}^\infty \frac{G_{2k}(\tau)}{2k} 
      \left(\frac{\bbeta h}{2\pi\sqrt{-1}}\right)^{2k}.
\end{equation}%
\begin{NB}
The sign in the first term is different from Kota's note.
Please check it. Aug. 31, H.N. I corrected my note. 2006/9/1. K.Y.
\end{NB}
Combining with \eqref{eq:h}, we get an explicit formula of
$\frac{\partial^2\mathcal F_0}
      {(\partial\log\Lambda)^2}$
in terms of $q^{1/8}$. For example, we have
\begin{equation*}
    \frac{\bbeta^2}{32} \frac{\partial^2\mathcal F_0}
      {(\partial\log\Lambda)^2}
      = -\frac{u}6 + \frac{h^2}{24} E_2 \bbeta^2 + \cdots.
\end{equation*}

\begin{NB}
It is better to study the behaviour of $h$ as in Kota's note (at
least in \verb+NB+). Aug. 31, H.N.

I add a formula. Sep. 3, H.N.
\end{NB}

\subsection{Explicit computations: the case of $\P^2$}\label{subsec:P2}
Let $H$ be the hyperplane bundle on $\P^2$, we denote by the same
letter its first Chern class. 
As an illustration of our results we  compute the holomorphic Euler characteristics of determinant line bundles on 
$M^{\P_2}_H(0,d)$ and $M^{\P_2}_H(H,d)$, and write the corresponding 
Hilbert series explicitly  for small $d$.

The determinant line bundles $\mu(H^{\otimes n})$ are by \eqref{eq:uL} defined on 
$M^{\P_2}_H(0,d)$ for all $n$ and on $M^{\P_2}_H(H,d)$ for $n$ even.
Let $Y$ be the blowup of $\P^2$ in a point, and let $E$
be the exceptional divisor. Denote by $H$ also its pullback to $Y$, and
write $F=H-E$. Then 
for $\epsilon$ sufficiently small 
$M^{Y}_{F+\epsilon H}(E,d+1)=\emptyset$, 
$M^{Y}_{F+\epsilon H}(H,d)=\emptyset$, and thus 
$\chi(M^{Y}_{F+\epsilon H}(E,d+1),\oo(\mu(H^{\otimes n})))=0$, 
$\chi(M^{Y}_{F+\epsilon H}(H,d),\oo(\mu(H^{\otimes 2n})))=0$ for all $n$.
On the other hand we get by  \corref{cor:blow} 
$\chi(M_{H}^{\P^2}(0,d),\oo(\mu(H^{\otimes n})))=\chi(M^{Y}_{H-\epsilon E}(E,d+1),\oo(\mu(H^{\otimes n})))$, 
$\chi(M_{H}^{\P^2}(H,d),\oo(\mu(H^{\otimes 2n})))=\chi(M^{Y}_{H-\epsilon E}(H,d),\oo(\mu(H^{\otimes 2n})))$. Thus we only have to sum the  wallcrossing over all the classes $\xi$ of type $E$ (respectively of type $H$)
with $\<\xi H\> >0>\<\xi F\>$.
These are $\big\{2m H-(2l+1)E\bigm| l\ge m>0\big\}$ for type $E$ and
$\big\{(2m-1) H-2lE\bigm| l\ge m>0\big\}$ for type $H$.

Putting this into \thmref{thm:q-expansion} and using the results of subsection
\secref{subsec:explicit}, and putting $\bbeta=1$, we obtain the following.
\begin{align*}
\sum_{d\ge 0}
\chi(M^{\P^2}_H(0,d),\oo(\mu(H^{\otimes n})))\Lambda^d&=
\Coeff_{q^0}\Bigg[
\sum_{l\ge m>0}
(-1)^{l+m+1} q^{\frac{1}{2}((l+\frac{1}{2})^2-m^2)}
e^{(m(n+3)-l-1/2)h}\\&\qquad\Big(-\frac{\theta_{11}(\frac{h}{2\pi \sqrt{-1}})}{\Lambda \theta_{01}}\Big)^{n^2+6n+8}
\frac{8\theta_{01}^8}{\Lambda\theta_{00}^3\theta_{10}^3}\frac{1}{\sqrt{1+u+\Lambda^4}}\Bigg],\\
\sum_{d\ge 0} \chi(M^{\P^2}_H(H,d),\oo(\mu(H^{\otimes 2n})))\Lambda^d&=
\Coeff_{q^0}\Bigg[
\sum_{l\ge m>0}
(-1)^{l+m} q^{\frac{1}{2}(l^2-(m-\frac{1}{2})^2)}
e^{((m-\frac{1}{2})(2n+3)-l)h}\\&\qquad\Big(-\frac{\theta_{11}(\frac{h}{2\pi \sqrt{-1}})}{\Lambda \theta_{01}}\Big)^{4n^2+12n+8}
\frac{8\theta_{01}^8}{\theta_{00}^3\theta_{10}^3}\frac{1}{\sqrt{1+u+\Lambda^4}}\Bigg].
\end{align*}
It is straightforward to write a maple program which computes the lower order terms in $\Lambda$. This computation can be extended to much  higher degrees in $\Lambda$,  in principle up to any given power. 
We get
\begin{align*}
\sum_{n\ge 0} \chi(M_H^{\P^2}(0,d),\oo(\mu(H^{\otimes n}))t^n&=
\frac{P_d(t)}{(1-t)^{d+1}}, \quad 5\le d\le 21,\\
\sum_{n\ge 0} \chi(M_H^{\P^2}(H,d),\oo(\mu(H^{\otimes 2n}))t^n&=
\frac{Q_d(t)}{(1-t)^{d+1}}, \quad 0\le d\le 24,
\end{align*}
with $P_d(t)\in \Z[t]$ of degree $d-5$ with $t^{d-5}P_d(1/t)=P_d$
and $Q_d(t)\in \Z[t]$ of degree $d-2$ with $t^{d-2}Q_d(1/t)=Q_d$ for $d\ge 4$.
In particular
{\allowdisplaybreaks\small 
\begin{align*}P_5&=1,\ P_9=1+t^2+t^4,\ P_{13}=1+t+7t^2+7t^3+22t^4+7t^5+7t^6+t^7+t^8,\\
P_{17}&=1+3t+27t^2+83t^3+312t^4+504t^5+680 t^6+504t^7+312t^8+83t^9+27t^{10}+3t^{11}+t^{12},\\
P_{21}&=1+6t+77t^2+484t^3+2877t^4+10374t^5+27027t^6
+46992t^7+57532t^8\\
&\qquad+ 46992t^9+27027t^{10}+10374t^{11}+2877t^{12}+484t^{13}
+77t^{14}+6t^{15}+t^{16};\\
Q_0&=1,\ Q_4=1+t+t^2,\ Q_8=1+12t+57t^2+92t^3+57t^4+12t^5+t^6,\\
Q_{12}&=1+43t+751t^2+5301t^3+16598 t^4+24137t^5+16598 t^6+5301t^7
+751t^8+43t^9+t^{10},\\
Q_{16}&=1+109t+5149t^2+103820t^3+
976685t^4+4609643t^5+11476395t^6+15506676t^7+\ldots\\
Q_{20}&=1+231t+25026t^2+1189860t^3+26750979t^4+308439936t^5+1946037411t^6\\
&\qquad+7038264246t^7+15046564512 t^8+19347012191t^9+\ldots\\
Q_{24}&=1+437t+97958t^2+9845240t^3+467190310t^4+
11368550417t^5+152640855877t^6\\
&\qquad+1196951395072t^7+5716465354180t^8
+17128652740280t^9+32841892687972t^{10}\\
&\qquad+40750517543272t^{11}+\ldots,
\end{align*}
\normalsize where} $\ldots$ stands for terms of degree larger than $\deg(Q_d)/2$.%
\begin{NB} I added the higher degree polynomials 10.11.LG\end{NB}
One checks that $P_d(1)=\Phi_0^{\P^2}(H^d)$, $Q_d(1)=2^d\Phi_H^{\P^2}(H^d)$, by comparing with \cite{EG2}, as required by the Hirzebruch-Riemann-Roch theorem.
In \cite{Da1},\cite{Da2} the $\chi(M_H^{\P^2}(0,d),\oo(\mu(H^{\otimes n}))$
were determined for $d\le 13$ and all $n$ and for $d=17, n=2,3$.

\subsection{Generalization to non-toric surfaces}\label{subsec:nontoric}
In this section we will generalize our results to arbitrary simply connected surfaces.  We extend \corref{cor:formula1} and 
\thmref{thm:q-expansion}
for the  wallcrossing terms
to any good wall $\xi$ on any simply connected projective surface $X$ with $p_g=0$. 
More generally let $X$ be a smooth  projective surface (not necessarily connected), and let $\xi\in\Pic(X)$ and $v\in K(X)$. We define the wallcrossing terms 
$\Delta_{\xi,T}^X(v,\Lambda)$, $\Delta_{\xi}^X(v,\Lambda)$ by the formulas
\eqref{wallct},\eqref{delT},\eqref{delT1}, where we replace 
in the summation index $d=4(n+m)-\xi^2-3$ by $d=4(n+m)-\xi^2-3\chi(\oo_X)$.
Then we show that these are computed by a suitable generalization of 
\corref{cor:formula1} and \thmref{thm:q-expansion}. 
This is done by adapting  the corresponding argument of \cite{GNY}
for the wallcrossing of the usual Donaldson  invariants, which is based on the  fact that  intersection numbers  on Hilbert schemes of points on $X$ are given by universal formulas in terms of intersection numbers on $X$.

If $X$ is a simply connected with $p_g=0$ and $\xi$ is a {\it good\/} class, then  \propref{wallcr}
shows that the wallcrossing of the $K$-theoretic Donaldson invariants for the wall defined by $\xi$ is given by the wallcrossing terms, thus we get a formula for the wallcrossing in terms of modular forms and elliptic functions. 
In the future we plan to adapt the arguments of \cite{Moc} to show that \propref{wallcr} and thus our wallcrossing formula also holds in case 
 $\xi$ is not good.

We start by sketching a proof of the  following result: 
\begin{Lemma}\label{lem:univ} Fix $r\in \Z$. There exist universal power series $A_i\in \Lambda\Q((T))[[\Lambda]]$, ($i=1,\ldots,7$), such that for all projective surfaces $X$, $\xi\in \Pic(X)$ and all 
$v\in K(X)$ of rank $r$ 
\begin{align*}
&\frac{(-T)^{\xi(\xi-K_X)/2+\chi(\oo_X)}\Lambda^{\xi^2+3\chi(\oo_X)}}{T^{{\xi}v^{(1)}/2}(1-T)^{\xi^2+2\chi(\oo_X)}}\Delta^X_{\xi,T}(v,\Lambda)=
\\
&\exp(\xi^2 A_1+\xi K_X A_2+ K_X^2 A_3+c_2(X) A_4+
\xi v^{(1)} A_5+K_X v^{(1)} A_6+ (v^{(1)})^2 A_7).
\end{align*}
Here, as before $v^{(1)}=c_1(v)+\frac{\rk(v)}{2}(c_1-K_X)$.
\end{Lemma}

A simple modification of the proof of \cite[Lemma 5.5]{GNY} shows
the following.
\begin{Lemma}\label{lem:unip}
Fix $n,m\ge 0$. Let $P$ be any polynomial in 
$\ch_{i_1}(\AA_+)$, $\ch_{i_2}(\AA_-)$, $\ch_{i_3}(\I_1)\xi^{i_4}/(v^{(1)})^{i_5}$
$\ch_{i_6}(\I_2)\xi^{i_7}/(v^{(1)})^{i_8},c_{i_9}(X^{[n]}\times X^{[n]})$
for $i_1,\ldots,i_9\in \Z_{\ge 0}$.
Then there exists a universal polynomial 
$Q$ (depending only on $P$) in $\xi^2$, $\xi K_X$, $K_X^2$, $c_2(X)$, 
$\xi v^{(1)}$, $K_X v^{(1)}$, $(v^{(1)})^2$, such that
$\int_{X^{[n]}\times X^{[m]}}P =Q$.
\end{Lemma}
The statement is very similar to \cite[Lemma 5.5]{GNY}. The only differences are
that  we replaced $X^{[l]}_2$ by  $X^{[n]}\times X^{[m]}$,
 and that we also allow the $c_{i}(X^{[n]}\times X^{[m]})$ in $P$. 
However looking at the proof of \cite[Lemma 5.5]{GNY} it obviously 
also works for $X^{[n]}\times X^{[m]}$, and in \cite{EGL} the argument is 
 also made for the  $c_{i}(X^{[n]})$. It readily generalizes to $X^{[n]}\times X^{[m]}$.
 
 Denote the left-hand-side of \lemref{lem:univ} by $\overline \Delta^X_{\xi,T}
 (v,\Lambda)$. By applying the Riemann-Roch theorem to   definition \eqref{delT}, we obtain that 
 $$\overline\Delta^X_{\xi,T}(v,\Lambda)=\sum_{n,m\ge 0}\sum_{i\in \Z}
 \Lambda^{4(n+m)} T^i \int_{X^{[n]}\times X^{[m]}}S_{n,m,i},$$ where $S_{n,m,i}$
 is a polynomial in the Chern characters of $\AA_+$, $\AA_-$, $\lambda_{\F_1}(v)$, $\lambda_{\F_2}(v)$ and the $c_j(X^{[n]}\times X^{[n]})$, which is zero for $i\ll 0$. By 
 \eqref{FC0} the Chern characters of the $\lambda_{\F_j}(v)$ are polynomials in the 
  $\ch_{i_1}(\I_j)\xi^{i_2}/(v^{(1)})^{i_3}$. Thus by
 \lemref{lem:unip} we see that $\overline \Delta^X_{\xi,T}
 (v,\Lambda)=\sum_{l\ge 0}
 \sum_{i\in \Z}\Lambda^{4l}P_{l,i} T^i$, where $P_{l,i}$
 is a universal polynomial in $\xi^2$, $\xi K_X$, $K_X^2$, $c_2(X)$, $\xi v^{(1)}$,  $K_X v^{(1)}$, $(v^{(1)})^2$, which is zero for $i\ll0$.
 From the definition \eqref{delT}, one readily computes that the coefficient of 
 $\Lambda^{0}$ of $ \overline \Delta^X_{\xi,T}
 (v,\Lambda)$ as a power series in $\Lambda$ is $1$.%
 \begin{NB} More detailed computation:
 When $n,m=0$, then 
 $\rk(\AA_-)=-\xi(\xi-K_X)/2-\chi(\oo_X)$, $\rk(\AA_+)=-\xi(\xi+K_X)/2-\chi(\oo_X)$ and
 $$\chi\big(X^{[0]}\times X^{[0]},\lambda_{\F_1}(v)\otimes \lambda_{\F_2}(v)
 \otimes \det(\AA_+)\otimes S_T(\AA_+^\vee)\otimes S_T(\AA_-)\big)=
 \frac{1}{(1-T)^{\rk(\AA_+)+\rk(\AA_-)}}.$$
 \end{NB}
 Now the  proof of \lemref{lem:univ} is finished by the same arguments as that of  \cite[Theorem 5.1]{GNY}.
  
  \begin{NB} There is a small mistake in the statement of \cite[Theorem 5.1]{GNY}. To prove \cite[Corollary 5.7]{GNY}, we need that the $A_i$ are
  in $\Lambda\Q((t^{-1}))[[\Lambda]]$, and this is indeed what the proof of 
  \cite[Theorem 5.1]{GNY} shows. The statement will have to be corrected when me make a revised version.
  \end{NB}
  
\begin{Corollary}\label{cor:unim}
\begin{enumerate}
\item
\corref{cor:formula1} and \thmref{thm:q-expansion}
hold for any simply connected smooth projective surface with $p_g=0$ and any 
$\xi\in \Pic(X)$. 
\item More generally for any smooth projective surface $X$ and any 
$\xi\in \Pic(X)$ we have
\begin{align*}
\Delta^X_{\xi,e^{-2\bbeta a}}(v,\bbeta\Lambda)&= 
\sqrt{-1}^{\<\xi,K_X\>}\frac{q^{-\frac{1}{2}(\frac{\xi}{2})^2}}{(\bbeta\Lambda)^{\chi(\oo_X)}}
\exp\Bigg(
\frac{\bbeta}{8}\frac{\partial^2\F_0}{\partial a\partial \log\Lambda}\<\xi,v^{(1)}+K_X\>\\&\qquad
+\frac{\bbeta^2}{32}\frac{\partial^2\F_0}{(\partial \log\Lambda)^2}
\<(v^{(1)}+K_X)^2\>\Bigg)\exp(A)^{4\chi(\oo_X)}\exp(B-A)^{\sigma}
\end{align*}
\end{enumerate}
\end{Corollary}
\begin{proof}
It is enough to show part (2). 
We put $T^{\frac{1}{2}}:=e^{-\bbeta a}$, and as above write $z=\frac{-\sqrt{-1}\bbeta\Lambda}{e^{-\bbeta a}-e^{\bbeta a}}=\frac{\sqrt{-1}\bbeta\Lambda T^{\frac{1}{2}}}{1-T}$.
Then by \eqref{eq:q_expansion} we have
$q^{\frac{1}{8}}=z \exp(l_1)$, with $l_1\in z^2\C[\bbeta,T^{\pm\rk(v)/2},\Lambda][[z^2]]\subset
\Lambda^2 \C[\bbeta,\Lambda]((T^{\frac{1}{2}}))$.
Similarly \eqref{eq:zeta''} implies $\frac{\partial^2 \F_0}{(\partial \log\Lambda)^2}\in \Lambda^2 \C[\bbeta,\Lambda]((T^{\frac{1}{2}}))$,  $\frac{\partial^2 \F_0^{inst}}{\partial a\partial \log\Lambda}\in\Lambda^2 \C[\bbeta,\Lambda]((T^{\frac{1}{2}}))$, and from the definition we see that
$\frac{\partial^2 \F_0^{pert}}{\partial a\partial \log\Lambda}=-8a.$
Finally by \eqref{eq:AB} we have
$\exp(A)=q^{\frac{1}{16}}\exp(l_2)$, $\exp(B-A)=\exp(l_3)$, with $l_2,l_3\in 
\Lambda\C[\bbeta,\Lambda]((T))$.
Thus we see that the left hand side of Corollary \ref{cor:unim} can be rewritten 
as 
$M\exp(\xi^2 B_1+\xi (v^{(1)}+K_X) B_2 +(v^{(1)}+K_X)^2 B_3
+c_2(X) B_4 + K_X^2 B_5)$, with  $B_i\in \Lambda\C((T^{\frac{1}{2}}))[[\Lambda]]$ and
\begin{align*}
M&=\sqrt{-1}^{\<\xi,K_X\>}
\Big(\frac{\sqrt{-1}\bbeta\Lambda T^{\frac{1}{2}}}{1-T}\Big)^{-\xi^2-2\chi(\oo_X)}
\frac{T^{\<\xi(v^{(1)}+K_X)\>/2}}{\Lambda^{\chi(\oo_X)}}=\frac{T^{\xi v^{(1)}/2}(1-T)^{\xi^2+2\chi(\oo_X)}}{(-T)^{\xi(\xi-K_X)/2+\chi(\oo_X)}\Lambda^{3\chi(\oo_X)}}.
\end{align*}
As the $A_i$, $i=1,\ldots,7$ of \lemref{lem:univ} are determined by 
the $\Delta_{\xi,T}^X(v,\Lambda)$ for toric surfaces, \corref{cor:formula1} implies the result.
\end{proof}

\begin{NB} Note that in the proof of \corref{cor:unim} the structure of 
the functions and of the formulas in \corref{cor:formula1} played an important 
role. This can be viewed as an additional check for the correctness of our
formulas. Basically we see that much of  general structure of  \corref{cor:formula1}
is determined by the perturbation part of $\F_0$.
\end{NB}

When $v = v(2L)$, we have $v^{(1)} + K_X = K_X - 2L$, which is equal
to the negative of the characteristic line bundle $\det W^\pm$ of the
$Spin^c$ structure $W^\pm$ induced from the complex structure of $X$
and the line bundle $L$ (see \subsecref{subsec:digress}).
Then
\(
  \sqrt{-1}^{-\<\xi,K_X\>}
  \left[\Delta^X_{\xi}(v,\bbeta\Lambda)\right]_{\Lambda^d}
\)
is a polynomial in $\<\xi,c_1(\det W^\pm)\>$ and $\< c_1(\det
W^\pm)^2\>$ whose coefficients depend only on $\<\xi^2\>$, $d$ and the
homotopy type of $X$. This statement is a natural analogue of the
Kotschick-Morgan conjecture \cite{KM} in the context of the $K$-theoretic
Donaldson invariants. Thus our formula above supports our belief that
the $K$-theoretic Donaldson invariants have a gauge theoretic
definition.
\begin{NB}
A paragraph added. Nov. 6, HN. 
\end{NB}

\appendix
\section{Seiberg-Witten curves for $K$-theoretic version}
\label{sec:SWcurve}

The purpose of this appendix is to prove some results on
Seiberg-Witten curves for the $K$-theoretic version with Chern-Simons
terms. In particular, we show 
\begin{aenume}
\item the perturbation part of the
Seiberg-Witten prepotential coincides with the genus $0$ part of the
perturbation part introduced in \S\ref{subsubsec:pert},
\item the Seiberg-Witten prepotential satisfies the contact term
  equation in \propref{prop:contact}.
\end{aenume}
The corresponding results of the Seiberg-Witten curves for the
homological version have been known \cite{HKP,Matone,STY,GM3}, and
were reproduced in \cite[\S2]{NY3}. Our proofs go along the same line,
while we need to consider the cases $r+m$ even and odd
separately. The adaptation might be standard to experts, but we cannot
find the statements or proofs in the literature.

\subsection{Seiberg-Witten curves}

We consider a family of curves parametrized by
$\vec{U} = (U_1,\dots,\linebreak[0]U_{r-1})$:
\begin{equation*}
\begin{split}
  & C_{\vec{U},m} : (-\sqrt{-1}\bbeta\Lambda)^rX^{(r+m)/2}
  \left(w + \frac1w\right) = P(X),
\\
  & \qquad P(X) = 
     X^r + U_1 X^{r-1} + U_2 X^{r-2} + \cdots + U_{r-1}X + (-1)^r
\end{split}
\end{equation*}
for $|m|\le r$, $m\in\Z$.%
\begin{NB}
$\bbeta$ is multiplied by $-\sqrt{-1}$. May 18, 2006.
\end{NB}
We call them {\it Seiberg-Witten curves}. When $r+m$ is odd, we should
understand this expression formally, and the rigorous definition will
be given soon below.
The projection $C_{\vec{U},m}\ni (w,X)\mapsto X\in \proj^1$ gives a
structure of hyperelliptic curves. The hyperelliptic involution
$\iota$ is given by $\iota(w) = 1/w$.

We introduce a new variable 
\(  
   Y = (-\sqrt{-1}\bbeta\Lambda)^r X^{(r+m)/2} (w - \frac1w).
\)
Thus we have 
\[
   Y^2 = P(X)^2 - 4 (-X)^{r+m} (\bbeta\Lambda)^{2r} 
   .
\]
This {\it does\/} make sense for $r+m$ odd also.

Note that $|m|\le r$ guarantees that the curve has genus $r-1$. Later
we further assume $|m|\neq r$.
\begin{NB}
  Added on Sep.\  7, 2006, H.N.
\end{NB}

If we replace the coordinate $X$ (near $0$) by $1/X$ (near $\infty$),
then the equation of the curve becomes
\begin{equation*}
  Y^2 = X^{2r}\big(P(1/X)^2 - 4(-1/X)^{r+m}(\bbeta\Lambda)^{2r}\big)
  = \widetilde P(X)^2-4(-X)^{r-m}(\bbeta\Lambda)^{2r},
\end{equation*}
where $\widetilde P(X) = X^r + (-1)^r U_{r-1} X^{r-1} + \cdots +
(-1)^r$.  Thus the curves for $m$ and $-m$ are essentially the same
(exactly the same when $r=2$), once written in a coordinate near $0$
and once in a coordinate near infinity.
\begin{NB}
I took remark from Lothar's note. Sep.\  7, 2006, H.N.  
\end{NB}

Let us define the {\it Seiberg-Witten differential\/} by 
\begin{equation*}
\begin{split}
dS &= \frac1{2\pi\sqrt{-1}
    \bbeta}\log X\frac{dw}w
  = \frac1{2\pi \sqrt{-1}
    \bbeta} {\log X}\, \frac{X^{(r+m)/2} (X^{-(r+m)/2}P(X))' dX}Y
\\
  &= \frac1{2\pi \sqrt{-1}
    \bbeta} {\log X}\,  \frac{2XP'(X) - (r+m) P(X)}{2XY} dX,
\end{split}
\end{equation*}
where we have used
\begin{equation*}
    X^{-(r+m)/2} Y \frac{dw}w
   = (-\sqrt{-1}\bbeta\Lambda)^r\left(w - \frac1w\right) \frac{dw}w
   = \left( X^{-(r+m)/2} P(X) \right)' dX.
\end{equation*}
This is a multi-valued meromorphic differential on $C_{\vec{U},m}$. The
last expression makes sense even in the case $r+m$ odd.

Let $X_1$,\dots, $X_r$ be the zeroes of $P(X) = 0$. We have $\prod X_i
= 1$.

\subsection{Homological limit $\bbeta\to 0$}\label{subsec:limit}

We move $\bbeta$ in a small disk around the origin. We see that the
Seiberg-Witten curve becomes the Seiberg-Witten curve for the
homological version (i.e.\ the $4$-dimensional gauge theory in the
physics terminology) at $\bbeta = 0$.

We choose $z_i$ with $X_i = e^{-\sqrt{-1}\bbeta z_i}$. We consider
$\vec{z} = (z_i)$ is a parameter for the curve. Let $X = \frac{2
-\sqrt{-1}\bbeta z}{2 +\sqrt{-1}\bbeta z}$. Then
\begin{equation*}
\begin{split}
   & (-\sqrt{-1}\bbeta)^{-r}  X^{-(r+m)/2} P(X) 
 \\
  & = 
  \left(1 + \frac{\bbeta^2 z^2}4 \right)^{-\frac{r+m}2}
  \left(1+\frac{\sqrt{-1}\bbeta}2\right)^{m}
  \prod_{i=1}^r \left[
     \frac{e^{-\sqrt{-1}\bbeta z_i} + 1}2 z
       - \frac{e^{-\sqrt{-1}\bbeta z_i} - 1}{-\sqrt{-1}\bbeta}
    \right].
\end{split}
\end{equation*}
If we introduce a new variable $y = (-\sqrt{-1}\bbeta)^{-r} Y (1 +
\frac{\sqrt{-1}}2\bbeta z)^{r}$, we have
\begin{equation*}
    y^2 = \prod_{i=1}^r \left[
     \frac{e^{-\sqrt{-1}\bbeta z_i} + 1}2 z
       - \frac{e^{-\sqrt{-1}\bbeta z_i} - 1}{-\sqrt{-1}\bbeta}
    \right]^2 - 4\Lambda^{2r}\left(1 + \frac{\bbeta^2
        z^2}4\right)^{r-m}
    \left(1-\frac{\sqrt{-1}\bbeta}2 z\right)^{2m}.
\end{equation*}
Therefore in the limit $\bbeta\to 0$, the Seiberg-Witten curve
converges to
\[
   \Lambda^r(w + \frac1w) = \prod_{i=1}^r (z - z_i)
   \text{ or }
   y^2 = \prod_{i=1}^r (z - z_i)^2 - 4\Lambda^{2r}.
\]
This is the Seiberg-Witten curve for the homological version. (The
variable $w$ is the same.) The Seiberg-Witten differential converges
to that of the homological version, i.e.\ $-\frac1{2\pi}z \frac{dw}w$.
\begin{NB}
$\frac1\bbeta\log\frac{2-\sqrt{-1}\bbeta z}{2+\sqrt{-1}\bbeta z}
\approx -\sqrt{-1}\bbeta z$.
\end{NB}

The points $X=0$, $\infty$ corresponds to $z =
\frac{2\sqrt{-1}}\bbeta$, $-\frac{2\sqrt{-1}}\bbeta$. Therefore in the
limit $\bbeta\to 0$, both points go to a common point $z=\infty$.

We find $X_i^\pm$ near $X_i$ such that
\begin{equation*}
   P(X_i^\pm) = \pm 2 (-\sqrt{-1}\bbeta\Lambda)^r (X_i^\pm)^{(r+m)/2}.
\end{equation*}
When $r+m$ is odd, we take the branch of $(X_i^\pm)^{1/2}$ so that
it is the same branch as $X_i^{1/2} = e^{-\sqrt{-1}\bbeta z_i/2}$.%
\begin{NB}
I changed the explanation of the choice of the branch. This
guarantees to fix our choice so that $X_i^\pm$ converges to $z_i^\pm$.
\end{NB}
Let us choose $z_i^\pm$ so that $\prod (z_i^\pm - z_i) = \pm
2\Lambda^r$. Then $X_i^\pm\to z_i^\pm$ (more precisely after moving to
the $z$-coordinates).
\begin{NB}
$(-\sqrt{-1}\bbeta)^{-r}  X^{-(r+m)/2} P(X)$ converges to
$\prod_i (z - z_i)$. The left hand side is equal to $\pm 2\Lambda^r$
at $X = X_i^\pm$. If we take another branch (for all $i$), then
it converges to $-\prod_i (z - z_i)$. So $X_i^\pm\to z_i^\mp$.
\end{NB}

The correspondence between the coefficients is%
\begin{NB} changed 30.10 LG\end{NB}
more tricky, as $U_i$ is
the $i^{\mathrm{th}}$ elementary symmetric function in
$e^{-\sqrt{-1}\bbeta z_i}$ while $u_i$ is the $i^{\mathrm{th}}$
elementary symmetric function in $z_i$, up to sign. For example, $r=2$
\begin{equation*}
  U_1 = - (e^{-\sqrt{-1}\bbeta z_1} + e^{-\sqrt{-1}\bbeta z_2})
  \approx - 2 + \frac{\bbeta^2}2 (z_1^2 + z_2^2)
  = - 2 - \bbeta^2 u_2.
\end{equation*}

\begin{NB}
The sign was corrected. Sep.\  7, 2006, H.N.
\end{NB}

\subsection{$a_i$, $a_i^D$ and the prepotential $\mathcal F_0$}

We first work in the region containing $z_1,\dots,z_r\in\R$ and $z_1 >
z_2 > \dots > z_r$. Then we will analytically continue to the whole
region.%
\begin{NB} which region? 30.10 LG\end{NB}
The curve itself is parametrized by $\vec{U}$, but its
homology basis introduced below depends on $\vec{z} = (z_i)$.
We also first suppose that $\Lambda$ is a sufficiently small positive
real number and then will analytically continue to a small punctured
disk.

We take cycles $A_i$, $B_j$ ($i=1,\dots,r$, $j=2,\dots,r$) so that it
gives the cycles for the Seiberg-Witten curves for the homological
version given in \cite[\S2]{NY2} at $\bbeta = 0$. Let us explain a
little bit more precisely: Our curve $C_{\vec{U}}$ is hyperelliptic
and is made up of two copies of the Riemann sphere, glued along the
$r$-cuts between $X_i^-$ and $X_i^+$. We then define $A_i$ as the
cycle encircling the cut between $X_i^-$ and $X_i^+$. Note that we
have $\sum_i A_i = 0$.
We choose cycles $B_j$ ($j=2,\dots,r$) as in \cite[Figure~1]{NY2},
i.e.\ $B_j$ is the sum $\sum_{k=2}^j C_k$ where $C_k$ is a cycle
starting from $X_{k-1}^\pm$, passing through $X_{k}^\pm$, and then
returning back to $X_{k-1}^\pm$ in the another sheet. Here the sign is
$+$ for $i$ odd, $-$ for $i$ even. Then $A_i$, $B_i$ ($i=2,\dots,r$)
form a symplectic basis of $H_1(C_{\vec{U}},\Z)$.

We define $a_{i}$, $a_j^D$ by
\begin{equation*}
   a_i = \int_{A_i} dS, \qquad
   a^D_j = \int_{B_j} dS,
   \qquad i=1,\dots,r,\ j=2,\dots,r.
\end{equation*}%
\begin{NB}
We drop the factor $2\pi \sqrt{-1}$ from $a^D_j$.
\end{NB}
We consider a region disjoint from a segment from $\infty$ to $0$
which does not pass $e^{-\sqrt{-1}\bbeta z_i}$. Therefore $\log X$ is
single-valued in the region.
We take the branch of $\log X$ so that it is given by $-\sqrt{-1} \bbeta z_i$
at $X_i = e^{-\sqrt{-1}\bbeta z_i}$.
The $A_i$, $B_i$ cycles are taken from the region.

\begin{NB}
Kota's more precise explanation:

We shall consider points $(z_1,z_2,...,z_r) \in {\Bbb C}^r$
which are sufficiently close to
$(z_1,z_2,...,z_r) \in {\Bbb R}^r$ with 
$z_r<z_{r-1}<\cdots <z_1$.
For simplicity, we assume that
$(z_1,z_2,...,z_r) \in {\Bbb R}^r$ with 
$z_r<z_{r-1}<\cdots <z_1$.
Then $0<e^{\bbeta z_r}<e^{\bbeta z_{r-1}}<\cdots<e^{\bbeta z_1}<\infty$.
Let ${\Bbb P}_{\Bbb R}^1 \subset {\Bbb P}_{\Bbb C}^1$ 
be the real projective line
which is a circle in ${\Bbb P}_{\Bbb C}^1$.
Let $l \subset {\Bbb P}_{\Bbb R}^1$ 
be the segment from $\infty$ to $0$ which does not pass
$e^{\bbeta z_i}$.
We take a compact neighborhood $V$ of $l$ which does not
meet the segment $[e^{\bbeta z_r},e^{\bbeta z_1}]$. 
Then ${\Bbb P}_{\Bbb C}^1 \setminus l$ is simply connected and 
we have a single valued holomorphic function $\log X$
on ${\Bbb P}_{\Bbb C}^1 \setminus V$.
Thus $\log X$ is a single valued holomorphic function on
$C_{\vec{U}} \setminus \pi^{-1}(V)$.
We take $A,B$ cycles in $C_{\vec{U}} \setminus \pi^{-1}(V)$.
\end{NB}

We have the following expansion:
\begin{equation}\label{eq:za}
\begin{split}
a_i&=
\frac{1}{2\pi \sqrt{-1}\bbeta} 
\int_{A_i} \log X
d \left[\log \left(\prod_{j=1}^r 
(X^{\frac{1}{2}}-e^{-\sqrt{-1}\bbeta z_j} X^{-\frac{1}{2}})\right)
   - \frac{m}2\log X\right]
 + O(\Lambda)\\
&= -\sqrt{-1} z_i+ O(\Lambda).
\end{split}
\end{equation}
\begin{NB}
Please check that the 5D Chern-Simons term ($m$) does not contribute
here.
Sep.\  7, 2006, H.N. 
\end{NB}

\begin{NB}
As explained in \cite{NY2}, we first move $A_i$ to a sheet where
$Y = \sqrt{P(X)^2 - 4(-X)^{r+m}(\bbeta\Lambda)^{2r}}$ is single-valued and
is approximated by $P(X)$ for small $\Lambda$. Then
$\frac{X^{(r+m)/2} (X^{-(r+m)/2}P(X))' dX}Y \approx 
\frac{(X^{-(r+m)/2}P(X))' dX}{X^{-(r+m)/2} P(X)}$.
\end{NB}

We invert the roles of $a_i$ and $U_p$, so we consider
$a_i$ as variables and $U_p$ are functions in $a_i$.

Let us differentiate the defining equation of $C_{\vec{U}}$ with
respect to $U_p$ by setting $w$ to be constant:
\begin{equation}\label{eq:dXdU}
   0 = \left(X^{-(r+m)/2}P(X)\right)' \frac{\partial X}{\partial U_p}
   + X^{(r-m)/2-p}.
\end{equation}
Therefore the differential of the Seiberg-Witten differential $dS$ is
\begin{equation}\label{eq:dSdU}
  \left.\frac{\partial}{\partial U_p} dS\right|_{w=\mathrm{const}}
  = - \frac1{2\pi \sqrt{-1}
    \bbeta} \frac{X^{(r-m)/2-p-1}}{(X^{-(r+m)/2}P(X))'}\frac{dw}{w}
  = - \frac1{2\pi \sqrt{-1}
    \bbeta} \frac{X^{r-p-1} dX}{Y}.
\end{equation}
It is well-known that these form a basis of holomorphic differentials
on $C_{\vec{U}}$ for $p=1,\dots,r-1$ (see e.g., \cite[\S2.3]{GH}). In
other words, the Seiberg-Witten differential is a `{\it potential\/}'
for holomorphic differentials.

Let $(\sigma_{i p})$ be the matrix given by
\begin{equation*}
   \sigma_{i p}
   = \pd{a_i}{U_p} 
   = - \frac1{2\pi \sqrt{-1}
     \bbeta}
    \int_{A_i} \frac{X^{r-p-1} dX}{Y}
   \qquad i =2,\dots,r, \ p=1,\dots,r-1.
\end{equation*}
If $(\sigma^{p i})$ is the inverse matrix, the normalized
holomorphic $1$-forms 
\[
  \omega_j
  = - \frac1{2\pi \sqrt{-1}
        \bbeta} \sum_p \sigma^{pj}
    \frac{X^{r-p-1}dX}{Y}
  = \left.\pd{}{a_j} dS\right|_{w=\mathrm{const}}
\]
satisfies $\int_{A_i} \omega_j = \delta_{ij}$. Therefore the {\it period matrix\/} $\tau =
(\tau_{ij})$ of the curve $C_{\vec{U}}$ is given by
\begin{equation}\label{eq:period}
  \tau_{ij}
  = \int_{B_i} \omega_j
  = 
  \pd{a^D_i}{a_j}.
\end{equation}
Since $(\tau_{ij})$ is symmetric (see e.g., \cite[\S2.2]{GH}), there
exists a locally defined function $\mathcal F_0$ such that
\begin{equation}\label{eq:SWprep}
   a_i^D = - \frac1{2\pi\sqrt{-1}}
   \frac{\partial\mathcal F_0}{\partial a_j}.
\end{equation}
It is unique up to a function independent of $a_i$.  The ambiguity
will be fixed later.%
\begin{NB}
In the homological version:

We fix the constant so that $\mathcal F_0$ is homogeneous of degree
$2$:
\begin{equation*}\label{eq:Euler}
   \left(\sum a_\alpha \frac{\partial}{\partial a_\alpha}
   + \Lambda \pd{}{\Lambda} \right) \mathcal F_0
 = 2\mathcal F_0.
\end{equation*}%
\end{NB}
This function $\mathcal F_0$ is called the {\it Seiberg-Witten
prepotential}. We may also write $\mathcal F_0(\vec{a})$ or $\mathcal
F_0(\vec{a};\Lambda)$.

\subsection{Perturbative part}

We determine the perturbative part of the prepotential $\mathcal F_0$
in this subsection.

Let
\begin{equation*}
\begin{split}
\overline\gamma_0(x|\bbeta;\Lambda)&=
2\left(
\frac{1}{\bbeta^2}(\Li_3(e^{-\bbeta x})-\zeta(3))+
\frac{x^2}{2}\log(\bbeta \Lambda)+
\frac{\pi^2}{6\bbeta}x \right)-
\frac{x^2\pi \sqrt{-1}}{2}-\frac{\bbeta x^3}{6},
\end{split}
\end{equation*}
where $\Li_3$ is the trilogarithm. See \cite[App.~B]{NY2} for the
definition and properties of polylogarithms. The relation to the
perturbative part $\widetilde\gamma_{\ve_1,\ve_2}(x|\bbeta;\Lambda)$
in \secref{subsec:Nek} is the following: We have defined
$\widetilde\gamma_{\ve_1,\ve_2}(x|\bbeta;\Lambda)$ first when $\bbeta
x >0$ and then analytically continued it to the whole plane. Then we
considered
\(
   \widetilde\gamma_{\ve_1,\ve_2}(x|\bbeta;\Lambda) + 
   \widetilde\gamma_{\ve_1,\ve_2}(-x|\bbeta;\Lambda).
\)
The coefficient of $1/\ve_1\ve_2$ is equal to
$\overline\gamma_0(x|\bbeta;\Lambda)$. See \cite[p.~510, the second
displayed formula from the bottom]{NY3}.
This becomes regular and its value is
\(
   -x^2 \left(\log\frac{\sqrt{-1}x}\Lambda\right) + \frac32 x^2
\)
at $\bbeta = 0$ ([loc.cit., p.510, the last displayed formula]).

\begin{Proposition}\label{prop:perturb}
\begin{equation*}
   \mathcal F_0(\vec{a};\Lambda)
   = -\sum_{i<j} \overline\gamma_0(a_i - a_j|\bbeta;\Lambda)
    - \frac{m\bbeta}6\sum_{i=1}^r a_i^3
   + O(\Lambda).
\end{equation*}
\end{Proposition}

The term $-\sum_{i<j} \overline\gamma_0(a_i - a_j|\bbeta;\Lambda)- \frac{m\bbeta}6\sum_{i=1}^r a_i^3$ is
called the {\it perturbative part\/} of $\mathcal F_0$. Recall that
$\mathcal F_0$ was defined up to a function (in $\Lambda$) independent
of $a_i$. We, in fact, prove
\begin{equation}\label{eq:to_be_checked}
\begin{split}
  & - \partial \mathcal F_0/\partial a_i = 2\pi\sqrt{-1} a^D_i
\\
  =\; & 
   - \sum_{j > 1} \overline\gamma_0'(a_1 - a_j|\bbeta;\Lambda)
   + \sum_{j : i < j} \overline\gamma_0'(a_i - a_j|\bbeta;\Lambda)
   - \sum_{j : j < i} \overline\gamma_0'(a_j - a_i|\bbeta;\Lambda)
\\
   & \qquad\qquad
   + \frac{m\bbeta}2\left(a_i^2 - a_1^2\right)
   + O(\Lambda).
\end{split}
\end{equation}
Then we take a function so that the above formula holds. The remained
ambiguity in $O(\Lambda)$ will be fixed later.

Note that $r=2$ case the term $\frac{m\bbeta}6\sum_{i=1}^r a_i^3$
vanishes as $a_1 + a_2 = 0$. Therefore this does not show up in
\subsecref{subsec:Nek}.

Let us describe the branch of $\overline\gamma_0$. As our $\bbeta$ is in a
small disk around the origin, it is enough for us to fix the branch
at $\bbeta = 0$. Then the ambiguity occurs only at
\(
   \log\left(\frac{\sqrt{-1}x}\Lambda\right).
\)
When $z_1 > \cdots > z_r$ and $\Lambda\in \mathbb R_{>0}$, $a_i$ is
pure imaginary and $\sqrt{-1}(a_i - a_j)\in \mathbb R_{>0}$ for $i <
j$. We then choose $\log\left(\frac{\sqrt{-1}x}\Lambda\right)\in
\mathbb R$. Therefore we have $a^D_i\in \R$.
\begin{NB}
This choice of the branch is the same as our original choice of
the branch of polylogarithms in $\overline\gamma_0$ explained above.
In the original choice we first choose the branch of the logarithm so
that it is real-valued on positive real numbers, and then we integrate
it to fix the branch of the polylogarithms. In the limit $\bbeta\to
0$, this choice goes to the above choice.
\end{NB}

Note that
\begin{equation}\label{eq:gamma'}
\overline\gamma_0'(x|\bbeta;\Lambda)=-2
\left(
\frac{1}{\bbeta}(\Li_2(e^{-\bbeta x})-\frac{\pi^2}{6})-
x\log(\bbeta \Lambda)
 \right)-
x\pi \sqrt{-1}-\frac{\bbeta x^2}{2}.
\end{equation}
We denote $\Li_2(e^{-\bbeta x}) - \frac{\pi^2}6 - \bbeta x
\log(\bbeta\Lambda)$ by $\widehat{\Li}_2(e^{-\bbeta x})$ for brevity. 
\begin{NB}
$\lim_{\bbeta\to 0} \frac1{\bbeta}\widehat{\Li}_2(e^{-\bbeta x})
= x \log\left(\frac{x}\Lambda\right) - x$ (see \cite[App.~B.3]{NY3}).
\end{NB}

Our proof is given so that it reduces to the proof of
\cite[Prop.~2.2]{NY2} when $\bbeta\to 0$. 
(The proof of \cite[Prop.~2.2]{NY2} was based on \cite{HKP} in turn.)

\begin{NB}
We use the following elementary formulas:
\begin{equation*}
\begin{split}
&-\sum_{1<j}\frac{\bbeta (a_1-a_j)^2}{2}+
\sum_{j : i<j}\frac{\bbeta (a_i-a_j)^2}{2}-
\sum_{j: j<i}\frac{\bbeta (a_j-a_i)^2}{2}\\
=& -\sum_{j : \textrm{any}}\frac{\bbeta (a_1-a_j)^2}{2}+
\sum_{j : \textrm{any}}\frac{\bbeta (a_i-a_j)^2}{2}-
\sum_{j : j < i}\bbeta (a_j-a_i)^2\\
=& \frac{r\bbeta}{2}(a_i^2-a_1^2)
-\sum_{j:j < i}\bbeta (a_j-a_i)^2.
\end{split}
\end{equation*}
And
\begin{equation*}
\sum_{1<j} (a_1-a_j)
- \sum_{j : i < j} (a_i-a_j)+
\sum_{j : j<i} (a_j-a_i)= r(a_1-a_i).
\end{equation*}
Hence
\begin{equation*}\label{eq:pert}
\begin{split}
-\pd{\mathcal F_0^{\mathrm{pert}}}{a_i} &=
\frac{2}{\bbeta}\sum_{1<j}\widehat{\Li}_2(e^{-\bbeta(a_1-a_j)})-
\frac{2}{\bbeta}\sum_{j : i<j}\widehat{\Li}_2(e^{-\bbeta(a_i-a_j)})+
\frac{2}{\bbeta}\sum_{j : i>j}\widehat{\Li}_2(e^{-\bbeta(a_j-a_i)})\\
&\quad -\frac{(r-m)\bbeta}{2}(a_i^2-a_1^2)+\sum_{i>j}\bbeta(a_j-a_i)^2
+ r\pi \sqrt{-1}(a_1-a_i) 
.
\end{split}
\end{equation*}%
\end{NB}

\begin{proof}[Proof of \propref{prop:perturb}]
In the proof we move $\Lambda$ in a punctured disk by the analytic
continuation, starting from positive real numbers. Then $a^D_i$ is a
multi-valued holomorphic function in $\Lambda$.

Let $C_i$ be a cycle starting from 
$e^{-\sqrt{-1}\bbeta z_{i-1}^{\pm}}$, passing through 
$e^{-\sqrt{-1}\bbeta z_i^{\pm}}$, and then returning back to
$e^{-\sqrt{-1}\bbeta z_{i-1}^{\pm}}$ in the another sheet.
Here the sign is $+$ for $i$ odd, $-$ for $i$ even.
Then $B_i=\sum_{k=2}^i C_k$.

We note that $\int_{C_i} dS$ is a local function of 
$\Lambda^{2r}$.
Since $C_i$ change to $C_i+A_i-A_{i-1}$ under 
the analytic continuation along $\Lambda^{2r} \to 
e^{2\pi \sqrt{-1}} \Lambda^{2r}$,
 $\int_{C_i} dS -(a_i-a_{i-1}) \log \Lambda^{2r}$
is a single valued function on the punctured disk
$0<|\Lambda^{2r}|\ll 1$.
\begin{NB}
In our lecture note, we write that 
$\int_{C_i} dS + \sqrt{-1}(z_i-z_{i-1}) \log \Lambda^{2r}$
is a single valued function on the punctured disk
$0<|\Lambda^{2r}|\ll 1$. This is not necessarily correct.
\end{NB}

We take a small positive real number $\delta$ with $|\Lambda| \ll
\delta$ and rewrite the integral as
\begin{equation*}
\begin{split}
\int_{C_i} dS &=2\int_{e^{-\sqrt{-1}\bbeta z_{i-1}^{\pm}}}^{e^{-\sqrt{-1}\bbeta z_i^{\pm}}}dS\\
&=2 \int_{e^{-\sqrt{-1}\bbeta z_{i-1}^{\pm}}}^{e^{-\sqrt{-1}\bbeta (z_{i-1}-\delta)}} dS
+2\int_{e^{-\sqrt{-1}\bbeta (z_{i-1}-\delta)}}^{e^{-\sqrt{-1}\bbeta (z_i +\delta)}} dS
+2\int_{e^{-\sqrt{-1}\bbeta (z_i+\delta)}}^{e^{-\sqrt{-1}\bbeta z_i^{\pm}}} dS.
\end{split}
\end{equation*}

\begin{NB}
  In \cite{NY2} we take $z_i-\delta$ and $z_{i-1}+\delta$. But
  this was wrong.
\end{NB}

We first compute the second term.%
\begin{NB}
This is regular at $\Lambda=0$ and we can approximate $Y$ by $P(X)$
on the interval. See \cite[p.42]{NY2}. Sep.\ 7, 2006, H.N.
\end{NB}
Let us write $\bbeta' =
-\sqrt{-1}\bbeta$ for brevity. Then
{\allowdisplaybreaks
\begin{equation*}
\begin{split}
&
  - 2\pi 
\int_{e^{\bbeta'(z_{i-1}-\delta)}}^{e^{\bbeta'(z_i +\delta)}} dS\\
=\; &
\begin{aligned}[t]
  & 
  - \int_{z_{i-1}-\delta}^{z_i+\delta} \frac{m\bbeta' t}2 dt 
  +\left[\sum_j \frac{\log X}{\bbeta'} \log \left( \frac{X^{\frac{1}{2}}-
  e^{\bbeta' z_j}X^{-\frac{1}{2}}}{\bbeta'\Lambda}\right)
   \right]_{e^{\bbeta'(z_{i-1}-\delta)}}^{e^{\bbeta'(z_i +\delta)}}
\\
  & \qquad\qquad 
-\int_{z_{i-1}-\delta}^{z_i +\delta}
\sum_j \left(\frac{\bbeta^{\prime} t}{2}+
\log\left(\frac{1-e^{-\bbeta'(t-z_j)}}{\bbeta'\Lambda}\right) 
\right)dt + O(\delta)
\end{aligned}
\\
=\; & \Biggl[
  - \frac{m\bbeta'}4 z_i^2
  + \sum_{j \ne i} z_i 
  \log \left( \frac{1-e^{-\bbeta'(z_i-z_j)}}{\bbeta'\Lambda} \right)
  +\frac{r}{2}\bbeta' z_i^2
  +z_i \log\left(\frac{1-e^{-\bbeta'\delta}}{\bbeta'\Lambda}\right) 
\\
&\quad \quad - \frac{r}{4}\bbeta' z_i^2
-\frac1{\bbeta'}\sum_{j > i}\widehat{\Li}_2(e^{-\bbeta'(z_i-z_j)})
+\frac1{\bbeta'}\sum_{j < i}\widehat{\Li}_2(e^{-\bbeta'(z_j-z_i)})+
\sum_{j<i} \frac{\bbeta'(z_j-z_i)^2}{2}\Biggr]
\\
& \quad
  - \Biggl[ \text{the same term with $z_i \to z_{i-1}$} \Biggr]
+ O(\delta).
\end{split}
\end{equation*}
Here} we have determined the branch of $\log$ so that this is
real-valued when $\bbeta = 0$ and $z_j$'s are all real with $z_1 >
\cdots > z_r$.
As $\bbeta$ is small, we have 
\(
   \log\left(\frac{1-e^{-\bbeta'(t-z_j)}}{\bbeta'\Lambda}\right)
   \approx \log(t - z_j).
\)
We may also suppose $t$ is real. Then the branch of $\log (t-z_j)$ is
given so that it is a real number, i.e.  $\log|t-z_j|$.
Therefore when $t < z_j$ (i.e.\ when we are integrating the summand $j
< i$), we have
\(
   \log\left(\frac{1-e^{-\bbeta'(t-z_j)}}{\bbeta'\Lambda}\right)
   = \log\left(\frac{1-e^{\bbeta'(t-z_j)}}{\bbeta'\Lambda}\right)
      - \bbeta'(t - z_j),
\)
with the branch of $\log$ in the right hand is determined so that it
is approximated by $\log |z_j - t| = \log (t - z_j)$.
Similarly we have
\(
   \int \log\left(\frac{1-e^{-\bbeta'(t-z_j)}}{\bbeta'\Lambda}\right) dt
   = -\frac1{\bbeta'}\widehat{\Li}_2(e^{\bbeta'(t-z_j)})
   -\frac{\bbeta'(t-z_j)^2}{2},
\)
and the branch of $\widehat{\Li}$ is given by the same way.

Let us turn to the third term:
\begin{equation*}
\begin{split}
-2\pi
\int_{e^{\bbeta'(z_i+\delta)}}^{e^{\bbeta' z_i^{\pm}}} dS &=
\frac{1}{\bbeta'}
\int_{e^{\bbeta'(z_i+\delta)}}^{e^{\bbeta' z_i^{\pm}}}\log X \frac{dw}{w}\\
& = \frac{1}{\bbeta'}
\int_{e^{\bbeta'(z_i+\delta)}}^{e^{\bbeta' z_i^{\pm}}}
\log e^{\bbeta' z_i} \frac{dw}{w}
+ \frac{1}{\bbeta'}
\int_{e^{\bbeta'(z_i+\delta)}}^{e^{\bbeta' z_i^{\pm}}}
(\log X-\log e^{\bbeta' z_i}) \frac{dw}{w}.
\end{split}
\end{equation*}

We take a positive number $N_{\delta} <\delta$ 
such that 
\begin{equation*}
N_{\delta}^r \left(\frac{P(e^{\bbeta'(z_i+\delta)})}
{\bbeta^{\prime r} e^{r \bbeta'(z_i+\delta)/2}}\right)^{-1} \ll \delta.
\end{equation*}
Then for $|\Lambda|<N_{\delta}$, we have
\begin{equation*}
\sqrt{1-4\Lambda^{2r}\left(\frac{P(e^{\bbeta'(z_i+\delta)})}
{\bbeta^{\prime r} e^{r \bbeta'(z_i+\delta)/2}}\right)^{-2}}=
1+O(\delta).
\end{equation*}

\begin{NB}
If $0<\bbeta' \ll 1$, then
$$
\frac{P(e^{\bbeta'(z_i+\delta)})}{\bbeta^{\prime r}} \approx
\prod_{i=1}^r(z-(z_i+\delta)),\;e^{\bbeta' z}=X.
$$
\end{NB}

We note that $w_{|e^{\bbeta' z_i^{\pm}}}=\pm 1$.
{\allowdisplaybreaks
\begin{equation*}
\begin{split}
& \frac{1}{\bbeta'}
\int_{e^{\bbeta'(z_i+\delta)}}^{e^{\bbeta' z_i^{\pm}}}
\log e^{\bbeta' z_i} \frac{dw}{w}
= 
z_i [\log w]_{e^{\bbeta'(z_i+\delta)}}^{e^{\bbeta' z_i^{\pm}}}
\\
=\; & 
z_i \left[ \log 
\frac{Y+P(X)}{2 (\bbeta' \Lambda)^r X^{(r+m)/2}} 
\right]_{e^{\bbeta'(z_i+\delta)}}^{e^{\bbeta' z_i^{\pm}}}\\
=\; & - 
z_i \log \left[\frac{1}{2(\bbeta' \Lambda)^r}
\left(\frac{P(e^{\bbeta'(z_i+\delta)})}
{w_{|e^{\bbeta' z_i^{\pm}}} e^{(r+m) \bbeta'(z_i+\delta)/2}}\right)
\left(1+\sqrt{1-4\Lambda^{2r}\left(\frac{P(e^{\bbeta'(z_i+\delta)})}
{\bbeta^{\prime r} e^{r \bbeta'(z_i+\delta)/2}}\right)^{-2}} \right)
 \right]\\
=\; & - 
z_i
\left(
\sum_{j \ne i} \log \left(\frac{1-e^{-\bbeta'(z_i-z_j)}}{\bbeta'\Lambda}\right)
  +\log\left(\frac{1-e^{-\bbeta' \delta}}{\bbeta'\Lambda}\right)
  +\frac{(r-m)\bbeta'}{2}z_i \right)+O(\delta), 
\end{split}
\end{equation*}
where} the branch of $\log$ is the same as before.
\begin{NB}
Let us explain more precisely. When $\bbeta = 0$,
\[
  \frac{P(e^{\bbeta'(z_i+\delta)})}
  {w_{|e^{\bbeta' z_i^{\pm}}} e^{r \bbeta'(z_i+\delta)/2}}
\]
converges to
$
   w|_{z_i^\pm} P_{\bbeta=0}(z_i+\delta)
$
where $P_{\bbeta=0}(z) = \prod_j (z - z_j)$. When all $z_j$ are real
and $z_1 > \cdots > z_r$, $P_{\bbeta=0}(z_i+\delta)$ and
$P_{\bbeta=0}(z_i^\pm)$ have the same sign. (This assertion is false
for $z_i-\delta$ as in \cite{NY2}.) Therefore 
\[
  w|_{z_i^\pm} P_{\bbeta=0}(z_i+\delta) = | P_{\bbeta=0}(z_i+\delta)|
 = \prod_j |z_i + \delta - z_j|
 \approx \delta \prod_{j\neq i} |z_i-z_j| .
\]
Therefore we should take the branch so that
\[
   \log \left(\frac{1-e^{-\bbeta'(z_i-z_j)}}{\bbeta'\Lambda}\right)
   \approx \log |z_i - z_j|.
\]
This is the same choice as before, i.e.\ real valued when $\bbeta=0$
and $z_1 > \cdots > z_r$.
\end{NB} 

\begin{Claim}
\begin{equation*}
\frac{1}{\bbeta'}
\int_{e^{\bbeta'(z_i+\delta)}}^{e^{\bbeta' z_i^{\pm}}}
(\log X-\log e^{\bbeta' z_i}) \frac{dw}{w}
=O(\delta).
\end{equation*}
\end{Claim}

\begin{proof}
If $X=e^{\bbeta' t}$ and $|t-z_i| \approx \delta$, we have
\begin{equation*}
\begin{split}
\frac{\log X-\log e^{\bbeta' z_i}}{X-e^{\bbeta' z_i}}&=
e^{-\bbeta' z_i}+O(\delta),\\
\frac{X-e^{\bbeta' z_i}}
{\prod_{j \ne i}(e^{\bbeta' z_i}-e^{\bbeta' z_j})^{-1}
P(X)}&=1+O(\delta).
\end{split}
\end{equation*}
Thus we get
\begin{equation*}
\log X-\log e^{\bbeta' z_i}=
e^{-\bbeta' z_i+\frac{r}{2}\bbeta' z_i}
\frac{\prod_{j \ne i}(e^{\bbeta' z_i}-e^{\bbeta' z_j})^{-1}
P(X)}{X^{\frac{r}{2}}}
+ E(X)
\end{equation*}
with $E(X) = O(\delta^2)$. The integration of $E(X)$ yields
$O(\delta^2) O(\log\delta) = O(\delta)$. 

For the main part we have
\begin{equation*}
\begin{split}
\int_{e^{\bbeta'(z_i+\delta)}}^{e^{\bbeta' z_i^{\pm}}}
X^{-\frac{r+m}{2}}P(X) \frac{dw}{w}
&=
\int_{e^{\bbeta'(z_i+\delta)}}^{e^{\bbeta' z_i^{\pm}}}
X^{-\frac{r+m}{2}}P(X) \frac{(X^{-\frac{r+m}{2}}P(X))'}{X^{-\frac{r+m}{2}}Y}dX \\
&= 
\left[X^{-\frac{r+m}{2}}Y\right]_{e^{\bbeta'(z_i+\delta)}}
^{e^{\bbeta' z_i^{\pm}}}\\
&= \bbeta^{\prime r} O(\delta).
\end{split}
\end{equation*}
Since 
$$
\prod_{j \ne i}(e^{\bbeta' z_i}-e^{\bbeta' z_j}) \approx
\bbeta^{\prime r-1} \prod_{j \ne i}(z_i-z_j),
$$
we get the assertion.
\end{proof}

The computation of the first term is similar. 
Since $O(\Lambda)\log \Lambda=O(\delta)$ for $\Lambda \ll \delta$,%
\begin{NB}
I do not understand why we need this. Could you check please ?
May. 29, HN.
\end{NB}
we have the following:
\begin{equation}\label{eq:C_i}
\begin{split}
& - 2\pi 
\int_{C_i} dS
-2\left[-\frac{r\bbeta^{\prime}}{4}(z_i^2-z_{i-1}^2)-
\frac{1}{\bbeta'}
\sum_{j > i}\widehat{\Li}_2(e^{-\bbeta'(z_i-z_j)})
+\frac{1}{\bbeta'}
\sum_{j < i}\widehat{\Li}_2(e^{-\bbeta'(z_j-z_i)})\right.\\
&\quad \quad 
+\frac{1}{\bbeta'}
\sum_{j > i-1}\widehat{\Li}_2(e^{-\bbeta'(z_{i-1}-z_j)})-
\frac{1}{\bbeta'}
\sum_{j < i-1}\widehat{\Li}_2(e^{-\bbeta'(z_j - z_{i-1})})\\
&\quad \quad \left. +\sum_{j<i} \frac{\bbeta^{\prime}(z_j-z_i)^2}{2}-
\sum_{j<i-1} \frac{\bbeta^{\prime}(z_j-z_{i-1})^2}{2}
\right]
=O(\delta).
\end{split}
\end{equation}
We now replace $\bbeta' z_i$ by $\bbeta a_i$. As $a_i + \sqrt{-1} z_i
= O(\Lambda)$ by \eqref{eq:za}, the left hand side is still
$O(\delta)$ after the replacement.%
\begin{NB}
We note that $\log (x+O(\Lambda))=\log x+O(\Lambda)$, $x>0$.
Since  
\begin{equation*}
\widehat{\Li}_2(e^{-x})=(x \log x-x)+O(x),
\end{equation*}
we have 
\begin{equation*}
\begin{split}
\widehat{\Li}_2(e^{-(x+O(\Lambda))})-\Li_2(e^{-x})
=&(x+O(\Lambda))\log(x+O(\Lambda))-x\log x+O(\Lambda)\\
=&(x+O(\Lambda))\log(1+\frac{1}{x}O(\Lambda))=O(\Lambda).
\end{split}
\end{equation*}
\end{NB}
Since the LHS is a single valued holomorphic function of $\Lambda$
on $0<|\Lambda|<N_{\delta}$, it is extended to a holomorphic function 
on $|\Lambda|<N_{\delta}$.
Since the LHS does not depend on $\delta$, 
it is 0 at $\Lambda=0$. 
Thus the left hand side of \eqref{eq:C_i} is $O(\Lambda)$. 

\begin{NB}
If $\Lambda \ne 0$, then
$O(\Lambda)$ is a (multi-valued) holomorphic function 
of $\Lambda,\bbeta,\vec{a}$.
By the above proof, there is a positive number $\delta$
which does not depend on the choice of $\bbeta$, $|\bbeta| \ll 1$
such that $|O(\Lambda)|<1$ on $|\Lambda|\leq \delta$.
Then 
$$
a_n:=\int_{|\Lambda|=\delta}
\frac{O(\Lambda)}{\Lambda^{n+1}}d\Lambda
$$
satisfies
$|a_n|<2\pi \delta^{-n}$.
Thus 
$$
O(\Lambda)=\sum_{n \geq 0} a_n \Lambda^n 
$$
is a holomorphic function of $\Lambda,\bbeta,\vec{a}$,
if $|\Lambda|<\delta$.
\end{NB}

Therefore we have
\begin{equation}
\begin{split}
  &{2\pi \sqrt{-1}} a^D_i = {2\pi \sqrt{-1}} \sum_{k=2}^i \int_{C_k} dS\\
= \; &
\left[-\frac{(r-m)\bbeta}{2}(a_i^2-a_{1}^2)-
\frac{2}{\bbeta}\sum_{j > i}
   \left( \Li_2(e^{-\bbeta(a_i-a_j)}) - \frac{\pi^2}6\right)
+\frac{2}{\bbeta}\sum_{j < i}
   \left( \Li_2(e^{-\bbeta(a_j-a_i)}) - \frac{\pi^2}6\right)\right.\\
&\quad \quad
+\frac{2}{\bbeta}\sum_{j > 1}
   \left({\Li}_2(e^{-\bbeta(a_{1}-a_j)}) - \frac{\pi^2}6\right)
+\sum_{j<i} {\bbeta(a_j-a_i)^2}
\\
&\qquad\qquad
+
2 r (a_i-a_1)\log (\bbeta \Lambda)
- r \pi \sqrt{-1} (a_i - a_1)
\Biggr]+O(\Lambda).
\end{split}
\end{equation}
By \eqref{eq:gamma'} we get \eqref{eq:to_be_checked}.
\end{proof}

\subsection{A renormalization group equation}

We assume $m\neq \pm r$ hereafter.

We give an analogue of the renormalization group equation for the
homological version (see \cite[\S2.4]{NY2}).

We set $w$ to be constant and differentiate the defining equation of
$C_{\vec{U}}$ with respect to $\log\Lambda$ to get
\begin{equation*}
  \frac{\partial X}{\partial \log\Lambda}
  = \frac{r X^{-(r+m)/2} P(X)}{(X^{-(r+m)/2}P(X))'}
  - \sum_{p=1}^{r-1} \frac{\partial U_p}{\partial\log\Lambda}
  \frac{X^{(r+m)/2-p}}{(X^{-(r+m)/2}P(X))'}.
\end{equation*}
Therefore
\begin{equation}\label{eq:dSdL}
  \begin{split}
  \left.\pd{}{\log\Lambda} dS\right|_{w=\mathrm{const}}
  & = \frac1{2\pi\sqrt{-1}\bbeta}\pd{X}{\log\Lambda} \frac{dw}{Xw}
\\
  & = \frac{1}{2\pi\sqrt{-1}\bbeta}
  \left[ \frac{rX^{-(r+m)/2}P(X)}{(X^{-(r+m)/2}P(X))'}
    - \sum_{p=1}^{r-1} \frac{\partial U_p}{\partial\log\Lambda}
  \frac{X^{(r+m)/2-p}}{(X^{-(r+m)/2}P(X))'}\right]
  \frac{dw}{Xw}
\\
  & = \frac{1}{2\pi\sqrt{-1}\bbeta} \left[
    \frac{r P(X) dX}{XY}
    - \sum_{p=1}^{r-1} \frac{\partial U_p}{\partial\log\Lambda}
    \frac{X^{r-p-1}dX}Y\right].
  \end{split}
\end{equation}
We thus have
\begin{align}
  0 = \pd{a_i}{\log\Lambda} & = 
  \frac{1}{2\pi\sqrt{-1}\bbeta} \int_{A_i}
    \frac{r P(X) dX}{XY}
    + \sum_{p=1}^{r-1} \frac{\partial U_p}{\partial\log\Lambda}
    \pd{a_i}{U_p}, \label{eq:periodA}
\\
    \pd{a_i^D}{\log\Lambda} & = 
  \frac{1}{2\pi\sqrt{-1}\bbeta} \int_{B_i}
    \frac{r P(X) dX}{XY}
    + \sum_{p=1}^{r-1} \frac{\partial U_p}{\partial\log\Lambda}
    \pd{a^D_i}{U_p}.
\end{align}
Combining these equalities we get
\begin{equation*}
  \pd{a_i^D}{\log\Lambda} = 
  \frac{1}{2\pi\sqrt{-1}\bbeta} \left[ \int_{B_i}
    \frac{r P(X) dX}{XY}
    - \sum_{j=2}^{r} \pd{a^D_i}{a_j}
     \int_{A_j} \frac{r P(X) dX}{XY}
    \right].
\end{equation*}

 From \eqref{eq:periodA} the meromorphic differential 
\(
  \frac{2\pi\sqrt{-1}\bbeta}r
  \left.\pd{}{\log\Lambda} dS\right|_{w=\mathrm{const}}
\)
has vanishing $A$-periods. Its poles are inverse images of $X =
0,\infty$. As they are not branched points, we have four points.
Let us denote them by $0_+$, $\infty_-$ ($w = \infty$), $0_-$,
$\infty_+$ ($w = 0$).
This convention is taken so that their residues are given by
\begin{equation*}
   0_\pm : \pm 1, \quad
   \infty_\pm : \pm 1.
\end{equation*}%
\begin{NB}
There was a mistake. Corrected on June 10. Let us check !
Since $P(X)/Y = (w + 1/w) / (w - 1/w)$, we have $P(X)/Y = 1$ at $w =
\infty$ and $= -1$ at $w = 0$. Therefore 
\(
  \left. \Res \left( \frac{2\pi\sqrt{-1}\bbeta}r
  \left.\pd{}{\log\Lambda} dS\right|_{w=\mathrm{const}}\right)
  \right/ \Res \frac{dX}X
\)
is $+1$ at $w= \infty$ and $-1$ at $w= 0$. As $\Res\frac{dX}X$ is $1$
at $X=0$ and $-1$ at $X=\infty$, we get the assertion.
\end{NB}
The assumption $m\neq \pm r$ is used here, otherwise $X=0$, $\infty$
may not correspond to $w=0,\infty$.
\begin{NB}
Added by H.N., Sep.\ 7, 2006   
\end{NB}

By the Riemann bilinear relation (see e.g., \cite[\S2.2]{GH}) we have
\begin{equation}\label{eq:dad}
\begin{split}
  -\frac1{2\pi\sqrt{-1}} \frac{\partial^2 \mathcal F_0}
  {\partial a_i\partial\log\Lambda}
  =
  \pd{a_i^D}{\log\Lambda} &= 
  \frac{r}{\bbeta}
  \int_{0_-+\infty_-}^{0_++\infty_+} \omega_i
  = 
  \frac{2r}{\bbeta}
  \int^{0_+}_{\infty_-} \omega_i
  ,
\end{split}
\end{equation}
where we have used the hyperelliptic involution $\iota$ in the second
equality. The path of the integral is taken disjoint from the cycles
$A_i$, $B_i$.

\begin{NB}
  This equation suggests that it is more natural to replace $a^D$ by
  $\sqrt{-1}a^D$, as we put
  $\frac{\bbeta}{2r}\pd{a_\alpha^D}{\log\Lambda}$ into theta functions
  in the contact term equations.

  I have changed $a^D$. (June 7).
\end{NB}

When $\bbeta\to 0$, two points $X = 0$, $\infty$ converge to a single
point $z = \infty$ as we observed in \subsecref{subsec:limit}. Here
more precisely, $0_+$, $\infty_-$ go to $z=\infty$, $w=\infty$ and 
$0_-$, $\infty_+$ goes to $z=\infty$, $w=0$.

As $\omega_i = \pd{}{a_i} dS$, this equation suggests 
\(
  \pd{\mathcal F_0}{\log\Lambda}
  = -\frac{4\pi\sqrt{-1} r}{\bbeta}
  \int^{0_+}_{\infty_-}
   dS
   .
\)
However the integral does not make sense as $dS$ has singularities at
$0_+$ and $\infty_-$. We overcome the difficulty by introducing a new
differential
\begin{equation*}
   dS' = \frac{Y}{P(X)} dS 
   = \frac1{2\pi \sqrt{-1}
    \bbeta} {\log X}\, \frac{X^{(r+m)/2} (X^{-(r+m)/2}P(X))' dX}{P(X)}
    .
\end{equation*}
Then $dS - dS'$ can be integrated from $0_+$ to $\infty_-$.%
\begin{NB}
Recall $P(X)/Y = 1$ at $w = \infty$, i.e., $0_+$, $\infty_-$. As
$P(X)/Y = -1$ at $w=0$, $dS + dS'$ is regular at $0_-$, $\infty_+$.  
\end{NB}
From (\ref{eq:dSdU}, \ref{eq:dSdL}) we have
\begin{equation}
\begin{split}
  \left.\frac{\partial}{\partial a_i} dS'\right|_{w=\mathrm{const}}
  & = - \frac1{2\pi \sqrt{-1}\bbeta}
    \sum_p \pd{U_p}{a_i} \frac{X^{r-p-1} dX}{P(X)}
    ,
\\
  \left.\pd{}{\log\Lambda} dS'\right|_{w=\mathrm{const}}
  & = \frac{1}{2\pi\sqrt{-1}\bbeta} \left[
    \frac{r dX}{X}
    - \sum_{p=1}^{r-1} \frac{\partial U_p}{\partial\log\Lambda}
    \frac{X^{r-p-1}dX}{P(X)}\right].
\end{split}
\end{equation}
Differentiating $P(X) = \prod (X - e^{-\sqrt{-1}\bbeta z_i})$ by
$U_p$, we get
\begin{equation*}
   \frac{X^{r-p-1}}{P(X)}
   = \sqrt{-1}\bbeta \sum_i
   \frac{e^{-\sqrt{-1}\bbeta z_i}}{X(X - e^{-\sqrt{-1}\bbeta z_i})}
   \pd{z_i}{U_p}
   = \sqrt{-1}\bbeta \sum_i
   \frac1{X - e^{-\sqrt{-1}\bbeta z_i}}
   \pd{z_i}{U_p},
\end{equation*}
where we have used $\sum_i z_i = 0$. Therefore we have
\begin{equation}\label{eq:int'}
\begin{split}
   \int_{\infty_-}^{0_+}\frac{X^{r-p-1}}{P(X)} dX
   & = \sqrt{-}\bbeta
   \left[\sum_i \pd{z_i}{U_p} \log(X-e^{-\sqrt{-1}\bbeta z_i})
   \right]_{X=\infty}^{X=0}
\\
   &= \sqrt{-1} \bbeta \sum_i \pd{z_i}{U_p}\log(-e^{-\sqrt{-1}\bbeta z_i})-
   \sum_i \pd{z_i}{U_p} \log(1-e^{-\sqrt{-1}\bbeta z_i}/X)_{|X=\infty}
\\
   &= \bbeta^2 \sum_i \pd{z_i}{U_p} z_i
   = \frac{\bbeta^2}2 \pd{}{U_p} \sum_i z_i^2
   = \frac{\bbeta^2}{2r} \pd{}{U_p} \sum_{i<j} (z_i-z_j)^2
   ,
\end{split}
\end{equation}
where we take a path in the upper half plane and we also used
$\sum_i z_i=0$.%
\begin{NB}
Let us expand as in Kota's note:
\[
   \frac{X^{r-p-1}}{P(X)}
    = \sum_i \frac{\lambda_i^p}{X-e^{-\sqrt{-1}\bbeta z_\alpha}},
\]
where
\[
\lambda_i^p=\frac{e^{-\sqrt{-1}\bbeta z_i (r-p-1)}}
{\prod_{j : j \ne i}(e^{-\sqrt{-1}\bbeta z_i}-e^{-\sqrt{-1}\bbeta z_j})}.
\]
Then $\lambda^p_i = -\bbeta \pd{z_i}{U_p}$.
\end{NB}
Therefore
\begin{equation*}
  \int_{\infty_-}^{0_+}
  \left.\frac{\partial}{\partial a_i} dS'\right|_{w=\mathrm{const}}
   = - \frac{\bbeta}{4\pi\sqrt{-1}r} \pd{}{a_i}
    \sum_{j < k} 
    (z_j - z_k)^2 
\end{equation*}
Combining with \eqref{eq:dad}, we get
\begin{equation*}
\begin{split}
   -\frac1{2\pi\sqrt{-1}}
   \frac{\partial^2 \mathcal F_0}{\partial a_i\partial\log\Lambda}
  &= \frac{2r}{\bbeta} \pd{}{a_i} \left[ 
     \int_{\infty_-}^{0_+} \left( dS - dS' \right)
    - \frac{\bbeta}{4\pi\sqrt{-1} r} \sum_{j < k} 
    (z_j - z_k)^2
  \right].
\end{split}
\end{equation*}
Therefore we have
\begin{equation}\label{eq:FL}
   -\frac1{2\pi\sqrt{-1}}
   \pd{\mathcal F_0}{\log\Lambda}
   = \frac{2r}{\bbeta} \int_{\infty_-}^{0_+} \left( dS - dS' \right)
    - \frac1{2\pi\sqrt{-1}} \sum_{i < j} 
    (z_i - z_j)^2
\end{equation}
up to a function of $\Lambda$ independent of $a_\alpha$. The right
hand side has a perturbative expansion as
\begin{equation*}
   \frac1{2\pi\sqrt{-1}}\sum_{i < j} (a_i - a_j)^2 + O(\Lambda).
\end{equation*}
This is exactly equal to the one given in
\propref{prop:perturb}. Therefore we finally fix the ambiguity of
$\mathcal F_0$ in $O(\Lambda)$ so that \eqref{eq:FL} holds.

When $\bbeta\to 0$, both points $0_+$, $\infty_-$ converges to
$z=\infty$, $w=\infty$. We have $dS = dS'$ at the limit
point. Therefore the first integral disappears in the limit and we get
\[
   \left.\pd{\mathcal F_0}{\log\Lambda}\right|_{\bbeta=0}
   = \sum_{i<j} (z_i - z_j)^2.
\]
This is nothing but the renormalization group equation
\cite[2.3]{NY2} in the homological version. On the other hand, if
$\bbeta$ stays nonzero, $\pd{\mathcal F_0}{\log\Lambda}$ could not be
expressed as a simple function in $U_p$.

We differentiate \eqref{eq:FL} by $\log\Lambda$:
{\allowdisplaybreaks
\begin{equation*}
\begin{split}
   & -\frac1{2\pi\sqrt{-1}}
   \frac{\partial^2 \mathcal F_0}{(\partial \log\Lambda)^2}
   = \frac{2r}{\bbeta} \int_{\infty_-}^{0_+} 
   \pd{}{\log\Lambda} \left( dS - dS' \right)
    - \frac1{2\pi\sqrt{-1}} \pd{}{\log\Lambda}
    \sum_{i < j} (z_i - z_j)^2
\\
  =\; &
  \frac{r}{\pi\sqrt{-1}\bbeta^2} \int_{\infty_-}^{0_+} 
  \left[
    \frac{r (P(X) - Y) dX}{XY}
    - \sum_{p=1}^{r-1} \frac{\partial U_p}{\partial\log\Lambda}
    \left(\frac{X^{r-p-1}}Y - \frac{X^{r-p-1}}{P(X)}\right)dX\right]
\\
   & \qquad\qquad - 
   \frac1{2\pi\sqrt{-1}}
    \sum_{i < j} \sum_{p=1}^{r-1}
     \pd{U_p}{\log\Lambda} \pd{}{U_p}
    (z_i - z_j)^2
\\
   =\; &
  \frac{r}{\pi\sqrt{-1}\bbeta^2} \int_{\infty_-}^{0_+} 
  \left[
    \frac{r (P(X) - Y) dX}{XY}
    - \sum_{p=1}^{r-1} \frac{\partial U_p}{\partial\log\Lambda}
    \frac{X^{r-p-1}dX}Y \right],
\end{split}
\end{equation*}
where} we have used \eqref{eq:int'} in the last equality.

Let us consider
\begin{equation}\label{eq:omega}
  \frac{(P(X) - Y) dX}{2XY}
    - \frac1{2r}
    \sum_{p=1}^{r-1} \frac{\partial U_p}{\partial\log\Lambda}
    \frac{X^{r-p-1}dX}Y.
\end{equation}
From \eqref{eq:periodA} and $\int_{A_\alpha} \frac{dX}X = 0$,  it also
has the vanishing $A$-periods. Its poles are $0_-$ and $\infty_+$
with residues $-1$ and $1$ respectively. These properties characterize
the meromorphic differential form uniquely. Let us denote it by
$\omega_{\infty+ - 0_-}$ as customary.%
\begin{NB}
There was a mistake: I corrected $\omega_{0_--\infty_+} \to 
\omega_{0_+-\infty_-}$ on June 10.
\end{NB}
Substituting this into above, we get
\begin{equation}\label{eq:d^2F_0}
   \frac{\partial^2\mathcal F_0}{(\partial\log\Lambda)^2}
   = \frac{4r^2}{\bbeta^2}
   \int_{0_+}^{\infty_-} \omega_{\infty_+ - 0_-}.
\end{equation}

\begin{NB}
We have
\begin{equation*}
   \int_{0_-}^{\infty_+} \omega_{0_+--\infty_-}
   = \int_{\iota(0_-)}^{\iota(\infty_+)} \omega_{\iota(0_+)-\iota(\infty_-)}
   = \int_{0_+}^{\infty_-} \omega_{0_- - \infty_+}.
\end{equation*}
In the version May 23, the sign was wrong. (June 10).
\end{NB}

\begin{NB} 
The following is the original proof:
\begin{proof}
We have
\begin{equation*}
   \left.\pd{}{a_\alpha} \frac{dX}{X} \right|_{w=\mathrm{const}}
   = \left.\pd{}{a_\alpha} d \log X \right|_{w=\mathrm{const}}
   = d\left( \frac1X \pd{X}{a_\alpha}\right)
   = \sum_p d\left(\frac1X \pd{X}{U_p} \right) \pd{U_p}{a_\alpha}
     .
\end{equation*}
Thanks to \eqref{eq:dXdU},
\(
   \frac1X \pd{X}{U_p}
\)
vanishes at $X = 0_-$, $\infty_+$.
\verb+\begin{NB}+
In fact,
\(
   \frac1X \pd{X}{U_p} = - \frac{X^{r-p}}{X^{r/2+1}(X^{-r/2} P(X))'}.
\)
The denominator is $(-1)^{r+1} r/2 \neq 0$ at $X = 0$. It is
approximately $X^r$ at $X = \infty$. Therefore we have
\(
   \left.\frac1X \pd{X}{U_p}\right|_{X = 0_-\ \mathrm{or}\ \infty_+}
   = 0.
\)
\verb+end{NB}+
Thus
\begin{equation*}
   \pd{}{a_\alpha} \int_{0_-}^{\infty_+} \omega_{0_+-\infty_-}
   = \frac{\pi\sqrt{-1}\bbeta}r \int_{0_-}^{\infty_+}
    \left.\frac{\partial^2}{\partial a_\alpha\partial\log\Lambda} dS
        \right|_{w=\mathrm{const}}
   = \frac{\pi\sqrt{-1}\bbeta}r \pd{}{\log\Lambda}
     \int_{0_-}^{\infty_+}\omega_\alpha
     .
\end{equation*}
This is equal to 
\(
   \frac{\pi\sqrt{-1}\bbeta^2}{2r^2}
   \frac{\partial^2 a_\alpha^D}{(\partial\log\Lambda)^2}
   =
   - \frac{\bbeta^2}{4r^2} \frac{\partial^3 \mathcal F_0}
    {\partial a_\alpha(\partial\log\Lambda)^2}
\)
by \eqref{eq:dad}. Thus we have the assertion modulo a constant
independent of $a_\alpha$.
\begin{verbatim}\begin{NB}\end{verbatim}
  We can probably prove that the constant is $0$ comparing the
  perturbative term.
\begin{verbatim}end{NB}\end{verbatim}
\end{proof}
\end{NB}

\begin{NB}
I look for an explicit formula for $\pd{\mathcal F_0}{\log\Lambda}$...
But I cannot find it so far.
\end{NB}

\subsection{Case $r+m$ even}\label{subsec:even}

We assume that $r+m$ is even in this subsection.
\begin{NB}
I added this assumption. Sep.\  7, 2006, H.N.  
\end{NB}

Recall that we set $X_1$, \dots, $X_r$ be the zeroes of $P(X) = 0$.
For small $\Lambda$, we can find $X_i^\pm$ near $X_i$ such
that $P(X_i^\pm) = \pm 2 (-\sqrt{-1}\beta\Lambda)^r(X_i^\pm)^{(r+m)/2}$.
These are branch points of the Seiberg-Witten curve $C_{\vec{U},m}$. We
have a natural partition of them as
\(
    \{ X_i^+\} \sqcup \{ X_i^-\},
\)
which corresponds to the even half-integer characteristic $E$. It is
the same as one in the homological version, i.e.\ 
\(
   {}^t(\frac12,\frac12,\frac12,\cdots).
\)
This is true regardless of the parity of $r$.

Recall that the Szeg\"o kernel of the hyperelliptic curve is
explicitly given by
\begin{equation*}
\begin{split}
\Psi_E(X_1,X_2)
  &= \frac{\Theta_E(\int^{X_2}_{X_1}\vec{\omega}|\tau)}
  {\Theta_E(0)E(X_1,X_2)}
  = \frac12 \left(
     \sqrt[4]{\frac{\psi_E(X_1)}{\psi_E(X_2)}}
     +
     \sqrt[4]{\frac{\psi_E(X_2)}{\psi_E(X_1)}}
     \right)
     \frac{\sqrt{dX_1 dX_2}}{X_2-X_1}
\\
 &= 
\frac{Y_2 \prod (X_1 - X_\alpha^+) + Y_1 \prod (X_2 - X_\alpha^+)}
   {2(X_2 - X_1)}
   \sqrt{\frac{dX_1 dX_2}{Y_1 Y_2 \prod(X_1 - X_\alpha^+)(X_2 - X_\alpha^+)}},
\end{split}
\end{equation*}
where $E$ is the prime form and
\begin{equation*}
  \psi_E(X) = \frac{\prod (X - X_\alpha^+)}{\prod (X - X_\alpha^-)}
  = \frac{P(X) - 2(-\sqrt{-1}\bbeta\Lambda)^rX^{(r+m)/2}}
  {P(X) + 2(-\sqrt{-1}\bbeta\Lambda)^rX^{(r+m)/2}}.
\end{equation*}
See \cite[p.12 Example]{Fay}.
\begin{NB}
The integral $\int_{X_2}^{X_1}$ was corrected. (June 10).
\end{NB}
We have $\psi_E(0_\pm) = \psi_E(\infty_\pm) = 1$.%
\begin{NB}
I have used $r+m\neq 0$, $2r$ also here.
Sep.\ 7, 2006, H.N.
\end{NB}
Therefore
\begin{equation}\label{eq:Szego}
    E(0_-,\infty_+)^2
    \left.{dX_1}\right|_{X_1=0_-}
      \left.\left(\frac{{dX_2}}{X_2^2}\right)\right|_{X_2=\infty_+}
    =  \frac{\Theta_E(\int_{0_-}^{\infty_+} \vec{\omega}|\tau)^2}
      {\Theta_E(0)^2}.
\end{equation}
On the other hand, \cite[p.17, Remark v)]{Fay} we have
\begin{equation}\label{eq:Fay}
   E(0_-,\infty_+)^2
    \left. dX_1\right|_{X_1=0_-}
      \left.\left(\frac{{dX_2}}{X_2^2}\right)\right|_{X_2=\infty_+}
   = \exp\left\{ \int_{0_+}^{\infty_-} \omega_{\infty_+ - 0_-}
     + \sum_{i=2}^{r} m_i
       \int_{0_+}^{\infty_-} \omega_i,
   \right\}
\end{equation}
where 
\(
  m_i = \frac1{2\pi}\int_{A_i}
  d\arg \frac{X - X(\infty_+)}{X - X(0_-)}.
\)
In our situation, this is equal to $0$. We thus get
\begin{equation}\label{eq:ThetaID}
  \frac{\Theta_E(\int_{0_-}^{\infty_+}\vec{\omega})^2}{\Theta_E(0)^2}
  = \exp\left\{
    \int_{0_+}^{\infty_-}
    \omega_{\infty_+ - 0_-}
  \right\}
  .
\end{equation}
By \eqref{eq:d^2F_0} and \eqref{eq:dad} we get
\begin{equation}\label{eq:contact}
   \frac{\Theta_E(\frac{\bbeta}{2r} \pd{a^D_\alpha}{\log\Lambda}|\tau)}
      {\Theta_E(0|\tau)}
  = \exp\left\{
    \frac{\bbeta^2}{8r^2}
    \frac{\partial^2\mathcal F_0}{(\partial\log\Lambda)^2}
  \right\}.
\end{equation}
Thus we get the contact term equation \cite[\protect (4.12) with $d=\frac{r}2$]{NY3}.%
\begin{NB}
The sign in $\exp$ was corrected. (June 10).
\end{NB}
More precisely, the above holds up to sign. However both sides go to
$1$ when $\bbeta\to 0$, so the above holds without the sign
ambiguity.
\begin{NB}
Originally I thought we can use the perturbative expansion: 

We remove the sign
umbiguity by considering the perturbative part.

I need to check this argument.
\end{NB}

\begin{NB}
As the contact term equation characterizes the prepotential $\mathcal
F_0$ up to the perturbation part, we conclude that the limit of
Nekrasov partition function is equal to the Seiberg-Witten
prepotential at least when $m=0$.
\end{NB}

\subsubsection{A differential equation for $U_p$}\label{subsubsec:diffU}

By \cite[Prop.~2.10 (38)]{Fay} we have
\begin{equation}\label{eq:Fay2.10}
   \frac{\Psi_E(X,0_-) \Psi_E(X,\infty_+)}{\Psi_E(0_-,\infty_+)}
   = \omega_{\infty_+ - 0_-} + \sum_{i=2}^r \left[
     \pd{\log\Theta_E}{\xi_i}(\int_{0_-}^{\infty_+} \vec{\omega})
     - \pd{\log\Theta_E}{\xi_i}(0)
     \right] \omega_i(X).
\end{equation}
The left hand side is equal to
\begin{equation*}
   \frac{(P(X) - Y) dX}{2XY}.
\end{equation*}%
\begin{NB}
More detail:
\begin{equation*}
  \mathrm{LHS} =
  \frac{
   \frac{Y|_{X=0_-}\prod (X - X_\alpha^+) + Y \prod (-X_\alpha^+)}
   {-2X}
   \frac{\frac{Y}{X^r}|_{X=\infty_+}\prod (X - X_\alpha^+)
     + Y \left.\prod (1-\frac{X_\alpha^+}X)\right|_{X=\infty_+}}2}
  {\frac{\frac{Y}{X^r}|_{X=\infty_+}\prod (- X_\alpha^+)
     + Y|_{X=0_-} \left.\prod (1-\frac{X_\alpha^+}X)\right|_{X=\infty_+}}2}
  \frac{dX}{Y\prod (X - X_\alpha^+)}.
\end{equation*}
We substitute $\frac{Y}{X^r}|_{X=\infty_+} = -1$,
$Y|_{X=0_-} = -\prod(-X_\alpha^+)$. Then
\begin{equation*}
\begin{split}
  & \mathrm{LHS} = \frac{(Y - \prod(X - X_\alpha^+))^2} 
  {4XY \prod (X- X_\alpha^+)}dX
  = \frac{Y^2 - 2 Y\prod (X - X_\alpha^+) + \prod(X - X_\alpha^+)^2}
  {4XY \prod (X- X_\alpha^+)}dX
\\
  =\; &
  \frac{\prod (X - X_\alpha^-) + \prod(X - X_\alpha^+) - 2Y}
  {4XY}dX
  = \frac{P(X) - Y}{2XY} dX.
  \end{split}
\end{equation*}

\end{NB}
As $E$ is an even characteristic, 
\(
  \pd{\log\Theta_E}{\xi_\alpha}(0) = 0.
\)
Looking at \eqref{eq:omega} we have
\begin{equation*}
    \frac1{2r}
    \sum_{p=1}^{r-1} \frac{\partial U_p}{\partial\log\Lambda}
    \frac{X^{r-p-1}dX}Y
    = 
    \sum_{i} 
     \pd{\log\Theta_E}{\xi_i}(\int_{0_-}^{\infty_+} \vec{\omega})
     \omega_i(X).
\end{equation*}
In other words,
\begin{equation}
    \frac1{2r}
    \frac{\partial U_p}{\partial\log\Lambda}
    = 
    - \frac1{2\pi\sqrt{-1}\bbeta}\sum_{i} 
      \left.\pd{\log\Theta_E}{\xi_i}\right|_{
      \vec{\xi} = -\frac{\bbeta}{4\pi\sqrt{-1}r}
      \frac{\partial^2 \mathcal F_0}{\partial \log\Lambda\partial \vec{a}}}
    \pd{U_p}{a_i}.
\end{equation}
This is an analog of the equation in \cite[Th.~2.4]{NY2}.
This equation suggests that it is possible to define $U_p$ in terms of
the instanton counting as in the homological version.

\begin{NB}
Let us take a limit $\bbeta\to 0$. Notice that we can replace $U_p$ by
any symmetric polynomial in $e^{-\sqrt{-1}\bbeta z_i}$. Therefore we can
replace $U_p$ by $u_p$ in the limit and get
\begin{equation*}
   \frac1{2r} \pd{u_p}{\log\Lambda}
   = - \frac1{8\pi^2 r} \sum_{i,j}
   \frac{\partial^2\log\Theta_E}{\partial\xi_i\partial \xi_j}(0)
   \frac{\partial^2 \mathcal F_0}
   {\partial\log\Lambda\partial a_j}
   \pd{u_p}{a_i}.
\end{equation*}
Note that the sign was wrong in \cite{NY2}.
\end{NB}

\subsubsection{Higher order equations}

By \cite[Cor.~2.19 (43)]{Fay} we have
\begin{equation*}
\begin{split}
  & \frac{\Theta_E(\sum_{i=1}^d y_i - \sum_{i=1}^d x_i)}
   {\Theta_E(0)}
   \frac{\prod_{i<j} E(x_i,x_j) E(y_j,y_i)}{\prod_{i,j} E(x_i,y_j)}
   = 
   \det\left(\frac{\Theta_E(y_j - x_i)}{\Theta_E(0) E(x_i,x_j)}\right)
\\
  =\; &
  \det\left(\Psi_E(x_i,y_j)\right).
\end{split}
\end{equation*}
Let us study the limit of this equation when all $x_i$ (resp.\ $y_j$)
goes to $0_-$ (resp.\ $\infty_+$).
As 
\(
   E(x_i,x_j) = \frac{(x_i - x_j)}{\sqrt{dx_i}\sqrt{dx_j}}
   \left(1 + O(x_i - x_j)^2\right),
\)
we have
\begin{equation*}
   \frac{\det\left(\Psi_E(x_i,y_j)\right)}
   {\prod_{i<j} E(x_i,x_j) E(y_j,y_i)}
   \to
   (-1)^{d(d-1)/2} \det
   \left( \left.\frac1{i! j!}\partial_x^i \partial_y^j
   (\Psi_E)(x,y)\right|_{\substack{x=0_-\\ y=\infty_+}}
   \right)_{0\le i,j\le d-1}
   .
\end{equation*}
Therefore the answer depends only on the differentials of $\Psi_E$ up
to order $d-1$.%
\begin{NB}
Corrected and commented added by H.N. Sep.\ 7, 2006.  
\end{NB}
Note that
\begin{equation*}
\begin{split}
   & \psi_E(X) 
   = \frac{P(X) - 2X^{(r+m)/2}(-\sqrt{-1}\bbeta\Lambda)^r}
      {P(X) + 2X^{(r+m)/2}(-\sqrt{-1}\bbeta\Lambda)^r}
\\
   =\; & 1 - \frac{4X^{(r+m)/2}(-\sqrt{-1}\bbeta\Lambda)^r}{P(X) +
   2X^{(r+m)/2}(-\sqrt{-1}\bbeta\Lambda)^r}
   = 
   \begin{cases}
   1 + O(X^{(r+m)/2}) & \text{as $X\to 0$},
   \\
   1 + O(X^{-(r-m)/2}) & \text{as $X\to \infty$}.
   \end{cases}
\end{split}
\end{equation*}
Therefore we can replace either $\psi_E(x_i)$ or $\psi_E(y_i)$ by $1$
when we compute the limit if $0\le d\le \max(r+m,r-m)/2$. We may
assume $m\le 0$ without loss of generality. Then we can replace
$\psi_E(y_i)$ by $1$. Thus
{\allowdisplaybreaks
\begin{equation}\label{eq:ThetaID2}
\begin{split}
   & \frac{\Theta_E(d \int_{0_-}^{\infty_+} \vec{\omega})}
   {\Theta_E(0)}
   = \left.\frac{\det \Psi_E(x_i,y_j)\prod_{i,j} E(x_i,y_j)}
   {\prod_{i<j} E(x_i,x_j) E(y_j,y_i)}\right|_{\substack{
    x_i = 0_- \\ y_j=\infty_+}}
\\
   =\; & 
   \left.\det\left(
       \frac12\left(
         \sqrt[4]{\psi_E(x_i)} + \frac1{\sqrt[4]{\psi_E(x_i)}}\right)
       \frac{\sqrt{dx_i}\sqrt{dy_j}}{y_j - x_i}\right)
     \frac{\prod_{i,j} E(x_i,y_j)}
   {\prod_{i<j} E(x_i,x_j) E(y_j,y_i)}\right|_{\substack{
    x_i = 0_- \\ y_j=\infty_+}}
\\
   =\; & 
   \begin{aligned}[t]
     & \prod_{i=1}^d\left(
     \frac12\left(
         \sqrt[4]{\psi_E(x_i)} + \frac1{\sqrt[4]{\psi_E(x_i)}}\right)
       \sqrt{dx_i} \right) \prod_{j=1}^d \sqrt{dy_j}
    \\
       & \qquad \times
       \left.
      \det\left(
       \frac{1}{y_j - x_i}\right)
     \frac{\prod_{i,j} E(x_i,y_j)}
   {\prod_{i<j} E(x_i,x_j) E(y_j,y_i)}\right|_{\substack{
    x_i = 0_- \\ y_j=\infty_+}}
   \end{aligned}
\\
   = \; &  \left(E(0_-,\infty_+) 
    \left.\sqrt{dX_1}\right|_{X_1=0_-}
      \left.\left(\frac{\sqrt{dX_2}}{X_2}\right)\right|_{X_2=\infty_+}
     \right)^{d^2}
\\
   = \; & \exp\left(\frac{d^2}2 \int_{0_+}^{\infty_-}
     \omega_{\infty_+ - 0_-}\right),
\end{split}
\end{equation}
where} we have used (\ref{eq:Szego}, \ref{eq:ThetaID}) in the last equality.
Hence we get
\begin{equation}
   \frac{\Theta_E(\frac{d \bbeta}{2r} \pd{a^D_\alpha}{\log\Lambda}|\tau)}
      {\Theta_E(0|\tau)}
  = \exp\left\{
    \frac{d^2 \bbeta^2}{8r^2}
    \frac{\partial^2\mathcal F_0}{(\partial\log\Lambda)^2}
  \right\}
\end{equation}
for $0\le d\le \max(r+m,r-m)/2$.
This is the same equation derived in \propref{prop:contact} under the
assumption \eqref{eq:conjvanish}.

\begin{NB}
It was not clear why we need the condition $d\le r/2$ in the previous
argument. 
\end{NB}

\subsection{Case $r+m$ odd}\label{subsec:rodd}

We assume that $r+m$ is odd in this subsection.

Let us introduce a new variable $W = \sqrt{X}$ and consider the
branched double covering 
\(
   p\colon \widehat C_{\vec{U},m} \to C_{\vec{U},m}
\)
given by
\begin{equation}\label{eq:2SW}
\begin{split}
  Y^2 &= P(W^2)^2 - 4(-W^2)^r (\bbeta\Lambda)^{2r}
\\
      &= \left(P(W^2) - 2 (\sqrt{-1}W \bbeta\Lambda)^{r}\right)
        \left(P(W^2) + 2 (\sqrt{-1} W \bbeta\Lambda)^{r}\right).
  \end{split}
\end{equation}
The branched points are $X = 0_\pm$, $\infty_\pm$. The genus of
$\widehat C_{\vec{U}}$ is $2r-1$.

For the new curve $\widehat C_{\vec{U},m}$ the calculation of the
previous section can be applied. We then use formulas in \cite[\S
5]{Fay} relating the theta functions for $\widehat C_{\vec{U},m}$ and
those for $C_{\vec{U},m}$. This is our strategy to prove the contact
term equation for the $r+m$ odd case.

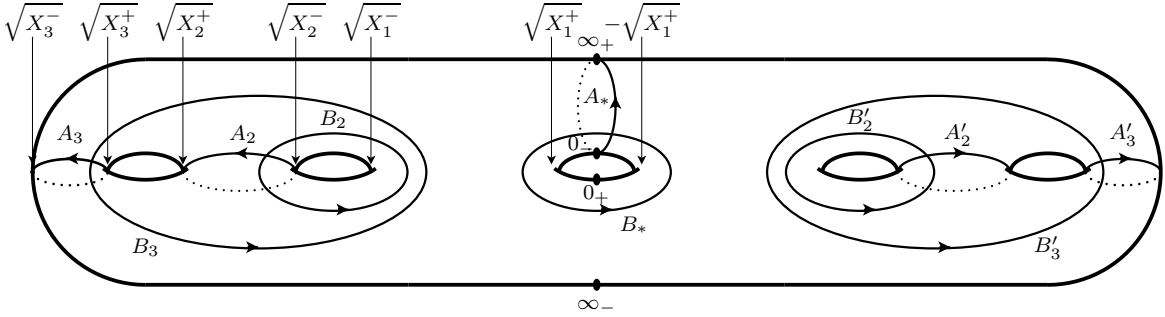
\begin{figure}[htbp]
\tiny
\centering
\psset{xunit=1mm,yunit=1mm,runit=1mm}
\psset{linewidth=0.3,dotsep=1,hatchwidth=0.3,hatchsep=1.5,shadowsize=1}
\psset{dotsize=0.7 2.5,dotscale=1 1,fillcolor=black}
\psset{arrowsize=1 2,arrowlength=1,arrowinset=0.25,tbarsize=0.7 5,bracketlength=0.15,rbracketlength=0.15}
\begin{pspicture}(0,0)(150,35)
\rput{0}(135,15){\parametricplot[linewidth=0.5,arrows=-]{-90}{90}{ t cos 15 mul t sin 15 mul }}
\rput{0}(15,15){\parametricplot[linewidth=0.5,arrows=-]{90}{270}{ t cos 15 mul t sin 15 mul }}
\rput{-0}(15,15){\parametricplot[linewidth=0.5,arrows=-]{-0}{180}{ t cos 5 mul t sin 2.5 mul }}
\rput{0}(15,17){\parametricplot[linewidth=0.5,arrows=-]{30}{150}{ t cos 6.5 mul t sin -3 mul }}
\rput{-0}(40,15){\parametricplot[linewidth=0.5,arrows=-]{-0}{180}{ t cos 5 mul t sin 2.5 mul }}
\rput{0}(40,17){\parametricplot[linewidth=0.5,arrows=-]{30}{150}{ t cos 6.5 mul t sin -3 mul }}
\rput{-0}(110,15){\parametricplot[linewidth=0.5,arrows=-]{-0}{180}{ t cos 5 mul t sin 2.5 mul }}
\rput{0}(110,17){\parametricplot[linewidth=0.5,arrows=-]{30}{150}{ t cos 6.5 mul t sin -3 mul }}
\rput{-0}(135,15){\parametricplot[linewidth=0.5,arrows=-]{-0}{180}{ t cos 5 mul t sin 2.5 mul }}
\rput{0}(135,17){\parametricplot[linewidth=0.5,arrows=-]{30}{150}{ t cos 6.5 mul t sin -3 mul }}
\rput{0}(0,15){\psellipse[](0,0)(0,0)}
\rput{0}(5,15){\parametricplot[linewidth=0.35,linestyle=dotted,dotsep=0.7,arrows=-]{-0}{180}{ t cos 5 mul t sin -1.75 mul }}
\rput{0}(5,15){\parametricplot[linewidth=0.35,arrows=-]{180}{360}{ t cos 5 mul t sin -1.75 mul }}
\rput{0}(145,15){\parametricplot[linewidth=0.35,linestyle=dotted,dotsep=0.7,arrows=-]{-0}{180}{ t cos 5 mul t sin -1.75 mul }}
\rput{0}(145,15){\parametricplot[arrows=-]{180}{360}{ t cos 5 mul t sin -1.75 mul }}
\rput{0}(27.5,15){\parametricplot[arrows=-]{180}{360}{ t cos 7.5 mul t sin -2.5 mul }}
\rput{0}(27.5,15){\parametricplot[linestyle=dotted,dotsep=0.7,arrows=-]{-0}{180}{ t cos 7.5 mul t sin -2.5 mul }}
\rput{0}(122.5,15){\parametricplot[arrows=-]{180}{360}{ t cos 7.5 mul t sin -2.5 mul }}
\rput{0}(122.5,15){\parametricplot[linestyle=dotted,dotsep=0.7,arrows=-]{-0}{180}{ t cos 7.5 mul t sin -2.5 mul }}
\rput{0}(110,15){\psellipse[](0,0)(10,-5)}
\rput{0}(120,15){\psellipse[](0,0)(22.5,-10)}
\psbezier{<-}(124,17.55)(122,17.55)(124,17.55)(124,17.55)
\psbezier{<-}(146,16.8)(144,16.8)(146,16.8)(146,16.8)
\psbezier{<-}(122,5)(120,5)(122,5)(122,5)
\psbezier{<-}(112,10)(110,10)(112,10)(112,10)
\rput{0}(40,15){\psellipse[](0,0)(10,-5)}
\psbezier{<-}(42,10)(40,10)(42,10)(42,10)
\rput{0}(30,15){\psellipse[](0,0)(22.5,-10)}
\psbezier{<-}(30,5)(28,5)(30,5)(30,5)
\psbezier{<-}(4.5,16.8)(6.5,16.8)(4.5,16.8)(4.5,16.8)
\psbezier{<-}(27,17.5)(29,17.5)(27,17.5)(27,17.5)
\rput(28,20){$A_2$}
\rput(25,20){}
\rput(5,20){$A_3$}
\rput(15,5){$B_3$}
\rput(40,22){$B_2$}
\rput(123,20){$A_2'$}
\rput(145,20){$A_3'$}
\rput(110,22){$B_2'$}
\rput(135,5){$B_3'$}
\psline[linewidth=0.5](100,30)(135,30)
\psline[linewidth=0.5](100,0)(135,0)
\psline[linewidth=0.5](15,30)(50,30)
\psline[linewidth=0.5](15,0)(50,0)
\rput{-0}(75,15){\parametricplot[linewidth=0.5,arrows=-]{-0}{180}{ t cos 5 mul t sin 2.5 mul }}
\rput{0}(75,17){\parametricplot[linewidth=0.5,arrows=-]{30}{150}{ t cos 6.5 mul t sin -3 mul }}
\rput{0}(75,15){\psellipse[](0,0)(10,-5)}
\psbezier{<-}(77,10)(75,10)(77,10)(77,10)
\rput{90}(75,23.75){\parametricplot[arrows=-]{-0}{180}{ t cos 6.25 mul t sin -2.5 mul }}
\rput{90}(75,23.75){\parametricplot[linestyle=dotted,dotsep=0.7,arrows=-]{180}{360}{ t cos 6.25 mul t sin -2.5 mul }}
\rput(75,32){$\infty_+$}
\rput(75,-3){$\infty_-$}
\rput(75,12){$0_+$}
\rput{0}(75,30){\psellipse[fillstyle=solid](0,0)(0.5,-0.5)}
\rput{0}(75,17.5){\psellipse[fillstyle=solid](0,0)(0.5,-0.5)}
\rput{0}(75,14){\psellipse[fillstyle=solid](0,0)(0.5,-0.5)}
\rput{0}(75,0){\psellipse[fillstyle=solid](0,0)(0.5,-0.5)}
\rput(73,18.5){$0_-$}
\rput(75,25){$A_*$}
\rput(80,8){$B_*$}
\psbezier{<-}(77.5,25)(77.5,23)(77.5,25)(77.5,25)
\rput(69,35){$\sqrt{X_1^+}$}
\rput(81,35){$-\sqrt{X_1^+}$}
\rput(65,15){}
\psline[linewidth=0.5](50,30)(100,30)
\psline[linewidth=0.5](50,0)(100,0)
\rput(65,15){}
\psline[linewidth=0.1]{->}(81,33)(81,16)
\psline[linewidth=0.1]{->}(69,33)(69,16)
\psline[linewidth=0.1]{->}(45,33)(45,16)
\rput(45,35){$\sqrt{X_1^-}$}
\psline[linewidth=0.1]{->}(35,33)(35,16)
\psline[linewidth=0.1]{->}(20,33)(20,16)
\psline[linewidth=0.1]{->}(10,33)(10,16)
\psline[linewidth=0.1]{->}(0,33)(0,16)
\rput(35,35){$\sqrt{X_2^-}$}
\rput(20,35){$\sqrt{X_2^+}$}
\rput(10,35){$\sqrt{X_3^+}$}
\rput(0,35){$\sqrt{X_3^-}$}
\rput(62,32){}
\end{pspicture}

\caption{Double cover of the Seiberg-Witten curve for $r=3$, $m$: even}
\label{fig:odd}
\end{figure}

Let us fix notations. See [loc.~cit.] for more detail.
Let $\phi$ be the involution $W\mapsto -W$ corresponding to the
projection $p$.
We choose a symplectic basis 
$A_2$, $B_2$,\dots, $A_r$, $B_r$, $A_*$, $B_*$, $A_2'$, $B_2'$,\dots,
$A_r'$, $B_r'$
of $H_1(\widehat C_{\vec{U}},\Z)$ as in Figure~\ref{fig:odd}, where
the involution $\phi$ is the rotation by $\pi$ about the vertical axis
passing through $0_\pm$, $\infty_\pm$.
They satisfy
\begin{enumerate}
\item $A_i$, $B_i$
$i=2,\dots r$ are taken so that they are in a single sheet of $p$
and mapped to the corresponding cycles in the original curve $C_{\vec{U}}$,
\item $A'_i = -\phi(A_i)$, $B'_i = -\phi(B_i)$,
\item $A_* + \phi(A_*) = 0 = B_* + \phi(B_*)$.
\end{enumerate}

The normalized holomorphic differentials $\Hat\omega_i$,
$\Hat\omega_*$, $\Hat\omega'_i$ on $\widehat C_{\vec{U}}$ 
satisfy
\begin{equation*}
  \phi^*\Hat\omega_i = - \Hat\omega'_i,\qquad
  \phi^*\Hat\omega_* = -\Hat\omega_*
\end{equation*}
and are related to those on $C_{\vec{U}}$ as
\begin{equation*}
   p^* \omega_i = \Hat\omega_i - \Hat\omega'_i.
\end{equation*}
%
We denote a vector in 
\( 
   \C^{2r-1}
\)
by $\left[\xi,\eta,\xi'\right]$ with $\xi,\xi'\in 
\C^{r-1}$,
$\eta\in\C$.
Let $\pi^*\colon J_0(C_{\vec{U}})\to J_0(\widehat C_{\vec{U}})$ be the
pull-back homomorphism of the divisor classes. It lifts to a map
$\C^{r-1}\to \C^{2r-1}$ by
\begin{equation*}
   \pi^*(\xi) = \left[\xi, 0, -\xi\right].
\end{equation*}
Let us choose two points $S, T$ from four branched points $0_\pm$,
$\infty_\pm$. Let $S'$, $T'$ be the remaining two points. Let $\xi_0 =
\frac14 \int_{S+T}^{S'+T'} \vec{\omega}$ where $\vec\omega =
(\omega_2,\dots,\omega_r)$ is the vector of the normalized holomorphic
differentials.
Then \cite[p.91 (102)]{Fay} says that there exists a unique half-period
$\left[0,c_*,0\right]\in J_0(\widehat C_{\vec{U}})$ such that
\begin{equation}\label{eq:ThetaCover}
  k_0 \defeq
  \frac{\widehat \Theta_{\left[c,c_*,-c\right]}(\pi^*\xi)}
    {\Theta_c(\xi+\xi_0)\Theta_c(\xi-\xi_0)} 
\end{equation}
is independent of $\xi\in\C^{r-1}$ and a half-integer characteristic $c$
for the curve $C_{\vec{U}}$. We choose $S, T = 0_-, \infty_-$, so
\begin{equation}
  \xi_0 = \frac14 \int_{0_-+\infty_-}^{0_++\infty_+}\vec\omega
  = \frac12 \int_{0_-}^{\infty_+}\vec\omega.
\end{equation}

The double cover $\widehat C_{\vec{U}}$ is also a hyperelliptic curve
by the involution $\widehat\iota\colon Y\mapsto -Y$. In
Figure~\ref{fig:odd} the involution $\widehat\iota$ is the rotation by
$\pi$ about the horizontal axis.
Note that $0_-$ and $\infty_+$ lie in the {\it same\/} sheet of the
covering $\widehat C_{\vec{U}}\to \widehat C_{\vec{U}}/\widehat\iota =
\proj^1$ as we have $P(X) \approx Y$ at both points. (We have $P(X)
\approx -Y$ in another sheet.) The sheet is the upper part of
$\widehat C_{\vec{U}}$ in Figure~\ref{fig:odd}.

The branched points are $W = \sqrt{X_i^\pm},
-\sqrt{X_i^\pm}$. (Recall that we have fixed the branch of
$\sqrt{X_i^\pm}$ so that $\sqrt{X_i^\pm} \approx \sqrt{X_i} =
e^{-\sqrt{-1}\bbeta z_i/2}$. 
We have a natural partition of them as
\(
   \{ \sqrt{X_i^+}, -\sqrt{X_i^-} \} \sqcup
   \{ \sqrt{X_i^-}, -\sqrt{X_i^+} \}.
\)
It corresponds to the factorization of the right hand side of
\eqref{eq:2SW}.%
\begin{NB}
Recall $\sqrt{X_i^\pm}$ is a solution of $P((\sqrt{X_i^\pm})^2)
=\pm 2(-\sqrt{-1}\bbeta\Lambda)^r (\sqrt{X_i^\pm})^r$. We have
$P((-\sqrt{X_i^\pm})^2)
=\mp 2(-\sqrt{-1}\bbeta\Lambda)^r (-\sqrt{X_i^\pm})^r$.
\end{NB}
Let $\widehat E$ be the corresponding even theta characteristic. We
now repeat the argument in \subsecref{subsec:even}. We do not
determine the characteristic $\widehat E$ explicitly at this moment,
as the argument goes through if $\widehat E$ corresponds to the above
partition.
We need to take the path $0_+\to\infty_-$ disjoint from $A$,
$B$-cycles. This can be accomplished if we shift $A_*$ a little
bit. For this choice, $m_i$ appeared in \eqref{eq:Fay} is also
$0$. The remaining arguments are unchanged, and by
\eqref{eq:ThetaID2} we get
\begin{equation}\label{eq:1}
    \frac{\widehat \Theta_{\widehat E}(
    2d\int^{\infty_+}_{0_-} \vechatom )}
  {\widehat \Theta_{\widehat E}(0)}
  = \exp\left\{ 2d^2\int_{0_+}^{\infty_-} \widehat \omega_{\infty_+-0_-}
  \right\},
\end{equation}
for $2d\le \max(r+m,r-m)$ (i.e.\ $d\le (\max(r+m,r-m)-1)/2$), where
$\vechatom$ is the vector of the normalized holomorphic differentials
of $\widehat C_{\vec{U}}$, and $\widehat \omega_{\infty_+-0_-}$ is the
meromorphic differential with $\Res_{\infty_+} = +1$, $\Res_{0_-}= -1$
having the vanishing $A$-periods.

\begin{Lemma}
The characteristic $\widehat E$ is of the form $\left[E, c_*, -E\right]$
where the half-period $\left[0, c_*, 0\right]$ corresponds to the
partition $\{ 0_+, \infty_+\}\sqcup \{ 0_-, \infty_-\}$ as above.
\end{Lemma}

\begin{proof}
We took the idea of proof from that of \cite[Prop.~5.3]{Fay}. We pinch
two cycles in $\widehat C_{\vec{U}}$ as in Figure~\ref{fig:degenerate}.
The limit is the union of a genus $1$ curve $C_*$ (containing $A_*$,
$B_*$) and two copies of $C_{\vec{U}}$.
These curves are glued at $P$ and $Q$ as in
Figure~\ref{fig:degenerate}, i.e., two points $P$, $Q$ in $C_*$ are
identified with a point in $C_{\vec{U}}$ and its copy in another
$C_{\vec{U}}$ respectively.

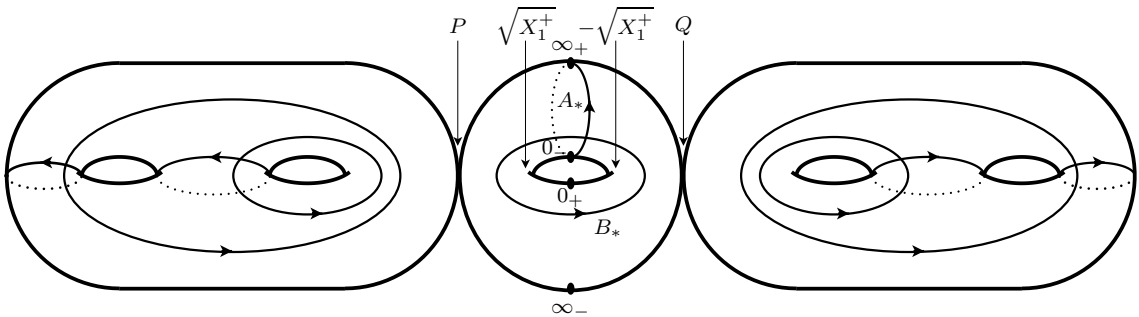
\begin{figure}[htbp]
\tiny
\centering
\psset{xunit=1mm,yunit=1mm,runit=1mm}
\psset{linewidth=0.3,dotsep=1,hatchwidth=0.3,hatchsep=1.5,shadowsize=1}
\psset{dotsize=0.7 2.5,dotscale=1 1,fillcolor=black}
\psset{arrowsize=1 2,arrowlength=1,arrowinset=0.25,tbarsize=0.7 5,bracketlength=0.15,rbracketlength=0.15}
\begin{pspicture}(0,0)(150,35)
\rput{0}(135,15){\parametricplot[linewidth=0.5,arrows=-]{-90}{90}{ t cos 15 mul t sin 15 mul }}
\rput{0}(15,15){\parametricplot[linewidth=0.5,arrows=-]{90}{270}{ t cos 15 mul t sin 15 mul }}
\rput{-0}(15,15){\parametricplot[linewidth=0.5,arrows=-]{-0}{180}{ t cos 5 mul t sin 2.5 mul }}
\rput{0}(15,17){\parametricplot[linewidth=0.5,arrows=-]{30}{150}{ t cos 6.5 mul t sin -3 mul }}
\rput{-0}(40,15){\parametricplot[linewidth=0.5,arrows=-]{-0}{180}{ t cos 5 mul t sin 2.5 mul }}
\rput{0}(40,17){\parametricplot[linewidth=0.5,arrows=-]{30}{150}{ t cos 6.5 mul t sin -3 mul }}
\rput{-0}(110,15){\parametricplot[linewidth=0.5,arrows=-]{-0}{180}{ t cos 5 mul t sin 2.5 mul }}
\rput{0}(110,17){\parametricplot[linewidth=0.5,arrows=-]{30}{150}{ t cos 6.5 mul t sin -3 mul }}
\rput{-0}(135,15){\parametricplot[linewidth=0.5,arrows=-]{-0}{180}{ t cos 5 mul t sin 2.5 mul }}
\rput{0}(135,17){\parametricplot[linewidth=0.5,arrows=-]{30}{150}{ t cos 6.5 mul t sin -3 mul }}
\rput{0}(0,15){\psellipse[](0,0)(0,0)}
\rput{0}(5,15){\parametricplot[linewidth=0.35,linestyle=dotted,dotsep=0.7,arrows=-]{-0}{180}{ t cos 5 mul t sin -1.75 mul }}
\rput{0}(5,15){\parametricplot[linewidth=0.35,arrows=-]{180}{360}{ t cos 5 mul t sin -1.75 mul }}
\rput{0}(145,15){\parametricplot[linewidth=0.35,linestyle=dotted,dotsep=0.7,arrows=-]{-0}{180}{ t cos 5 mul t sin -1.75 mul }}
\rput{0}(145,15){\parametricplot[arrows=-]{180}{360}{ t cos 5 mul t sin -1.75 mul }}
\rput{0}(27.5,15){\parametricplot[arrows=-]{180}{360}{ t cos 7.5 mul t sin -2.5 mul }}
\rput{0}(27.5,15){\parametricplot[linestyle=dotted,dotsep=0.7,arrows=-]{-0}{180}{ t cos 7.5 mul t sin -2.5 mul }}
\rput{0}(122.5,15){\parametricplot[arrows=-]{180}{360}{ t cos 7.5 mul t sin -2.5 mul }}
\rput{0}(122.5,15){\parametricplot[linestyle=dotted,dotsep=0.7,arrows=-]{-0}{180}{ t cos 7.5 mul t sin -2.5 mul }}
\rput{0}(110,15){\psellipse[](0,0)(10,-5)}
\rput{0}(120,15){\psellipse[](0,0)(22.5,-10)}
\psbezier{<-}(124,17.55)(122,17.55)(124,17.55)(124,17.55)
\psbezier{<-}(146,16.8)(144,16.8)(146,16.8)(146,16.8)
\psbezier{<-}(122,5)(120,5)(122,5)(122,5)
\psbezier{<-}(112,10)(110,10)(112,10)(112,10)
\rput{0}(40,15){\psellipse[](0,0)(10,-5)}
\psbezier{<-}(42,10)(40,10)(42,10)(42,10)
\rput{0}(30,15){\psellipse[](0,0)(22.5,-10)}
\psbezier{<-}(30,5)(28,5)(30,5)(30,5)
\psbezier{<-}(4.5,16.8)(6.5,16.8)(4.5,16.8)(4.5,16.8)
\psbezier{<-}(27,17.5)(29,17.5)(27,17.5)(27,17.5)
\rput(25,20){}
\psline[linewidth=0.5](105,30)(135,30)
\psline[linewidth=0.5](105,0)(135,0)
\psline[linewidth=0.5](15,30)(45,30)
\psline[linewidth=0.5](15,0)(45,0)
\rput{-0}(75,15){\parametricplot[linewidth=0.5,arrows=-]{-0}{180}{ t cos 5 mul t sin 2.5 mul }}
\rput{0}(75,17){\parametricplot[linewidth=0.5,arrows=-]{30}{150}{ t cos 6.5 mul t sin -3 mul }}
\rput{0}(75,15){\psellipse[](0,0)(10,-5)}
\psbezier{<-}(77,10)(75,10)(77,10)(77,10)
\rput{90}(75,23.75){\parametricplot[arrows=-]{-0}{180}{ t cos 6.25 mul t sin -2.5 mul }}
\rput{90}(75,23.75){\parametricplot[linestyle=dotted,dotsep=0.7,arrows=-]{180}{360}{ t cos 6.25 mul t sin -2.5 mul }}
\rput(75,32){$\infty_+$}
\rput(75,-3){$\infty_-$}
\rput(75,12){$0_+$}
\rput{0}(75,30){\psellipse[fillstyle=solid](0,0)(0.5,-0.5)}
\rput{0}(75,17.5){\psellipse[fillstyle=solid](0,0)(0.5,-0.5)}
\rput{0}(75,14){\psellipse[fillstyle=solid](0,0)(0.5,-0.5)}
\rput{0}(75,0){\psellipse[fillstyle=solid](0,0)(0.5,-0.5)}
\rput(73,18.5){$0_-$}
\rput(75,25){$A_*$}
\rput(80,8){$B_*$}
\psbezier{<-}(77.5,25)(77.5,23)(77.5,25)(77.5,25)
\rput(69,35){$\sqrt{X_1^+}$}
\rput(81,35){$-\sqrt{X_1^+}$}
\rput(65,15){}
\rput(65,15){}
\psline[linewidth=0.1]{->}(81,33)(81,16)
\psline[linewidth=0.1]{->}(69,33)(69,16)
\rput{0}(45,15){\parametricplot[linewidth=0.5,arrows=-]{-90}{90}{ t cos 15 mul t sin 15 mul }}
\rput{0}(105,15){\parametricplot[linewidth=0.5,arrows=-]{90}{270}{ t cos 15 mul t sin 15 mul }}
\rput{0}(75,15){\psellipse[linewidth=0.5](0,0)(15,-15)}
\psline[linewidth=0.1]{->}(90,33)(90,19)
\psline[linewidth=0.1]{->}(60,33)(60,19)
\rput(90,35){$Q$}
\rput(60,35){$P$}
\end{pspicture}

\caption{Degenerate curve}
\label{fig:degenerate}
\end{figure}

Then it is enough to calculate the characteristic in the limit. It is
clear that the $C_{\vec{U}}$-parts have characteristic $E$ and $-E$
respectively. 

Let us concentrate on the genus $1$ part. Among the
original branched points, $\pm \sqrt{X_1^+}$ are contained in $C_*$,
and $P$, $Q$ are new branched points. As the limit of the partition
corresponding to $\widehat E$, we get the partition
\(
   \{ \sqrt{X_1^+}, Q \} \sqcup \{ -\sqrt{X_1^+}, P \}.
\)
This can be seen by pinching only one of the two cycles, say one
corresponding to $Q$. As each part has the equal number of branched
points, we must have
\(
   \{ Q, \sqrt{X_1^+}, \sqrt{X_2^+}, \dots \}
   \sqcup \{ -\sqrt{X_1^+}, \sqrt{X_1^-}, \sqrt{X_2^-}, \dots \}.
\)
Pinching the remaining cycle corresponding to $P$, we get the
assertion.
On the other hand, the partition
\(
   \{ 0_+, \infty_+ \} \sqcup \{ 0_-, \infty_-\}
\)
of the branched points of $\phi$ is clearly preserved under the
degeneration.

Thus the elliptic curve $C_*$ has two hyperelliptic involutions
$\widehat\iota$ and $\phi$, and we have the corresponding partitions
of branched points
\(
   \{ \sqrt{X_1^+}, Q \} \sqcup \{ -\sqrt{X_1^+}, P \}
\)
and 
\(
   \{ 0_+, \infty_+ \} \sqcup \{ 0_-, \infty_-\}.
\)
It is clear from the picture that  both give rise to the same
characteristic of the theta function (in fact, it is $\theta_{00}$).
\end{proof}

The denominator of the left hand side of \eqref{eq:1} is
\begin{equation*}
   k_0 \Theta_E(\xi_0)^2
   = k_0 \Theta_E(\frac12 \int^{\infty_+}_{0_-} \vec\omega)^2.
\end{equation*}
On the other hand, we have
\begin{equation*}
   2d \int^{\infty_+}_{0_-} \vechatom
   = \left[ d \int^{\infty_+}_{0_-} \vec\omega,
     d, 
     -d \int^{\infty_+}_{0_-} \vec\omega\right].
\end{equation*}
To evaluating the value of the theta function at this point, we can
replace $d$ by $0$ as $d$ is an integer.
Therefore the numerator of the left hand side of \eqref{eq:1} is equal to
\begin{equation*}
  \widehat \Theta_{\widehat E}\left(
    d \int^{\infty_+}_{0_-} \vec\omega,  0, 
     -d \int^{\infty_+}_{0_-} \vec\omega\right)
  = k_0
  \Theta_E((d + \frac12) \int^{\infty_+}_{0_-} \vec\omega)
  \Theta_E((d - \frac12) \int^{\infty_+}_{0_-} \vec\omega).
\end{equation*}

On the other hand, we have $\widehat\omega_{\infty_+-0_-} = \frac12
p^*(\omega_{\infty_+-0_-})$.%
\begin{NB}
Let $p(W) = W^2 = X$. Then $\frac{dW}W = \frac12 \frac{dX}X$.
\end{NB}
Therefore the right hand side of \eqref{eq:1} is
\begin{equation*}
  \exp\left\{ d^2 \int_{0_+}^{\infty_-} \omega_{\infty_+-0_-}
  \right\}.
\end{equation*}
Thus we have
\begin{equation*}
   \frac{\Theta_E((d + \frac12) \int^{\infty_+}_{0_-} \vec\omega)}
  {\Theta_E(\frac12 \int^{\infty_+}_{0_-} \vec\omega)}
  = \exp\left\{ \frac{d(d+1)}2 \int_{0_+}^{\infty_-} \omega_{\infty_+-0_-}
  \right\},
\end{equation*}
i.e.\ 
\begin{equation}
   \frac{\Theta_E((d+\frac12)\frac{\bbeta}{2r}
   \pd{a^D_\alpha}{\log\Lambda}|\tau)}
      {\Theta_E(\frac{\bbeta}{4r}\pd{a^D_\alpha}{\log\Lambda}|\tau)}
  = \exp\left\{ \frac{d(d+1)}2
    \frac{\bbeta^2}{4r^2}
    \frac{\partial^2\mathcal F_0}{(\partial\log\Lambda)^2}
  \right\}
\end{equation}
for $0\le d\le (\max(r+m,r-m)-1)/2$.
This is the same equation derived in \propref{prop:contact} under the
assumption \eqref{eq:conjvanish}.

\begin{NB}
The following subsubsection is incomplete, and is comment out.
Nov. 9, HN
\subsubsection{A differential equation for $U_p$}

From \eqref{eq:Fay2.10} and the subsequent two sentences we have
\begin{equation*}
   \frac{(P(W^2) - Y) dW}{2WY}
   = 
   \begin{aligned}[t]
   & \widehat\omega_{\infty_+-0_-} + 
     \pd{\log\widehat\Theta_{\widehat E}}{\eta}
     (\int_{0_-}^{\infty_+} \vechatom)
     \hat\omega_{*}(X)
\\
   & \qquad + \sum_{i}
     \pd{\log\widehat\Theta_{\widehat E}}{\xi_i}
     (\int_{0_-}^{\infty_+} \vechatom)
     \hat\omega_{i}(X)
     + 
     \pd{\log\widehat\Theta_{\widehat E}}{\xi'_i}
     (\int_{0_-}^{\infty_+} \vechatom)
     \hat\omega_{i}'(X)
     ,
   \end{aligned}
\end{equation*}
where $[(\xi_i), \eta, (\xi'_i)]$ is the coordinate system of
$\C^{2r-1}$. Taking the integration over the cycles $A_i$, $A_i'$
($i=2,\dots,r$), we get
\begin{equation*}
\begin{split}
    \int_{A_i} \frac{P(X) dX}{4XY}
   &= \int_{A_i} \frac{P(W^2) dW}{2WY}
   = \pd{\log\widehat\Theta_{\widehat E}}{\xi_i}
     (\int_{0_-}^{\infty_+} \vechatom),
\\
   &= -\int_{A'_i} \frac{P(W^2) dW}{2WY}
   = -\pd{\log\widehat\Theta_{\widehat E}}{\xi'_i}
     (\int_{0_-}^{\infty_+} \vechatom).
\end{split}
\end{equation*}
From \eqref{eq:ThetaCover} we have
\begin{equation*}
   \pd{\log\widehat\Theta_{\widehat E}}{\xi_i}(\pi^*\xi)
   - \pd{\log\widehat\Theta_{\widehat E}}{\xi_i'}(\pi^*\xi)
   = \pd{\log\Theta_E}{\xi}(\xi+\xi_0) + 
   \pd{\log\Theta_E}{\xi}(\xi-\xi_0).
\end{equation*}
Therefore
\begin{equation*}
   \int_{A_i}\frac{P(X) dX}{2XY}
   = 
\end{equation*}
\end{NB}

\begin{NB}
****************

The following subsection was moved from \verb+change_var.tex+
curve.
Sep. 7, 2006, H.N.  

*******************
\end{NB}

\subsection{rank $2$ case}
We assume $r=2$, $m=0$ in this subsection.

We have $P(X) = X^2 + U_1 X + 1$. Then
\begin{equation*}
\begin{split}
Y^2 = P(X)^2 - 4X^2\bbeta^4\Lambda^4
    & = \left\{ X^2 + U_1 X + 1 + 2X\bbeta^2\Lambda^2\right\}
    \left\{ X^2 + U_1 X + 1 - 2X\bbeta^2\Lambda^2\right\}
\\
    & = \left\{\alpha_+ (X+1)^2 - \beta_+ (X-1)^2\right\}
    \left\{ \alpha_- (X+1)^2 - \beta_- (X-1)^2\right\}
\end{split}
\end{equation*}
where 
\[ 
   \alpha_\pm = \frac12 + \frac{U_1}4 \pm \frac{\bbeta^2\Lambda^2}2,
   \quad
\beta_\pm = -\frac12 + \frac{U_1}4 \pm \frac{\bbeta^2\Lambda^2}2.
\]%
\begin{NB}
I have exchanged $\alpha_\pm\leftrightarrow \alpha_\mp$.
This is compatible with our change of the Seiberg-Witten
curve in May, 2006.
Aug. 31, H.N. 
\end{NB}
Then the solutions of $P(X)^2 - 4 X^2\bbeta^4\Lambda^4 = 0$ are
\begin{equation*}
\frac{-\sqrt{\frac{\beta_+}{\alpha_+}}+1}{-\sqrt{\frac{\beta_+}{\alpha_+}}-1},
\quad
\frac{-\sqrt{\frac{\beta_-}{\alpha_-}}+1}{-\sqrt{\frac{\beta_-}{\alpha_-}}-1},
\quad
\frac{\sqrt{\frac{\beta_-}{\alpha_-}}+1}{\sqrt{\frac{\beta_-}{\alpha_-}}-1},
\quad
\frac{\sqrt{\frac{\beta_+}{\alpha_+}}+1}{\sqrt{\frac{\beta_+}{\alpha_+}}-1}.
\end{equation*}
Here we choose the branch of $\sqrt{\beta_{\pm}/\alpha_{\pm}}$ so that
the above are $X_1^+$, $X_1^-$, $X_2^-$, $X_2^+$ in
sequence. (Recall the $A$-cycle encircles $X_2^-$, $X_2^+$, and
$B$-cycles encircles $X_1^-$, $X_2^-$.)
We introduce new variables
\[
  x = \sqrt{\frac{\alpha_+}{\beta_+}} \frac{X+1}{X-1},
 \quad
  y =  \frac1{\sqrt{\beta_+\beta_-}} \frac{Y}{(X-1)^2}
\]
Then the Seiberg-Witten curve is
\begin{equation*}
    y^2 = (1 - x^2)(1 - \kappa^2 x^2),
\end{equation*}
where
\begin{equation*}
  \kappa = \sqrt{\frac{\alpha_-\beta_+}{\alpha_+\beta_-}}
  = \sqrt{1 + \frac{\bbeta^2\Lambda^2}
    {\frac{U_1^2}{16} - \left(\frac{\bbeta^2\Lambda^2}2 + \frac12\right)^2}}.
\end{equation*}

In the $x$-coordinates, the
$A$-cycle encircles $1$, $1/\kappa$, and the $B$-cycle encircles $\pm
1/\kappa$.
Note that the $A$-cycle encircles $\pm 1$ usually, so $A$, $B$-cycles
are interchanged in our convention. Note also that the curve has
period $2\tau$ instead of $\tau$ usually.
Therefore when we use various formulas in textbooks (e.g.\
\cite{W-W}), we need to replace $\tau$ by $-2/\tau$.%
\begin{NB}
In \cite{W-W}, the symbol $\theta_i$, $i=1,2,3,4$ are used:
$\theta_1=-\theta_{11}, \theta_2=\theta_{10},\theta_3=\theta_{00},
\theta_4=\theta_{01}$.
\end{NB}
From \cite[22 $\cdot$ 11]{W-W} we have 
\begin{equation*}
   \sqrt{\frac{\alpha_-\beta_+}{\alpha_+\beta_-}} =
   \kappa = \frac{\theta_{10}(-2/\tau)^2}{\theta_{00}(-2/\tau)^2}
   .
\end{equation*}%
\begin{NB}
If $0<\kappa<1$, then there is a unique 
$\widetilde{\tau} \in {\Bbb R}_{>0}\sqrt{-1}$ such that
$\kappa=
\frac{\theta_{10}(\widetilde{\tau})^2}{\theta_{00}(\widetilde{\tau})^2}$.
$\widetilde{\tau}$ is given by
$$
\int_1^{1/\kappa^2}\frac{dx}{\sqrt{(1-x^2)(1-\kappa^2 x^2)}}/
\int_0^1\frac{dx}{\sqrt{(1-x^2)(1-\kappa^2 x^2)}}
$$
where the path $[1, 1/\kappa]$ is (almost) contained in the upper half plane
and $\sqrt{(1-x^2)(1-\kappa^2 x^2)}_{|x=0}=1$.
In our case, $0<\kappa(\tau)<1$ and
$(-2/\tau) \in {\Bbb R}_{>0}\sqrt{-1}$ if $-\sqrt{-1}\bbeta>0$,
$1 \gg \Lambda>0$ and $z_1>0$.
\begin{equation}\label{eq:U_1}
\begin{split}
 \frac{\partial a}{\partial U_1}
  & = - \frac1{2\pi\sqrt{-1}\bbeta} \int_A \frac{dX}Y
  = \frac1{2\pi\sqrt{-1}\bbeta\sqrt{\alpha_+\beta_-}}
   \int_1^{1/\kappa} \frac{dx}y\\
\frac{\partial a^D}{\partial U_1}
  & = - \frac1{2\pi\sqrt{-1}\bbeta} \int_B \frac{dX}Y
  = \frac2{2\pi\sqrt{-1}\bbeta\sqrt{\alpha_+\beta_-}}
   \int_1^{0} \frac{dx}y
\end{split}
\end{equation}
where we have to take a branch of $y$ such that $y_{|x=0}=1$.
Indeed it is easy to see that $(X,Y)=(-1,4\sqrt{\beta_+\beta_-})$
correspond to $(x,y)=(0,1)$
Therefore $\widetilde{\tau}=(-2/\tau)$. 
2006/9/5, K.Y.
\end{NB}
Therefore
\begin{equation}\label{eq:U_1}
\begin{split}
  \frac{U_1^2}{16}
  & = \frac{\bbeta^2\Lambda^2}
  {\kappa^2-1} + \left(\frac{\bbeta^2\Lambda^2}2 + \frac12 \right)^2
  = - \bbeta^2\Lambda^2
  \frac{\theta_{00}(-2/\tau)^4}{\theta_{01}(-2/\tau)^4}
   + \left(\frac{\bbeta^2\Lambda^2}2 + \frac12 \right)^2
\\
  &= \frac14\left(
    1 - \frac{\theta_{00}(\tau)^4 + \theta_{10}(\tau)^4}
    {\theta_{00}(\tau)^2\theta_{10}(\tau)^2}\bbeta^2\Lambda^2 + \bbeta^4\Lambda^4
  \right).
\end{split}
\end{equation}

We also have
\begin{equation*}
  \begin{split}
  \frac{\partial a}{\partial U_1}
  & = - \frac1{2\pi\sqrt{-1}\bbeta} \int_A \frac{dX}Y
  = \frac1{2\pi\sqrt{-1}\bbeta\sqrt{\alpha_+\beta_-}}
   \int_1^{1/\kappa} \frac{dx}y
  = \frac{K'(-2/\tau)}{2\pi\bbeta\sqrt{\alpha_+\beta_-}}
  .
  \end{split}
\end{equation*}%
\begin{NB}
Use \cite[22 $\cdot$ 32]{W-W} to relate $\int_1^{1/\kappa}dx/y$
to $K$ with $\kappa'$. 2006/9/5, K.Y. 
\end{NB}
Note that $\alpha_-\beta_+ = \alpha_+\beta_-
+\bbeta^2\Lambda^2$. Therefore
\begin{equation*}
    \alpha_+\beta_- = - \bbeta^2\Lambda^2
    \frac{\theta_{00}(-2/\tau)^4}{\theta_{01}(-2/\tau)^4},
    \quad
    \alpha_-\beta_+ = - \bbeta^2\Lambda^2
    \frac{\theta_{10}(-2/\tau)^4}{\theta_{01}(-2/\tau)^4}.
\end{equation*}
Substituting $K'(-2/\tau) = \sqrt{-1}\pi\theta_{00}(-2/\tau)^2/\tau$ 
(\cite[22 $\cdot$ 32]{W-W}) we get
\begin{equation}\label{eq:dadU}
  \frac{\partial a}{\partial U_1}
  = \frac{\theta_{01}(-2/\tau)^2}{2\bbeta^2\Lambda\tau}
  = \sqrt{-1}\frac{\theta_{00}(\tau)\theta_{10}(\tau)}{2\bbeta^2\Lambda}.
\end{equation}
Here we fix the sign so that it coincides with the formula
for the homological version when $\bbeta \to 0$, i.e.\
\(
  {da}/{du} 
  = -\sqrt{-1}{\theta_{00}(\tau)\theta_{10}(\tau)}/{2\Lambda}.
\)%
\begin{NB}
I think that 
$$
\sqrt{\alpha_+ \beta_-}=(-\sqrt{-1})\Lambda 
\frac{\theta_{00}(-2/\tau)^2}{\theta_{01}(-2/\tau)^2}$$
if $-\sqrt{-1},\Lambda>0$.
It seems that $U_1=-2\sqrt{1+u+\bbeta^4 \Lambda^4}$.
So I think that 
\begin{equation}
  \frac{\partial a}{\partial U_1}
  = \frac{\theta_{01}(-2/\tau)^2}{2\bbeta^2\Lambda\tau}
  = \sqrt{-1}\frac{\theta_{00}(\tau)\theta_{10}(\tau)}{2\bbeta^2\Lambda}.
\end{equation}
(I think that $\frac{\partial a}{\partial U_1} \in {\Bbb R}_{<0}\sqrt{-1}$).
2006/9/5, K.Y.

The sign is corrected. Sep.\ 7, 2006, H.N. The mistake was come from
the mistake in the equation at the end of \subsecref{subsec:limit}.
Sep.\ 7, 2006, H.N.
\end{NB}

Let $\sn(\bullet,\kappa(-2/\tau))$, $\cn(\bullet,\kappa(-2/\tau))$,
$\dn(\bullet,\kappa(-2/\tau))$ be Jacobi's elliptic functions for the
period $-2/\tau$.
From \eqref{eq:dad} we have
\begin{equation*}
  \bbeta \pd{a^D}{\log\Lambda}
  = {4 
    }{}
  \int_{0_-}^{\infty_+} \omega
  = - \frac{2\sqrt{-1}}{K'(-2/\tau)}
  \int_{-\sqrt{\frac{\alpha_+}{\beta_+}}}^{\sqrt{\frac{\alpha_+}{\beta_+}}}
  \frac{dx}y
  = - \frac{4\sqrt{-1}}{K'(-2/\tau)}
  \sn^{-1} (\sqrt{\frac{\alpha_+}{\beta_+}}).
\end{equation*}
Here we have used that $\omega$ is normalized so that
\(
  \int_A \omega = 2 \int_{1}^{1/\kappa} \omega = 1,
\)
and hence $\omega = \frac{dx}{2\sqrt{-1}K'y}$.

Let 
\(
   h := - \frac14 \frac{\partial^2 \mathcal F_0}
                       {\partial a\partial\log\Lambda}
      = \frac{\pi\sqrt{-1}}2\pd{a^D}{\log\Lambda}.
\)
Then by using addition theorem for theta functions and the definition 
of Jacobi's elliptic functions, 
\begin{equation*}
  \begin{split}
    \frac{\theta_{11}(\frac{\bbeta h}{2\pi\sqrt{-1}},\tau)}
    {\theta_{01}(\frac{\bbeta h}{2\pi\sqrt{-1}},\tau)}
  &= \frac{\theta_{10}(\frac{\bbeta h}{4\pi\sqrt{-1}},\frac{\tau}2)
    \theta_{11}(\frac{\bbeta h}{4\pi\sqrt{-1}},\frac{\tau}2)}
  {\theta_{00}(\frac{\bbeta h}{4\pi\sqrt{-1}},\frac{\tau}2)
    \theta_{01}(\frac{\bbeta h}{4\pi\sqrt{-1}},\frac{\tau}2)}
  = \sqrt{-1}\frac{\theta_{01}(\frac{\bbeta h}{2\pi\sqrt{-1}\tau},-\frac2{\tau})
    \theta_{11}(\frac{\bbeta h}{2\pi\sqrt{-1}\tau},-\frac2{\tau})}
  {\theta_{00}(\frac{\bbeta h}{2\pi\sqrt{-1}\tau},-\frac2{\tau})
    \theta_{10}(\frac{\bbeta h}{2\pi\sqrt{-1}\tau},-\frac2{\tau})}
\\
  &= -\sqrt{-1} \kappa'(-\frac2{\tau})
  \frac{\sn(\frac{K\bbeta h}{\pi\sqrt{-1}\tau},\kappa(-\frac2{\tau}))}
  {\cn(\frac{K\bbeta h}{\pi\sqrt{-1}\tau},\kappa(-\frac2{\tau}))
    \dn(\frac{K \bbeta h}{\pi\sqrt{-1}\tau},\kappa(-\frac2{\tau}))},
  \end{split}
\end{equation*}
where $K = K(-2/\tau)$. As 
\(
  \frac{K\bbeta h}{\pi\sqrt{-1}\tau}
  = - \frac{\bbeta h}{2\pi\sqrt{-1}} \sqrt{-1} K' =
-\sn^{-1}(\sqrt{\frac{\alpha_+}{\beta_+}}),
\)
the above is equal to
\begin{equation*}
\setbox0=\hbox{$1 - \dfrac{\alpha_-\beta_+}{\alpha_+\beta_-}$}
  -\sqrt{-1}
  \sqrt{\copy0}
  \sqrt{\frac{\alpha_+}{\beta_+}}
  \sqrt{\frac{\beta_+}{\beta_+-\alpha_+}}
  \sqrt{\frac{\beta_-}{\beta_--\alpha_-}}
  = \pm \bbeta\Lambda,
\end{equation*}
where we have used $\alpha_\pm - \beta_\pm = 1$. Hence we get
\begin{equation}\label{eq:sn}
    \frac{\theta_{11}(\frac{\bbeta h}{2\pi\sqrt{-1}},\tau)}
    {\theta_{01}(\frac{\bbeta h}{2\pi\sqrt{-1}},\tau)}
   = - \bbeta\Lambda.  
\end{equation}
Here the sign was fixed by considering the limit $\bbeta\to 0$:
\begin{equation*}
   \frac1\bbeta\frac{\theta_{11}(\frac{\bbeta h}{2\pi\sqrt{-1}},\tau)}
   {\theta_{01}(\frac{\bbeta h}{2\pi\sqrt{-1}},\tau)}
   \xrightarrow{\bbeta\to 0}
   - \frac{\theta_{11}'(0,\tau)}{\theta_{01}} \frac1{8\pi\sqrt{-1}}
   \frac{\partial^2 \mathcal F_0}{\partial a\partial\log\Lambda}
   = - \Lambda
\end{equation*}%
\begin{NB}
In more detail:
\begin{equation*}
- \frac{\theta_{11}'(0,\tau)}{\theta_{01}} \frac1{8\pi\sqrt{-1}}
   \frac{\partial^2 \mathcal F_0}{\partial a\partial\log\Lambda}
   = -\pi\theta_{00}\theta_{10} \frac1{2\pi\sqrt{-1}}\frac{du}{da}
   = -\pi\theta_{00}\theta_{10}
   \frac\Lambda{\pi}\frac1{\theta_{00}\theta_{10}}
\end{equation*}
\end{NB}
The equation~\eqref{eq:sn} can be also derived from the blowup formula
\cite[Prop.~3.2(1)]{NY3} for $c_1 = \mathrm{odd}$ together with the
argument in \cite[\S4.3]{NY3}.

\end{document}